\pdfoutput=1

\documentclass[12pt,reqno]{amsart}

\usepackage{pdflscape}
\usepackage{cancel}
\usepackage[dvipsnames]{xcolor}
\usepackage{mathrsfs}
\usepackage{amsmath,accents}
\usepackage{amsfonts}
\usepackage{amssymb}
\usepackage[all,cmtip]{xy}
\usepackage{hyperref}
\usepackage{enumerate}  
\usepackage{enumitem}
\usepackage{ulem}
\usepackage[all]{xy}
\usepackage[toc,page]{appendix}
\usepackage{tikz-cd}
\usepackage[all,cmtip]{xy}
\usepackage[toc,page]{appendix}
\usepackage{amsmath, amsthm, amsfonts, amssymb}
\usepackage{graphicx}
\usepackage{subcaption}
\usepackage{tikz-cd}
\usepackage[utf8]{inputenc}
\usepackage[english]{babel}
\usepackage{dsfont}
\usepackage{xcolor}
\usepackage{arydshln}
\usepackage{faktor}
\usepackage{mathtools}
\usepackage{enumitem}
\usepackage{scalerel}

\def\endpf{\hbox{\vrule height1.5ex width.5em}}

\usepackage{dsfont}

\makeatletter
\newcommand*{\da@rightarrow}{\mathchar"0\hexnumber@\symAMSa 4B }
\newcommand*{\da@leftarrow}{\mathchar"0\hexnumber@\symAMSa 4C }
\newcommand*{\xdashrightarrow}[2][]{%
	\mathrel{%
		\mathpalette{\da@xarrow{#1}{#2}{}\da@rightarrow{\,}{}}{}%
	}%
}
\newcommand{\xdashleftarrow}[2][]{%
	\mathrel{%
		\mathpalette{\da@xarrow{#1}{#2}\da@leftarrow{}{}{\,}}{}%
	}%
}
\newcommand{\xdashdownarrow}[2][]{%
	\mathrel{%
		\mathpalette{\da@xarrow{#1}{#2}\da@downarrow{}{}{\,}}{}%
	}%
}
\newcommand*{\da@xarrow}[7]{%
	\sbox0{$\ifx#7\scriptstyle\scriptscriptstyle\else\scriptstyle\fi#5#1#6\m@th$}%
	\sbox2{$\ifx#7\scriptstyle\scriptscriptstyle\else\scriptstyle\fi#5#2#6\m@th$}%
	\sbox4{$#7\dabar@\m@th$}%
	\dimen@=\wd0 %
	\ifdim\wd2 >\dimen@
	\dimen@=\wd2 %
	\fi
	\count@=2 %
	\def\da@bars{\dabar@\dabar@}%
	\@whiledim\count@\wd4<\dimen@\do{%
		\advance\count@\@ne
		\expandafter\def\expandafter\da@bars\expandafter{%
			\da@bars
			\dabar@ 
		}%
	}%
	\mathrel{#3}%
	\mathrel{%
		\mathop{\da@bars}\limits
		\ifx\\#1\\%
		\else
		_{\copy0}%
		\fi
		\ifx\\#2\\%
		\else
		^{\copy2}%
		\fi
	}%
	\mathrel{#4}%
}
\makeatother

\usepackage{graphicx}

\newcommand\tikzmark[1]{%
  \tikz[remember picture,overlay]\coordinate (#1);}

\let\amp=&

\DeclareSymbolFont{extraup}{U}{zavm}{m}{n}
\DeclareMathSymbol{\varheart}{\mathalpha}{extraup}{86}
\DeclareMathSymbol{\vardiamond}{\mathalpha}{extraup}{87}

\newtheorem{theorem}{Theorem}[section]
\newtheorem{lemma}[theorem]{Lemma}
\newtheorem{corollary}[theorem]{Corollary}
\newtheorem{proposition}[theorem]{Proposition}
\newtheorem{claim}{Claim}[theorem]

\theoremstyle{definition}
\newtheorem{question}[theorem]{Question}
\newtheorem{definition}[theorem]{Definition}

\newtheorem{remark}[theorem]{Remark}

\newtheorem{definitionlemma}[theorem]{Definition-Lemma}

\date{}

\makeatletter
\let\@wraptoccontribs\wraptoccontribs
\makeatother

\begin{document}

\title[Canonical blow-ups of Grassmannians II]{\bf  Canonical blow-ups of Grassmannians II}

\author[Hanlong Fang and Mingyi Zhang]{Hanlong Fang and Mingyi Zhang}

\vspace{3cm} \maketitle

\begin{abstract}
We give a linear algebraic construction of the Lafforgue spaces associated to the Grassmannians $G(2,n)$ by blowing up certain explicitly defined monomial ideals, which sharpens and generalizes a result of Faltings. As an application, we provide a family of  homogeneous varieties with high complexity and with nice compactifications, which exhibits the notion of homeward compactification introduced in our previous work in a  non-spherical setting.
\end{abstract}

\tableofcontents

\section{Introduction}

A celebrated theorem of Kapranov \cite{Ka1,Ka2} 
shows that the ubiquitous Grothendieck-Knudsen compactification $\overline M_{0,n}$ of stable $n$-pointed projective lines is isomorphic to the Chow quotient of the Grassmannian $G(2,n)$ by the maximal torus.  Faltings \cite{Fal1,Fal2} and Lafforgue \cite{L1,L2}  related the sub-torus actions on Grassmannians to compacifications in algebraic geometry and number theory from different perspectives: Faltings investigated  compactifications of quotients of Cartesian products of projective linear groups by the diagonal subgroups, in order to resolve singularities of local models of Shimura varieties in places of bad reduction; Lafforgue's original motivation is to compactify the truncated moduli stacks of the Drinfeld shtukas with arbitrary level structures.   Their constructions coincide surprisingly though the approaches are distinct. Later on, Lafforgue established a general theory  beyond the Cartesian product case. From a geometric point of view,  one of the most interesting and beautiful part of the theory is the case of $r=2$, where the compactification is smooth and the boundary is a simple normal crossing divisor, which extends De Concini-Procesi compactification to a non-spherical setting. 

A major technique used in \cite{L1,L2} is a refined notion of the secondary polytopes introduced by \cite{KSZ} where the subdivisions are further required to be compatible with the thin Schubert cells (see \cite{GS,GGMS}). Compared with the case of the Grothendieck-Knudsen compactification where cross-ratios are exploited, the approach is of a complicated combinatorial nature. 
Faltings \cite{Fal2} was concerned with a similar question to derive normal forms in terms of purely linear algebra data beyond the case of De Concini-Procesi compactification of symmetric spaces \cite{Fal1}. 
Alternatively, Cartwright-Sturmfels \cite{CS} introduced a (singular) compactification as multigraded Hilbert schemes. 

In this paper, we are concerned with the following natural question.

\begin{question}\label{quest}
How to construct the Faltings-Lafforgue compactification by linear algebraic data when $r=2$ as that for $\overline M_{0,n}$ by sequential blow-ups \cite{Ka1,Ka2,GHP,Ke}?
\end{question}

Before proceeding, we first introduce  notations following \cite{L2}, \cite{Fal1},  \cite{Fan}.  

Fix a vector $\underline s=(s_1,s_2,\cdots,s_N)$ such that each component is a positive integer, and denote $n:=\sum_{t=1}^Ns_t$. Let $E:=E_1\oplus E_2\oplus\cdots\oplus E_N$ be a sum of free $\mathbb Z$-modules, where $E_t$ is of rank $s_t$ for $1\leq t\leq N$.  Denote by ${\rm Gr}^{2,E}$ the Grassmannian scheme over ${\rm Spec}\,\mathbb Z$ which represents the functor parametrizing rank $(n-2)$ free quotients of $E$ (see \S 3.2.7 of \cite{EH} for instance). 

Denote by $\mathbb V^{\underline s}$ the set consisting of vectors $\underline v=(v_1,v_2,\cdots, v_N)\in\mathbb Z
^{N}$ such that $\sum_{t=1}\nolimits^Nv_{t}=2$ and $0\leq v_{t}\leq s_{t}$ for  $1\leq t\leq N$. For each $\underline v\in\mathbb V^{\underline s}$, denote  $\wedge^{\underline v}E_{\bullet}:=\wedge^{v_1}E_1\otimes\wedge^{v_2}E_2\otimes\cdots\otimes\wedge^{v_N}E_N$. We have the Pl\"ucker embedding \begin{equation*}
e:{\rm Gr}^{2,E}\hookrightarrow\mathbb P\left(\wedge^2E\right)=\left\{(x_{\underline v})_{\underline v\in\mathbb V^{\underline s}}\in \left(\prod\nolimits_{\underline v}\wedge^{\underline v}E_{\bullet}-\{0\}\right)\Big/\mathbb G_m\right\}.
\end{equation*}
We can define an open subscheme of ${\rm Gr}^{2,E}$ by
\begin{equation*}
\begin{split}
{\rm Gr}^{2,E}_{\mathbb V^{\underline s}}:&=\left\{F\in E\left|\,\dim\left(F\cap \left(\bigoplus\nolimits_{t\in I}E_{t}\right)\right)=\min\left\{2,\sum\nolimits_{t\in I}s_{t}\right\},\,\,\forall\, I\subset\{1,2,\cdots,N\}\right.\right\}\\
&=\left\{\left.(x_{\underline v})_{\underline v\in\mathbb V^{\underline s}}\in\left(\prod\nolimits_{\underline v}\wedge^{\underline v}E_{\bullet}-\{0\}\right)\Big/\mathbb G_m\right|x_{\underline v}\neq 0,\,\,\forall\, \underline v\in\mathbb V^{\underline s}\right\}.\\
\end{split}  
\end{equation*} 
We define a diagonal sub-torus of the maximal torus $T_{\rm max}$ by
\begin{equation*}
T_{\underline s}:= \left\{\left.(\underbrace{\delta_1,\ldots,\delta_1}_{s_1},\underbrace{\delta_2,\ldots,\delta_2}_{s_2},\ldots,\underbrace{\delta_N,\ldots,\delta_N}_{s_N}) \right|\,\delta_{t} \in \mathbb G_m,\,\forall\, 1\leq t\leq N\right\},
\end{equation*}
 which acts 
on $G(2,n)$ by right matrix multiplication. Since ${\rm Gr}^{2,E}_{\mathbb V^{\underline s}}$ is preserved by the action of the group ${\rm Aut}(E_1)\times\cdots\times{\rm Aut}(E_N)$, and in particular its center $T_{\underline s}\cong\mathbb G^N_m$, we can
define a quasi-projective scheme
\begin{equation*}
\overline {\rm Gr}^{2,E}_{\mathbb V^{\underline s}}:=  {\rm Gr}^{2,E}_{\mathbb V^{\underline s}}\big/\mathbb G^N_{m}.    
\end{equation*}

Lafforgue constructed a scheme $\overline{\Omega}^{\mathbb V^{\underline s},E}$ over ${\rm Spec}\,\mathbb Z$, and, among other things, proved that
\begin{theorem}[Theorem 3.6 of \cite{L2}]
$\overline{\Omega}^{\mathbb V^{\underline s},E}$ is an equivariant smooth compactification of $\overline {\rm Gr}^{2,E}_{\mathbb V^{\underline s}}$, and the boundary $\overline{\Omega}^{\mathbb V^{\underline s},E}\mathbin{\scaleobj{1.5}{\backslash}}\overline {\rm Gr}^{2,E}_{\mathbb V^{\underline s}}$ is a simple normal crossing divisor.
\end{theorem}

For convenience, we write the Pl\"ucker embedding $e:{\rm Gr}^{2,E}\rightarrow\mathbb P\left(\wedge^2E\right)$ as $e:G(2,n)\rightarrow\mathbb P^{N_{2,n}}$, where $N_{2,n}:=\frac{n(n-1)}{2}-1$.  Denote by \begin{equation*}\left[\cdots ,z_I,\cdots\right]_{\,I\in\mathbb I_{2,n}}
\end{equation*} the homogeneous coordinates for the ambient projective space $\mathbb P^{N_{2,n}}$, where \begin{equation*}
\mathbb I_{2,n}:=\left\{(i_1,i_2)\in\mathbb Z^2\big|1\leq i_1<i_2\leq n\right\}.
\end{equation*}

In the following, we  introduce certain special sets of monomials in $\mathbb Z\left[\cdots ,z_I,\cdots\right]_{\,I\in\mathbb I_{2,n}}$, which will play a fundamental role in the paper. Note that $\underline s$ divides $\{1,\cdots, n\}$ into $N$  blocks
\begin{equation*}
\Delta_t(\underline s):=\left\{1+\sum\nolimits^{t-1}_{i=1}s_i,2+\sum\nolimits^{t-1}_{i=1}s_i,\cdots,\sum\nolimits^t_{i=1}s_i\right\},\,\,\,\,1\leq t\leq N.
\end{equation*}
For each $\underline v\in\mathbb V^{\underline s}$ , we define 
\begin{equation*}
\mathbb I^{\underline s}_{\underline v}:=\left\{(i_1,i_2)\in\mathbb I_{2,n}\left|{\rm \,the\,\, cardinality\,\, of\,\,}\{i_1,i_2\}\cap\Delta_t(\underline s)\,\,{\rm is\,\,}v_t,\,\,\forall\,1\leq t\leq N\right.\right\}.  
\end{equation*}
Now, for each $\underline w\in\mathbb Z
^{N}$, we define  
$G^{\underline s}_{\underline w}$ to be the set consisting of all monomials $z_{I_1}z_{I_2}\cdots z_{I_k}$ such that there is a decomposition $\underline w=\underline v_1+\underline v_2+\cdots+\underline v_k$, where $\underline v_i\in\mathbb V^{\underline s}$ and $I_i\in \mathbb I^{\underline s}_{\underline v_i}$ for all $1\leq i\leq k$. 
We denote by ${\rm ht}_{\underline s}(\underline w)$ the number of vectors $\underline v_i\in\mathbb V^{\underline s}$ appearing in each decomposition when $G^{\underline s}_{\underline w}$ is nonempty, and adopt the convention that ${\rm ht}_{\underline s}(\underline w)=0$ when $G^{\underline s}_{\underline w}$ is empty.

Denote by $\mathcal G^{\underline s}_{\underline w}$ the ideal sheaf of the closed subscheme of $\mathbb P^{N_{2,n}}$
determined by the homogeneous ideal generated by $G^{\underline s}_{\underline w}$. Note that when $G^{\underline s}_{\underline w}$ is empty $\mathcal G^{\underline s}_{\underline w}$ is the structure sheaf. We then define an  ideal sheaf on $G(2,n)$ by
\begin{equation*}
\mathscr I^{\underline s}_{\underline w}:=e^{-1}\mathcal G^{\underline s}_{\underline w}\cdot\mathcal O_{G(2,n)}, 
\end{equation*}
which is the inverse image ideal sheaf under the Pl\"ucker embedding.

Faltings proved that

\begin{theorem}[Theorem 11 of \cite{Fal2}]\label{fall} Let $\underline s=(2,2,\cdots,2)$. There exists a positive integer $\mathcal B$ depending on the combinotorics of $\mathbb V^{\underline s}$, such that the blow-up of $G(2,n)$ with respect to the product of all ideal sheaves $\mathscr I^{\underline s}_{\underline w}$ with  ${\rm ht}_{\underline s}(\underline w)\leq\mathcal B$ can reproduce the compactification $\overline{\Omega}^{\mathbb V^{\underline s},E}$.
\end{theorem}

We observe that a combination of Theorem \ref{fall} and Lemma 3.1 of \cite{L2} yields that the bound $\mathcal B$ can be simply taken as $2$, and hence no combinatorical information is required at all.  By a far more delicate method, we generalize it to sub-torus actions and  show that among all the vectors $\underline w$ with ${\rm ht}_{\underline s}(\underline w)=2$ only the ones similar to the cross-ratios for $\overline M_{0,n}$ are essential. 

More precisely,  denote by $C^{\underline s}$ the set consisting of the following vectors $\underline w=(w_t)_{1\leq t\leq N}$.
\begin{enumerate}
\item $w_{i_1}=w_{i_2}=1$ for certain $1\leq i_1<i_2\leq N$, and $w_t=0$ otherwise.
    
\item $w_{i_1}=2$ for a certain $1\leq i_1\leq N$ such that $s_{i_1}\geq 2$, and $w_t=0$ otherwise.

\item $w_{i_1}=w_{i_2}=w_{i_3}=w_{i_4}=1$ for certain $1\leq i_1<i_2<i_3<i_4\leq N$, and $w_t=0$ otherwise.
\end{enumerate}
For each $\underline w\in C^{\underline s}$, denote by $\mathbb P^{N^{\underline s}_{\underline w}}$ the projective space with dimension $N^{\underline s}_{\underline w}:={\rm card}\left(G^{\underline s}_{\underline w}\right)-1$. Notice that the blow-up of $G(2,n)$ with respect to the product of ideal sheaves $\prod\nolimits_{\underline w\in C^{\underline s}}\mathscr I^{\underline s}_{\underline w}$ can be canonically  identified with a closed subscheme of $\mathbb P^{N_{2,n}}\times\prod\nolimits_{\underline w\in C^{\underline s}}\mathbb {P}^{N^{\underline s}_{\underline w}}$  (see \S \ref{pre} for details). Denote by $\mathcal M^{\underline s}_{n}$ the  image of the blow-up under the projection from $\mathbb P^{N_{2,n}}\times\prod\nolimits_{\underline w\in C^{\underline s}}\mathbb {P}^{N^{\underline s}_{\underline w}}$ to $\prod\nolimits_{\underline w\in C^{\underline s}}\mathbb {P}^{N^{\underline s}_{\underline w}}$.

We answer Question \ref{quest} as follows.
\begin{theorem}\label{CR}
$\mathcal M^{\underline s}_{n}$ is isomorphic to $\overline{\Omega}^{\mathbb V^{\underline s},E}$ over ${\rm Spec}\,\mathbb Z$.
\end{theorem}

Note that when $\underline s=(1,1,\cdots,1)$, the  nontrivial terms in $\prod\nolimits_{\underline w\in C^{\underline s}}\mathscr I^{\underline s}_{\underline w}$ correspond to the cross ratios for $\overline M_{0,n}$. Also, similar to the transformation law among cross ratios (see \cite{GHP}), our method leads to transformation laws among local coordinate charts.

\medskip 

As an application, we show that a simple resolution of the singularities of the blow-up defined above indeed yields a nice compactification of  ${\rm Gr}^{2,E}_{\mathbb V^{\underline s}}$ that is compatible with the compactification $\overline{\Omega}^{\mathbb V^{\underline s},E}$ of $\overline {\rm Gr}^{2,E}_{\mathbb V^{\underline s}}$ in the sense of \cite{Fan}.  

Let us be more precise. For each $1\leq t\leq N$, denote by $\mathcal G^{\underline s}_t$ the ideal sheaf of the closed subscheme of $\mathbb P^{N_{2,n}}$
determined by the homogeneous ideal generated by $\left\{z_I\right\}_{I\in\mathbb I^{\underline s}_{t}}$, where
\begin{equation*}
\mathbb I^{\underline s}_{t}:=\left\{(i_1,i_2)\in\mathbb I_{2,n}\big|\,\{i_1,i_2\}\cap\Delta_t(\underline s)\neq\emptyset\right\}.
\end{equation*}
Let $\mathscr I^{\underline s}_{t}$ be the ideal sheaf on $G(2,n)$ given by
\begin{equation*}
\mathscr I^{\underline s}_t:=e^{-1}\mathcal G^{\underline s}_t\cdot\mathcal O_{G(2,n)}.
\end{equation*}
Denote by $\mathbb P^{N^{\underline s}_{t}}$ the projective space with dimension $N^{\underline s}_{t}:={\rm card}\left(\mathbb I^{\underline s}_{t}\right)-1$. Denote by $\mathcal T^{\underline s}_n$ the blow-up of $G(2,n)$ with respect to the product of ideal sheaves $\left(\prod\nolimits_{\underline w\in C^{\underline s}}\mathscr I^{\underline s}_{\underline w}\right)\cdot\left(\prod\nolimits_{t=1}^N\mathscr I^{\underline s}_t\right)$. Since $\mathcal T^{\underline s}_n$ can be canonically  identified with a closed subscheme of $\mathbb P^{N_{2,n}}\times\prod\nolimits_{t=1}^N\mathbb P^{N^{\underline s}_{t}}\times\prod\nolimits_{\underline w\in C^{\underline s}}\mathbb {P}^{N^{\underline s}_{\underline w}}$ (see \S \ref{pre} for details),  the projection from 
$\mathbb P^{N_{2,n}}\times\prod\nolimits_{t=1}^N\mathbb P^{N^{\underline s}_{t}}\times\prod\nolimits_{\underline w\in C^{\underline s}}\mathbb {P}^{N^{\underline s}_{\underline w}}$ to $\prod\nolimits_{\underline w\in C^{\underline s}}\mathbb {P}^{N^{\underline s}_{\underline w}}$ induces a morphism $\mathcal P^{\underline s}_n:\mathcal T^{\underline s}_n\rightarrow\mathcal M^{\underline s}_n$.

Then, we have
\begin{theorem}\label{HCR}
$\mathcal T^{\underline s}_n$ is a smooth, projective scheme over ${\rm Spec}\,\mathbb Z$ such that the following holds. 
\begin{enumerate}[label={\rm(\Alph*)}]

\item ${\rm Gr}^{2,E}_{\mathbb V^{\underline s}}$ is an open subscheme of $\mathcal T^{\underline s}_n$, and the embedding ${\rm Gr}^{2,E}_{\mathbb V^{\underline s}}\hookrightarrow\mathcal T^{\underline s}_n$ is ${\rm Aut}(E_1)\times\cdots\times{\rm Aut}(E_N)$-equivariant.  The boundary $\mathcal T^{\underline s}_n\mathbin{\scaleobj{1.5}{\backslash}}{\rm Gr}^{2,E}_{\mathbb V^{\underline s}}$ is a simple normal crossing divisor. 

\item The morphism $\mathcal P^{\underline s}_n:\mathcal T^{\underline s}_n\rightarrow\mathcal M^{\underline s}_n$ is flat and ${\rm Aut}(E_1)\times\cdots\times{\rm Aut}(E_N)$-equivariant. 

\item There is a closed subscheme $\mathfrak M^{\underline s}_n$ of $\mathcal T^{\underline s}_n$ such that:
\begin{enumerate}[label={\rm(\alph*)}]
\item The restriction of $\mathcal P^{\underline s}_n$ induces an isomorphism from $\mathfrak M^{\underline s}_n$ to $\mathcal M^{\underline s}_n$.
\item $\mathfrak M^{\underline s}_n$ is a connected component of the fixed point scheme of $\mathcal T^{\underline s}_n$ under the sub-torus action $T_{\underline s}$ induced from  ${\rm Gr}^{2,E}$.


\end{enumerate}

\end{enumerate}

\end{theorem}

Theorem \ref{HCR} shows that the notion of homeward compactification introduced in \cite{Fan} exists for non-spherical homogeneous varieties of arbitrarily high complexity.\medskip

We now briefly describe the organization of the paper and the basic ideas for the proof. A large part of the paper is devoted to the smoothness of $\mathcal M^{\underline s}_n$, $\mathcal T^{\underline s}_n$, and the flatness of  $\mathcal P^{\underline s}_n$. For that purpose, we introduce an algorithm to produce inductively explicit coordinate charts for $\mathcal T^{\underline s}_n$: each time we perform a subsequence of smooth, local blow-ups by bubble sort algorithm, such that the resulting scheme is a product of certain affine spaces and  $\mathcal T^{\underline s^{\prime}}_{n^{\prime}}$ with $n^{\prime}<n$. To prove that $\mathcal M^{\underline s}_n$
is isomorphic to the Lafforgue space, we first deal with the lowest-dimensional nontrivial case $N=4$ by the explicit local coordinate charts for $\mathcal M^{\underline s}_{n}$ and $\overline{\Omega}^{\mathbb V^{\underline s},E}$ to construct a suitable embedding of a torsor over $\mathcal M^{\underline s}_n$ to the torsor ${\Omega}^{\mathbb V^{\underline s},E}$ defined by Lafforgue. For the general case, we modify Lemma 3.3 of \cite{ST} to establish the finiteness of the functorial face morphism to the product of the Lafforgue spaces with $N=4$, and then complete the proof by the smoothness of $\mathcal M^{\underline s}_n$, which bypass the formidable computation.  

The organization of the paper is as follows. In \S\ref{pre}, after introducing notations and basic properties, we give a precise definition of $\mathcal M^{\underline s}_n$ by iterated blow-ups, and discuss certain geometry of Grassmannians. To  illustrate the idea of the proof, we first investigate the well-studied case of $\overline M_{0,n}$ in \S \ref{mta}; \S \ref{sta} is devoted to the general case by the same idea. In  \S \ref{ahom}, we establish Theorem \ref{HCR} by studying the geometry of $\mathcal T^{\underline s}_{n^{\prime}}$.  We prove Theorem \ref{CR} when $N=4$ and when $N\geq 5$ respectively in \S \ref{n4} and \S \ref{n5}. 

For convenience of the reader, we attach two appendices to the paper. In Appendix \ref{ccm4}, we provide an explicit open cover of $\mathcal M^{\underline s}_n$ when $N=4$, up to permutations of blocks and that of columns within a block. In Appendix \ref{le4}, we prove Claim \ref{torsor} for the remaining cases.
\medskip

{\noindent\bf Acknowledgement.}  The first author is very grateful to
Prof.~Xiaowen Hu, Prof.~Zhan Li, and Prof.~Enlin Yang for answering questions on algebraic geometry with patience, and Dr.~Xian Wu for helpful discussions on Chow quotients. This research is supported by National Key R\&D Program of China (No.~2022YFA1006700) and NSFC grant (No.~12201012).

\section{Preliminaries}\label{pre}

Throughout the paper, we will call $\underline s=(s_1,s_2,\cdots,s_N)$ a {\it size vector}, if each $s_i$ is a positive integer. For brevity, we will write $\mathbb A^{m}$ for the affine spaces ${\rm Spec}\,\mathbb Z\left[x_1,x_2,\cdots,x_m\right]$ over ${\rm Spec}\,\mathbb Z$, and $\mathbb P^{m}$ for the projective spaces ${\rm Proj}\,\mathbb Z\left[x_0,x_1,\cdots,x_m\right]$ over ${\rm Spec}\,\mathbb Z$.


\subsection{Notations and basic properties}

We recall that
\begin{definition}[Definition 6.14 of \cite{GW}] Let $f : X \rightarrow Y$ be a morphism of schemes.
We say that $f$ is smooth (of relative dimension $d$) at $x\in X$, if there exist affine open neighborhoods $U$ of $x$ and $V={\rm Spec}\, R$ of $f(x)$ such that $f(U)\subset V$, and an open
immersion
\begin{equation*}
j:U\hookrightarrow {\rm Spec}\,R[T_1,\cdots,T_n]/(f_1,\cdots,f_{n-d})
\end{equation*}
of $R$-schemes for suitable $n$ and $f_i$, such that the Jacobian matrix
\begin{equation*}
\left(\frac{\partial f_i}{\partial T_j}(x)\right)\in   M_{(n-d)\times n}\left(\mathcal O_{X,x}/\mathfrak m_x\right)  
\end{equation*}
has rank $n-d$. Here $\mathcal O_{X,x}$ is the local ring
of $X$ at $x$, $\mathfrak m_x$ is the maximal ideal of $\mathcal O_{X,x}$, and $\mathcal O_{X,x}/\mathfrak m_x$ is the residue
field. We say that $f:X\rightarrow Y$ is smooth (of relative dimension
$d$), if it is smooth (of relative dimension $d$) at all points $x\in X$.
\end{definition}

When investigating images under projective morphisms, the following lemma is important.
\begin{lemma}[Theorem 13.40 of \cite{GW}]\label{proper} Let $X$ be a scheme. Then, the structure morphism ${\rm Proj}\,\mathcal O_X[x_0,\cdots,x_n] \rightarrow X$ is
separated, quasi-compact,  and
universally closed.
\end{lemma}

Let $X$ and $Y$ be integral, smooth, quasi-projctive schemes over ${\rm Spec}\,\mathbb Z$. Let $\mathcal R_{{\rm Spec}\,\mathbb Z}(X,Y)$ be the set of pairs $(U,\widetilde f)$ where $U\subset X$ is an open subscheme of $X$ and $\widetilde f:U\rightarrow Y$ is a ${\rm Spec}\,\mathbb Z$-morphism. We call $(U,\widetilde f)$, $(V,\widetilde g)\in\mathcal R_{{\rm Spec}\,\mathbb Z}$ equivalent if there exists an open  subscheme $W\subset U\cap V$ such that $\widetilde f|_W=\widetilde g|_W$. Clearly this is an equivalence relation
on $\mathcal R_{{\rm Spec}\,\mathbb Z}$. Let ${\rm dom}(f)= {\rm dom}_{{\rm Spec}\,\mathbb Z}(f)$ be the set of points of $x\in X$
such that there exists a representative $(U,\widetilde f)$ of $f$ such that $x\in U$ and $\widetilde f$ is an ${\rm Spec}\,\mathbb Z$-morphism.
\begin{definitionlemma}[Definition 9.26 and Proposition 9.27 of \cite{GW}]
An equivalence class in $\mathcal R_{{\rm Spec}\,\mathbb Z}(X,Y)$ is called a ${{\rm Spec}\,\mathbb Z}$-rational map from 
$X$ to $Y$. For a ${{\rm Spec}\,\mathbb Z}$-rational map $f$ from $X$ to $Y$, we denote it by $f:X\dashrightarrow Y$, and call ${\rm dom}(f)$ the domain of definition of $f$. Then there exists a unique ${{\rm Spec}\,\mathbb Z}$-morphism $f_0:{\rm dom}(f)\rightarrow Y$ in the class $f$; we denote $f_0$ by $f$ by a slight abuse of notation.
\end{definitionlemma}

The following well-known lemma provides a paradigm to establish flatness in the paper.

\begin{lemma}\label{zmodel}
Let $f:\mathbb Z[t]\rightarrow\mathbb Z[x,y]$ be the ring homomophism defined by $t\mapsto xy$. Then $f$ is flat.
\end{lemma}
{\noindent\bf Proof of Lemma \ref{zmodel}.}
According to Theorem 1 in Chapter 2 of \cite{Mat}, to show that $f$ is flat, it suffices to show that the following holds. If $a_i\in \mathbb Z[t]$, $x_i\in \mathbb Z[x,y]$, $1\leq i\leq r$, and $\sum_{i=1}^rf(a_i)x_i=0$, then there exists an integer $s$, and elements $b_{ij}\in \mathbb Z[t]$ and $z_j\in \mathbb Z[x,y]$, $1\leq j\leq s$, such that     $\sum_{i=1}^r{a_i}b_{ij}=0$ for all $j$ and $x_i = \sum_{j=1}^s f(b_{ij})z_j$ for all $i$.

Write 
$x_i= \sum_{m,l\ge 0}B^i_{m,l}x^{m}y^{l}\in\mathbb Z[x,y]$ for $1\leq i\leq r$. Since there are only finite many $B^i_{m,l}\neq 0$, we can define
\begin{equation*}
s:=\max\{i+m+l\left|\right.B_{m,l}^i\neq 0\}.    
\end{equation*} 
For  $1\leq i\leq r$ and $1\leq j\leq s$, define
\begin{equation*}
\left\{\begin{aligned}
&b_{ij} := \sum\nolimits_{m\ge 0}B^i_{m,m+j}t^m \in \mathbb Z[t]\\
&z_j := y^j\in\mathbb Z[x,y]\\
\end{aligned}\right.. 
\end{equation*}

Then Lemma \ref{zmodel} follows.\,\,\,\,\,$\endpf$.
\medskip

Recall that
\begin{lemma}\label{prin}
Let $X$ be a scheme of finite type, and $\mathscr{I}_1$, $\mathscr{I}_{2}$ be non-zero ideal sheaves on $X$. Let $\varphi: {\rm{Bl}}_{\mathscr{I}_1\cdot\mathscr{I}_2}X\to X$ be the blow-up of $X$ with respect to the product of ideal sheaves $\mathscr{I}_1\cdot\mathscr{I}_2$. Then $\varphi^{-1}(\mathscr{I}_1)\cdot \mathcal O_{{\rm{Bl}}_{\mathscr{I}_1\cdot\mathscr{I}_2}X}$ is an invertible ideal sheaf.
\end{lemma}
{\noindent\bf  Proof of Lemma \ref{prin}.} Since the property is local, we may assume that $X={\rm Spec}\, A$, and the ideals $I_1$, $I_2$, associated to $\mathscr{I}_1$, $\mathscr{I}_2$ are given by $I_1=(a_1,\ldots,a_r)$, $I_2=(b_1,\ldots, b_s)$, respectively. 

Recall that ${\rm{Bl}}_{\mathscr{I}_1\cdot\mathscr{I}_2}X = \textrm{Proj}(R)$, where $R:=A\oplus J\oplus J^2 \oplus\cdots$ is the Rees algebra of $J:=I_1\cdot I_2$ in $A$. Then ${\rm{Bl}}_{\mathscr{I}_1\cdot\mathscr{I}_2}X$ has an open cover $\{D_+(a_ib_j)\}_{}$. Without loss of generality, it suffices to show that the ideal $\varphi^{-1}(I_1)\cdot R_{(a_1b_1)}$ is principal. Computation yields that
\begin{equation*}
\begin{split}
&R_{(a_1b_1)}=A\left[\frac{a_1b_1}{a_1b_1},\cdots,\frac{a_ib_j}{a_1b_1},\cdots,\frac{a_rb_s}{a_1b_1}\right],\\
&\varphi^{-1}(I)\cdot R_{(a_1b_1)}= (a_1,\cdots,a_r)\cdot A\left[\frac{a_1b_1}{a_1b_1},\cdots,\frac{a_ib_j}{a_1b_1},\cdots,\frac{a_rb_s}{a_1b_1}\right].\\  
\end{split}  
\end{equation*}
Notice that for any $1\leq i\leq r$,
\begin{equation*}
a_i\cdot\frac{a_1b_1}{a_1b_1}-a_1\cdot\frac{a_ib_1}{a_1b_1}=0,    
\end{equation*}
hence
\begin{equation*}
a_i=a_1\frac{a_ib_1}{a_1b_1}\in (a_1) \cdot R_{(a_1b_1)}.   
\end{equation*}

The proof is complete.\,\,\,\,$\endpf$
\medskip

The following lemma plays an important role in our inductive argument. 

\begin{lemma}\label{stur}Let $X_0$ be an integral scheme of finite type, and $\mathscr{I}_1$, $\mathscr{I}_{2}$, $\cdots, \mathscr{I}_m$ be non-zero ideal sheaves on $X_0$. For $1\leq i\leq m$, inductively define  $f_i:X_i\rightarrow X_{i-1}$ to be the blow-up of $X_{i-1}$ with respect to the ideal sheaf $\left(f_1\circ f_2\circ\cdots\circ f_{i-1}\right)^{-1}(\mathscr{I}_i)\cdot \mathcal O_{X_{i-1}}$. Then, $X_m$ is isomorphic to the blow-up of $X_0$ with respect to the product of ideal sheaves $\mathscr I_1\cdot\mathscr I_2\cdots\mathscr I_m$. \end{lemma} {\bf\noindent Proof of Lemma \ref{stur}.} It suffices to prove for the case of $m=2$.

Let $\rho:Y\rightarrow X_0$ be the blow-up of $X_0$ with respect to the product of ideal sheaves $\mathscr I_1\cdot\mathscr I_2$. By the universal property of blowing up, we can show that $(f_1\circ f_2)^{-1}(\mathscr I_1)\cdot\mathcal O_{X_2}$, $(f_1\circ f_2)^{-1}(\mathscr I_2)\cdot\mathcal O_{X_2}$ are invertible, hence the product of ideal sheaves \begin{equation*}
\left((f_1\circ f_2)^{-1}(\mathscr I_1)\cdot\mathcal O_{X_2}\right)\cdot\left((f_1\circ f_2)^{-1}(\mathscr I_2)\cdot\mathcal O_{X_2}\right)=(f_1\circ f_2)^{-1}(\mathscr I_1\cdot\mathscr I_2)\cdot\mathcal O_{X_2}    
\end{equation*} is invertible. Then there is a unique morphism  $G:X_2\rightarrow Y$ factoring $\rho$.

By Lemma \ref{prin}, we can conclude that $\rho^{-1}(\mathscr I_1)\cdot\mathcal O_Y$, $\rho^{-1}(\mathscr I_2)\cdot\mathcal O_Y$ are invertible. Thus, by the universal property of blowing up, $Y\rightarrow X_{i-1}$ inductively lifts to $Y\rightarrow X_i$ for $i=1,2$. This defines a morphism $F:Y\rightarrow X_2$ factoring $f_1\circ f_2$. Since $Y$ and $X_2$ are integral, and $G$, $F$ are isomorphisms on open subschemes, they must be isomorphisms. 

We complete the proof of Lemma \ref{stur}. \,\,\,$\endpf$ 

\subsection{Precise  definition of \texorpdfstring{$\mathcal M_n^{\underline s}$}{dd} and an interpretation by iterated blow-ups}\label{pdm} 
\begin{definition}\label{rspn} Denote by $\mathcal Q^{\underline s}_n$ the blow-up of $G(2,n)$ with respect to the product of ideal sheaves $\prod\nolimits_{t=1}^N{\mathscr I}^{\underline s}_t$. Denote the blow-up by $\widetilde{\mathcal R}^{\underline s}:\mathcal Q_{n}^{\underline s}\rightarrow G(2,n)$.
\end{definition}

By Lemma \ref{stur}, we can conclude that $\mathcal T_{n}^{\underline s}$ is the blow-up of $\mathcal Q_{n}^{\underline s}$
with respect to the product of the inverse image ideal sheaves $\prod\nolimits_{\underline w\in C^{\underline s}}\left(\widetilde{\mathcal R}^{\underline s}\right)^{-1}\mathscr I^{\underline s}_{\underline w}\cdot\mathcal O_{\mathcal Q_{n}^{\underline s}}$.

\begin{definition}\label{ret}
Denote by $\widehat{\mathcal R}^{\underline s}_n:\mathcal T_{n}^{\underline s}\rightarrow G(2,n)$, $\overline{\mathcal R}^{\underline s}_n:\mathcal T_{n}^{\underline s}\rightarrow\mathcal Q_{n}^{\underline s}$ the corresponding blow-ups. 
\end{definition}

Next, we identify $\mathcal Q_{n}^{\underline s}$ with a closed subscheme of $\mathbb P^{N_{2,n}}\times\prod\nolimits_{t=1}^N\mathbb {P}^{N^{\underline s}_{t}}$ via a sequence of blow-ups. Let $\sigma$ be a permutation of $\{1,2,\cdots,N\}$. Let $\phi_1^{\sigma}:Y^{\sigma}_1\rightarrow G(2,n)$ be the blow-up of $G(2,n)=:Y^{\sigma}_0$ with respect to the ideal sheaf $\mathscr I^{\underline s}_{\sigma(1)}$. For  $2\leq t\leq N$, inductively define $\phi^{\sigma}_{t}:Y^{\sigma}_{t}\rightarrow Y^{\sigma}_{t-1}$ to be the blow-up of $Y^{\sigma}_{t-1}$ with respect to the inverse image ideal sheaf  $(\phi^{\sigma}_1\circ \phi^{\sigma}_{2}\circ\cdots\circ \phi^{\sigma}_{t-1})^{-1}\mathscr I^{\underline s}_{\sigma(t)}\cdot\mathcal O_{Y^{\sigma}_{t-1}}$. According to Lemma \ref{stur}, we have $Y^{\sigma} _N\cong\mathcal Q_{n}^{\underline s}$. 

Define a rational map $\widetilde{\mathcal K}^{{\underline s}}:G(2,n)\dashrightarrow\mathbb P^{N_{2,n}}\times\prod\nolimits_{t=1}^N\mathbb {P}^{N^{\underline s}_{t}}$ by 
\begin{equation}\label{awhk}
\widetilde{\mathcal K}^{{\underline s}}:=\left(e,\left(\cdots,F^{\underline s}_{t}\circ e,\cdots\right)_{1\leq t\leq N}\right),  
\end{equation}
where rational maps $F^{\underline s}_{t}:\mathbb P^{N_{2,n}}\dashrightarrow \mathbb P^{N^{\underline s}_{t}}$ are given by
\begin{equation}\label{ftr}
F^{\underline s}_{t}\left([\cdots ,z_I,\cdots]_{I\in\mathbb I_{2,n}}\right):=\left[\cdots, z_{I},\cdots\right]_{{I}\in\mathbb I^{\underline s}_{t}}.
\end{equation}
Regarding the blow-ups  $\phi^{\sigma}_{t}:Y^{\sigma}_{t}\rightarrow Y^{\sigma}_{t-1}$,
$1\leq t\leq N$, we have by definition that
\begin{equation*}
Y^{\sigma}_{t}={\rm Proj}_{Y^{\sigma}_{t}}  \left(\bigoplus\nolimits_{d\geq 0}\left((\phi^{\sigma}_1\circ \phi^{\sigma}_{2}\circ\cdots\circ \phi^{\sigma}_{t-1})^{-1}\mathscr I^{\underline s}_{\sigma(t)}\cdot\mathcal O_{Y^{\sigma}_{t-1}}\right)^d\right). 
\end{equation*}
Therefore, we can prove by induction that $Y^{\sigma}_t$ is a closed subscheme of $\mathbb P^{N_{2,n}}\times\prod\nolimits_{i=1}^t\mathbb {P}^{N^{\underline s}_{i}}$ for $1\leq t\leq N$.
Moreover,  as a topological space, $Y^{\sigma}_N$ is the closure of the image $\widetilde{\mathcal K}^{{\underline s}}\left(U\right)$ in $\mathbb P^{N_{2,n}}\times\prod\nolimits_{t=1}^N\mathbb {P}^{N^{\underline s}_{t}}$, where $\widetilde{\mathcal K}^{{\underline s}}$ is defined by (\ref{awhk}) and $U$ is the domain of  definition of $\widetilde{\mathcal K}^{{\underline s}}$. 

Similarly, we have the following natural identification of $\mathcal T_{n}^{\underline s}$ with a closed subscheme of $\mathbb P^{N_{2,n}}\times\prod\nolimits_{t=1}^N\mathbb {P}^{N^{\underline s}_{t}}\times\prod\nolimits_{\underline w\in C^{\underline s}}\mathbb {P}^{N_{\underline w}^{\underline s}}$ via sequential blow-ups. Denote by $N_{\underline s}$ the cardinality of $C^{\underline s}$. Take a total order  $\sigma$ on $C^{\underline s}$, or equivalently, a bijection $\sigma:\left\{1,\cdots,N_{\underline s}\right\}\rightarrow C^{\underline s}$. Let $\psi^{\sigma}_{1}:W^{\sigma}_{1}\rightarrow \mathcal Q_{n}^{\underline s}$ be the blow-up of $\mathcal Q_{n}^{\underline s}=:W^{\sigma}_{0}$ with respect to the inverse image ideal sheaf  $\left(\widetilde{\mathcal R}^{\underline s}\right)^{-1}\mathscr I^{\underline s}_{\sigma(1)}\cdot\mathcal O_{\mathcal Q_{n}^{\underline s}}$. For $2\leq k\leq N_{\underline s}$, inductively  define $\psi^{\sigma}_{k}:W^{\sigma}_{k}\rightarrow W^{\sigma}_{k-1}$ to be the blow-up of $W^{\sigma}_{k-1}$ with respect to  $\left(\widetilde{\mathcal R}^{\underline s}\circ\psi^{\sigma}_{1}\circ\psi^{\sigma}_{2}\circ\cdots\circ\psi^{\sigma}_{k-1}\right)^{-1}\mathscr I^{\underline s}_{\sigma(k)}\cdot\mathcal O_{W^{\sigma}_{k-1}}$. Then, we have
$W^{\sigma}_{N_{\underline s}}\cong\mathcal T^{\underline s}_{n}$.

Define a rational map $\mathcal K_{n}^{{\underline s}}:G(2,n)\dashrightarrow\prod\nolimits_{\underline w\in C^{\underline s}}\mathbb {P}^{N^{\underline s}_{\underline w}}$ by
\begin{equation}\label{ks}
\mathcal K_{n}^{{\underline s}}:=\left(\cdots,F^{\underline s}_{\underline w}\circ e,\cdots\right)_{\underline w\in C^{\underline s}},  
\end{equation}
where rational maps $F^{\underline s}_{\underline w}:\mathbb P^{N_{2,n}}\dashrightarrow \mathbb P^{N^{\underline s}_{\underline w}}$ are given by
\begin{equation}\label{Fw}
F^{\underline s}_{\underline w}\left([\cdots ,z_I,\cdots]_{I\in\mathbb I_{2,n}}\right):=\left[\cdots, \prod\nolimits_{i=1}^{{\rm ht}_{\underline s}\left(\underline w\right)}z_{I_i},\cdots\right]_{\prod\nolimits_{i=1}^{{\rm ht}_{\underline s}\left(\underline w\right)}z_{I_i}\in G^{\underline s}_{\underline w}}.
\end{equation}
Define a rational map $\widehat{\mathcal K}_{n}^{\underline s}:G(2,n)\dashrightarrow \mathbb {P}^{N_{2,n}}\times\prod\nolimits_{t=1}^N\mathbb {P}^{N_{t}^{\underline s}}\times\prod\nolimits_{\underline w\in C^{\underline s}}\mathbb {P}^{N_{\underline w}^{\underline s}}$  by 
\begin{equation}\label{whk}
\widehat{\mathcal K}_{n}^{\underline s}:=\left(\widetilde{\mathcal K}^{\underline s},\,\mathcal K_{n}^{\underline s}\right)=\left(e,\left(\cdots,F^{\underline s}_{t}\circ e,\cdots\right)_{1\leq t\leq N},\left(\cdots,F^{\underline s}_{\underline w}\circ e,\cdots\right)_{\underline w\in C^{\underline s}}\right).  
\end{equation}
We can prove by induction that  $W^{\sigma}_{k}$, $1\leq k\leq N_{\underline s}$, is a closed subscheme of \begin{equation*}
\mathbb P^{N_{2,n}}\times\prod\nolimits_{t=1}^N\mathbb {P}^{N^{\underline s}_{t}}\times\prod\nolimits_{i=1}^k\mathbb {P}^{N_{\sigma(i)}^{\underline s}}.
\end{equation*} 
Moreover,  as a topological space, $\mathcal T_{n}^{\underline s}$ is the closure of the image $\widehat{\mathcal K}^{{\underline s}}_n\left(U\right)$ in $\mathbb P^{N_{2,n}}\times\prod\nolimits_{t=1}^N\mathbb {P}^{N^{\underline s}_{t}}\times\prod\nolimits_{\underline w\in C^{\underline s}}\mathbb {P}^{N_{\underline w}^{\underline s}}$, where $\widehat{\mathcal K}^{{\underline s}}_n$ is defined by (\ref{whk}) and $U$ is the domain of  definition of $\widehat{\mathcal K}^{{\underline s}}_n$.

Now, we have by Lemma \ref{proper} that
\begin{definition}\label{msss}
Denote by $\mathcal M_{n}^{\underline s}$ the image of 
$\mathcal T^{\underline s}_n$  under the natural projection from
$\mathbb {P}^{N_{2,n}}\times\prod\nolimits_{t=1}^N\mathbb {P}^{N^{\underline s}_{t}}\times\prod\nolimits_{\underline w\in C^{\underline s}}\mathbb {P}^{N^{\underline s}_{\underline w}}$ to $\prod\nolimits_{\underline w\in C^{\underline s}}\mathbb {P}^{N^{\underline s}_{\underline w}}$, which is equipped with a reduced, closed subscheme structure of $\prod\nolimits_{\underline w\in C^{\underline s}}\mathbb {P}^{N^{\underline s}_{\underline w}}$. Denote by $\mathcal P_{n}^{\underline s}:\mathcal T_{n}^{\underline s}\rightarrow \mathcal M_{n}^{\underline s}$ the induced morphism.
\end{definition}

\begin{remark}\label{last}
It is clear that as a topological space $\mathcal M_{n}^{\underline s}$ is the closure of the image ${\mathcal K}^{{\underline s}}_n\left(U\right)$ in $\prod\nolimits_{\underline w\in C^{\underline s}}\mathbb {P}^{N_{\underline w}^{\underline s}}$, where ${\mathcal K}^{{\underline s}}_n$ is defined by (\ref{ks}) and $U$ is the domain of  definition of ${\mathcal K}^{{\underline s}}_n$.
\end{remark}

\subsection{Geometry of Grassmannians} 
We first recall an open cover of $G(2,n)$ following \cite{EH}.

For $1\leq j_1<j_2\leq N$, define a closed subscheme of \begin{equation*}\mathbb A^{2(n-2)}:={\rm Spec}\left(\mathbb Z\left[x_{11},x_{12},\cdots,x_{1n},x_{21},x_{22},\cdots,x_{2n}\right]\right)
\end{equation*} by
\begin{equation}\label{u12}
U_{j_1j_2}:={\rm Spec}\left(\mathbb Z\left[\cdots,x_{ij},\cdots\right]_{\substack{i=1,2,\\1\leq j\leq n}}\Big/(x_{1j_1}-1,x_{1j_2},x_{2j_1},x_{2j_2}-1)\right).
\end{equation} 
Arranged as a matrix, the corresponding embedding  $e_{j_1j_2}:U_{j_1j_2}\hookrightarrow \mathbb A^{2(n-2)}$ takes the form 
\begin{equation}\label{ej1j2}\small(\cdots,x_{ij},\cdots)\mapsto\left(\begin{matrix}x_{11}&\cdots& x_{1(j_1-1)}\\x_{21}&\cdots& x_{2(j_1-1)}\end{matrix}  \hspace{-.1in}\begin{matrix}  &\hfill\tikzmark{a}\\  &\hfill\tikzmark{b}  \end{matrix} \,\,\,\,  \begin{matrix}  1\,\\0\,\\\end{matrix}\hspace{-.11in}\begin{matrix}  &\hfill\tikzmark{c}\\  &\hfill\tikzmark{d}  \end{matrix}\hspace{-.11in}\begin{matrix}\,\,\,\,\,\,\,\,x_{1(j_1+1)}&\cdots& x_{1(j_2-1)}\\\,\,\,\,\,\,\,\,x_{2(j_1+1)}&\cdots&x_{2(j_2-1)}\end{matrix}\hspace{-.11in}\begin{matrix}  &\hfill\tikzmark{g}\\  &\hfill\tikzmark{h}  \end{matrix}\,\,\,\,\begin{matrix}0\,\\1\,\\\end{matrix}\hspace{-.11in}\begin{matrix}  &\hfill\tikzmark{e}\\  &\hfill\tikzmark{f}\end{matrix}\hspace{-.1in}\begin{matrix}\,\,\,\,\,\,\,\,x_{1(j_2+1)}&\cdots&x_{1n}\\\,\,\,\,\,\,\,\,x_{2(j_2+1)}&\cdots& x_{2n}\end{matrix}\right).  \tikz[remember picture,overlay]   \draw[dashed,dash pattern={on 4pt off 2pt}] ([xshift=0.5\tabcolsep,yshift=7pt]a.north) -- ([xshift=0.5\tabcolsep,yshift=-2pt]b.south);\tikz[remember picture,overlay]   \draw[dashed,dash pattern={on 4pt off 2pt}] ([xshift=0.5\tabcolsep,yshift=7pt]c.north) -- ([xshift=0.5\tabcolsep,yshift=-2pt]d.south);\tikz[remember picture,overlay]   \draw[dashed,dash pattern={on 4pt off 2pt}] ([xshift=0.5\tabcolsep,yshift=7pt]e.north) -- ([xshift=0.5\tabcolsep,yshift=-2pt]f.south);\tikz[remember picture,overlay]   \draw[dashed,dash pattern={on 4pt off 2pt}] ([xshift=0.5\tabcolsep,yshift=7pt]g.north) -- ([xshift=0.5\tabcolsep,yshift=-2pt]h.south);
\end{equation}

For each $(i_1,i_2)\in\mathbb I_{2,n}$, define a regular function $P_{(i_1,i_2)}$  on $\mathbb A^{2(n-2)}$ by 
\begin{equation*}
P_{(i_1,i_2)}:=x_{1i_1}\cdot x_{2i_2}-x_{1i_2}\cdot x_{2i_1}.
\end{equation*}
It is clear that $P_{(i_1,i_2)}\circ e_{j_1j_2}$ defines an embedding from $U_{j_1j_2}$ to $G(2,n)\subset\mathbb P^{N_{2,n}}$ by
\begin{equation}\label{pplu}
(\cdots,x_{ij},\cdots)\mapsto\left[\cdots,P_{(i_1,i_2)}\left(e_{j_1j_2}\left(\cdots,x_{ij},\cdots\right)\right),\cdots\right]_{(i_1,i_2)\in\mathbb I_{2,n}},
\end{equation}
which is indeed the Pl\"ucker embedding. By identifying $U_{j_1j_2}$ with the images under (\ref{pplu}), we derive an open cover $\{U_{j_1j_2}\}$ of $G(2,n)$.  


\begin{definition}\label{pluc}
We call $P_{(i_1,i_2)}\circ e_{j_1j_2}$ Pl\"ucker coordinate functions on  $G(2,n)$. 
\end{definition}

We introduce a systematic way to construct blow-ups as follows. Define an open subscheme  $\mathcal {G}(2,n)$ of 
$\mathbb A^{2(n-2)}$
by
\begin{equation}\label{gltor}
\big\{\left.\mathfrak p\in\mathbb A^{2(n-2)}\right|\,x_{1j_1}x_{2j_2}-x_{2j_1}x_{1j_2}\notin\mathfrak p,\,{\rm for\,\,certain\,\,}1\leq j_1<j_2\leq n\big\}.
\end{equation} 
It is clear that
$\mathcal {G}(2,n)$ is a ${\rm GL}_2$-torsor over $G(2,n)$. We denote the corresponding projection by \begin{equation}\label{pi}
\pi:\mathcal {G}(2,n)\rightarrow G(2,n).    
\end{equation}
Let $\underline s=(s_1,\cdots,s_N)$ be a size vector and denote $n=\sum_{t=1}^Ns_t$. For any morphism $L:U\rightarrow\mathcal {G}(2,n)$ such that
for each ${\underline w}\in C^{\underline s}$ the ideal sheaf $\left(\pi\circ L\right)^{-1}\mathscr I^{\underline s}_{\underline w}\cdot\mathcal O_U$ is non-zero, we have 
\begin{definition}\label{mat}
Denote by $\mathcal R^{\underline s}_L:\mathbb M^{\underline s}_L\rightarrow U$ the blow-up of $U$ with respect to the product of inverse image ideal sheaves \begin{equation*}
\prod\nolimits_{\underline w\in C^{\underline s}}\left(\pi\circ L\right)^{-1}\mathscr I^{\underline s}_{\underline w}\cdot\mathcal O_U.    
\end{equation*} 
\end{definition}


To end this section, we give a geometric interpretation of the ideal sheaves $\mathscr I^{\underline s}_t$. Note that the blow-up with respect to $\mathscr I^{\underline s}_{\underline w}$ produces singularities corresponding to fixed points of generic $T_{\underline s}$-orbits on $G(2,n)$. The intuition of the blow-up with respect to $\mathscr I^{\underline s}_t$ is to solve such singularities.

For $1\leq t\leq N$, denote by $S_t$ the subscheme of $G(2,n)$ defined by the ideal sheaf $\mathscr I^{\underline s}_t$.

\begin{lemma}\label{subg}
$S_t$ is isomorphic to a sub-Grassmannian $G(2,n-s_t)$ of $G(2,n)$ for $1\leq t\leq N$. For each nonempty subset $\mathcal P\subset\{1,2,\cdots,N\}$, the  intersection $\bigcap\nolimits_{t\in\mathcal P}S_t$ is smooth. 
\end{lemma}
{\bf\noindent Proof of Lemma \ref{subg}.} Fix $1\leq t\leq N$. Let $U_{j_1j_2}$ be defined by (\ref{u12}). If $j_1$ or $j_2$ is in the sub-block $\Delta_t(\underline s)=\left[\sum_{i=1}^{t-1}s_{i}+1, \sum_{i=1}^{t}s_{i}\right]$, then $(j_1,j_2)\in\mathbb I^{\underline s}_{t}$ and hence 
$S_t\cap U_{j_1j_2}=\emptyset$.   

Assume that $j_1,j_2\notin\Delta_t(\underline s)$. It is clear that $S_t\cap U_{j_1j_2}$ is defined by the ideal generated by \begin{equation*}
\left\{P_{I_{i\alpha}}\circ e_{j_1j_2}\left| \,i=1,2,\,\,{\rm and}\,\,\alpha\in\Delta_t(\underline s)\right.\right\},
\end{equation*}
where $I_{i\alpha}=\left(j_i,\alpha\right)$ if $j_i<\alpha$ and $I_{i\alpha}=\left(\alpha,j_i\right)$ if $\alpha<j_i$.
Noticing that $P_{I_{i\alpha}}\circ e_{j_1j_2}=x_{i\alpha}$, we can conclude that $S_t$ is the sub-Grassmannian of $G(2,n)$ parametrizing rank $(n-2-s_t)$ free quotients of $\bigoplus\nolimits_{\substack{1\leq\alpha\leq N,\,\,{\rm and}\,\,\alpha\neq t}}E_{\alpha}$.

The proof of Lemma \ref{subg} is complete.\,\,\,$\endpf$


\section{The Case of Maximal Torus Actions}\label{mta}

In this section, we are concerned with the following proposition.

\begin{proposition}\label{ms}
Let $N=n\geq 3$ and $\underline s=(1,1,\cdots,1)\in\mathbb Z^N$. Then, $\mathcal T_{n}^{\underline s}$ and $\mathcal M_{n}^{\underline s}$ are smooth over ${\rm Spec}\,\mathbb Z$ of relative dimensions $2n-4$ and $n-3$, respectively.
Moreover, the morphism $\mathcal P^{\underline s}_{n}:\mathcal T_{n}^{\underline s}\rightarrow\mathcal M_{n}^{\underline s}$ is flat.
\end{proposition}

We adopt the convention that $\underline s=(1,1,\cdots,1)\in\mathbb Z^N$ and $N=n$ throughout this section.  Note that our approach does not rely on the  geometry of stable curves, and hence it can be adapted to the case of sub-torus actions.

\subsection{Parametrizations for \texorpdfstring{$\mathcal Q_{n}^{\underline s}$}{dd}}\label{blf}

We first introduce the following notation.
\begin{definition}
Let $\mathbb J^{\underline s}$ be the set of indices $\left(j_1,j_2,(j^+_1,\cdots,j^+_l),(j^-_1,\cdots,j^-_{m})\right)$ such that:
\begin{enumerate}[label=(\arabic*)]
\item  $l$, $m$ are non-negative integers, and $l+m=N-2$;

\item $\{j_1,j_2,j^+_1,\cdots,j^+_l,j^-_1,\cdots,j^-_{m}\}=\{1,2,\cdots,N\}$;

\item $j_1<j_2$, $j^+_1<j^+_2<\cdots<j^+_l$, and $j^-_1<j^-_2<\cdots<j^-_m$.
\end{enumerate}
\end{definition}

We now associate to each $\tau=\left(j_1,j_2,(j^+_1,\cdots,j^+_l),(j^-_1,\cdots,j^-_{m})\right)\in\mathbb J^{\underline s}$  an affine space \begin{equation*}
{\rm Spec}\,\mathbb Z\left[\overrightarrow  A,\overrightarrow  H,\overrightarrow  \Xi\right], 
\end{equation*}
where 
\begin{equation}\label{bp}
\begin{split}
&\overrightarrow  A:=\left(a_{j^+_1},a_{j^+_2},\cdots,a_{j^+_l},a_{j^-_1},a_{j^-_2},\cdots,a_{j^-_m}\right),\\
&\overrightarrow  H:=\left(\eta_{j^+_1},\eta_{j^+_2},\cdots,\eta_{j^+_l}\right),\,\,\,\overrightarrow \Xi:=\left(\xi_{j^-_1},\xi_{j^-_2},\cdots,\xi_{j^-_m}\right).\\
\end{split}
\end{equation}
Define a morphism $\Gamma^{\tau}:{\rm Spec}\,\mathbb Z\left[\overrightarrow  A,\overrightarrow  H,\overrightarrow  \Xi\right]\rightarrow U_{j_1j_2}$  by the following homomorphism between the coordinate rings
\begin{equation}\label{ngamma}
\left\{\begin{aligned}
&\,x_{1j^+_{\alpha}}\mapsto a_{j^+_{\alpha}},\,\,\,\,\,\,\,\,\,\,\,\,\,\,\,\,\,\,\,\,\,\,x_{2j^+_{\alpha}}\mapsto a_{j^+_{\alpha}}\cdot\eta_{j^+_{\alpha}},\,\,\,\,\,\,\,\,\,\,{\rm for}\,\,1\leq\alpha\leq l\\
&\,x_{1j^-_{\beta}}\mapsto a_{j^-_{\beta}}\cdot\xi_{j^-_{\beta}},\,\,\,\,\,\,\,\,\,\,x_{2j^-_{\beta}}\mapsto a_{j^-_{\beta}},\,\,\,\,\,\,\,\,\,\,\,\,\,\,\,\,\,\,\,\,\,\,{\rm for}\,\,1\leq\beta\leq m\\
\end{aligned}\right.,
\end{equation}
where $U_{j_1j_2}$ is the open subscheme of $G(2,n)$ defined by (\ref{u12}).

Define a rational map  $J^{\tau}:{\rm Spec}\,\mathbb Z\left[\overrightarrow  A,\overrightarrow  H,\overrightarrow  \Xi\right]\dashrightarrow\mathbb P^{N_{2,n}}\times\prod\nolimits_{t=1}^N\mathbb P^{N^{\underline s}_{t}}$ by $J^{\tau}:=\widetilde{\mathcal K}^{\underline s}\circ \Gamma^{\tau}$.  
\begin{lemma}\label{em}
The rational map $J^{\tau}$ defined above extends to  a locally closed embedding from ${\rm Spec}\,\mathbb Z\left[\overrightarrow  A,\overrightarrow  H,\overrightarrow  \Xi\right]$ to $\mathbb P^{N_{2,n}}\times\prod\nolimits_{t=1}^N\mathbb P^{N^{\underline s}_{t}}$.
\end{lemma} 
{\noindent\bf Proof of Lemma \ref{em}.} Without loss of generality, we can assume that 
\begin{equation*}
\tau=\left(j_1,j_2,(j^+_1,j^+_2,\cdots,j^+_l),(j^-_1,j^-_2,\cdots,j^-_{m})\right)= \left(1,2,(3,4,\cdots,l+2),(l+3,l+4,\cdots,N)\right).
\end{equation*}

We first show that the rational maps $F^{\underline s}_{t}\circ e\circ \Gamma^{\tau}$, $1\leq t\leq N$, extends to  a morphism on ${\rm Spec}\,\mathbb Z\left[\overrightarrow  A,\overrightarrow  H,\overrightarrow  \Xi\right]$. 
When $t=1$,
\begin{equation*}
\begin{split}
\left(F^{\underline s}_{1}\circ e\circ  \Gamma^{\tau}\right)\left(\overrightarrow A,\overrightarrow H,\overrightarrow \Xi\right)&=\bigg[\cdots,\left(P_{(1,i)}\circ e_{j_1j_2}\circ\Gamma^{\tau}\right)\left(\overrightarrow A,\overrightarrow H,\overrightarrow \Xi\right),\cdots\bigg]_{2\leq i\leq n}\\
=&\left[1,\,a_{3}\eta_{3},\,a_{4}\eta_{4},\,\cdots,\,a_{l+2}\eta_{l+2},\, a_{l+3},\, a_{l+4},\,\cdots,\,a_{N}\right],\\
\end{split}
\end{equation*}
where  $P_{(1,i)}\circ e_{j_1j_2}$ are the Pl\"ucker coordinate functions on $G(2,n)$. Similarly, when $t=2$,
\begin{equation*}
\begin{split}
\left(F^{\underline s}_{2}\circ e\circ \Gamma^{\tau}\right)&\left(\overrightarrow A,\overrightarrow H,\overrightarrow \Xi\right)=\left[1,-a_{3},-a_{4},\cdots,-a_{l+2},-a_{l+3}\xi_{l+3},-a_{l+4}\xi_{l+4},\cdots,-a_{N}\xi_{N}\right].\\
\end{split}
\end{equation*}

When $m=0$,  the rational map $F^{\underline s}_{t}\circ e\circ \Gamma^{\tau}$, $3\leq t\leq N$, takes the form 
\begin{equation*}
\begin{split}
&\left(F^{\underline s}_{t}\circ e\circ \Gamma^{\tau}\right)\left(\overrightarrow A,\overrightarrow H,\overrightarrow \Xi\right)=\left[a_{t}\eta_{t},\,-a_{t},\,a_3a_t(\eta_{t}-\eta_3),\,a_4a_t(\eta_{t}-\eta_4),\cdots,\right.\\
&\,\,\,\,\,\,\,\,\,\,\,\,\,\left.a_{t-1}a_t(\eta_{t}-\eta_{t-1}),\,a_ta_{t+1}(\eta_{t+1}-\eta_t),\,a_ta_{t+2}(\eta_{t+2}-\eta_t),\,\cdots,\,a_ta_{N}(\eta_{N}-\eta_t)\right].\\
\end{split}
\end{equation*}
Then, we can define the extension by
\begin{equation*}
\begin{split}
&\left(F^{\underline s}_{t}\circ e\circ \Gamma^{\tau}\right)\left(\overrightarrow A,\overrightarrow H,\overrightarrow \Xi\right):=\left[-\eta_{t},\,1,\,-a_3(\eta_{t}-\eta_3),\,-a_4(\eta_{t}-\eta_4),\cdots,\right.\\
&\,\,\,\,\,\,\,\,\,\,\,\,\left.-a_{t-1}(\eta_{t}-\eta_{t-1}),-a_{t+1}(\eta_{t+1}-\eta_t),\,-a_{t+2}(\eta_{t+2}-\eta_t),\,\cdots,\,-a_{N}(\eta_{N}-\eta_t)\right].\\
\end{split}
\end{equation*}
When $l=0$, the rational map $F^{\underline s}_{t}\circ e\circ \Gamma^{\tau}$, $3\leq t\leq N$, takes the form 
\begin{equation*}
\begin{split}
&\left(F^{\underline s}_{t}\circ e\circ \Gamma^{\tau}\right)\left(\overrightarrow A,\overrightarrow H,\overrightarrow \Xi\right)=\left[a_{t},\,-a_{t}\xi_t,\,a_3a_t(\xi_{3}-\xi_t),\,a_4a_t(\xi_{4}-\xi_t),\,\cdots,\right.\\
&\,\,\,\,\,\,\,\,\,\,\,\left.a_{t-1}a_t(\xi_{t-1}-\xi_{t}),\,a_ta_{t+1}(\xi_{t}-\xi_{t+1}),\,a_ta_{t+2}(\xi_{t}-\xi_{t+2}),\,\cdots,\,a_ta_{N}(\xi_{t}-\xi_N)\right].\\
\end{split}
\end{equation*}
We can define the extension by
\begin{equation*}
\begin{split}
&\left(F^{\underline s}_{t}\circ e\circ \Gamma^{\tau}\right)\left(\overrightarrow A,\overrightarrow H,\overrightarrow \Xi\right)=\left[1,\,-\xi_{t},\,a_3(\xi_{3}-\xi_t),\,a_4(\xi_{4}-\xi_t),\right.\\
&\left.\,\cdots,\,a_{t-1}(\xi_{t-1}-\xi_{t-1}),\,a_{t+1}(\xi_{t}-\xi_{t+1}),a_{t+2}(\xi_{t}-\xi_{t+2}),\,\cdots,\,a_{N}(\xi_{t}-\xi_N)\right].\\
\end{split}
\end{equation*}

When $l,m\geq1$, the rational map $F^{\underline s}_{t}\circ e\circ \Gamma^{\tau}$, $3\leq t\leq l+2$, takes the form  
\begin{equation*}
\begin{split}
&\left(F^{\underline s}_{t}\circ e\circ \Gamma^{\tau}\right)\left(\overrightarrow A,\overrightarrow H,\overrightarrow \Xi\right)=\left[a_{t}\eta_{t},\,-a_{t},\,a_3a_t(\eta_{t}-\eta_3),\,\cdots,\,a_{t-1}a_t(\eta_{t}-\eta_{t-1}),\right.\\
&\,\,\,\,\,\,\,\,\,\,\,\,\left.a_ta_{t+1}(\eta_{t+1}-\eta_t),\,\cdots,\,a_ta_{l+2}(\eta_{l+2}-\eta_t),\,a_ta_{l+3}(1-\eta_t\xi_{l+3}),\,\cdots,\,a_ta_{N}(1-\eta_t\xi_{N})\right],\\
\end{split}
\end{equation*}
and the extension is given by
\begin{equation*}
\begin{split}
&\left(F^{\underline s}_{t}\circ e\circ \Gamma^{\tau}\right)\left(\overrightarrow A,\overrightarrow H,\overrightarrow \Xi\right)=\left[-\eta_{t},\,1,\,-a_3(\eta_{t}-\eta_3),\,\cdots,\,-a_{t-1}(\eta_{t}-\eta_{t-1}),\right.\\
&\,\,\,\,\,\,\,\,\,\,\,\,\left.-a_{t+1}(\eta_{t+1}-\eta_t),\,\cdots,-a_{l+2}(\eta_{l+2}-\eta_t),-a_{l+3}(1-\eta_t\xi_{l+3}),\,\cdots,\,-a_{N}(1-\eta_t\xi_{N})\right].\\
\end{split}
\end{equation*}
When $l,m\geq1$, the rational map $F^{\underline s}_{t}\circ e\circ \Gamma^{\tau}$, $l+3\leq t\leq N$, takes the form  
\begin{equation*}
\begin{split}
&\left(F^{\underline s}_{t}\circ e\circ \Gamma^{\tau}\right)\left(\overrightarrow A,\overrightarrow H,\overrightarrow \Xi\right)=\left[a_{t},\,-a_{t}\xi_t,\,a_3a_t(1-\eta_{3}\xi_t),\,\cdots,\,a_{l+2}a_t(1-\eta_{l+2}\xi_{t}),\right.\\
&\,\,\,\,\,\,\,\,\,\,\,\,\,\,\,\left.a_{l+3}a_t(\xi_{l+3}-\xi_t),\,\cdots,a_{t-1}a_t(\xi_{t-1}-\xi_t),\,a_ta_{t+1}(\xi_{t}-\xi_{t+1}),\,\cdots,\,a_ta_{N}(\xi_{t}-\xi_N)\right],\\
\end{split}
\end{equation*}
and the extension is given by
\begin{equation*}
\begin{split}
&\left(F^{\underline s}_{t}\circ e\circ \Gamma^{\tau}\right)\left(\overrightarrow A,\overrightarrow H,\overrightarrow \Xi\right)=\left[1,\,-\xi_{t},\,a_3(1-\eta_{3}\xi_t),\,\cdots,\,a_{l+2}(1-\eta_{l+2}\xi_t),\right.\\
&\,\,\,\,\,\,\,\,\,\,\,\,\,\,\,\left.a_{l+3}(\xi_{l+3}-\xi_t),\,\cdots,\,a_{t-1}(\xi_{t-1}-\xi_{t}),a_{t+1}(\xi_{t}-\xi_{t+1}),\,\cdots,\,a_{N}(\xi_{t}-\xi_N)\right].\\
\end{split}
\end{equation*}

For brevity, we denote the extension by $J^{\tau}$ as well. To prove that the morphism $J^{\tau}$
is a locally closed embedding, it suffices to  show that the corresponding ring homomorphism is surjective. Notice that for each coordinate $x$ in $\overrightarrow  A$, $\overrightarrow  H$,  $\overrightarrow  \Xi$, there is an inhomogeneous coordinate of a certain projective space maps to $\pm x$ by the above formulas. The proof of Lemma \ref{em} is complete.
\,\,\,\,$\endpf$
\medskip

Denote by $A^{\tau}$ the image of ${\rm Spec}\,\mathbb Z\left[\overrightarrow  A,\overrightarrow  H,\overrightarrow  \Xi\right]$ under  $J^{\tau}$, $\tau\in\mathbb J^{\underline s}$, which is equipped with the reduced scheme structure as a locally closed subscheme of $\mathbb P^{N_{2,n}}\times\prod\nolimits_{t=1}^N\mathbb {P}^{N^{\underline s}_{t}}$. We then have 
\begin{lemma}\label{coor2} $\mathcal Q_{n}^{\underline s}$ is a union of $A^{\tau}$ for all $\tau\in\mathbb J^{\underline s}$.
\end{lemma}

{\bf\noindent Proof of Lemma \ref{coor2}.} It suffices to prove that $\bigcup_{\tau\in\mathbb J^{\underline s}}A^{\tau}\supset\left( \widetilde{\mathcal R}^{\underline s}\right)^{-1}\left(U_{12}\right)$, where $U_{12}$ is defined by (\ref{u12}). 
Motivated by Lemma \ref{subg}, we shall define a closed embedding \begin{equation*}
E:U_{12}\times\prod\nolimits_{i=3}^{N}\mathbb P^1\hookrightarrow \left(\mathbb P^{N_{2,n}}\mathbin{\scaleobj{1.5}{\backslash}}\{z_{(1,2)}\neq 0\}\right)\times\prod\nolimits_{t=1}^{N}\mathbb P^{N^{\underline s}_{t}}.    
\end{equation*}

Denote by $[u_t,v_t]$ the homogeneous coordinates for the $t^{\rm th}$ copy of $\mathbb P^1$ in $\prod\nolimits_{t=3}^{N}\mathbb P^1$, $3\leq t\leq N$. Denote by $\left[z^t_{(1,t)},z^t_{(2,t)},\cdots, z^t_{(t-1,t)},z^t_{(t,t+1)},z^t_{(t,t+2)},\cdots,z^t_{(t,n)}\right]$ the homogeneous coordinates for $\mathbb P^{N^{\underline s}_{t}}$, $1\leq t\leq N$. Define
\begin{equation*}
E:=\left(e\circ{\rm Pr},f_1,f_2,\cdots,f_N\right),
\end{equation*}
where ${\rm Pr}:U_{12}\times\prod\nolimits_{t=3}^{N}\mathbb P^1\rightarrow U_{12}$ is the projection, and $f_t:U_{12}\times\prod\nolimits_{t=3}^{N}\mathbb P^1\rightarrow\mathbb P^{N^{\underline s}_{t}}$ are defined as follows. $f_1$ is defined by
\begin{equation*}
\left[z^1_{(1,2)},z^1_{(1,3)},z^1_{(1,4)},\cdots, z^1_{(1,n)}\right]\mapsto\left[1,x_{23},x_{24},\cdots, x_{2n}\right],    
\end{equation*}
and $f_2$ is defined by
\begin{equation*}
\left[z^2_{(1,2)},z^2_{(2,3)},z^2_{(2,4)},\cdots, z^1_{(2,n)}\right]\mapsto\left[1,-x_{13},-x_{14},\cdots, -x_{1n}\right];   
\end{equation*}
for $3\leq t\leq N$, $f_t$ is defined by
\begin{equation*}
\begin{split}
&\,\,\,\,\,\,\,\,\,\,\,\,\,\,\,\,\,\,\,\,\,\,\,\,\,\,\,\,\,\,\left[z^t_{(1,t)},z^t_{(2,t)},z^t_{(3,t)},z^t_{(4,t)},\cdots, z^t_{(t-1,t)},z^t_{(t,t+1)},z^t_{(t,t+2)},\cdots, z^t_{(t,N)}\right]\mapsto\\
&\left[v_t,-u_t,\left|\begin{matrix}
 x_{13}&u_t\\
 x_{23}&v_t 
\end{matrix}\right|,\,\left|\begin{matrix}
 x_{14}&u_t\\
 x_{24}&v_t 
\end{matrix}\right|,\,\cdots,\,\left|\begin{matrix}
 x_{1(t-1)}&u_t\\
 x_{2(t-1)}&v_t 
\end{matrix}\right|,\,\left|\begin{matrix}
u_t&x_{1(t+1)}\\
v_t&x_{2(t+1)}\\
\end{matrix}\right|,\,\left|\begin{matrix}
u_t&x_{1(t+2)}\\
v_t&x_{2(t+2)}\\
\end{matrix}\right|,\,\cdots,\,\left|\begin{matrix}
u_t&x_{1N}\\
v_t&x_{2N}\\
\end{matrix}\right|\right].\\
\end{split}  
\end{equation*}

Now, it is clear that the restriction of the rational map $\widetilde{\mathcal K}^{\underline s}|_{U_{12}}$ equals to $E\circ\check{\mathcal K}$, where the rational map $\check{\mathcal K}:U_{12}\dashrightarrow U_{12}\times\prod\nolimits_{t=3}^{N}\mathbb P^1$ is defined by
\begin{equation*}
\left\{\begin{aligned}
&\,x_{it}\mapsto x_{it},\,\,\,\,\,\,\,\,\,\,\,\,\,\,\,\,\,\,\,\,\,\,\,\,\,\,\,\,\,\,\,\,\,\,\,\,\,{\rm for}\,\,1\leq i\leq 2\,\,{\rm and}\,\,3\leq t\leq n\\
&\,[u_t,v_t]\mapsto [x_{1t},x_{2t}],\,\,\,\,\,\,\,\,\,\,\,\,\,\,{\rm for}\,\,3\leq t\leq N\\
\end{aligned}\right.\,.
\end{equation*}
Hence, by the same argument as that in \S \ref{pdm}, we can identify  $\left(\widetilde{\mathcal R}^{\underline s}\right)^{-1}\left(U_{12}\right)$
with a reduced, closed subscheme of  $U_{12}\times\prod\nolimits_{t=3}^{N}\mathbb P^1$, which as a topological space is the closure of the image of $U_{12}$ under the rational map $\check{\mathcal K}$.

Notice that for each $\tau\in\mathbb J^{\underline s}$ such that $j_1=1$ and $j_2=2$,  $\check {\mathcal K}\circ\Gamma^{\tau}$ extends to a locally closed embedding from ${\rm Spec}\,\mathbb Z\left[\overrightarrow  A,\overrightarrow  H,\overrightarrow  \Xi\right]$ to $U_{12}\times\prod\nolimits_{t=3}^{N}\mathbb P^1$. Denote by $\check A^{\tau}$ the image of ${\rm Spec}\,\mathbb Z\left[\overrightarrow  A,\overrightarrow  H,\overrightarrow  \Xi\right]$ under this morphism.
Then, to prove Lemma \ref{coor2}, it suffices to show that $\left(\widetilde{\mathcal R}^{\underline s}\right)^{-1}\left(U_{12}\right)$ is covered by $\check A^{\tau}$ where $\tau\in\mathbb J^{\underline s}$ run over all the indices such that $j_1=1$ and $j_2=2$. 

Recall that $U_{12}\times\prod\nolimits_{i=3}^{N}\mathbb P^1$ is covered by all affine spaces
\begin{equation*}
\check U^{\tau}:=U_{12}\times\prod\nolimits_{\alpha=1}^l\left\{\left.[u_{j^+_{\alpha}},v_{j^+_{\alpha}}]\right|u_{j^+_{\alpha}}\neq 0\right\}\times\prod\nolimits_{\beta=1}^m\left\{\left.[u_{j^-_{\beta}},v_{j^-_{\beta}}]\right|v_{j^-_{\beta}}\neq 0\right\},   
\end{equation*}
where $\tau=\left(1,2,(j^+_1,\cdots,j^+_l),(j^-_1,\cdots,j^-_{m})\right)\in\mathbb J^{\underline s}$. Moreover, we can show that  $\check A^{\tau}$ is a closed subscheme of $\check U^{\tau}$,  hence $\left(\widetilde{\mathcal R}^{\underline s}\right)^{-1}\left(U_{12}\right)\cap \check U^{\tau}=\check A^{\tau}$.

We complete the proof of Lemma \ref{coor2}.\,\,\,\,$\endpf$
\medskip

According to the above proof, it is clear that $A^{\tau}$ are open subschemes of $\mathcal Q_{n}^{\underline s}$ for all $\tau\in\mathbb J^{\underline s}$. Then, we can conclude that
\begin{corollary}\label{coor3}
$\mathcal Q_{n}^{\underline s}$ is smooth over ${\rm Spec}\,\mathbb Z$. 
\end{corollary}

In the remaining of this subsection, we will prove
\begin{lemma}\label{wn4}
Proposition \ref{ms} holds when $n=4$.
\end{lemma}
\begin{remark}
When $n=3$, $\mathcal T_{n}^{\underline s}$ is isomorphic to $\mathcal Q_{n}^{\underline s}$, which is a blow-up of $\mathbb P^2$ along $3$ points, and $\mathcal M_{n}^{\underline s}$ is isomorphic to a single point.
\end{remark}
{\noindent\bf Proof of Lemma \ref{wn4}.} When $n=4$, $\overline{\mathcal R}^{\underline s}_4:\mathcal T_{4}^{\underline s}\rightarrow\mathcal Q_{4}^{\underline s}$ is the blow-up of $\mathcal Q_{4}^{\underline s}$ with respect to $\left(\widetilde {\mathcal R}^{\underline s}\right)^{-1}\mathscr I^{\underline s}_{\underline w}\cdot\mathcal O_{\mathcal Q^{\underline s}_4}$ with $\underline w=(1,1,1,1)$.

Thanks to the Pl\"ucker relation $z_{(1,2)}\cdot z_{(3,4)}-z_{(1,3)}\cdot z_{(2,4)}+z_{(1,4)}\cdot z_{(2,3)}=0$, we can conclude that $\mathcal M_{4}^{\underline s}\subset\mathbb P^{N^{\underline s}_{\underline w}}\cong\mathbb P^2$ is contained in a certain $\mathbb P^1$. Define an open subscheme of $U_{12}$ by 
\begin{equation*}
V:={\rm Spec}\left(\mathbb Z[y]\left[\cdots,x_{ij},\cdots\right]_{\substack{i=1,2,\\1\leq j\leq 4}}\Big/(x_{11}-1,x_{12},x_{21},x_{22}-1,y(x_{13}\cdot x_{24}-x_{23}\cdot x_{14})=1)\right).
\end{equation*}
Then the rational map $\mathcal P^{\underline s}_{4}\circ\widehat{\mathcal K}^{\underline s}_4$ is well-defined as a morphism on $V$, which takes the form
\begin{equation*} \left(\cdots,x_{ij},\cdots\right)\mapsto[x_{13}\cdot x_{24}-x_{23}\cdot x_{14},\,\,-x_{23}\cdot x_{14},\,\,-x_{13}\cdot x_{24}].   
\end{equation*}
We can thus conclude  $\mathcal M^{\underline s}_{4}\cong\mathbb P^1$. 

Hence, to prove Lemma \ref{wn4}, it suffices to show that  $\left(\widehat{\mathcal R}_{4}^{\underline s}\right)^{-1}(U_{12})$ is smooth, and the restriction of  $\mathcal P^{\underline s}_{4}$ to $\left(\widehat{\mathcal R}_{4}^{\underline s}\right)^{-1}(U_{12})$ is flat. 

Recall Lemma \ref{coor2} that $\left(\widetilde{\mathcal R}^{\underline s}\right)^{-1}(U_{12})$  is covered by the following open subschemes:

\begin{enumerate}[label={(\arabic*)}]
\item $\tau=\left(1,2,(j_1^+,j_2^+)\right)=(1,2,(3,4))$.
$A^{\tau}\cong{\rm Spec}\,\mathbb Z\left[\overrightarrow  A,\overrightarrow  H\right]={\rm Spec}\,\mathbb Z\left[a_{3},a_{4},\eta_{3},\eta_{4}\right]$.

\item $\tau=\left(1,2,(j_1^+),(j_1^-)\right)=(1,2,(3,4))$. $A^{\tau}\cong{\rm Spec}\,\mathbb Z\left[\overrightarrow  A,\overrightarrow  H,\overrightarrow  \Xi\right]={\rm Spec}\,\mathbb Z\left[a_{3},a_4,\eta_3,\xi_{4}\right]$.

\item $\tau=\left(1,2,(j_1^-,j_2^-)\right)=(1,2,(3,4))$. $A^{\tau}\cong{\rm Spec}\,\mathbb Z\left[\overrightarrow  A,\overrightarrow  \Xi\right]={\rm Spec}\,\mathbb Z\left[a_{3},a_{4},\xi_{3},\xi_{4}\right]$.

\item $\tau=\left(1,2,(j_1^+),(j_1^-)\right)=(1,2,(4,3))$. $A^{\tau}\cong{\rm Spec}\,\mathbb Z\left[\overrightarrow  A,\overrightarrow  H,\overrightarrow  \Xi\right]={\rm Spec}\,\mathbb Z\left[a_{4},a_3,\eta_4,\xi_{3}\right]$.

\end{enumerate}
We next study the blow-up $\overline{\mathcal R}_{4}^{{\underline s}}:\mathcal T_{4}^{\underline s}\rightarrow\mathcal Q_{4}^{\underline s}$ on a case by case basis.
\smallskip

{\bf\noindent Case (1).}  
The restriction to $A^{\tau}$ of the ideal sheaf $\left(\widetilde {\mathcal R}^{\underline s}\right)^{-1}\mathscr I^{\underline s}_{\underline w}\cdot\mathcal O_{\mathcal Q^{\underline s}_4}$ is generated by
\begin{equation*}
\left|\begin{matrix}1&0\\0&1\end{matrix}\right|\cdot\left|\begin{matrix}a_3&a_{4}\\a_{3}\cdot\eta_3&a_{4}\cdot\eta_4\end{matrix}\right|,\,\,\,\left|\begin{matrix}1&a_{3}\\0&a_{3}\cdot\eta_3\end{matrix}\right|\cdot\left|\begin{matrix}0&a_{4}\\1&a_{4}\cdot\eta_4\end{matrix}\right|,\,\,\,\left|\begin{matrix}1&a_{4}\\0&a_{4}\cdot\eta_4\end{matrix}\right|\cdot\left|\begin{matrix}0&a_{3}\\1&a_{3}\cdot\eta_3\end{matrix}\right|.
\end{equation*}
Then,  $\left(\widetilde {\mathcal R}^{\underline s}\right)^{-1}\mathscr I^{\underline s}_{\underline w}\cdot\mathcal O_{\mathcal Q^{\underline s}_4}$ is a product of the invertible ideal sheaf generated by $a_3a_4$, and the ideal sheaf generated by $\{\eta_3,\eta_4\}$. Hence, 
locally in $A^{\tau}$ $\overline{\mathcal R}_{4}^{{\underline s}}$ is  the blow-up of ${\rm Spec}\,\mathbb Z\left[\overrightarrow  A,\overrightarrow  H\right]$ along the coordinate plane $\{\eta_3=\eta_4=0\}$, which is smooth over ${\rm Spec}\,\mathbb Z$.

Next, we consider the morphism ${\mathcal P}_{4}^{{\underline s}}$. Recall that $\left(\overline{\mathcal R}_{4}^{{\underline s}}\right)^{-1}\left(A^{\tau}\right)$ is isomorphic to the closed subscheme of ${\rm Spec}\,\mathbb Z\left[\overrightarrow  A,\overrightarrow  H\right]\times{\rm Proj}\,\mathbb Z\left[u_1,v_1\right]$ defined by
\begin{equation*}
u_1\cdot\eta_4=v_1\cdot\eta_3. 
\end{equation*}
Then $\left(\overline{\mathcal R}_{4}^{{\underline s}}\right)^{-1}\left(A^{\tau}\right)$ is covered by the following two open subschemes:
\begin{enumerate}[label={$\bullet$}]
\item $V_1\cong{\rm Spec}\,\mathbb Z\left[a_3,a_4,\eta_3,x\right]$, where the embedding $V_1\hookrightarrow\left(\overline{\mathcal R}_{4}^{{\underline s}}\right)^{-1}\left(A^{\tau}\right)$ is induced by the homomorphism between coordinate rings
\begin{equation*}
\frac{v_1}{u_1}\mapsto x\,\,\,\,{\rm and}\,\,\,\,\eta_4\mapsto\eta_3\cdot x.
\end{equation*}
\item  $V_2\cong{\rm Spec}\,\mathbb Z\left[a_3,a_4,\eta_4,y\right]$, where the embedding $V_2\hookrightarrow\left(\overline{\mathcal R}_{4}^{{\underline s}}\right)^{-1}\left(A^{\tau}\right)$ is induced by the homomorphism between coordinate rings
\begin{equation*}
\frac{u_1}{v_1}\mapsto y\,\,\,\,{\rm and}\,\,\,\,\eta_3\mapsto\eta_4\cdot y.
\end{equation*}
\end{enumerate}
Resticted to $V_1$, the morphism ${\mathcal P}_{4}^{{\underline s}}:\mathcal T_{4}^{\underline s}\rightarrow\mathcal M_{4}^{\underline s}\subset\mathbb P^{N^{\underline s}_{\underline w}}$ takes the form
\begin{equation*}\begin{split}
{\mathcal P}_{4}^{{\underline s}}\left(a_{3},a_{4},\eta_3,x\right)&=[a_3 a_4\eta_3 x-a_3\eta_3a_4,-a_3\eta_3 a_4,-a_3 a_4\eta_3x] \\
&=[x-1,-1,x].
\end{split}
\end{equation*} Restricted to $V_2$,  ${\mathcal P}_{4}^{{\underline s}}$ takes the form 
\begin{equation*}\begin{split}
{\mathcal P}_{4}^{{\underline s}}\left(a_{3},a_{4},\eta_4,y\right)&=[a_3a_4\eta_4-a_3\eta_4y a_4,-a_3\eta_4 y a_4,-a_3 a_4\eta_4] \\
&=[1-y,-y,-1].
\end{split}
\end{equation*}  Hence  ${\mathcal P}_{4}^{{\underline s}}$ is flat in $\left(\overline{\mathcal R}_{4}^{{\underline s}}\right)^{-1}\left(A^{\tau}\right)$. 
\smallskip

{\bf\noindent Case (2).} Computation yields that  the restriction to $A^{\tau}$ of $\left(\widetilde {\mathcal R}^{\underline s}\right)^{-1}\mathscr I^{\underline s}_{\underline w}\cdot\mathcal O_{\mathcal Q^{\underline s}_4}$  is generated by
 \begin{equation*}
\left|\begin{matrix}1&0\\0&1\end{matrix}\right|\cdot\left|\begin{matrix}a_3&a_{4}\cdot\xi_4\\a_{3}\cdot\eta_3&a_{4}\end{matrix}\right|,\,\,\,\left|\begin{matrix}1&a_{3}\\0&a_{3}\cdot\eta_3\end{matrix}\right|\cdot\left|\begin{matrix}0&a_{4}\cdot\xi_4\\1&a_4\end{matrix}\right|,\,\,\,\left|\begin{matrix}1&a_{4}\cdot\xi_4\\0&a_{4}\end{matrix}\right|\cdot\left|\begin{matrix}0&a_{3}\\1&a_{3}\cdot\eta_3\end{matrix}\right|,
\end{equation*}
which is the invertible ideal sheaf generated by $a_3a_4$. Then, locally in $A^{\tau}$, $\overline{\mathcal R}_{4}^{{\underline s}}$ is an isomorphism. The morphism ${\mathcal P}_{4}^{{\underline s}}$ takes the form
\begin{equation*}
\begin{split}
{\mathcal P}_{4}^{{\underline s}}\left(a_{3},a_{4},\eta_3,\xi_4\right)&=[a_3 a_4-a_3 a_4\eta_3\xi_4,-a_3 a_4\eta_3\xi_4,-a_3a_4] \\
&=[1-\eta_3\xi_4,-\eta_3\xi_4,-1].
\end{split}
\end{equation*} 
By Lemma \ref{zmodel}, we conclude that ${\mathcal P}_{4}^{{\underline s}}$ is flat in $\left(\overline{\mathcal R}_{4}^{{\underline s}}\right)^{-1}\left(A^{\tau}\right)$.

\smallskip

{\bf\noindent Cases (3)} and {\bf (4).} By permuting the first two columns, we can reduce them to  {Cases (1) and (2)}. We omit the details for brevity.
\smallskip

The proof is complete.\,\,\,\,$\endpf$

\subsection{Factorization}\label{FAC} In this subsection, we will introduce the method of factorization by considering a special subsequence of the blow-ups. We first fix certain notations.

Notice that when $\underline s=(1,1,\cdots,1)$, $G^{\underline s}_{\underline w}$ consists of at least two elements if and only if $w_{i_1}=w_{i_2}=w_{i_3}=w_{i_4}=1$,  and $w_t=0$ otherwise, for certain $1\leq i_1<i_2<i_3<i_4\leq N$. For such $\underline w$, $\mathbb P^{N^{\underline s}_{\underline w}}\cong\mathbb P^2$, and the rational map \begin{equation}\label{rrpe}
F^{\underline s}_{\underline w}:\mathbb P^{N_{2,n}}\dashrightarrow \mathbb P^{N^{\underline s}_{\underline w}}   
\end{equation} takes the form 
\begin{equation}\label{Fw2}
F^{\underline s}_{\underline w}\left([\cdots ,z_I,\cdots]_{I\in\mathbb I_{2,n}}\right)=\left[z_{(i_1,i_2)}\cdot z_{(i_3,i_4)},\,\,\,\,z_{(i_1,i_3)}\cdot z_{(i_2,i_4)},\,\,\,\,z_{(i_1,i_4)}\cdot z_{(i_2,i_3)}\right].
\end{equation}
Hence, in the remaining of this section, we denote (\ref{rrpe}) by $F_{i_1i_2i_3i_4}:\mathbb P^{N_{2,n}}\dashrightarrow \mathbb P^{N_{i_1i_2i_3i_4}}$, and the ideal sheaf $\mathscr I^{\underline s}_{\underline w}$ by  $\mathscr I^{\underline s}_{i_1i_2i_3i_4}$, to have more explicit referential indices.
By omitting the trivial projective spaces, we identify $\prod\nolimits_{\underline w\in C^{\underline s}}\mathbb {P}^{N_{\underline w}^{\underline s}}$ with $\prod\nolimits_{1\leq i_1<i_2<i_3<i_4\leq n}\mathbb P^{N_{i_1i_2i_3i_4}}$, and use them interchangeably for convenience. 


We now fix $\tau_+:=\left(1,2,(j_1^+,j_2^+,\cdots,j^+_{N-2})\right)=\left(1,2,(3,4,\cdots,n)\right)\in\mathbb J^{\underline s}$.
We take a total order $\sigma$ on the set  \begin{equation}\label{i1234}
\{(i_1,i_2,i_3,i_4)\in\mathbb Z^4|1\leq i_1<i_2<i_3<i_4\leq n\}    
\end{equation} so that
\begin{equation*}
\sigma(\alpha)\in\{(1,2,i_3,i_4)|\,\,3\leq i_3<i_4\leq n\}\,\,\,\,{\rm for}\,\,\,\, 1\leq\alpha\leq\frac{1}{2}(n-2)(n-3).
\end{equation*}
Define blow-ups
$\check{\mathcal R}_n:W^{\sigma}_{\frac{1}{2}(n-2)(n-3)}\rightarrow \mathcal Q_{n}^{\underline s}$, $\ddot{\mathcal R}_n:\mathcal T_{n}^{\underline s}\rightarrow W^{\sigma}_{\frac{1}{2}(n-2)(n-3)}$ respectively by
\begin{equation}\label{12a1}
\begin{split}
&\check {\mathcal R}_n:=\psi^{\sigma}_{1}\circ\psi^{\sigma}_{2}\circ\cdots\circ\psi^{\sigma}_{k}\circ\cdots\circ\psi^{\sigma}_{\frac{1}{2}(n-2)(n-3)},\\
&\ddot{\mathcal R}_n:=\psi^{\sigma}_{\frac{1}{2}(n-2)(n-3)+1}\circ\psi^{\sigma}_{\frac{1}{2}(n-2)(n-3)+2}\circ\cdots\circ\psi^{\sigma}_{\frac{1}{24}n(n-1)(n-2)(n-3)}.\\    
\end{split}
\end{equation}
(see \S \ref{pdm} for the definitions of $\psi^{\sigma}_{k}$, $W^{\sigma}_{k}$). All blow-ups fit into the following commutative diagram.
\vspace{-0.07in}
\begin{equation*}
\small
\begin{tikzcd}
\mathcal T_{n}^{\underline s}\arrow[rrr,"\ddot{\mathcal R}_n"]\arrow[rrrrrr, "\overline{\mathcal R}^{\underline s}_n"',bend right=15]\arrow[rrrrrrrrr, "\widehat{\mathcal R}_{n}^{{\underline s}}",bend left=15]&&&W^{\sigma}_{\frac{1}{2}(n-2)(n-3)}\arrow[rrr, "\check{\mathcal R}_n"]&&&\mathcal Q_{n}^{\underline s}\arrow[rrr, "\widetilde{\mathcal R}^{\underline s}"]&&&G(2,n).
\end{tikzcd}\vspace{-3pt}
\end{equation*}

Associate to each permutation  $\lambda$ of $\{1,2,\cdots,N-2\}$ an affine space ${\rm Spec}\,\mathbb Z\left[\overrightarrow  A,\overrightarrow  E\right]$, where
\begin{equation}\label{bp1}
\begin{split}&\overrightarrow  A:=\left(a_{3},a_{4},\cdots,a_{N}\right),\,\,\overrightarrow  E:=\left(\epsilon^+_{1},\epsilon^+_{2},\cdots,\epsilon^+_{N-2}\right).\\
\end{split}
\end{equation}
Define a morphism $\Sigma^{\lambda}:{\rm Spec}\,\mathbb Z\left[\overrightarrow  A,\overrightarrow  E\right]\rightarrow U_{12}$  by
\begin{equation}\label{ngamma1}
\left\{
\begin{aligned}
&x_{1j}\mapsto a_{j},\,\,\,\,\,\,\,\,\,\,\,\,\,\,\,\,\,\,\,\,\,\,\,\,\,\,\,\,\,\,\,\,\,\,\,\,\,\,\,\,\,\,{\rm for}\,\,3\leq j\leq n\\
&x_{2j}\mapsto a_j\cdot\prod\nolimits_{\gamma=1}^{\lambda(j-2)}\epsilon^+_{\gamma},\,\,\,\,\,\,\,\,\,\,{\rm for}\,\,3\leq j\leq n\\
\end{aligned}\right.\,\,.
\end{equation}
Similarly to Lemma \ref{em}, we can derive a locally closed embedding
\begin{equation*}
\Omega^{\lambda}:{\rm Spec}\,\mathbb Z\left[\overrightarrow  A,\overrightarrow  E\right]\longrightarrow\mathbb P^{N_{2,n}}\times\prod\nolimits_{t=1}^N\mathbb P^{N^{\underline s}_{t}}\times\prod\nolimits_{\substack{3\leq i_3<i_4\leq n}}\mathbb P^{N_{12i_3i_4}}    
\end{equation*}
by extending the rational map
\begin{equation*}
\left(\widetilde{\mathcal K}^{\underline s}\circ\Sigma^{\lambda},\left(\cdots,F_{12i_3i_4}\circ e\circ \Sigma^{\lambda},\cdots\right)_{3\leq i_3<i_4\leq n}\right).
\end{equation*} 
Denote by  $B^{\lambda}$ the image of ${\rm Spec}\,\mathbb Z\left[\overrightarrow  A,\overrightarrow  E\right]$ under $\Omega^{\lambda}$, which is equipped with the reduced scheme structure as a locally closed subscheme of $\mathbb P^{N_{2,n}}\times\prod\nolimits_{t=1}^N\mathbb {P}^{N^{\underline s}_{t}}\times\prod\nolimits_{\substack{3\leq i_3<i_4\leq n}}\mathbb P^{N_{12i_3i_4}}$. 

\begin{lemma}\label{fac} $B^{\lambda}$ is an open subscheme of $\left(\check{\mathcal R}_n\right)^{-1}\left(A^{\tau_+}\right)$, and
\begin{equation*}\bigcup\nolimits_{{\rm all\,\,permutations\,\,}\lambda\,\,{\rm of\,\,}\{1,2,\cdots,N-2\}}B^{\lambda}=\left(\check{\mathcal R}_n\right)^{-1}\left(A^{\tau_+}\right).
\end{equation*}
In particular, $\left(\check{\mathcal R}_n\right)^{-1}\left(A^{\tau_+}\right)$ is smooth over ${\rm Spec}\,\mathbb Z$.
\end{lemma}

{\bf\noindent Proof of Lemma \ref{fac}.}
We shall prove by induction on $n$. The cases $n=2,3$ are trivial, and the case $n=4$ follows from the proof of Lemma \ref{wn4}. Suppose that Lemma \ref{fac} holds for all $n$ such that $2\leq n\leq k$ where $k\geq 4$. We proceed to prove that it holds for $n=k+1$.

Take a total order $\sigma$ on  $\{(i_1,i_2,i_3,i_4)|1\leq i_1<i_2<i_3<i_4\leq k+1\}$ such that
\begin{equation}\label{ngamma3}
\small
\left\{\begin{aligned}
&\,\sigma(\alpha)\in\{(1,2,i_3,i_4)|\,\,3\leq i_3<i_4\leq k\},\,\,\,\,{\rm for}\,\,\,\, 1\leq\alpha\leq\frac{1}{2}(k-2)(k-3),\\
&\,\sigma(\alpha)\in\{(1,2,i_3,k+1)|\,\,3\leq i_3\leq k\},\,\,\,\,{\rm for}\,\,\,\, \frac{1}{2}(k-2)(k-3)+1\leq\alpha\leq\frac{1}{2}(k-1)(k-2).\\
\end{aligned}\right.\,\,\,
\end{equation}
Define blow-ups
$\check{\mathcal R}:W^{\sigma}_{\frac{1}{2}(k-2)(k-3)}\rightarrow \mathcal Q_{k+1}^{\underline s}$ and $\check{\mathcal R}^0:W^{\sigma}_{\frac{1}{2}(k-1)(k-2)}\rightarrow W^{\sigma}_{\frac{1}{2}(k-2)(k-3)}$ respectively by
$\check {\mathcal R}:=\psi^{\sigma}_{1}\circ\psi^{\sigma}_{2}\circ\cdots\circ\psi^{\sigma}_{\frac{1}{2}(k-2)(k-3)}$ and $\check{\mathcal R}^0:=\psi^{\sigma}_{\frac{1}{2}(k-2)(k-3)+1}\circ\psi^{\sigma}_{\frac{1}{2}(k-2)(k-3)+2}\circ\cdots\circ\psi^{\sigma}_{\frac{1}{2}(k-1)(k-2)}$.

Associate to each permutation  $\iota$ of $\{1,2,\cdots,k-2\}$ an affine space ${\rm Spec}\,\mathbb Z\left[\overrightarrow  A,\overrightarrow  E,\overrightarrow  H\right]$, where
$\overrightarrow  A:=\left(a_{3},a_{4},\cdots,a_{k+1}\right)$, $\overrightarrow  E:=\left(\epsilon^+_{1},\epsilon^+_{2},\cdots,\epsilon^+_{k-2}\right)$, $\overrightarrow  H:=\left(\eta_{k+1}\right)$. Define a morphism $\Sigma_{k}^{\iota}:{\rm Spec}\,\mathbb Z\left[\overrightarrow  A,\overrightarrow  E,\overrightarrow  H\right]\rightarrow U_{12}$ by
\begin{equation*}
\left\{\begin{aligned}
&\,x_{1j}\mapsto a_{j},\,\,\,\,\,\,\,\,\,\,\,\,\,\,\,\,\,\,\,\,\,\,\,\,\,\,\,\,\,\,\,\,\,\,\,\,\,\,\,\,\,{\rm for}\,\,3\leq j\leq k+1\\
&\,x_{2j}\mapsto a_{j}\cdot\prod\nolimits_{\gamma=1}^{\iota(j-2)}\epsilon^+_{\gamma},\,\,\,\,\,\,\,\,\,\,{\rm for}\,\,3\leq j\leq k\\
&\,x_{2(k+1)}\mapsto a_{k+1}\cdot\eta_{k+1}\\
\end{aligned}\right.\,\,\,.
\end{equation*}
By the same token, we can define a locally closed embedding 
\begin{equation*}
\Omega_{k}^{\iota}:{\rm Spec}\,\mathbb Z\left[\overrightarrow  A,\overrightarrow  E,\overrightarrow  H\right]\longrightarrow\mathbb P^{N_{2,k+1}}\times\prod\nolimits_{t=1}^{k+1}\mathbb P^{N^{\underline s}_{t}}\times\prod\nolimits_{\substack{3\leq i_3<i_4\leq k}}\mathbb P^{N_{12i_3i_4}}
\end{equation*}  by extending the rational map
$\left(\widetilde{\mathcal K}^{\underline s}\circ \Sigma^{\iota}_{k},\left(\cdots,F_{12i_3i_4}\circ e\circ\Sigma^{\iota}_{k},\cdots\right)_{3\leq i_3<i_4\leq k}\right)$. Denote by $B^{\iota}_{k}$ the image of ${\rm Spec}\,\mathbb Z\left[\overrightarrow  A,\overrightarrow  E,\overrightarrow  H\right]$ under the morphism $\Omega^{\iota}_{k}$, which is equipped with the reduced scheme structure as a locally closed subscheme of $\mathbb P^{N_{2,k+1}}\times\prod\nolimits_{t=1}^{k+1}\mathbb {P}^{N^{\underline s}_{t}}\times\prod\nolimits_{\substack{3\leq i_3<i_4\leq k}}\mathbb P^{N_{12i_3i_4}}$. 

Notice that the blow-up $\check {\mathcal R}$ is with respect to an ideal sheaf that is generated by polynomials independent of the $a_{k+1},\eta_{k+1}$ variables. Then, we can conclude by the induction hypothesis  that each  $B^{\iota}_{k}$ is an open subscheme of $\left(\check{\mathcal R}\right)^{-1}\left(A^{\tau_+}\right)$, and
\begin{equation*}
\bigcup\nolimits_{{\rm all\,\,permutations\,\,}\iota\,\,{\rm of\,\,}\{1,2,\cdots,k-2\}}B^{\iota}_{k}=\left(\check{\mathcal R}\right)^{-1}\left(A^{\tau_+}\right).  
\end{equation*}

Fix a permutation $\iota$ of $\{1,2,\cdots,k-2\}$. Denote by $\Lambda$ the set consisting of all permutations $\lambda$ of $\{1,2,\cdots,k-1\}$ such that the following holds.
\begin{enumerate}[label={$\bullet$}]
\item For any $1\leq i_1,i_2\leq k-2$, $\lambda(i_1)<\lambda(i_2)$ if and only if
$\iota(i_1)<\iota(i_2)$.
\end{enumerate}
To prove Lemma \ref{fac} for $n=k+1$, it suffices to show that for each $\lambda\in\Lambda$, $B^{\lambda}$ is an open subscheme of $\left(\check{\mathcal R}^0\right)^{-1}\left(B^{\iota}_{k}\right)$, and 
\begin{equation*}
\bigcup\nolimits_{{\rm all\,\,permutations\,\,}\lambda\in\Lambda}B^{\lambda}=\left(\check{\mathcal R}^0\right)^{-1}\left(B^{\iota}_{k}\right).
\end{equation*}
Without loss of generality, we can assume that $\iota$ is the identity permutation in the following.

We rearrange the order $\sigma$ so that (\ref{ngamma3}) holds and moreover
\begin{equation*}
\small
\sigma(\alpha)=\left(1,2,\alpha+2-\frac{1}{2}(k-2)(k-3),k+1\right),
\end{equation*}
for $\frac{1}{2}(k-2)(k-3)+1\leq\alpha\leq\frac{1}{2}(k-1)(k-2)$. For simplicity, we denote by $\psi_j$ the blow-up  $\psi^{\sigma}_{j-2+\frac{1}{2}(k-2)(k-3)}$ for $3\leq j\leq k$. 

Computation yields that  the ideal sheaf $\left(\widetilde{\mathcal R}^{{\underline s}}\circ\check{\mathcal R}\right)^{-1}\mathscr I^{\underline s}_{123(k+1)}\cdot\mathcal O_{B^{\iota}_{k}}$ is generated in $B^{\iota}_{k}$ by
$\eta_{k+1}$ and $\epsilon^+_{1}$. Then, $\left(\psi_3\right)^{-1}\left(B^{\iota}_{k}\right)$ is isomorphic to the closed subscheme of ${\rm Spec}\,\mathbb Z\left[\overrightarrow  A,\overrightarrow  E,\overrightarrow  H\right]\times{\rm Proj}\,\mathbb Z\left[u_1,v_1\right]$ defined by
\begin{equation*}
u_1\cdot\eta_{k+1}=v_1\cdot\epsilon^+_{1}.
\end{equation*} 
It is clear that $\left(\psi_3\right)^{-1}\left(B^{\iota}_{k}\right)$ is covered by the following two open subschemes:
\begin{enumerate}[label={$\bullet$}]
\item $V_1\cong{\rm Spec}\,\mathbb Z\left[a_{3},\cdots,a_{k+1},\epsilon^+_{1},\epsilon^+_{2},\cdots,\epsilon^+_{k-2},\nu_1\right]$, where the embedding $V_1\hookrightarrow\left(\psi_3\right)^{-1}\left(B^{\iota}_{k}\right)$ is induced by
\begin{equation*}
\frac{v_1}{u_1}\mapsto \nu_1\,\,\,\,{\rm and}\,\,\,\,\eta_{k+1}\mapsto\epsilon^+_1\cdot \nu_1.
\end{equation*}

\item $V_2\cong{\rm Spec}\,\mathbb Z\left[a_{3},\cdots,a_{k+1},\epsilon^+_{2},\cdots,\epsilon^+_{k-2},\eta_{k+1},\mu_1\right]$, where the embedding $V_1\hookrightarrow\left(\psi_3\right)^{-1}\left(B^{\iota}_{k}\right)$ is induced by \begin{equation*}
\frac{u_1}{v_1}\mapsto \mu_1\,\,\,\,{\rm and}\,\,\,\,\epsilon^+_1\mapsto\eta_{k+1}\cdot \mu_1.
\end{equation*}\label{112b}
\end{enumerate}
The blow-ups $\psi_j$ are trivial in $V_2$ for $4\leq j\leq k$, hence after renaming the variables we can conclude that $\left(\psi_4\circ\psi_5\circ\cdots\circ\psi_{k}\right)^{-1}\left(V_2\right)=B^{\lambda}$, where $\lambda$ is the permutation defined by
\begin{equation*}
\lambda(j)=\left\{\begin{aligned}
&\,j+1\,\,\,\,{\rm for}\,\,\,\, 1\leq j\leq k-2\\
&\,1\,\,\,\,\,\,\,\,\,\,\,\,\,\,{\rm for}\,\,\,\, j=k-1\\
\end{aligned}\right.\,\,\,.
\end{equation*}

Computation yields that the ideal sheaf $\left(\widetilde{\mathcal R}^{{\underline s}}\circ\check{\mathcal R}\circ\psi_3\right)^{-1}\mathscr I^{\underline s}_{124(k+1)}\cdot\mathcal O_{V_1}$ is generated in $V_1$ by
$\epsilon^+_1\cdot\nu_1$ and $\epsilon^+_{1}\cdot\epsilon^+_{2}$. Then,  $\left(\psi_4\right)^{-1}\left(V_1\right)$ is isomorphic to the closed subscheme of $V_1\times{\rm Proj}\,\mathbb Z\left[u_2,v_2\right]$ defined by
\begin{equation*}
u_2\cdot\nu_1=v_2\cdot\epsilon^+_{2}.
\end{equation*} 
Moreover, $\left(\psi_4\right)^{-1}\left(V_1\right)$ is covered by the following two open subschemes:
\begin{enumerate}[label={$\bullet$}]
\item $U_1\cong{\rm Spec}\,\mathbb Z\left[a_{3},\cdots,a_{k+1},\epsilon^+_{1},\epsilon^+_{2},\epsilon^+_{3},\cdots,\epsilon^+_{k-2},\nu_2\right]$, where the embedding $U_1\hookrightarrow\left(\psi_4\right)^{-1}\left(V_1\right)$ is induced by 
\begin{equation*}
\frac{v_2}{u_2}\mapsto\nu_2\,\,\,\,{\rm and}\,\,\,\,\nu_1\mapsto\epsilon^+_2\cdot \nu_2.
\end{equation*}

\item $U_2\cong{\rm Spec}\,\mathbb Z\left[a_{3},\cdots,a_{k+1},\epsilon^+_{1},\epsilon^+_{3},\cdots,\epsilon^+_{k-2},\nu_1,\mu_2\right]$, where the embedding $U_2\hookrightarrow\left(\psi_4\right)^{-1}\left(V_1\right)$ is induced by 
\begin{equation*}
\frac{u_2}{v_2}\mapsto\mu_2\,\,\,\,{\rm and}\,\,\,\,\epsilon^+_2\mapsto\nu_1\cdot \mu_2.
\end{equation*}
\end{enumerate}
Since $\psi_j$ are trivial in $U_2$ for $5\leq j\leq k$, after renaming the variables we can conclude that $\left(\psi_5\circ\cdots\circ\psi_{k}\right)^{-1}\left(U_2\right)=B^{\lambda}$, where $\lambda$ is the permutation defined by
\begin{equation*}
\lambda(j)=\left\{\begin{aligned}
&\,1\,\,\,\,\,\,\,\,\,\,\,\,\,\,{\rm for}\,\,\,\, j=1\\
&\,j+1\,\,\,\,{\rm for}\,\,\,\, 2\leq j\leq k-2\\
&\,2\,\,\,\,\,\,\,\,\,\,\,\,\,\,{\rm for}\,\,\,\, j=k-1\\
\end{aligned}\right.\,\,\,.
\end{equation*}

Computation yields that the ideal sheaf $\left(\widetilde{\mathcal R}^{{\underline s}}\circ\check{\mathcal R}\circ\psi_3\circ\psi_4\right)^{-1}\mathscr I^{\underline s}_{125(k+1)}\cdot\mathcal O_{U_1}$ is generated in $U_1$ by $\epsilon^+_1\cdot\epsilon^+_{2}\cdot\nu_2$ and $\epsilon^+_{1}\cdot\epsilon^+_{2}\cdot\epsilon^+_{3}\}$, hence  $\left(\psi_5\right)^{-1}\left(U_1\right)$ is isomorphic to the closed subscheme of $U_1\times{\rm Proj}\,\mathbb Z\left[u_3,v_3\right]$ defined by
\begin{equation*}
u_3\cdot\nu_2=v_3\cdot\epsilon^+_{3}. 
\end{equation*}
By the same argument (similar to the bubble sort algorithm), we can continue the above procedure all the way up to $\psi_{k}$, and finally prove Lemma \ref{fac} hold for $n=k+1$.

We complete the proof of Lemma \ref{fac}.\,\,\,\,\,\,$\endpf$


\begin{remark}\label{fac2}
Let $\tau_-:=\left(1,2,(j_1^-,\cdots,j^-_{N-2})\right)=\left(1,2,(3,4,\cdots,n)\right)\in\mathbb J^{\underline s}$. We can construct an open cover $\{B^{\lambda}\}_{{\rm all\,\,permutations\,\,}\lambda\,\,{\rm of\,\,}\{1,2,\cdots,N-2\}}$ of 
$\left(\check{\mathcal R}_n\right)^{-1}\left(A^{\tau_-}\right)$ in the same way.
\end{remark}

\subsection{Reduction to a local form}
We first modify the open cover $\left\{A^{\tau}\right\}_{\tau\in\mathbb J^{\underline s}}$ of $\mathcal Q_{n}^{\underline s}$ as follows. For each $\tau=\left(j_1,j_2,(j^+_1,\cdots,j^+_l),(j^-_1,\cdots,j^-_{m})\right)\in\mathbb J^{\underline s}$, we can define an open subscheme $\mathring A^{\tau}$ of $A^{\tau}\cong{\rm Spec}\,\mathbb Z\left[\overrightarrow  A,\overrightarrow  H,\overrightarrow  \Xi\right]$
by
\begin{equation}\label{punc}
\mathring A^{\tau}:=\left\{\,\mathfrak p\in A^{\tau}\left|\,1-\eta_{j^+_{\alpha}}\cdot\xi_{j^-_{\beta}}\notin\mathfrak p,\,\,\forall\,1\leq\alpha\leq l\,\,{\rm and}\,\,1\leq\beta\leq m\right.\right\}.
\end{equation}
Note that if $l=0$ or $m=0$,  $\mathring A^{\tau}=A^{\tau}$.

\begin{definitionlemma}\label{trun}
$\left\{\mathring A^{\tau}\right\}_{\tau\in\mathbb J^{\underline s}}$ is a finite open cover of  $\mathcal Q_{n}^{\underline s}$. We call $\mathring A^{\tau}$ the (truncated) coordinate charts of  $\mathcal Q_{n}^{\underline s}$.
\end{definitionlemma}
{\bf\noindent Proof of Definition-Lemma \ref{trun}.} Fix an arbitrary $\mathfrak p\in \mathcal Q_{n}^{\underline s}$. Since $\left\{A^{\tau}\right\}$ is an open cover of $\mathcal Q_{n}^{\underline s}$, there exists an index $\tau=\left(j_1,j_2,(j^+_1,\cdots,j^+_{l}),(j^-_1,\cdots,j^-_{m})\right)\in\mathbb J^{\underline s}$ such that \begin{enumerate}
    \item $\mathfrak p\in A^{\tau}$,
    
    \item $|l-m|$ is the maximum of all the numbers $\big|\widetilde l-\widetilde m\big|$ associated to the indices $\tilde{\tau}=\left(j_1,j_2,(\tilde j^+_1,\cdots,\tilde j^+_{\widetilde l}),(\tilde j^-_1,\cdots,\tilde j^-_{\widetilde m})\right)\in\mathbb J^{\underline s}$ such that $\mathfrak p\in A^{\tilde{\tau}}$.
\end{enumerate}

\begin{claim}\label{vani}
When  $m\geq l$, $\eta_{ j^+_1},\eta_{j^+_2},\cdots,\eta_{ j^+_{ l}}\in\mathfrak p$; when $ l\geq m$, $\xi_{ j^-_1},\xi_{j^-_2},\cdots,\xi_{j^-_{m}}\in\mathfrak p$. 
\end{claim} 
{\bf\noindent Proof of Claim \ref{vani}.} 
We may assume that
$m\geq l$. If there exists  $\eta_{j^+_{\alpha_0}}\notin\mathfrak p$, we can define a new index $\hat\tau:=\left(j_1,j_2,(\hat j^+_1,\cdots,\hat j^+_{l-1}),(\hat j^-_1,\cdots,\hat j^-_{m+1})\right)\in\mathbb J^{\underline s}$ by
\begin{equation*}
\begin{split}
&\hat j^+_1=j_1^+,\,\,\,\,\hat j^+_2=j_2^+,\,\,\,\,\cdots,\,\,\,\,\hat j^+_{\alpha_0-1}=j_{\alpha_0-1}^+,\,\,\,\,\hat j^+_{\alpha_0}=j_{\alpha_0+1}^+,\,\,\,\,\hat j^+_{\alpha_0+1}=j_{\alpha_0+2}^+,\,\,\,\,\cdots,\,\,\,\,\hat j^+_{l-1}=j_{l}^+,\\
&\hat j^-_1=j_1^-,\,\,\,\,\hat j^-_2=j_2^-,\,\,\,\,\cdots,\,\,\,\,\hat j^-_{\beta_0-1}=j_{\beta_0-1}^-,\,\,\,\,\,\,\hat j^+_{\beta_0}=j_{\alpha_0}^+,\,\,\,\,\,\,\hat j^-_{\beta_0+1}=j_{\beta_0}^-,\,\,\,\,\cdots,\,\,\,\,\hat j^-_{m+1}=j_{m}^-.\\
\end{split}    
\end{equation*}

Then the transition map between 
\begin{equation*}
A^{\tau}\cong{\rm Spec}\,\mathbb Z\left[\cdots, a_{ j^+_{\alpha}},\cdots, a_{j^-_{\beta}},\cdots,\eta_{j^+_{\alpha}},\cdots, \xi_{j^-_{\beta}},\cdots\right]    
\end{equation*}
and
\begin{equation*}
A^{\hat\tau}\cong{\rm Spec}\,\mathbb Z\left[\cdots,\hat a_{\hat j^+_{\alpha}},\cdots,\hat a_{\hat j^-_{\beta}},\cdots,\hat \eta_{\hat j^+_{\alpha}},\cdots,\hat \xi_{\hat j^-_{\beta}},\cdots\right]    
\end{equation*}
takes the form
\begin{equation*}
\left\{\begin{aligned}
&\,\hat a_{\hat j^+_{\alpha}}\mapsto a_{j^+_{\alpha}},\,\,\,\,\,\,\,\,\,\,\,\,\,\,\,\,\,\,\,\,\,\,\,\,\hat\eta_{\hat j^+_{\alpha}}\mapsto \eta_{j^+_{\alpha}},\,\,\,\,\,\,\,\,\,\,\,\,\,\,\,\,\,\,\,{\rm for}\,\,1\leq\alpha\leq \alpha_0-1\\
&\,\hat a_{\hat j^+_{\alpha}}\mapsto a_{j^+_{\alpha+1}},\,\,\,\,\,\,\,\,\,\,\,\,\,\,\,\,\,\,\,\,\hat\eta_{\hat j^+_{\alpha}}\mapsto \eta_{j^+_{\alpha+1}},\,\,\,\,\,\,\,\,\,\,\,\,\,\,\,{\rm for}\,\,\alpha_0\leq\alpha\leq l-1\\
&\,\hat a_{\hat j^-_{\beta_0}}\mapsto a_{j^+_{\alpha_0}}\cdot \eta_{j^+_{\alpha_0}},\,\,\,\,\,\,\,\,\hat\xi_{\hat j^-_{\beta_0}}\mapsto 1/\eta_{j^+_{\alpha_0}}\\
&\,\hat a_{\hat j^-_{\beta}}\mapsto a_{j^-_{\beta}},\,\,\,\,\,\,\,\,\,\,\,\,\,\,\,\,\,\,\,\,\,\,\,\,\hat\xi_{j^-_{\beta}}\mapsto \xi_{j^-_{\beta}},\,\,\,\,\,\,\,\,\,\,\,\,\,\,\,\,\,\,\,\,{\rm for}\,\,1\leq\beta\leq \beta_0-1,\\
&\,\hat a_{\hat j^-_{\beta}}\mapsto a_{j^-_{\beta-1}},\,\,\,\,\,\,\,\,\,\,\,\,\,\,\,\,\,\,\,\,\hat\xi_{j^-_{\beta}}\mapsto \xi_{j^-_{\beta-1}},\,\,\,\,\,\,\,\,\,\,\,\,\,\,\,\,{\rm for}\,\,\beta_0+1\leq\beta\leq m+1\\
\end{aligned}\right..
\end{equation*}
It is easy to verify that $\mathfrak p\in A^{\hat\tau}$ and hence $|m+1-(l-1)|=|m-l|+2$, which is a contradiction. We complete the proof of Claim \ref{vani}.\,\,\,\,$\endpf$\smallskip

Definition-Lemma \ref{trun} then follows by Claim.\,\,\,\,$\endpf$
\medskip

We state the following (equivalent) local form of Proposition \ref{ms}.

\begin{lemma}\label{loc}  Let $N\geq 4$ be an integer. Let $\underline s=(1,\cdots,1)\in\mathbb Z^N$. For each (truncated) coordinate chart ${\mathring A^{\tau}}\subset\mathcal Q_{n}^{\underline s}$ where $\tau=\left(1,2,(j^+_1,\cdots,j^+_l),(j^-_1,\cdots,j^-_{m})\right)\in\mathbb J^{\underline s}$, there is a finite open cover $\{T^{\tau}_{\alpha}\}_{\alpha}$ of  $\left(\overline{\mathcal R}^{\underline s}_n\right)^{-1}\left({\mathring A^{\tau}}\right)$ with the following properties.

\begin{enumerate}[label=(\arabic*),ref=\arabic*]

\item Each $T^{\tau}_{\alpha}$ is an open subscheme of $\left(\overline{\mathcal R}^{\underline s}_n\right)^{-1}\left({\mathring A^{\tau}}\right)$, and is smooth over ${\rm Spec}\,\mathbb Z$.
\label{loc1}    
\item There is a locally closed subscheme $M^{\tau}_{\alpha}$ of $\prod\nolimits_{\underline w\in C^{\underline s}}\mathbb {P}^{N^{\underline s}_{\underline w}}$ for each $\alpha$, such that as a topological space $M^{\tau}_{\alpha}$ is the image of  $T^{\tau}_{\alpha}$ under the projection from 
$\mathbb {P}^{N_{2,n}}\times\prod\nolimits_{t=1}^N\mathbb {P}^{N^{\underline s}_{t}}\times\prod\nolimits_{\underline w\in C^{\underline s}}\mathbb {P}^{N^{\underline s}_{\underline w}}$  to  $\prod\nolimits_{\underline w\in C^{\underline s}}\mathbb {P}^{N^{\underline s}_{\underline w}}$. Moreover, $M^{\tau}_{\alpha}$ is smooth over ${\rm Spec}\,\mathbb Z$ of relative dimension  $n-3$. 
\label{loc2}  

\item The projection from $T^{\tau}_{\alpha}$ to $M^{\tau}_{\alpha}$ is a flat morphism.
\label{loc3}  

\end{enumerate}
\end{lemma}

Notice that Lemma \ref{loc} follows by  Proposition \ref{ms} immediately, for a flat morphism of finite type of Noetherian schemes is  open. Conversely, we have
\begin{lemma}\label{red}
Proposition \ref{ms} holds for an integer $N\geq 4$, if Lemma \ref{loc} holds for such $N$.
\end{lemma}

{\noindent\bf Proof Lemma \ref{red}.} Notice that the constructions of $\widehat{\mathcal R}_{n}^{\underline s}:\mathcal T_{n}^{\underline s}\rightarrow G(2,n)={\rm Gr}^{2,E}$ and $\widetilde{\mathcal R}^{\underline s}:\mathcal Q_{n}^{\underline s}\rightarrow G(2,n)={\rm Gr}^{2,E}$ are invariant under the permutations of the factors $E_t$ of  $E=\oplus_{t=1}^NE_{t}$. Hence, to prove that $\mathcal T_{n}^{\underline s}$ is smooth, it suffices to show that  $\left(\widehat{\mathcal R}_{n}^{\underline s}\right)^{-1}(U_{12})$ is smooth. According to  Definition-Lemma \ref{trun}, we can derive that $\mathcal T_{n}^{\underline s}$ is smooth by property (\ref{loc1}) in Lemma \ref{loc}.

By  property (\ref{loc2}) in Lemma \ref{loc}, as topological spaces, we have 
\begin{equation}\label{beta}
\mathcal M_{n}^{\underline s}=\bigcup\, M_{\beta}, 
\end{equation}
where $M_{\beta}$ are the images of $M^{\tau}_{\alpha}$ under auotomorphisms of $\prod\nolimits_{\underline w\in C^{\underline s}}\mathbb {P}^{N^{\underline s}_{\underline w}}$ induced by the permutations of $E_t$.
Denote by $\overline{M_{\beta}}$ the closure of $M_{\beta}$. 
We claim that for each $M_{\beta}$ in (\ref{beta})
\begin{equation*}
\overline{M_{\beta}}=\mathcal M_{n}^{\underline s}.
\end{equation*}
Otherwise, there exists an $\overline{M_{\beta^{\prime}}}$ that is a proper subscheme of $\mathcal M_{n}^{\underline s}$. Since  there exists an $M_{\beta^{\prime\prime}}$ such that $\overline{M_{\beta^{\prime\prime}}}=\mathcal M_{n}^{\underline s}$, $\overline{M_{\beta^{\prime}}}$ is a proper subscheme of $\overline{M_{\beta^{\prime\prime}}}$, which contradicts property (2) that they have the same dimension. Then $\mathcal M_{n}^{\underline s}$ is smooth over ${\rm Spec}\,\mathbb Z$ of relative dimension  $n-3$.

Since we have shown that the image of $T^{\tau}_{\alpha}$ is an open subscheme of $\mathcal M_{n}^{\underline s}$, the flatness of $\mathcal P^{\underline s}_{n}$ follows by property (\ref{loc3}) in Lemma \ref{loc}.

The proof of Lemma \ref{red} is complete. \,\,\,\,$\endpf$

\subsection{Induction for Type (\ref{mix})}

For $\tau\in\mathbb J^{\underline s}$ with $j_1=1,j_2=2$, we define its type as follows.

\begin{enumerate}[label={\bf Type\,(\arabic*)},ref=\arabic*]

\item $\tau=\left(1,2,(j^+_1,\cdots,j^+_l),(j^-_1,\cdots,j^-_{m})\right)\in\mathbb J^{\underline s}$, where $l,m\geq1$. 
\label{mix}

\item  $\tau=\left(1,2,(j_1^+,\cdots,j^+_{N-2})\right)\in\mathbb J^{\underline s}$, or $\tau=\left(1,2,(j_1^-,\cdots,j^-_{N-2})\right)\in\mathbb J^{\underline s}$. 
\label{up}
\end{enumerate}

In this subsection, we shall prove Lemma \ref{loc} for the (truncated) coordinate charts of {\bf Type (\ref{mix})}, provided that Proposition \ref{ms} holds for all $N^{\prime}$ such that $3\leq N^{\prime}<N$. The idea is to decompose $\mathcal T_{n}^{\underline s}$ locally into products of open suschemes of $\mathcal T_{l+2}^{\underline s_+}$ and 
$\mathcal T_{m+2}^{\underline s_-}$, where \begin{equation}
\underline s_+:=(1,1,\cdots,1)\in\mathbb Z^{l+2},\,\,\,\,\,\, \underline s_-:=(1,1,\cdots,1)\in\mathbb Z^{m+2}. 
\end{equation}

\begin{lemma}\label{uv}
Suppose that  Proposition \ref{ms} holds for all integers $N^{\prime}$ such that $3\leq N^{\prime}<N$. 
Then Lemma \ref{loc} holds for all (truncated) coordinate charts of {\bf Type (\ref{mix})}.
\end{lemma}
{\bf\noindent Proof of Lemma \ref{uv}.} 
Fix $\tau=\left(1,2,(j^+_1,\cdots,j^+_l),(j^-_1,\cdots,j^-_{m})\right)\in\mathbb J^{\underline s}$, where $l,m\geq1$ and $l+m=N-2$. Partition set (\ref{i1234}) as follows.
\begin{enumerate}[label={(\Alph*0)},ref={\Alph*0}]
\item $\{(1,2,j^+_{\alpha_1},j^+_{\alpha_2})|1\leq\alpha_1<\alpha_2\leq l\}$.

\item $\{(1,2,j^-_{\beta_1},j^-_{\beta_2})|1\leq\beta_1<\beta_2\leq m\}$.

\end{enumerate}

\begin{enumerate}[label={(A\arabic*)},ref={A\arabic*}]\setcounter{enumi}{0}
\item $\{(1,j^+_{\alpha_1},j^+_{\alpha_2},j^+_{\alpha_3})|1\leq\alpha_1<\alpha_2<\alpha_3\leq l\}$.

\item $\{(2,j^+_{\alpha_1},j^+_{\alpha_2},j^+_{\alpha_3})|1\leq\alpha_1<\alpha_2<\alpha_3\leq l\}.$

\item $\{(j^+_{\alpha_1},j^+_{\alpha_2},j^+_{\alpha_3},j^+_{\alpha_4})|1\leq\alpha_1<\alpha_2<\alpha_3<\alpha_4\leq l\}$.

\end{enumerate}

\begin{enumerate}[label={(B\arabic*)},ref={B\arabic*}]\setcounter{enumi}{0}
\item $\{(1,j^-_{\beta_1},j^-_{\beta_2},j^-_{\beta_3})|1\leq\beta_1<\beta_2<\beta_3\leq m\}$.

\item $\{(2,j^-_{\beta_1},j^-_{\beta_2},j^-_{\beta_3})|1\leq\beta_1<\beta_2<\beta_3\leq m\}$.

\item $\{(j^-_{\beta_1},j^-_{\beta_2},j^-_{\beta_3},j^-_{\beta_4})|1\leq\beta_1<\beta_2<\beta_3<\beta_4\leq m\}.$

\end{enumerate}
\begin{enumerate}[label={(C\arabic*)},ref={C\arabic*}]
\item $\{(1,2,i_3,i_4)|\begin{matrix}\{i_3,i_4\}=\{j^+_{\alpha_1},j^-_{\beta_1}\}\,\,{\rm where}\,\,1\leq\alpha_1\leq l\,\,{\rm and}\,\,1\leq\beta_1\leq m 
\end{matrix}\}$.

\item $\{(1,i_2,i_3,i_4)|\{i_2,i_3,i_4\}=\{j^+_{\alpha_1},j^+_{\alpha_2},j^-_{\beta_1}\}{\,\,\rm where}\,\,1\leq\alpha_1<\alpha_2\leq l\,\,{\rm and}\,\,1\leq\beta_1\leq m\}$.

\item $\{(2,i_2,i_3,i_4)|\{i_2,i_3,i_4\}=\{j^+_{\alpha_1},j^-_{\beta_1},j^-_{\beta_2}\}{\,\,\rm where}\,\,1\leq\alpha_1\leq l\,\,{\rm and}\,\,1\leq\beta_1<\beta_2\leq m\}$.

\item $\{(1,i_2,i_3,i_4)|\{i_2,i_3,i_4\}=\{j^+_{\alpha_1},j^-_{\beta_1},j^-_{\beta_2}\}{\,\,\rm where}\,\,1\leq\alpha_1\leq l\,\,{\rm and}\,\,1\leq\beta_1<\beta_2\leq m\}$.

\item $\{(2,i_2,i_3,i_4)|\{i_2,i_3,i_4\}=\{j^+_{\alpha_1},j^+_{\alpha_2},j^-_{\beta_1}\}{\,\,\rm where}\,\,1\leq\alpha_1<\alpha_2\leq l\,\,{\rm and}\,\,1\leq\beta_1\leq m\}$.

\item $\left\{(i_1,i_2,i_3,i_4)\left|\begin{matrix}
\{i_1,i_2,i_3,i_4\}=\{j^+_{\alpha_1},j^+_{\alpha_2},j^-_{\beta_1},j^-_{\beta_2}\}\,\,{\rm where}\,\\
 1\leq\alpha_1<\alpha_2\leq l\,\,{\rm and}\,\,1\leq\beta_1<\beta_2\leq m
\end{matrix}\right.\right\}$.

\item $\left\{(i_1,i_2,i_3,i_4)\left|\begin{matrix}
\{i_1,i_2,i_3,i_4\}=\{j^+_{\alpha_1},j^+_{\alpha_2},j^+_{\alpha_3},j^-_{\beta_1}\}\,\,{\rm where}\,\\
 1\leq\alpha_1<\alpha_2<\alpha_3\leq l\,\,{\rm and}\,\,1\leq\beta_1\leq m
\end{matrix}\right.\right\}$.

\item $\left\{(i_1,i_2,i_3,i_4)\left|\begin{matrix}
\{i_1,i_2,i_3,i_4\}=\{j^+_{\alpha_1},j^-_{\beta_1},j^-_{\beta_2},j^-_{\beta_3}\}\,\,{\rm where}\,\\
 1\leq\alpha_1\leq l\,\,{\rm and}\,\,1\leq\beta_1<\beta_2<\beta_3\leq m
\end{matrix}\right.\right\}$.

\end{enumerate}
Define a total order $\sigma$ so that the $4$-tuples in an upper subset are strictly smaller than the ones in a lower subset. Equivalently, we first blow up with respect to the (pull-backs of) ideal sheaves $\mathscr I^{\underline s}_{i_1i_2i_3i_4}$ with indices $(i_1,i_2,i_3,i_4)$ in (A0), (B0), then (A1), (A2), (A3), (B1), (B2), (B3), and finally (C1)-(C8).

Let ${\mathcal R}_{A0}:=\psi^{\sigma}_{1}\circ\psi^{\sigma}_{2}\circ\cdots\circ\psi^{\sigma}_{\frac{1}{2}(l-2)(l-3)}$ be the sequence of blow-ups with respect to the (pull-backs of) ideal sheaves $\mathscr I^{\underline s}_{i_1i_2i_3i_4}$ with indices in (A0). Let $\mathcal R_{A0B0}:=\psi^{\sigma}_{\frac{1}{2}(l-2)(l-3)+1}\circ\psi^{\sigma}_{\frac{1}{2}(l-2)(l-3)+2}\circ\cdots\circ\psi^{\sigma}_{\frac{1}{2}(l-2)(l-3)+\frac{1}{2}(m-2)(m-3)}$ the sequence of blow-ups with respect to the (pull-backs of) ideal sheaves $\mathscr I^{\underline s}_{i_1i_2i_3i_4}$ with indices in (B0). Denote by $\mathcal R_{AB}$ the sequence of blow-ups with respect to the (pull-backs of) ideal sheaves $\mathscr I^{\underline s}_{i_1i_2i_3i_4}$ with indices in (A1), (A2), (A3), (B1), (B2), (B3);  denote by $\mathcal R_{ABC}$ the sequence  of blow-ups with respect to the (pull-backs of) ideal sheaves $\mathscr I^{\underline s}_{i_1i_2i_3i_4}$ with indices in (C1)-(C8).

We first investigate  $\mathcal R_{A0}$ and $\mathcal R_{A0B0}$ via Lemma \ref{fac}. Associate to any permutations  $\lambda^+$ of $\{1,2,\cdots,l\}$ and $\lambda^-$ of $\{1,2,\cdots,m\}$ an affine space ${\rm Spec}\,\mathbb Z\left[\overrightarrow  A,\overrightarrow  E_+,\overrightarrow  E_-\right]$, where
\begin{equation*}
\small
\begin{split}
&\overrightarrow  A:=\left(a_{j^+_1},\cdots,a_{j^+_l},a_{j^-_1},\cdots,a_{j^-_m}\right),\,\overrightarrow  E_+:=\left(\epsilon_1^+,\cdots,\epsilon^+_l\right),\,\overrightarrow E_-:=\left(\epsilon^-_1,\cdots,\epsilon^-_m\right).\\
\end{split}
\end{equation*}
Define a morphism $\Sigma^{\lambda^+,\lambda^-}:{\rm Spec}\,\mathbb Z\left[\overrightarrow  A,\overrightarrow  E_+,\overrightarrow  E_-\right]\rightarrow U_{12}$ by
\begin{equation*}
\left\{\begin{aligned}
&\,x_{1j^+_{\alpha}}\mapsto a_{j^+_{\alpha}},\,\,\,\,\,\,\,\,\,\,\,\,\,\,\,\,\,\,\,\,\,\,\,\,\,\,\,\,\,\,\,\,\,\,\,\,\,\,\,{\rm for}\,\,1\leq {\alpha}\leq l\\
&\,x_{2j^+_{\alpha}}\mapsto a_{j^+_{\alpha}}\cdot\prod\nolimits_{\gamma=1}^{\lambda^+(\alpha)}\epsilon^+_{\gamma},\,\,\,\,\,\,\,\,{\rm for}\,\,1\leq \alpha\leq l\\
&\,x_{1j^-_{\beta}}\mapsto a_{j^-_{\beta}}\cdot\prod\nolimits_{\gamma=1}^{\lambda^-(\beta)}\epsilon^-_{\gamma},\,\,\,\,\,\,\,\,{\rm for}\,\,1\leq {\beta}\leq m\\
&\,x_{2j^-_{\beta}}\mapsto a_{j^-_{\beta}},\,\,\,\,\,\,\,\,\,\,\,\,\,\,\,\,\,\,\,\,\,\,\,\,\,\,\,\,\,\,\,\,\,\,\,\,\,\,\,{\rm for}\,\,1\leq {\beta}\leq m\\
\end{aligned}\right..
\end{equation*}
We can define a locally closed embedding $\Omega^{\lambda^+,\lambda^-}$ from ${\rm Spec}\,\mathbb Z\left[\overrightarrow  A,\overrightarrow  E_+,\overrightarrow  E_-\right]$ to \begin{equation}\label{pabmath}
\begin{split}
&\mathbb P^{N_{2,n}}\times\prod\nolimits_{t=1}^N\mathbb P^{N^{\underline s}_{t}}\times\prod\nolimits_{\substack{1\leq\alpha_1<\alpha_2\leq l}}\mathbb P^{N_{12j^+_{\alpha_1}j^+_{\alpha_2}}}\times\prod\nolimits_{\substack{1\leq\beta_1<\beta_2\leq m}}\mathbb P^{N_{12j^-_{\beta_1}j^-_{\beta_2}}} \\
\end{split}
\end{equation} by extending the rational map $\left(\widetilde{\mathcal K}^{\underline s}\circ\Sigma^{\lambda^+,\lambda^-},\left(\cdots,F_{12i_3i_4}\circ e\circ \Sigma^{\lambda^+,\lambda^-},\cdots\right)\right)$.  
Denote by $B^{\lambda^+,\lambda^-}$ the image of ${\rm Spec}\,\mathbb Z\left[\overrightarrow  A,\overrightarrow  E_+,\overrightarrow  E_-\right]$ under $\Omega^{\lambda^+,\lambda^-}$, which is equipped with the reduced scheme structure as a locally closed subscheme of (\ref{pabmath}).

\begin{claim}\label{df} $B^{\lambda^+,\lambda^-}$ is an open subscheme of  $\left({\mathcal R}_{A0}\circ{\mathcal R}_{A0B0}\right)^{-1}\left(A^{\tau}\right)$,  and
\begin{equation*}
\bigcup\nolimits_{{\rm all\,\,permutations\,\,}\lambda^+\,\,{\rm of\,\,}\{1,2,\cdots,l\}\,\,{\rm and}\,\,\lambda^-\,\,{\rm of\,\,}\{1,2,\cdots,m\}}B^{\lambda^+,\lambda^-}=\left({\mathcal R}_{A0}\circ{\mathcal R}_{A0B0}\right)^{-1}\left(A^{\tau}\right).
\end{equation*}
\end{claim}
{\noindent\bf Proof of Claim \ref{df}.} It is clear that $U_{12}\cong U^{\tau}_{12,l}\times U^{\tau}_{12,m}$, where  $U^{\tau}_{12,l}$ (resp. $U^{\tau}_{12,m}$) is the closed subscheme of $U_{12}$ defined by
\begin{equation}\label{subgr}
x_{1j^-_{\beta}}=x_{2j^-_{\beta}}=0,\,\,\forall\,1\leq\beta\leq m\,\,\,\left({\rm resp.}\,\,x_{1j^+_{\alpha}}=x_{2j^+_{\alpha}}=0,\,\,\forall\,1\leq\alpha\leq l\right).
\end{equation}
Note that $U^{\tau}_{12,l}$ (resp. $U^{\tau}_{12,m}$) is an open subscheme of the sub-Grassmannian $G(2,l+2)$ (resp. $G(2,m+2)$) of $G(2,n)$ determined by (\ref{subgr}). 
Hence, $A^{\tau}\subset\mathcal Q_{n}^{\underline s}$ has the Cartesian product structure  $A^{\tau}\cong A^{\tau_+}\times A^{\tau_-}$. Here $A^{\tau_+}$ is the open subscheme of the preimage of $U^{\tau}_{12,l}$  under the blow-up \begin{equation*}
\widetilde{\mathcal R}^{\underline s_+}:\mathcal Q^{\underline s_+}_{l+2}\rightarrow G(2,l+2),
\end{equation*}
associated to the size vector $\underline s_+:=(1,1,\cdots,1)\in\mathbb Z^{l+2}$ and the index $\tau_+:=(1,2,(j^+_1,j^+_2,\cdots,j^+_l))$ $=(1,2,(3,4,\cdots,l+2))\in\mathbb J^{\underline s_+}$; $A^{\tau_-}$ is the open subscheme of the preimage of $U^{\tau}_{12,m}$  under \begin{equation*}
\widetilde{\mathcal R}^{\underline s_-}:\mathcal Q^{\underline s_-}_{m+2}\rightarrow G(2,m+2),
\end{equation*}
associated to the size vector $\underline s_-:=(1,1,\cdots,1)\in\mathbb Z^{m+2}$ and the index $\tau_-:=(1,2,(j^-_1,j^-_2,\cdots,j^-_m))$ $=(1,2,(3,4,\cdots,m+2)) \in\mathbb J^{\underline s_-}$.

Notice that the ideal sheaves
$\left(\widetilde{\mathcal R}^{{\underline s}}\right)^{-1}\mathscr I^{\underline s}_{i_1i_2i_3i_4}\cdot\mathcal O_{A^{\tau}}$ with indices in (A0)
are generated by certain polynomials in  $a_{j^+_1},a_{j^+_2},\cdots,a_{j^+_l}$ and $\overrightarrow H$. Applying Lemma \ref{fac} to $A^{\tau_+}$, we can conclude that $\left({\mathcal R}_{A0}\right)^{-1}\left(A^{\tau}\right)$ is covered by open subschemes $B^{\lambda^+}\cong{\rm Spec}\,\mathbb Z\left[\overrightarrow A,\overrightarrow  E_+,\overrightarrow  \Xi\right]$ as follows. Define a morphism $\Sigma^{\lambda^+}:{\rm Spec}\,\mathbb Z\left[\overrightarrow  A,\overrightarrow  E_+,\overrightarrow  \Xi\right]\rightarrow U_{12}$ by
\begin{equation*}
\left\{\begin{aligned}
&\,x_{1j^+_{\alpha}}\mapsto a_{j^+_{\alpha}},\,\,\,\,\,\,\,\,\,\,\,\,\,\,\,\,\,\,\,\,\,\,\,\,\,\,\,\,\,\,\,\,\,\,\,\,\,\,\,{\rm for}\,\,1\leq {\alpha}\leq l\\
&\,x_{2j^+_{\alpha}}\mapsto a_{j^+_{\alpha}}\cdot\prod\nolimits_{\gamma=1}^{\lambda^+(\alpha)}\epsilon^+_{\gamma},\,\,\,\,\,\,\,\,{\rm for}\,\,1\leq \alpha\leq l\\
&\,x_{1j^-_{\beta}}\mapsto a_{j^-_{\beta}}\cdot\xi_{j^-_{\beta}},\,\,\,\,\,\,\,\,\,\,\,\,\,\,\,\,\,\,\,\,\,\,\,\,\,\,\,{\rm for}\,\,1\leq {\beta}\leq m\\
&\,x_{2j^-_{\beta}}\mapsto a_{j^-_{\beta}},\,\,\,\,\,\,\,\,\,\,\,\,\,\,\,\,\,\,\,\,\,\,\,\,\,\,\,\,\,\,\,\,\,\,\,\,\,\,\,{\rm for}\,\,1\leq {\beta}\leq m\\
\end{aligned}\right.\,\,\,.
\end{equation*}
We can define a locally closed embedding \begin{equation*}
\Omega^{\lambda^+}:{\rm Spec}\,\mathbb Z\left[\overrightarrow  A,\overrightarrow  E_+,\overrightarrow  \Xi\right]\rightarrow \mathbb P^{N_{2,n}}\times\prod\nolimits_{t=1}^N\mathbb P^{N^{\underline s}_{t}}\times\prod\nolimits_{\substack{1\leq\alpha_1<\alpha_2\leq l}}\mathbb P^{N_{12j^+_{\alpha_1}j^+_{\alpha_2}}}    
\end{equation*} by extending the rational map $\left(\widetilde{\mathcal K}^{\underline s}\circ\Sigma^{\lambda^+},\left(\cdots,F_{12j^+_{\alpha_1}j^+_{\alpha_2}}\circ e\circ \Sigma^{\lambda^+},\cdots\right)\right)$. Let $B^{\lambda^+}$ be the image of ${\rm Spec}\,\mathbb Z\left[\overrightarrow  A,\overrightarrow  E_+,\overrightarrow  \Xi\right]$ under $\Omega^{\lambda^+}$, which is equipped with the reduced scheme structure as a locally closed subscheme of $\mathbb P^{N_{2,n}}\times\prod\nolimits_{t=1}^N\mathbb P^{N^{\underline s}_{t}}\times\prod\nolimits_{\substack{1\leq\alpha_1<\alpha_2\leq l}}\mathbb P^{N_{12j^+_{\alpha_1}j^+_{\alpha_2}}}$.  
 
Similarly, the ideal sheaves
$\left(\widetilde{\mathcal R}^{{\underline s}}\circ\mathcal R_{A0}\right)^{-1}\mathscr I^{\underline s}_{i_1i_2i_3i_4}\cdot\mathcal O_{B^{\lambda^+}}$ with indices in (B0)
are generated by certain polynomials in  $a_{j^-_1},\cdots,a_{j^-_m}$, $\overrightarrow \Xi$. Hence, $\left({\mathcal R}_{A0B0}\right)^{-1}\left(B^{\lambda^+}\right)\cong{\rm Spec}\,\mathbb Z\left[\overrightarrow  A,\overrightarrow  E_+\right]\times \mathbb V$, where $\mathbb V$ is the blow-up of $A^{\tau_-}\subset\mathcal Q^{\underline s_-}_{m+2}$ with respect to the product of ideal sheaves
\begin{equation*}
\prod\nolimits_{1\leq \beta_1\leq \beta_2\leq m}\left(\widetilde{\mathcal R}^{{\underline s}_-}\right)^{-1}\mathscr I^{\underline s_-}_{12\beta_1\beta_2}\cdot\mathcal O_{A^{\tau_-}}.
\end{equation*}
By Remark \ref{fac2}, we can conclude that $\left({\mathcal R}_{A0B0}\right)^{-1}\left(B^{\lambda^+}\right)$ is covered by open subschemes $B^{\lambda^+,\lambda^-}$ for all permutations $\lambda^-$ of $\{1,2,\cdots,m\}$.

We complete the proof of Claim \ref{df}.\,\,\,\,$\endpf$

\begin{claim}\label{gi}
For any permutations  $\lambda^+$ of $\{1,2,\cdots,l\}$ and $\lambda^-$ of $\{1,2,\cdots,m\}$, the restriction of  $\mathcal R_{ABC}$ to
$\left(\mathcal R_{AB}\circ\mathcal R_{ABC}\right)^{-1}\left(B^{\lambda^+,\lambda^-}\cap \left({\mathcal R}_{A0}\circ{\mathcal R}_{A0B0}\right)^{-1}\left({\mathring A^{\tau}}\right)\right)$ is an isomorphism.
\end{claim}

{\bf\noindent Proof of Claim \ref{gi}.} It suffices to show that, for any index $(i_1,i_2,i_3,i_4)$ in (C1)-(C8), the restriction  to 
$\left(\mathcal R_{AB}\right)^{-1}\left(B^{\lambda^+,\lambda^-}\cap \left({\mathcal R}_{A0}\circ{\mathcal R}_{A0B0}\right)^{-1}\left({\mathring A^{\tau}}\right)\right)$ of the inverse image ideal sheaf
\begin{equation*}
\left(\widetilde{\mathcal R}^{\underline s}\circ{\mathcal R}_{A0}\circ{\mathcal R}_{A0B0}\circ\mathcal R_{AB}\right)^{-1}\mathscr I^{\underline s}_{i_1i_2i_3i_4}\cdot\mathcal O_{\left(\mathcal R_{AB}\right)^{-1}\left(B^{\lambda^+,\lambda^-}\right)}    
\end{equation*} is invertible.
The argument is on a case by case basis.

For any index in  (C1), $\left(\widetilde{\mathcal R}^{\underline s}\circ{\mathcal R}_{A0}\circ{\mathcal R}_{A0B0}\right)^{-1}\mathscr I^{\underline s}_{12i_3i_4}\cdot\mathcal O_{B^{\lambda^+,\lambda^-}}$ is generated by $a_{j^+_{\alpha_1}}a_{j^-_{\beta_1}}$.

For any index in (C2), $\left(\widetilde{\mathcal R}^{\underline s}\circ{\mathcal R}_{A0}\circ{\mathcal R}_{A0B0}\right)^{-1}\mathscr I^{\underline s}_{1i_2i_3i_4}\cdot\mathcal O_{B^{\lambda^+,\lambda^-}}$  is generated by
\begin{equation*}
\begin{split}
&\left\{a_{j^+_{\alpha_1}}a_{j^+_{\alpha_2}}a_{j^-_{\beta_1}}\epsilon^+_{1}\left(\prod\nolimits_{\gamma=2}^{\lambda^+\left(\alpha_{1}\right)}\epsilon^+_{\gamma}\right)\left(1-\left(\epsilon^+_{1}\epsilon^-_1\right)\prod\nolimits_{\gamma=2}^{\lambda^+\left(\alpha_{2}\right)}\epsilon^+_{\gamma}\cdot\prod\nolimits_{\gamma=2}^{\lambda^-\left(\beta_{1}\right)}\epsilon^-_{\gamma}\right)\right.,\\
&\,\,\,\,\,\,\,\,\,\,\,\,\,\left.a_{j^+_{\alpha_1}}a_{j^+_{\alpha_2}}a_{j^-_{\beta_1}}\epsilon^+_1\left(\prod\nolimits_{\gamma=2}^{\lambda^+\left(\alpha_{2}\right)}\epsilon^+_{\gamma}\right)\left(1-\left(\epsilon^+_{1}\epsilon^-_1\right)\prod\nolimits_{\gamma=2}^{\lambda^+\left(\alpha_{1}\right)}\epsilon^+_{\gamma}\cdot\prod\nolimits_{\gamma=2}^{\lambda^-\left(\beta_{1}\right)}\epsilon^-_{\gamma}\right)\right\}.\\ 
\end{split}
\end{equation*}
Notice that the following factors are invertible in $\left({\mathcal R}_{A0}\circ{\mathcal R}_{A0B0}\right)^{-1}\left({\mathring A^{\tau}}\right)$ by (\ref{punc}),
\begin{equation*}
1-\left(\epsilon^+_{1}\epsilon^-_1\right)\prod\nolimits_{\gamma=2}^{\lambda^+\left(\alpha_{2}\right)}\epsilon^+_{\gamma}\cdot\prod\nolimits_{\gamma=2}^{\lambda^-\left(\beta_{1}\right)}\epsilon^-_{\gamma},\,\,\,\,\,\,\,\,\,\,\,\,1-\left(\epsilon^+_{1}\epsilon^-_1\right)\prod\nolimits_{\gamma=2}^{\lambda^+\left(\alpha_{1}\right)}\epsilon^+_{\gamma}\cdot\prod\nolimits_{\gamma=2}^{\lambda^-\left(\beta_{1}\right)}\epsilon^-_{\gamma}.    \end{equation*}
When restricted to 
$B^{\lambda^+,\lambda^-}\cap \left({\mathcal R}_{A0}\circ{\mathcal R}_{A0B0}\right)^{-1}\left({\mathring A^{\tau}}\right)$, the ideal sheaf is generated by
\begin{equation*}
a_{j^+_{\alpha_1}}a_{j^+_{\alpha_2}}a_{j^-_{\beta_1}}\epsilon^+_1\left(\prod\nolimits_{\gamma=1}^{\min\{\lambda^+\left(\alpha_{1}\right),\lambda^+\left(\alpha_{2}\right)\}}\epsilon^+_{\gamma}\right).
\end{equation*}
The same holds for indices in  (C3).  

For any index in  (C4), $\left(\widetilde{\mathcal R}^{\underline s}\circ{\mathcal R}_{A0}\circ{\mathcal R}_{A0B0}\right)^{-1}\mathscr I^{\underline s}_{1i_2i_3i_4}\cdot\mathcal O_{B^{\lambda^+,\lambda^-}}$  is generated by
\begin{equation*}
\left\{a_{j^+_{\alpha_1}}a_{j^-_{\beta_1}}a_{j^-_{\beta_2}}\left(\prod\nolimits_{\gamma=1}^{\lambda^+\left(\alpha_{1}\right)}\epsilon^+_{\gamma}\right)\left(\prod\nolimits_{\gamma=1}^{\lambda^-\left(\beta_{1}\right)}\epsilon^-_{\gamma}\right),\,\,\,a_{j^+_{\alpha_1}}a_{j^-_{\beta_1}}a_{j^-_{\beta_2}}\left(\prod\nolimits_{\gamma=1}^{\lambda^+\left(\alpha_{1}\right)}\epsilon^+_{\gamma}\right)\left(\prod\nolimits_{\gamma=1}^{\lambda^-\left(\beta_{2}\right)}\epsilon^-_{\gamma}\right)\right\}.  
\end{equation*}
Computation yields that the ideal sheaf is generated by
\begin{equation*}
a_{j^+_{\alpha_1}}a_{j^-_{\beta_1}}a_{j^-_{\beta_2}}\left(\prod\nolimits_{\gamma=1}^{\lambda^+\left(\alpha_{1}\right)}\epsilon^+_{\gamma}\right)\left(\prod\nolimits_{\gamma=1}^{\min\{\lambda^-\left(\beta_{1}\right),\lambda^-\left(\beta_{2}\right)\}}\epsilon^-_{\gamma}\right).
\end{equation*}
The same holds for indices in (C5).

Similarly for any index in  (C6), $\left(\widetilde{\mathcal R}^{\underline s}\circ{\mathcal R}_{A0}\circ{\mathcal R}_{A0B0}\right)^{-1}\mathscr I^{\underline s}_{i_1i_2i_3i_4}\cdot\mathcal O_{B^{\lambda^+,\lambda^-}}$ is generated by
\begin{equation*}
\small
a_{j^+_{\alpha_1}}a_{j^+_{\alpha_2}}a_{j^-_{\beta_1}}a_{j^-_{\beta_2}}\left(1-\left(\epsilon^+_{1}\epsilon^-_1\right)\prod_{\gamma=2}^{\lambda^+\left(\alpha_{\sigma(1)}\right)}\epsilon^+_{\gamma}\cdot\prod_{\gamma=2}^{\lambda^-\left(\beta_1\right)}\epsilon^-_{\gamma}\right)\left(1-\left(\epsilon^+_{1}\epsilon^-_1\right)\prod_{\gamma=2}^{\lambda^+\left(\alpha_{\sigma(2)}\right)}\epsilon^+_{\gamma}\cdot\prod_{\gamma=2}^{\lambda^-\left(\beta_1\right)}\epsilon^-_{\gamma}\right), 
\end{equation*}
for all permutations $\sigma$ of $\{1,2\}$. We can conclude by (\ref{punc}) that the ideal sheaf is invertible when restricted to 
$B^{\lambda^+,\lambda^-}\cap \left({\mathcal R}_{A0}\circ{\mathcal R}_{A0B0}\right)^{-1}\left({\mathring A^{\tau}}\right)$.

For any index in (C7), $\left(\widetilde{\mathcal R}^{\underline s}\circ{\mathcal R}_{A0}\circ{\mathcal R}_{A0B0}\right)^{-1}\mathscr I^{\underline s}_{i_1i_2i_3i_4}\cdot\mathcal O_{B^{\lambda^+,\lambda^-}}$ is generated by
\begin{equation*}
a_{j^+_{\alpha_1}}a_{j^+_{\alpha_2}}a_{j^+_{\alpha_3}}a_{j^-_{\beta_1}}\epsilon^+_1\left(\prod_{\gamma=2}^{\lambda^+\left(\alpha_{\sigma(1)}\right)}\epsilon^+_{\gamma}-\prod_{\gamma=2}^{\lambda^+\left(\alpha_{\sigma(2)}\right)}\epsilon^+_{\gamma}\right)\left(1-\left(\epsilon^+_{1}\epsilon^-_1\right)\prod_{\gamma=2}^{\lambda^+\left(\alpha_{\sigma(3)}\right)}\epsilon^+_{\gamma}\cdot\prod_{\gamma=2}^{\lambda^-\left(\beta_1\right)}\epsilon^-_{\gamma}\right),   
\end{equation*}
for all permutations $\sigma$ of $\{1,2,3\}$. When restricted to 
$B^{\lambda^+,\lambda^-}\cap \left({\mathcal R}_{A0}\circ{\mathcal R}_{A0B0}\right)^{-1}\left({\mathring A^{\tau}}\right)$, we have by (\ref{punc}) that the ideal sheaf is a product of the invertible ideal sheaf generated by $a_{j^+_{\alpha_1}}a_{j^+_{\alpha_2}}a_{j^+_{\alpha_3}}a_{j^-_{\beta_1}}\epsilon^+_{1}$, and the ideal sheaf, which we denote by $\mathcal I_{j^+_{\alpha_1}j^+_{\alpha_2}j^+_{\alpha_3}}$, generated by
\begin{equation*}
\left\{\prod\nolimits_{\gamma=2}^{\lambda^+\left(\alpha_1\right)}\epsilon^+_{\gamma}-\prod\nolimits_{\gamma=2}^{\lambda^+\left(\alpha_{2}\right)}\epsilon^+_{\gamma},\,\,\,\,\,\,\prod\nolimits_{\gamma=2}^{\lambda^+\left(\alpha_1\right)}\epsilon^+_{\gamma}-\prod\nolimits_{\gamma=2}^{\lambda^+\left(\alpha_{3}\right)}\epsilon^+_{\gamma},\,\,\,\,\,\,\prod\nolimits_{\gamma=2}^{\lambda^+\left(\alpha_2\right)}\epsilon^+_{\gamma}-\prod\nolimits_{\gamma=2}^{\lambda^+\left(\alpha_{3}\right)}\epsilon^+_{\gamma}\right\}.
\end{equation*}
Computation yields that,  for any index in (A2),   $\left(\widetilde{\mathcal R}^{\underline s}\circ{\mathcal R}_{A0}\circ{\mathcal R}_{A0B0}\right)^{-1}\mathscr I^{\underline s}_{2j^+_{\alpha_1}j^+_{\alpha_2}j^+_{\alpha_3}}\cdot\mathcal O_{B^{\lambda^+,\lambda^-}}$,  $1\leq\alpha_1<\alpha_2<\alpha_3\leq l$,  is a product of $\mathcal I_{j^+_{\alpha_1}j^+_{\alpha_2}j^+_{\alpha_3}}$ and the invertible ideal sheaf generated by $a_{j^+_{\alpha_1}}a_{j^+_{\alpha_2}}a_{j^+_{\alpha_3}}\epsilon^+_{1}$. Therefore, $\left(\mathcal R_{AB}\right)^{-1}\left(\mathcal I_{j^+_{\alpha_1}j^+_{\alpha_2}j^+_{\alpha_3}}\right)\cdot\mathcal O_{\left({\mathcal R}_{A0}\circ{\mathcal R}_{A0B0}\circ\mathcal R_{AB}\right)^{-1}\left({\mathring A^{\tau}}\right)}$ is invertible,  when restricted to  $\left(\mathcal R_{AB}\right)^{-1}\left(B^{\lambda^+,\lambda^-}\cap \left({\mathcal R}_{A0}\circ{\mathcal R}_{A0B0}\right)^{-1}\left({\mathring A^{\tau}}\right)\right)$.  The same holds for indices in  (C8). 

We complete the proof of Claim \ref{gi}.\,\,\,\,$\endpf$
\smallskip

Now,  to prove Lemma \ref{uv}, it suffices to show that the following holds for all permutations $\lambda^+$ of $\{1,2,\cdots,l\}$ and $\lambda^-$ of $\{1,2,\cdots,m\}$.
\begin{enumerate}[label={\bf(P\arabic*)},ref={\bf P\arabic*}]
\item $\left(\mathcal R_{AB}\right)^{-1}\left(B^{\lambda^+,\lambda^-}\cap \left({\mathcal R}_{A0}\circ{\mathcal R}_{A0B0}\right)^{-1}\left({\mathring A^{\tau}}\right)\right)$
is a locally  closed subscheme of 
\begin{equation}\label{dudu}
\mathbb {P}^{N_{2,n}}\times\prod\limits_{t=1}^N\mathbb {P}^{N^{\underline s}_{t}}\times\left(\prod_{\substack{(i_1,i_2,i_3,i_4)\,\,{\rm in\,\,}\\{\rm (A0),\,(A1),\,(A2),\,(A3)}}}\mathbb P^{N_{i_1i_2i_3i_4}}\times\prod_{\substack{(i_1,i_2,i_3,i_4)\,\,{\rm in\,\,}\\{\rm (B0),\,(B1),\,(B2),\,(B3)}}}\mathbb P^{N_{i_1i_2i_3i_4}}\right),
\end{equation}
and, moreover, it is smooth over ${\rm Spec}\,\mathbb Z$.

\label{one}

\item Under the projection from $\mathbb {P}^{N_{2,n}}\times\prod\nolimits_{t=1}^N\mathbb {P}^{N^{\underline s}_{t}}\times\prod\nolimits_{\underline w\in C^{\underline s}}\mathbb {P}^{N^{\underline s}_{\underline w}}$  to $\prod\nolimits_{\underline w\in C^{\underline s}}\mathbb {P}^{N^{\underline s}_{\underline w}}$ the image of \begin{equation}\label{rabc}
\left(\mathcal R_{AB}\circ\mathcal R_{ABC}\right)^{-1}\left(B^{\lambda^+,\lambda^-}\cap \left({\mathcal R}_{A0}\circ{\mathcal R}_{A0B0}\right)^{-1}\left({\mathring A^{\tau}}\right)\right)    
\end{equation} has a reduced scheme structure as a locally closed subscheme of $\prod\nolimits_{\underline w\in C^{\underline s}}\mathbb {P}^{N^{\underline s}_{\underline w}}$; the subscheme is smooth over ${\rm Spec}\,\mathbb Z$ of relative dimension  $n-3$; the projection from (\ref{rabc}) to its image is a flat morphism.\label{two}

\end{enumerate}
\smallskip

Fix permutations $\lambda^+$ of $\{1,2,\cdots,l\}$ and $\lambda^-$ of $\{1,2,\cdots,m\}$. We first prove that

\begin{claim}\label{dddr}
$\left(\mathcal R_{AB}\right)^{-1}\left(B^{\lambda^+,\lambda^-}\right)$  is smooth over ${\rm Spec}\,\mathbb Z$. 
\end{claim}

{\noindent\bf Proof of Claim \ref{dddr}.} Let $\check{\mathcal R}_{l+2}$, $\ddot{\mathcal R}_{l+2}$, and $\check{\mathcal R}_{m+2}$, $\ddot{\mathcal R}_{m+2}$ be the blow-ups defined by setting $n=l+2$ and $n=m+2$ in (\ref{12a1}),  respectively. Let $B^{\lambda^+}$ and $B^{\lambda^-}$ be the open subschemes of $\left(\check{\mathcal R}_{l+2}\right)^{-1}\left(\mathcal Q_{l+2}^{\underline s_+}\right)$ and $\left(\check{\mathcal R}_{m+2}\right)^{-1}\left(\mathcal Q_{m+2}^{\underline s_-}\right)$ given by Lemma \ref{fac} and Remark \ref{fac2}, respectively. By the assumption that Proposition \ref{ms} holds for all $3\leq N^{\prime}<N$, it suffices to show that 
\begin{equation}\label{rab}
\left(\mathcal R_{AB}\right)^{-1}\left(B^{\lambda^+,\lambda^-}\right)\cong\left(\ddot{\mathcal R}_{l+2}\right)^{-1}\left(B^{\lambda^+}\right)\times\left(\ddot{\mathcal R}_{m+2}\right)^{-1}\left(B^{\lambda^-}\right).
\end{equation}

We denote $\mathbb A^{l}:={\rm Spec}\,\mathbb Z\left[a_{j^+_1},\cdots,a_{j^+_l}\right]$, $\mathbb A^{m}:={\rm Spec}\,\mathbb Z\left[a_{j^-_1},\cdots,a_{j^-_m}\right]$, $\mathbb A_{+}:={\rm Spec}\,\mathbb Z\left[\epsilon^+_1\right]$, $\mathbb A_{-}:={\rm Spec}\,\mathbb Z\left[\epsilon^-_1\right]$,
$\mathbb A^{\lambda^+}:={\rm Spec}\,\mathbb Z\left[\epsilon^+_2, \epsilon^+_3, \cdots,\epsilon^+_l\right]$, and  $\mathbb A^{\lambda^-}:={\rm Spec}\,\mathbb Z\left[\epsilon^-_2, \epsilon^-_3, \cdots,\epsilon^-_m\right]$. Recall the isomorphism $\Omega^{\lambda^+,\lambda^-}$ from \begin{equation}\label{aa}
\left(\mathbb A^{l}\times\mathbb A_{+}\times \mathbb A^{\lambda^+}\right)\times\left(\mathbb A^{m}\times\mathbb A_-\times \mathbb A_-^{\lambda^-}\right)    
\end{equation} to $B^{\lambda^+,\lambda^-}$ defined by (\ref{pabmath}), and let $\left(\Omega^{\lambda^+,\lambda^-},{\rm Id}\right)$ be its base change from
\begin{equation}\label{aapp}
\small
\left(\mathbb A^{l}\times\mathbb A_{+}\times \mathbb A^{\lambda^+}\right)\times\left(\mathbb A^{m}\times\mathbb A_-\times \mathbb A_-^{\lambda^-}\right)\times \prod\nolimits_{\substack{(i_1,i_2,i_3,i_4)\,\,{\rm in}\\{\rm (A1),\,(A2),\,(A3)}}}\mathbb P^{N_{i_1i_2i_3i_4}}\times\prod\nolimits_{\substack{(i_1,i_2,i_3,i_4)\,\,{\rm in}\\{\rm (B1),\,(B2),\,(B3)}}}\mathbb P^{N_{i_1i_2i_3i_4}}  
\end{equation}
to
\begin{equation}\label{bbpp}
B^{\lambda^+,\lambda^-}\times \prod\nolimits_{\substack{(i_1,i_2,i_3,i_4)\,\,{\rm in}\\{\rm (A1),\,(A2),\,(A3)}}}\mathbb P^{N_{i_1i_2i_3i_4}}\times\prod\nolimits_{\substack{(i_1,i_2,i_3,i_4)\,\,{\rm in}\\{\rm (B1),\,(B2),\,(B3)}}}\mathbb P^{N_{i_1i_2i_3i_4}}.    
\end{equation} Denote by ${\rm pr}_A$ (resp. ${\rm pr}_B$) the projection from (\ref{aa}) to $\mathbb A^{l}\times\mathbb A_{+}\times \mathbb A^{\lambda^+}$ (resp. $\mathbb A^{m}\times\mathbb A_{-}\times \mathbb A^{\lambda^-}$).

We define a rational map $\mathcal K_1:=({\rm Id},\mathcal K_{1A}\circ{\rm pr}_A,\mathcal K_{1B}\circ{\rm pr}_B)$ from (\ref{aa}) to (\ref{aapp}) as follows. Define 
\begin{equation*}
\mathcal K_{1A}:\mathbb A^{l}\times\mathbb A_{+}\times \mathbb A^{\lambda^+}\dashrightarrow \prod\nolimits_{\substack{(i_1,i_2,i_3,i_4)\,\,{\rm in}\\{\rm (A1),\,(A2),\,(A3)}}}\mathbb P^{N_{i_1i_2i_3i_4}}   
\end{equation*} by $\mathcal K_{1A}:=\left(\cdots, F^A_{i_1i_2i_3i_4},\cdots\right)$, where for $(i_1,i_2,i_3,i_4)$ in  (A1), \begin{equation}\label{a1}
\begin{split}
F^A_{i_1i_2i_3i_4}:=&\left[\left(\prod_{\gamma=2}^{\lambda^+\left(\alpha_1\right)}\epsilon^+_{\gamma}\right)\left(\prod_{\gamma=2}^{\lambda^+\left(\alpha_2\right)}\epsilon^+_{\gamma}-\prod_{\gamma=2}^{\lambda^+\left(\alpha_3\right)}\epsilon^+_{\gamma}\right),\left(\prod_{\gamma=2}^{\lambda^+\left(\alpha_2\right)}\epsilon^+_{\gamma}\right)\left(\prod_{\gamma=2}^{\lambda^+\left(\alpha_1\right)}\epsilon^+_{\gamma}-\prod_{\gamma=2}^{\lambda^+\left(\alpha_3\right)}\epsilon^+_{\gamma}\right),\right.\\
&\,\,\,\,\,\,\,\,\,\,\,\,\,\,\,\,\,\,\left.\left(\prod_{\gamma=2}^{\lambda^+\left(\alpha_3\right)}\epsilon^+_{\gamma}\right)\left(\prod_{\gamma=2}^{\lambda^+\left(\alpha_1\right)}\epsilon^+_{\gamma}-\prod_{\gamma=2}^{\lambda^+\left(\alpha_2\right)}\epsilon^+_{\gamma}\right)\right],
\end{split} \end{equation}
for $(i_1,i_2,i_3,i_4)$ in (A2),
\begin{equation}\label{a2}
\begin{split}
F^A_{i_1i_2i_3i_4}:=&\left[\left(\prod_{\gamma=2}^{\lambda^+\left(\alpha_2\right)}\epsilon^+_{\gamma}-\prod_{\gamma=2}^{\lambda^+\left(\alpha_3\right)}\epsilon^+_{\gamma}\right),\left(\prod_{\gamma=2}^{\lambda^+\left(\alpha_1\right)}\epsilon^+_{\gamma}-\prod_{\gamma=2}^{\lambda^+\left(\alpha_3\right)}\epsilon^+_{\gamma}\right),\left(\prod_{\gamma=2}^{\lambda^+\left(\alpha_1\right)}\epsilon^+_{\gamma}-\prod_{\gamma=2}^{\lambda^+\left(\alpha_2\right)}\epsilon^+_{\gamma}\right)\right],
\end{split}
\end{equation}
for $(i_1,i_2,i_3,i_4)$ in (A3),
\begin{equation}\label{a3}
\footnotesize
\begin{split}
&F^A_{i_1i_2i_3i_4}:=\left[\left(\prod_{\gamma=2}^{\lambda^+\left(\alpha_1\right)}\epsilon^+_{\gamma}-\prod_{\gamma=2}^{\lambda^+\left(\alpha_2\right)}\epsilon^+_{\gamma}\right)\left(\prod_{\gamma=2}^{\lambda^+\left(\alpha_3\right)}\epsilon^+_{\gamma}-\prod_{\gamma=2}^{\lambda^+\left(\alpha_4\right)}\epsilon^+_{\gamma}\right),\right.\\
&\left.\left(\prod_{\gamma=2}^{\lambda^+\left(\alpha_1\right)}\epsilon^+_{\gamma}-\prod_{\gamma=2}^{\lambda^+\left(\alpha_3\right)}\epsilon^+_{\gamma}\right)\left(\prod_{\gamma=2}^{\lambda^+\left(\alpha_2\right)}\epsilon^+_{\gamma}-\prod_{\gamma=2}^{\lambda^+\left(\alpha_4\right)}\epsilon^+_{\gamma}\right),\left(\prod_{\gamma=2}^{\lambda^+\left(\alpha_1\right)}\epsilon^+_{\gamma}-\prod_{\gamma=2}^{\lambda^+\left(\alpha_4\right)}\epsilon^+_{\gamma}\right)\left(\prod_{\gamma=2}^{\lambda^+\left(\alpha_2\right)}\epsilon^+_{\gamma}-\prod_{\gamma=2}^{\lambda^+\left(\alpha_3\right)}\epsilon^+_{\gamma}\right)\right].
\end{split}
\end{equation}
We can define $\mathcal K_{1B}:    \mathbb A^{m}\times\mathbb A_{-}\times \mathbb A^{\lambda^-}\dashrightarrow \prod\nolimits_{\substack{(i_1,i_2,i_3,i_4)\,\,{\rm in}\\{\rm (B1),\,(B2),\,(B3)}}}\mathbb P^{N_{i_1i_2i_3i_4}}$ in a similar way by cancelling common monomials in the homogeneous coordinates (\ref{Fw2}) for $F_{i_1i_2i_3i_4}$. Denote by $\mathcal T_{\rm model}$ the  closure of the image $\mathcal K_1\left({\rm dom}\left(\mathcal K_1\right)\right)$ in (\ref{aapp}). Then $\mathcal T_{\rm model}$ has a natural reduced scheme structure, for we can show that it is a blow-up of (\ref{aa}) by the same argument as that in \S \ref{pdm}.

Denote by $\mathcal K_2$ the inverse rational map of the restriction of the blow-up $\mathcal R_{AB}\left|_{\left(\mathcal R_{AB}\right)^{-1}\left(B^{\lambda^+,\lambda^-}\right)}\right.$, which is from $B^{\lambda^+,\lambda^-}$ to (\ref{bbpp}). We can write $\mathcal K_2=\left({\rm Id}, \mathcal K_{2A}, \mathcal K_{2B} \right)$ where 
\begin{equation*}
\begin{split}
&\mathcal K_{2A}:=\left(\cdots,\left.\left(F_{i_1i_2i_3i_4}\circ e\circ\widetilde  {\mathcal R}^{\underline s}\circ{\mathcal R}_{A0}\circ\mathcal R_{A0B0}\right)\right|_{B^{\lambda^+,\lambda^-}},\cdots\right)_{(i_1,i_2,i_3,i_4)\,\,{\rm in\,\,}  (A1), (A2), (A3)},\\ 
&\mathcal K_{2B}:=\left(\cdots,\left.\left(F_{i_1i_2i_3i_4}\circ e\circ\widetilde  {\mathcal R}^{\underline s}\circ{\mathcal R}_{A0}\circ\mathcal R_{A0B0}\right)\right|_{B^{\lambda^+,\lambda^-}},\cdots\right)_{(i_1,i_2,i_3,i_4)\,\,{\rm in\,\,}  (B1), (B2), (B3)}.\\    
\end{split} 
\end{equation*}
It is easy to verify that 
\begin{equation}\label{cd}
\left(\Omega^{\lambda^+,\lambda^-},{\rm Id}\right)\circ\mathcal K_1=\mathcal K_2\circ\Omega^{\lambda^+,\lambda^-},
\end{equation}
and hence diagram (\ref{cps}) commutes. 
Therefore, $\left(\mathcal R_{AB}\right)^{-1}\left(B^{\lambda^+,\lambda^-}\right)\cong  \mathcal T_{\rm model}$.

In what follows, we shall show that \begin{equation}\label{prod1}\mathcal T_{\rm model}\cong\left(\ddot{\mathcal R}_{l+2}\right)^{-1}\left(B^{\lambda^+}\right)\times\left(\ddot{\mathcal R}_{m+2}\right)^{-1}\left(B^{\lambda^-}\right). \end{equation}

Define a rational map $\mathcal K_{3A}:B^{\lambda^+}\dashrightarrow\prod\nolimits_{\substack{1\leq i_1<i_2<i_3<i_4\leq l+2\\(i_1,i_2)\neq (1,2)}}\mathbb P^{N_{i_1i_2i_3i_4}}$ by
\begin{equation}\label{k3a}
\begin{split}
&\mathcal K_{3A}:=\left(\cdots,\left.\left(F^l_{i_1i_2i_3i_4}\circ e_l\circ\widetilde  {\mathcal R}^{\underline s_+}\circ{\check{\mathcal R}}_{l+2}\right)\right|_{B^{\lambda^+}},\cdots\right)_{\substack{1\leq i_1<i_2<i_3<i_4\leq l+2\\(i_1,i_2)\neq (1,2)}},\\
\end{split}
\end{equation}
where $e_l:G(2,l+2)\hookrightarrow\mathbb P^{N_{2,l+2}}$ is the Pl\"ucker embedding, and
$F^l_{i_1i_2i_3i_4}$ is defined as in (\ref{Fw2}) by
\begin{equation*}
F^l_{i_1i_2i_3i_4}\left([\cdots ,z_I,\cdots]_{I\in\mathbb I_{2,l+2}}\right):=\left[z_{(i_1,i_2)}\cdot z_{(i_3,i_4)},\,\,\,\,z_{(i_1,i_3)}\cdot z_{(i_2,i_4)},\,\,\,\,z_{(i_1,i_4)}\cdot z_{(i_2,i_3)}\right].
\end{equation*} Similarly we can define a rational map $\mathcal K_{3B}:B^{\lambda^-}\dashrightarrow\prod\nolimits_{\substack{1\leq i_1<i_2<i_3<i_4\leq m+2\\(i_1,i_2)\neq (1,2)}}\mathbb P^{N_{i_1i_2i_3i_4}}$. By the same argument as that in \S\ref{pdm}, we can show that, as topological spaces, $\left(\ddot{\mathcal R}_{l+2}\right)^{-1}\left(B^{\lambda^+}\right)$, $\left(\ddot{\mathcal R}_{m+2}\right)^{-1}\left(B^{\lambda^-}\right)$ are respectively  the closures of $B^{\lambda^+}$,  $B^{\lambda^-}$ under the rational maps
\begin{equation}\label{k3ab}
\begin{split}
&\left({\rm Id},\mathcal K_{3A}\right):B^{\lambda^+}\dashrightarrow B^{\lambda^+}\times\prod\nolimits_{\substack{1\leq i_1<i_2<i_3<i_4\leq l+2\\(i_1,i_2)\neq (1,2)}}\mathbb P^{N_{i_1i_2i_3i_4}},\\
&\left({\rm Id},\mathcal K_{3B}\right):B^{\lambda^-}\dashrightarrow B^{\lambda^-}\times\prod\nolimits_{\substack{1\leq i_1<i_2<i_3<i_4\leq m+2\\(i_1,i_2)\neq (1,2)}}\mathbb P^{N_{i_1i_2i_3i_4}}.\\
\end{split} 
\end{equation}
Identifying $j^+_{\alpha}$ with $\alpha+2$ for $1\leq\alpha\leq l$, and $j^-_{\beta}$ with $\beta+2$ for $1\leq\beta\leq m$,
we have isomorphisms
\begin{equation*}
\begin{split}
&j_A:\prod\nolimits_{\substack{(i_1,i_2,i_3,i_4){\,\,\rm in}\\{\rm (A1),(A2),(A3)}}}\mathbb P^{N_{i_1i_2i_3i_4}}\rightarrow\prod\nolimits_{\substack{1\leq i_1<i_2<i_3<i_4\leq l+2\\(i_1,i_2)\neq (1,2)}}\mathbb P^{N_{i_1i_2i_3i_4}},\\
&j_B:\prod\nolimits_{\substack{(i_1,i_2,i_3,i_4){\rm \,\,in}\\{\rm (B1),(B2),(B3)}}}\mathbb P^{N_{i_1i_2i_3i_4}}\rightarrow\prod\nolimits_{\substack{1\leq i_1<i_2<i_3<i_4\leq m+2\\(i_1,i_2)\neq (1,2)}}\mathbb P^{N_{i_1i_2i_3i_4}}.
\end{split}
\end{equation*}
By the same token, we can derive that $\mathcal K_{3A}\circ\Omega^{\lambda^+}=j_A\circ\mathcal K_{1A}$ and $\mathcal K_{3B}\circ\Omega^{\lambda^-}=j_B\circ\mathcal K_{1B}$, and 
thus establish (\ref{prod1}). (All fit into commutative diagram (\ref{deccc}), where ${\rm Pr}_A$ and ${\rm Pr}_B$ respectively denote the projections from $B^{\lambda^+}\times B^{\lambda^-}$ to $B^{\lambda^+}$ and $B^{\lambda^-}$.)

We complete the proof of Claim \ref{dddr}.\,\,\,\,$\endpf$\smallskip 

\begin{landscape}\label{tttt}

\begin{equation}\label{cps}
\scriptsize
\begin{tikzcd}
&&\\
\left(\mathcal R_{AB}\right)^{-1}\left(B^{\lambda^+,\lambda^-}\right)\subset B^{\lambda^+,\lambda^-}\times\prod\limits_{\substack{(i_1,i_2,i_3,i_4)\\{\rm in\,\,  (A1)-(A3)}}}\mathbb P^{N_{i_1i_2i_3i_4}}\times\prod\limits_{\substack{(i_1,i_2,i_3,i_4)\\{\rm in\,\,(B1)-(B3)}}}\mathbb P^{N_{i_1i_2i_3i_4}}\arrow[rr, "\Gamma_2"]&&B^{\lambda^+,\lambda^-}\arrow[ll,dashrightarrow,"{\mathcal K_2=\left({\rm Id}, {\textcolor{blue}{\mathcal K_{2A}}},{\textcolor{red}{\mathcal K_{2B}}}\right)}"',bend right=30]\\
&&\\
\left(\mathbb A^{l}\times\mathbb A_+\times \mathbb A_+^{\lambda^+}\right)\times\left(\mathbb A^{m}\times\mathbb A_-\times \mathbb A_-^{\lambda^-}\right)\times\prod\limits_{\substack{(i_1,i_2,i_3,i_4)\\{\rm in\,\, (A1)-(A3)}}}\mathbb P^{N_{i_1i_2i_3i_4}}\times\prod\limits_{\substack{(i_1,i_2,i_3,i_4)\\{\rm in\,\,(B1)-(B3)}}}\mathbb P^{N_{i_1i_2i_3i_4}}\arrow[uu, "{\Huge{\left(\Omega^{\lambda^+,\lambda^-},\,\mathrm{Id}\right)}}","{\large\cong}"']\arrow[rr, "\Gamma_{1}"]&&\left(\mathbb A^{l}\times\mathbb A_+\times \mathbb A_+^{\lambda^+}\right)\times\left(\mathbb A^{m}\times\mathbb A_-\times \mathbb A_-^{\lambda^-}\right)\arrow[uu,"{\Omega^{\lambda^+,\lambda^-}}","{\large\cong}"']\arrow[ll,dashrightarrow,"{{\mathcal K_1=\left({\rm Id},\,{\textcolor{blue}{\mathcal K_{1A}\circ{\rm pr}_A}},\,{\textcolor{red}{\mathcal K_{1B}\circ{\rm pr}_B}}\right)}}",bend left=35]\\
\end{tikzcd}.
\end{equation}

$\vspace{1em}$

\begin{equation}\label{c3}
\footnotesize
\begin{tikzcd} 
\left(\mathbb A^{l}\times\mathbb A_+\times\mathbb A^{\lambda^+}\right)\times\prod\limits_{\substack{1\leq i_1<i_2<i_3<i_4\leq l+2\\(i_1,i_2)\neq (1,2)}}\mathbb P^{N_{i_1i_2i_3i_4}}\times \prod\limits_{\substack{3\leq i_3<i_4\leq l+2}}\mathbb P^{N_{12i_3i_4}}\arrow[r,"{{\rm Pr}_+}"]
&\left(\mathbb A^{l}\times\mathbb A_+\right)\times\prod\limits_{\substack{1\leq i_1<i_2<i_3<i_4\leq l+2\\(i_1,i_2)\neq (1,2)}}\mathbb P^{N_{i_1i_2i_3i_4}}\times \prod\limits_{\substack{3\leq i_3<i_4\leq l+2}}\mathbb P^{N_{12i_3i_4}}\arrow[dd]\\
&\\
\left(\mathbb A^{l}\times\mathbb A_+\times\mathbb A^{\lambda^+}\right)\times\prod\limits_{\substack{1\leq i_1<i_2<i_3<i_4\leq l+2\\(i_1,i_2)\neq (1,2)}}\mathbb P^{N_{i_1i_2i_3i_4}}\arrow[uu,"{\mathcal J=\left({\rm Id},\,\mathcal K_{3A0}\circ\Omega^{\lambda^+}\circ {\rm Pr}_1\right)}",shift left=3] \arrow[r]
&\prod\limits_{\substack{1\leq i_1<i_2<i_3<i_4\leq l+2\\(i_1,i_2)\neq (1,2)}}\mathbb P^{N_{i_1i_2i_3i_4}}\times \prod\limits_{\substack{3\leq i_3<i_4\leq l+2}}\mathbb P^{N_{12i_3i_4}} \\
\end{tikzcd}
\vspace{-20pt} 
\end{equation}

\end{landscape}

\begin{landscape}

$\vspace{2em}$

\begin{equation}\label{deccc}
\scriptsize
\begin{tikzcd}
\left(\mathcal R_{AB}\right)^{-1}\left(B^{\lambda^+,\lambda^-}\right)\arrow[dd]\arrow[dd,"\mathcal R_{AB}"'] &\mathcal T_{\rm model}\arrow[l,"{\Large\cong}"']\arrow[dd,"\Gamma_{AB}"']\arrow[r,"{\Large{\cong}}"] &\left(\ddot{\mathcal R}_{l+2}\right)^{-1}\left(B^{\lambda^+}\right)\times\left(\ddot{\mathcal R}_{m+2}\right)^{-1}\left(B^{\lambda^-}\right)\arrow[dd,"{\left(\ddot{\mathcal R}_{l+2},\ddot{\mathcal R}_{m+2}\right)}"]\arrow[dd] &\\ 
&&&\\
B^{\lambda^+,\lambda^-}\arrow[dd,"{\mathcal K_{2B0}}"',green]&\left(\mathbb A^{l}\times\mathbb A_+\times \mathbb A_+^{\lambda^+}\right)\times\left(\mathbb A^{m}\times\mathbb A_-\times \mathbb A_-^{\lambda^-}\right)\arrow[l,"{\Large{\cong}}","{\Omega^{\lambda^+,\lambda^-}}"'] \arrow[r,"{\Large{\cong}}"',"{\left(\Omega^{\lambda^+},\,\Omega^{\lambda^-}\right)}"]&B^{\lambda^+}\times B^{\lambda^-}\arrow[dd,"{\mathcal K_{3B0}\circ{\rm Pr}_B}"',green]\arrow[ddddr,"{\mathcal K_{3B}\circ{\rm Pr}_B}",dashrightarrow,red]&\\
&&&\\
\prod\limits_{\substack{(i_1,i_2,i_3,i_4)\\{\rm in\,\,(B0)}}}\mathbb P^{N_{i_1i_2i_3i_4}}
\arrow[rr,"\Large{\cong}"',"j_{B0}"]&&\prod\limits_{\substack{3\leq i_3<i_4\leq m+2}}\mathbb P^{N_{12i_3i_4}}&\\
&&&\\
&\prod\limits_{\substack{(i_1,i_2,i_3,i_4)\\{\rm in\,\, (B1)-(B3)}}}\mathbb P^{N_{i_1i_2i_3i_4}}\arrow[uuuul,"{\mathcal K_{2B}}"'{yshift=10pt,xshift=-10pt},dashleftarrow,red]\arrow[uuuu,"{\mathcal K_{1B}\circ{\rm pr}_B}"'{yshift=27pt,xshift=1pt},dashleftarrow,crossing over,red]\arrow[rr,"\Large{\cong}"',"j_B"]&&\prod\limits_{\substack{1\leq i_1<i_2<i_3<i_4\leq m+2\\(i_1,i_2)\neq(1,2)}}\mathbb P^{N_{i_1i_2i_3i_4}}\\
&&&\\
\prod\limits_{\substack{(i_1,i_2,i_3,i_4)\\{\rm in\,\, (A0)}}}\mathbb P^{N_{i_1i_2i_3i_4}}\arrow[uuuuuu,"{\mathcal K_{2A0}}",leftarrow,xshift=1.3ex,orange]\arrow[rr,"\Large{\cong}"',"j_{A0}"]
&&\prod\limits_{\substack{3\leq i_3<i_4\leq l+2}}\mathbb P^{N_{12i_3i_4}}\arrow[uuuuuu,"{\mathcal K_{3A0}\circ{\rm Pr}_A}",leftarrow,xshift=1.8ex,orange]&\\
&&&\\
&\prod\limits_{\substack{(i_1,i_2,i_3,i_4)\\{\rm in\,\,(A1)-(A3)}}}\mathbb P^{N_{i_1i_2i_3i_4}}\arrow[rr,"\Large{\cong}"',"j_{A}"]\arrow[uuuuuuuul,"{\mathcal K_{2A}}",dashleftarrow,bend right=15,blue]\arrow[uuuuuuuu,"{\mathcal K_{1A}\circ{\rm pr}_A}"'{yshift=40pt},dashleftarrow,bend right=55,blue]&&\prod\limits_{\substack{1\leq i_1<i_2<i_3<i_4\leq l+2\\(i_1,i_2)\neq(1,2)}}\mathbb P^{N_{i_1i_2i_3i_4}}\arrow[uuuuuuuul,"{\mathcal K_{3A}\circ{\rm Pr}_A}"{yshift=25pt,xshift=-15pt},dashleftarrow,bend right=10,blue]&\\
\end{tikzcd} \end{equation}

\end{landscape}

We proceed to establish  (\ref{two}). 

It is clear that $\left(\ddot{\mathcal R}_{l+2}\right)^{-1}\left(B^{\lambda^+}\right)$ is isomorphic to a closed subscheme of \begin{equation}\label{small1}\left(\mathbb A^{l}\times\mathbb A_+\times\mathbb A^{\lambda^+}\right)\times\prod\nolimits_{\substack{1\leq i_1<i_2<i_3<i_4\leq l+2\\(i_1,i_2)\neq (1,2)}}\mathbb P^{N_{i_1i_2i_3i_4}},
\end{equation}
which as a topological space is  the  closure of the image in (\ref{small1}) under the rational map $\left({\rm Id},\mathcal K_{3A}\circ\Omega^{\lambda^+}\right)$ from $\mathbb A^{l}\times\mathbb A_+\times\mathbb A^{\lambda^+}$ to (\ref{small1}). 

Denote by ${\rm Pr}_1$ the projection from (\ref{small1}) to $\mathbb A^{l}\times\mathbb A_+\times\mathbb A^{\lambda^+}$. We can define a morphism
\begin{equation*}
\mathcal K_{3A0}:B^{\lambda^+}\longrightarrow\prod\nolimits_{\substack{3\leq i_3<i_4\leq l+2}}\mathbb P^{N_{12i_3i_4}}\cong \prod\nolimits_{\substack{(i_1,i_2,i_3,i_4)\\{\rm in\,\,(A0)}}}\mathbb P^{N_{i_1i_2i_3i_4}}    
\end{equation*}
by extending the rational map $\left(\cdots,\left.\left(F^l_{12i_3i_4}\circ e_l\circ\widetilde  {\mathcal R}^{\underline s_+}\circ{\check{\mathcal R}}_{l+2}\right)\right|_{B^{\lambda^+}},\cdots\right)_{\substack{3\leq i_3<i_4\leq l+2}}$. Then, there is an embedding from (\ref{small1}) to
\begin{equation}\label{big}
\left(\mathbb A^{l}\times\mathbb A_+\times\mathbb A^{\lambda^+}\right)\times\prod\nolimits_{\substack{1\leq i_1<i_2<i_3<i_4\leq l+2\\(i_1,i_2)\neq (1,2)}}\mathbb P^{N_{i_1i_2i_3i_4}}\times \prod\nolimits_{\substack{3\leq i_3<i_4\leq l+2}}\mathbb P^{N_{12i_3i_4}}
\end{equation}
defined by $\mathcal J:=\left({\rm Id},\mathcal K_{3A0}\circ\Omega^{\lambda^+}\circ {\rm Pr}_1\right)$.

Denote by ${\rm Pr}_{+}$ the projection from (\ref{big}) to
\begin{equation}\label{small2}
\left(\mathbb A^{l}\times\mathbb A_+\right)\times\prod\nolimits_{\substack{1\leq i_1<i_2<i_3<i_4\leq l+2\\(i_1,i_2)\neq (1,2)}}\mathbb P^{N_{i_1i_2i_3i_4}}\times \prod\nolimits_{\substack{3\leq i_3<i_4\leq l+2}}\mathbb P^{N_{12i_3i_4}}.  
\end{equation}
We claim that ${\rm Pr}_{+}\circ{\mathcal J}$ is a locally closed embedding. Without loss of generality, we may assume that $\lambda^+$ is the identity permutation. Write the homogeneous coordinates for $\mathbb P^{N_{i_1i_2i_3i_4}}$ as $\left[z_{i_1i_2i_3i_4},z_{i_1i_3i_2i_4},z_{i_1i_4i_2i_3}\right]$. Computation yields that
\begin{equation*}
\begin{split}
&\left({\rm Pr}_{+}\circ{\mathcal J}\right)\bigg(\left(a_3,a_4,\cdots,a_{l+2}\right),\left(\epsilon^+_1\right),\left(\epsilon^+_2,\epsilon^+_3,\cdots,\epsilon^+_{l}\right),\\
&\left.\,\,\,\,\,\,\,\,\,\,\,\,\,\,\,\,\,\,\,\,\,\,\,\,\,\,\,\,\,\,\,\,\,\,\,\,\,\,\,\,\,\,\,\,\,\,\,\,\,\,\,\,\,\,\,\,\,\left(\cdots,\left[z_{i_1i_2i_3i_4},z_{i_1i_3i_2i_4},z_{i_1i_4i_2i_3}\right],\cdots\right)_{\substack{1\leq i_1<i_2<i_3<i_4\leq l+2\\(i_1,i_2)\neq (1,2)}}\right)\\
&\,\,\,\,\,\,=\left(\left(a_3,a_4,\cdots,a_{l+2}\right),\left(\epsilon^+_1\right),\left(\cdots,\left[z_{i_1i_2i_3i_4},z_{i_1i_3i_2i_4},z_{i_1i_4i_2i_3}\right],\cdots\right)_{\substack{1\leq i_1<i_2<i_3<i_4\leq l+2\\(i_1,i_2)\neq (1,2)}},\right.\\
&\,\,\,\,\,\,\,\,\,\,\,\,\,\,\,\,\,\,\,\,\,\,\,\,\,\,\,\,\,\,\,\,\,\,\,\,\,\,\,\,\,\,\,\,\,\,\,\,\,\,\,\,\,\,\,\,\left.\left(\cdots,\left[\left(1-\prod\nolimits_{i_3-1}^{i_4-2}\epsilon^+_{\gamma}\right),1,\left(\prod\nolimits_{i_3-1}^{i_4-2}\epsilon^+_{\gamma}\right)\right],\cdots\right)_{3\leq i_3<i_4\leq l+2}\right).\\
\end{split}    
\end{equation*}
By taking $(i_3,i_4)=(3,4),(4,5),\cdots,(l+1,l+2)$ in \begin{equation}\label{i1i2i3}
\left[\left(1-\prod\nolimits_{i_3-1}^{i_4-2}\epsilon^+_{\gamma}\right),1,\left(\prod\nolimits_{i_3-1}^{i_4-2}\epsilon^+_{\gamma}\right)\right],    
\end{equation} we can conclude that ${\rm Pr}_{+}\circ{\mathcal J}$ is a locally closed embedding.

Now by projecting (\ref{small2}) to  \begin{equation}\label{small3}
\prod\nolimits_{\substack{1\leq i_1<i_2<i_3<i_4\leq l+2\\(i_1,i_2)\neq (1,2)}}\mathbb P^{N_{i_1i_2i_3i_4}}\times \prod\nolimits_{\substack{3\leq i_3<i_4\leq l+2}}\mathbb P^{N_{12i_3i_4}}\cong\prod\nolimits_{\substack{(i_1,i_2,i_3,i_4)\,\,{\rm in}\\{\rm (A0),\,(A1),\,(A2),\,(A3)}}}\mathbb P^{N_{i_1i_2i_3i_4}},   
\end{equation}
we can conclude that
\begin{equation}\label{rma}\left(\ddot{\mathcal R}_{l+2}\right)^{-1}\left(B^{\lambda^+}\right)\cong\mathbb A^{l}\times\mathbb A_+\times\mathbb M^{\lambda^+}.
\end{equation} Here $\mathbb M^{\lambda^+}$ is a locally closed subscheme of
(\ref{small3}) such that, as a topological space, it is the closure of the image in (\ref{small3}) under the rational map $\left(\mathcal K_{3A},\mathcal K_{3A0}\right)$. Similarly,  we can conclude that
\begin{equation}
\label{rmb} \left(\ddot{\mathcal R}_{m+2}\right)^{-1}\left(B^{\lambda^-}\right)\cong\mathbb A^{m}\times\mathbb A_-\times\mathbb M^{\lambda^-},   
\end{equation} where $\mathbb M^{\lambda^-}$ is a locally closed subscheme of 
\begin{equation}
\prod\nolimits_{\substack{1\leq i_1<i_2<i_3<i_4\leq m+2\\(i_1,i_2)\neq (1,2)}}\mathbb P^{N_{i_1i_2i_3i_4}}\times \prod\nolimits_{\substack{3\leq i_3<i_4\leq m+2}}\mathbb P^{N_{12i_3i_4}}\cong\prod\nolimits_{\substack{(i_1,i_2,i_3,i_4)\,\,{\rm in}\\{\rm (B0),\,(B1),\,(B2),\,(B3)}}}\mathbb P^{N_{i_1i_2i_3i_4}}.    
\end{equation} Moreover, according to the assumption that Proposition \ref{ms} holds for all $3\leq N^{\prime}<N$, we can derive that $\mathbb M^{\lambda^+}$,  $\mathbb M^{\lambda^-}$ are smooth over ${\rm Spec}\,\mathbb Z$ of relative dimension  $l-1$, $m-1$, respectively.

\medskip

Take integers $q_A\in\{1,2,\cdots,l\}$, $q_B\in\{1,2,\cdots,m\}$ such that $\lambda^+(q_A)=\lambda^-(q_B)=1$. Without loss of generality, we may assume that $j^+_{q_A}<j^-_{q_B}$. Computation yields that 
\begin{equation}\label{cf}
\begin{split}
&\left.\left(F_{12j^+_{q_A}j^-_{q_B}}\circ e\circ\widetilde  {\mathcal R}^{\underline s}\circ{\mathcal R}_{A0}\circ\mathcal R_{A0B0}\right)\right|_{B^{\lambda^+,\lambda^-}}\left(\overrightarrow  A,  \overrightarrow  E_+, \overrightarrow  E_-\right)\\ 
&\,\,\,=\left[a_{j^+_{q_A}}\cdot a_{j^-_{q_B}}-a_{j^+_{q_A}}\cdot a_{j^-_{q_B}}\cdot \epsilon_1^+\cdot \epsilon^-_1,\,\,\,\,-a_{j^+_{q_A}}\cdot a_{j^-_{q_B}}\cdot \epsilon_1^+\cdot \epsilon^-_1,\,\,\,\,-a_{j^+_{q_A}}\cdot a_{j^-_{q_B}}\right]\\
&\,\,\,=\left[1-\epsilon_1^+\cdot \epsilon^-_1,\,\,\,\,-\epsilon_1^+\cdot \epsilon^-_1,\,\,\,\,-1\right],\\
\end{split}
\end{equation}
which thus extends to a morphism. We can define a morphism $\mathcal L$ from $\left(\mathcal R_{AB}\right)^{-1}\left(B^{\lambda^+,\lambda^-}\right)$ to
\begin{equation}\label{1ab0}
U:=\prod\nolimits_{\substack{(i_1,i_2,i_3,i_4)\,\,{\rm in}\\{\rm (A0),\,(A1),\,(A2),\,(A3)}}}\mathbb P^{N_{i_1i_2i_3i_4}}\times\prod\nolimits_{\substack{(i_1,i_2,i_3,i_4)\,\,{\rm in}\\{\rm (B0),\,(B1),\,(B2),\,(B3)}}}\mathbb P^{N_{i_1i_2i_3i_4}}\times\mathbb P^{12j^+_{q_A}j^-_{q_B}}   
\end{equation}
by extending the rational maps $\left.\left(F_{i_1i_2i_3i_4}\circ e\circ\widetilde  {\mathcal R}^{\underline s}\circ{\mathcal R}_{A0}\circ\mathcal R_{A0B0}\circ\mathcal R_{AB}\right)\right|_{\left(\mathcal R_{AB}\right)^{-1}\left(B^{\lambda^+,\lambda^-}\right)}$. Composing isomorphisms (\ref{rab}), (\ref{rma}),  (\ref{rmb}), we can conclude by (\ref{cf}) that $\mathcal L$ is the same as the morphism $(l_1,l_2)$ from
\begin{equation}\label{pp+-}\left(\mathbb A^{l}\times\mathbb A_+\times \mathbb M^{\lambda^+}\right)\times\left(\mathbb A^{m}\times\mathbb A_-\times \mathbb M^{\lambda^-}\right)
\end{equation}
to $\mathbb M^{\lambda^+}\times\mathbb M^{\lambda^-}\times{\rm Spec}\,\mathbb Z[t]$, which is
defined as follows.
$l_1$ is the projection from (\ref{pp+-}) to $\mathbb M^{\lambda^+}\times\mathbb M^{\lambda^-}$; $l_2:=l_{t}\circ p$ where $p$ is the projection from (\ref{pp+-}) to $\mathbb A_+\times\mathbb A_-$, and  $l_{t}:\mathbb A_+\times\mathbb A_-\rightarrow{\rm Spec}\,\mathbb Z[t]$ is induced by the ring homomorphism $t\mapsto \epsilon_1^+\epsilon_1^-$. Then, it is clear that the following holds. 
\begin{enumerate}[label=$\bullet$]
\item The image of $\left(\mathcal R_{AB}\right)^{-1}\left(B^{\lambda^+,\lambda^-}\right)$ under the morphism $\mathcal L$ has a reduced scheme structure as a locally closed subscheme of (\ref{1ab0}), which we denote by $\mathbb M$.

\item $\mathbb M$ is isomorphic to $\mathbb M^{\lambda^+}\times\mathbb M^{\lambda^-}\times{\rm Spec}\,\mathbb Z[t]$, and hence it is smooth over ${\rm Spec}\,\mathbb Z$ of relative dimension  $n-3$.

\item The morphism induced by $\mathcal L$ from
$\left(\mathcal R_{AB}\right)^{-1}\left(B^{\lambda^+,\lambda^-}\right)$ to $\mathbb M$ is flat, which we also denote by $\mathcal L$ by a slight abuse of notation.
\end{enumerate}  

Recall that by Claim \ref{gi} the restrictions of the rational maps \begin{equation*}
\chi_{i_1i_2i_3i_4}:=\left.\left(F_{i_1i_2i_3i_4}\circ e\circ\widetilde  {\mathcal R}^{\underline s}\circ{\mathcal R}_{A0}\circ\mathcal R_{A0B0}\circ{\mathcal R}_{AB}\right)\right|_{\left(\mathcal R_{AB}\right)^{-1}\left(B^{\lambda^+,\lambda^-}\cap \left({\mathcal R}_{A0}\circ{\mathcal R}_{A0B0}\right)^{-1}\left({\mathring A^{\tau}}\right)\right)}    
\end{equation*} extends to a morphism for any $1\leq i_1<i_2<i_3<i_4\leq n$. Now, to prove (\ref{two}) it suffices to show that $\chi_{i_1i_2i_3i_4}$ for indices in (C1)-(C8) further factor through $\mathcal L:\left(\mathcal R_{AB}\right)^{-1}\left(B^{\lambda^+,\lambda^-}\right)\longrightarrow\mathbb M$.

We may identify $\left(\mathcal R_{AB}\right)^{-1}\left(B^{\lambda^+,\lambda^-}\right)$ with (\ref{pp+-}) via isomorphisms (\ref{rab}), (\ref{rma}), (\ref{rmb}), and identify $\mathbb M$ with $\mathbb M^{\lambda^+}\times\mathbb M^{\lambda^-}\times\mathbb A$. Without loss of generality, we assume that $\lambda^+$ and $\lambda^-$ are the identity permutations, and  $j^+_{\alpha_1}<j^+_{\alpha_2}<\cdots<j^+_{\alpha_l}<j^-_{\beta_1}<j^-_{\beta_2}<\cdots<j^-_{\beta_m}$. According to (\ref{i1i2i3}), $\chi_{12j^+_{\alpha_p}j^+_{\alpha_{p+1}}}$, $1\leq p\leq l-1$, is induced by the ring homomorphism
\begin{equation*}
z_{12j^+_{\alpha_p}j^+_{\alpha_{p+1}}}\left(z_{1j^+_{\alpha_p}2j^+_{\alpha_{p+1}}}\right)^{-1}\mapsto1-\epsilon^+_{p+1},\,\,\,\,z_{1j^+_{\alpha_{p+1}}2j^+_{\alpha_p}}\left(z_{1j^+_{\alpha_p}2j^+_{\alpha_{p+1}}}\right)^{-1}\mapsto\epsilon^+_{p+1},
\end{equation*}
and $\chi_{12j^-_{\beta_q}j^-_{\beta_{q+1}}}$,  $1\leq q\leq m-1$, is induced by the ring homomorphism
\begin{equation*}
z_{12j^-_{\beta_q}j^-_{\beta_{q+1}}}\left(z_{1j^-_{\beta_q}2j^-_{\beta_{q+1}}}\right)^{-1}\mapsto1-\epsilon^-_{q+1},\,\,\,\,z_{1j^-_{\beta_{q+1}}2j^-_{\beta_q}}\left(z_{1j^-_{\beta_q}2j^-_{\beta_{q+1}}}\right)^{-1}\mapsto\epsilon^-_{q+1},\\
\end{equation*}
where $\left[z_{i_1i_2i_3i_4},z_{i_1i_3i_2i_4},z_{i_1i_4i_2i_3}\right]$ are homogeneous coordinates for $\mathbb P^{N_{i_1i_2i_3i_4}}$. Then it suffices to show that for indices in (C1)-(C8), each component of $\chi_{i_1i_2i_3i_4}$ is a polynomial in  $\left(\epsilon^+_1\epsilon^-_1\right)$, $\epsilon^+_2,\epsilon^+_3,\cdots,\epsilon^+_l$, and $\epsilon^-_2,\epsilon^-_3,\cdots,\epsilon^-_m$. Computation yields that for $1\leq\alpha_1\leq l$ and $1\leq\beta_1\leq m$, $\chi_{12j^+_{\alpha_1}j^-_{\beta_1}}$ is induced by the ring homomorphism
\begin{equation*}
\left\{\begin{aligned}
&\frac{z_{12j^+_{\alpha_1}j^-_{\beta_{1}}}}{z_{1j^-_{\beta_{1}}2j^+_{\alpha_1}}}\mapsto\frac{a_{j^+_{\alpha_1}}a_{j^-_{\beta_1}}-a_{j^+_{\alpha_1}}a_{j^-_{\beta_1}}\left(\epsilon_1^+\epsilon^-_1\right)\prod\nolimits_{\gamma=2}^{\alpha_1}\epsilon^+_{\gamma}\prod\nolimits_{\gamma=2}^{\beta_1}\epsilon^-_{\gamma}}{-a_{j^+_{\alpha_1}}\cdot a_{j^-_{\beta_1}}}=-1+\left(\epsilon_1^+\epsilon^-_1\right)\prod\nolimits_{\gamma=2}^{\alpha_1}\epsilon^+_{\gamma}\prod\nolimits_{\gamma=2}^{\beta_1}\epsilon^-_{\gamma}\\ 
&\frac{z_{1j^+_{\alpha_1}2j^-_{\beta_{1}}}}{z_{1j^-_{\beta_{1}}2j^+_{\alpha_1}}}\mapsto\frac{-a_{j^+_{\alpha_1}}a_{j^-_{\beta_1}}\left(\epsilon_1^+\epsilon^-_1\right)\prod\nolimits_{\gamma=2}^{\alpha_1}\epsilon^+_{\gamma}\prod\nolimits_{\gamma=2}^{\beta_1}\epsilon^-_{\gamma}}{-a_{j^+_{\alpha_1}}\cdot a_{j^-_{\beta_1}}}=\left(\epsilon_1^+\epsilon^-_1\right)\prod\nolimits_{\gamma=2}^{\alpha_1}\epsilon^+_{\gamma}\prod\nolimits_{\gamma=2}^{\beta_1}\epsilon^-_{\gamma}\\
\end{aligned}.\right.
\end{equation*}
The same result holds for indices in (C2)-(C7), and we omit the details here for brevity.

We complete the proof of Lemma \ref{red}.\,\,\,\,\,\,$\endpf$

\subsection{Induction for Type (\ref{up})}

In this subsection, we shall prove that

\begin{lemma}\label{ud}
Suppose that  Proposition \ref{ms} holds for all integers $N^{\prime}$ such that $3\leq N^{\prime}<N$. Then Lemma \ref{loc} holds for all (truncated) coordinate charts of {\bf Type (\ref{up})}.
\end{lemma}
{\bf\noindent Proof of Lemma \ref{ud}.} Without loss of generality, we assume
$\tau=\left(1,2,(j_1^+,\cdots,j^+_{N-2})\right)$. 

Let $B^{\lambda}\subset\left(\check{\mathcal R}_N\right)^{-1}\left(A^{\tau}\right)$ be the open subscheme associated to any permutation $\lambda$ as in Lemma \ref{fac}.  Without loss of generality, we may assume that $\lambda$ is the identity permutation. 
Denote $\mathbb A_+:={\rm Spec}\,\mathbb Z\left[\epsilon^+_2,\epsilon^+_3,\cdots,\epsilon^+_{N-2}\right]$. Define a morphism $j:\mathbb A_+\rightarrow \mathcal {G}(2,n)$ by 
\begin{equation}\label{j0}
\left(\begin{matrix}
1&0&1&1&1&\cdots&1&\cdots&1\\
0&1&1&\epsilon^+_2&\epsilon^+_2\epsilon^+_3&\cdots&\epsilon^+_2\epsilon^+_3\cdots\epsilon^+_k&\cdots&\epsilon^+_2\epsilon^+_3\cdots\epsilon^+_{N-2}
\end{matrix}\right).  
\end{equation}
According to Definition \ref{mat}, we can define a blow-up $\mathcal R^{\underline s}_{j}:\mathbb M^{\underline s}_{j}\rightarrow \mathbb A_+$  with respect to the product of ideal sheaves $\prod\nolimits_{\underline w\in C^{\underline s}}\left(\pi\circ j\right)^{-1}\mathscr I^{\underline s}_{\underline w}\cdot\mathcal O_{\mathbb A_+}$. Similarly to  (\ref{rma}), we can prove that
\begin{equation}\label{prod31}
\left(\ddot{\mathcal R}_{N}\right)^{-1}\left(B^{\lambda}\right)\cong {\rm Spec}\,\mathbb Z\left[\overrightarrow A\right]\times{\rm Spec}\,\mathbb Z\left[\epsilon^+_1\right]\times\mathbb M^{\underline s}_{j}. 
\end{equation}

As in \S \ref{pdm}, we will interpret $\mathcal R^{\underline s}_{j}$ as a sequence of blow-ups. We partition (\ref{i1234}) as follows.
\begin{enumerate}[label={(E\arabic*)},ref={E\arabic*}]
\setcounter{enumi}{-1}
\item $\{(1,2,i_3,i_4)|\,\,3\leq i_3<i_4\leq n\}$.

\item $\{(1,i_2,i_3,i_4)|\,\,3\leq i_2<i_3<i_4\leq n\}$.

\item $\{(i_1,i_2,i_3,i_4)|\,\,3\leq i_1<i_2<i_3<i_4\leq n\}$.

\item $\{(2,i_2,i_3,i_4)|\,\,3\leq i_2<i_3<i_4\leq n\}$.

\end{enumerate}
Define a total order $\sigma$ so that the $4$-tuples in an upper subset are strictly smaller than the ones in a lower subset. 
For $0\leq k<\frac{n!}{24(n-4)!}$, we inductively define $\check\psi^{\sigma}_{k+1}:\check W^{\sigma}_{k+1}\rightarrow \check W^{\sigma}_{k}$ to be the blow-up of $\check W^{\sigma}_{k}$ with respect to $\left(\check\psi^{\sigma}_{1}\circ\check\psi^{\sigma}_{2}\circ\cdots\circ\check\psi^{\sigma}_{k}\right)^{-1}\left(\left(\pi\circ j\right)^{-1}\mathscr I^{\underline s}_{\underline w}\cdot\mathcal O_{\mathbb A_+}\right)\cdot\mathcal O_{\check W^{\sigma}_{k}}$, where $\check W^{\sigma}_{0}=\mathbb A_+$, and $\mathscr I^{\underline s}_{\underline w}$ are identified with  $\mathscr I^{\underline s}_{i_1i_2i_3i_4}$ as in \S \ref{FAC}.

For simplicity, we denote 
\begin{equation*}
\begin{split}
&\check {\mathcal R}_{\varheart}:=\check\psi^{\sigma}_{1}\circ\check\psi^{\sigma}_{2}\circ\cdots\check\psi^{\sigma}_{\frac{(n-2)!}{24\cdot(n-6)!}+\frac{(n-2)!}{6\cdot(n-5)!}+\frac{(n-2)!}{2\cdot(n-4)!}},\\ 
&\check{\mathcal R}_{\vardiamond}:=\check\psi^{\sigma}_{\frac{(n-2)!}{24\cdot(n-6)!}+\frac{(n-2)!}{6\cdot(n-5)!}+\frac{(n-2)!}{2\cdot(n-4)!}+1}\circ\check\psi^{\sigma}_{\frac{(n-2)!}{24\cdot(n-6)!}+\frac{(n-2)!}{6\cdot(n-5)!}+\frac{(n-2)!}{2\cdot(n-4)!}+2}\circ\cdots\circ\check\psi^{\sigma}_{\frac{n!}{24\cdot(n-4)!}}.\\
\end{split}
\end{equation*}
It is clear that $\check {\mathcal R}_{\vardiamond}$ is the sequence of blow-ups with respect to the (pull-backs of) ideal sheaves $\mathscr I^{\underline s}_{i_1i_2i_3i_4}$ with indices in (E3). Similarly to Claim \ref{gi}, we can show that $\check {\mathcal R}_{\vardiamond}$ is an isomorphism.

Define a morphism  $L:\mathbb A_+\rightarrow \mathcal {G}(2,n-1)$ by
\begin{equation}\label{j1}
\left(\begin{matrix}
1&1&1&1&\cdots&1\\
0&1&\epsilon^+_2&\epsilon^+_2\epsilon^+_3&\cdots&\epsilon^+_2\epsilon^+_3\cdots\epsilon^+_{N-2}
\end{matrix}\right).  
\end{equation}  
Denote
\begin{equation*}
\underline s^*=(1,1,\cdots,1)\in\mathbb Z^{N-1}. 
\end{equation*}
Then, it is clear that
\begin{equation}\label{prod32}\mathbb M^{\underline s}_{j}\cong\left(\check{\mathcal R}_{\varheart}\right)^{-1}\left(\mathbb A_+\right)\cong\mathbb M^{\underline s^*}_L,
\end{equation} where ${\mathcal R}^{\underline s^*}_{L}:\mathbb M^{\underline s^*}_{L}\rightarrow\mathbb A_+$ is the blow-up given by Definition \ref{mat}.

Define an open cover $\left\{C^k_N\right\}_{k=2}^{N-1}$ of $\mathbb A_+$ by 
\begin{equation}\label{c1n}
C^k_N:=\left\{\mathfrak p\in \mathbb A_+\left| \epsilon^+_{\gamma}\notin \mathfrak p,\,\forall\,2\leq\gamma\leq k-1,\,\,{\rm and}\,\,1-\prod\nolimits_{\gamma=j_1}^{j_2}\epsilon^+_{\gamma}\notin \mathfrak p,\,\forall\,2\leq j_1\leq k\leq j_2\leq N-2\right.\right\}. 
\end{equation}
Then, to prove Lemma \ref{ud} it suffices to show that
\begin{enumerate}[label={\bf(Q1)},ref={\bf Q1}]
\item  for each $2\leq k\leq N-1$,  $\left({\mathcal R}^{\underline s^*}_{L}\right)^{-1}\left(C^k_N\right)$ is a locally closed subscheme of $\prod\nolimits_{\underline w\in C^{\underline s}}\mathbb {P}^{N^{\underline s}_{\underline w}}$, and is smooth over ${\rm Spec}\,\mathbb Z$ of relative dimension  $n-3$.
\label{p11}
\end{enumerate}

We proceed to prove {\bf (\ref{p11})} holds on a case by case basis.

\medskip

{\noindent\bf Case I ($k=N-1$).} Define a morphism $J$ from ${\rm Spec}\,\mathbb Z\left[\xi_3,\xi_4,\cdots,\xi_{N-1}\right]$ to $\mathcal {G}(2,n-1)$ by
\begin{equation}\label{j2}
\left(\begin{matrix}
1&0&\xi_3&\xi_4&\cdots&\xi_{n-1}\\
0&1&1&1&\cdots&1
\end{matrix}\right).
\end{equation}
Define an open subscheme $D$ of ${\rm Spec}\,\mathbb Z\left[\xi_3,\xi_4,\cdots,\xi_{N-1}\right]$ by \begin{equation*}
D:=\left\{\mathfrak p\in {\rm Spec}\,\mathbb Z\left[\xi_3,\xi_4,\cdots,\xi_{N-1}\right]\big|\,\xi_k+1\notin\mathfrak p,\,\,\forall\, 3\leq k\leq N-1 \right\}.
\end{equation*}
Applying the change of coordinates
\begin{equation*}
\begin{split}
\xi_k:=\left(\epsilon_2\epsilon_3\cdots\epsilon_{k-1}\right)^{-1}-1,\,\,\,\,\,\,\forall\,3\leq k\leq N-1,
\end{split}   \end{equation*}
we can derive that
\begin{equation*}
\begin{split}
D&={\rm Spec}\,\mathbb Z\left[\xi_3+1,\cdots,\xi_{N-1}+1,\left(\xi_3+1\right)^{-1},\cdots,\left(\xi_{N-1}+1\right)^{-1}\right]\\
&\cong{\rm Spec}\,\mathbb Z\left[\epsilon_2,\cdots,\epsilon_{N-2},\epsilon_2^{-1},\cdots,\epsilon_{N-2}^{-1}\right]=C^{N-1}_N.\\
\end{split}   
\end{equation*}
Combining (\ref{j1}), (\ref{j2}) and the fact that
\begin{equation*}\begin{split}
&\left(\begin{matrix}1&1&1&\cdots&1\\0&1&\epsilon^+_2&\cdots&\epsilon^+_2\epsilon^+_3\cdots\epsilon^+_{N-2}\end{matrix}\right)=\left(\begin{matrix}1&1\\0&1\end{matrix}\right)\cdot\left(\begin{matrix}1&0&\frac{1}{\epsilon^+_2}-1&\cdots&\frac{1}{\epsilon^+_2\epsilon^+_3\cdots\epsilon^+_{N-2}}-1\\0&1&1&\cdots&1\\\end{matrix}\right)\\
&\,\,\,\,\,\,\,\,\,\,\,\,\,\,\,\,\,\,\,\,\,\,\,\,\cdot{\rm diag}\left(1,1,\epsilon^+_2,\epsilon^+_3,\cdots,\epsilon^+_2\epsilon^+_3\cdots\epsilon^+_{N-2}\right),\\
\end{split}\end{equation*}
we can conclude that
\begin{equation*}
\prod\nolimits_{\underline w\in C^{\underline s}}\left(\pi\circ L\right)^{-1}\mathscr I^{\underline s}_{\underline w}\cdot\mathcal O_{C^{N-1}_{N}}\cong a\cdot\prod\nolimits_{\underline w\in C^{\underline s}}\left(\pi\circ J\right)^{-1}\mathscr I^{\underline s}_{\underline w}\cdot\mathcal O_{D},
\end{equation*}
where $a$ is a non-zero element of $\Gamma\left(D,\mathcal O_{D}\right)$. Then, $\left({\mathcal R}^{\underline s^*}_{L}\right)^{-1}\left(C^{N-1}_N\right)$ is isomorphic to an open subscheme of $\mathbb M^{\underline s^*}_{J}$.

Now we consider the restriction of the blow-up $\overline{\mathcal R}_{n-1}^{\underline s^*}:\mathcal T_{n-1}^{\underline s^*}\rightarrow\mathcal Q_{n-1}^{\underline s^*}$ to $\left(\overline{\mathcal R}_{n-1}^{\underline s^*}\right)^{-1}\left(A^{\tau_-}\right)$, where $\tau_-=\left(1,2,(j_1^-,\cdots,j^-_{N-3})\right)\in\mathbb J^{\underline s^*}$. 
Similarly to the proof of (\ref{rma}), one can show that
\begin{equation*}
\left(\overline{\mathcal R}_{n-1}^{{\underline s^*}}\right)^{-1}\left(A^{\tau_-}\right)\cong{\rm Spec}\,\mathbb Z\left[\overrightarrow A\right]\times \mathbb M^{\underline s^*}_{\widetilde J}.
\end{equation*}
By the assumption that $\mathcal T_{n-1}^{\underline s^*}$ is smooth,
we can complete the proof for case I.
\smallskip

{\noindent\bf Case II ($k=2$).} 
Partition the set  
\begin{equation}\label{dudu1}
\{(1,i_2,i_3,i_4)\in\mathbb Z^4|\,3\leq i_2<i_3<i_4\leq n\}\cup\{(i_1,i_2,i_3,i_4)\in\mathbb Z^4|\,3\leq i_1<i_2<i_3<i_4\leq n\}   
\end{equation}
as follows.
\begin{enumerate}[label={(G\arabic*)},ref=G\arabic*]
\setcounter{enumi}{-1}
\item $\{(1,3,i_3,i_4)|4\leq i_3<i_4\leq n\}$.

\item $\{(1,i_2,i_3,i_4)|4\leq i_2< i_3<i_4\leq n\}$.

\item $\{(i_1,i_2,i_3,i_4)|4\leq i_1<i_2<i_3<i_4\leq n\}$.

\item $\{(3,i_2,i_3,i_4)|4\leq i_2<i_3<i_4\leq n\}$.

\end{enumerate}
Define a total order $\sigma$ so that the $4$-tuples in an upper subset are strictly smaller than the ones in a lower subset. 
Corresponding to indices in (G0), (G1), (G2), we define
\begin{equation*}
\check {\mathcal R}_{2,1}:=\check\psi^{\sigma}_1\circ\check\psi^{\sigma}_2\circ\cdots\circ\check\psi^{\sigma}_{\frac{(n-3)!}{24\cdot(n-7)!}+\frac{(n-3)!}{6\cdot(n-6)!}+\frac{(n-3)(n-4)}{2}},
\end{equation*}
and corresponding to indices in (G3), we define
\begin{equation*}
\check{\mathcal R}_{2,2}:=\check\psi^{\sigma}_{\frac{(n-3)!}{24\cdot(n-7)!}+\frac{(n-3)!}{6\cdot(n-6)!}+\frac{(n-3)(n-4)}{2}+1}\circ\check\psi^{\sigma}_{\frac{(n-3)!}{24\cdot(n-7)!}+\frac{(n-3)!}{6\cdot(n-6)!}+\frac{(n-3)(n-4)}{2}+2}\circ\cdots\circ\check\psi^{\sigma}_{\frac{(n-3)!}{24\cdot(n-7)!}+\frac{(n-3)!}{3\cdot(n-6)!}+\frac{(n-3)(n-4)}{2}}.
\end{equation*}

Computation yields that for $4\leq i_2<i_3<i_4\leq n$, $\left(\pi\circ j\right)^{-1}\mathscr I^{\underline s}_{1i_2i_3i_4}\cdot\mathcal O_{C^2_N}$ is a product of the invertible ideal sheaf generated by $\prod\nolimits_{\gamma=2}^{i_2}\epsilon^+_{\gamma}\cdot\prod\nolimits_{\gamma=2}^{i_3}\epsilon^+_{\gamma}$ and the ideal sheaf $\mathcal I_{1i_2i_3i_4}$ generated by
\begin{equation}\label{3A1}
\left\{\prod\nolimits_{\gamma=i_3+1}^{i_4}\epsilon^+_{\gamma}-1,\,\,\,\,\prod\nolimits_{\gamma=i_2+1}^{i_4}\epsilon^+_{\gamma}-1,\,\,\,\,\prod\nolimits_{\gamma=i_2+1}^{i_3}\epsilon^+_{\gamma}-1\right\}.
\end{equation}
For $4\leq i_2<i_3<i_4\leq n$, $\left(\pi\circ j\right)^{-1}\mathscr I^{\underline s}_{3i_2i_3i_4}\cdot\mathcal O_{C^2_N}$ is generated by
\begin{equation*}
\begin{split}
&\left\{\left(\prod\nolimits_{\gamma=2}^{i_2}\epsilon^+_{\gamma}-1\right)\left(\prod\nolimits_{\gamma=i_3+1}^{i_4}\epsilon^+_{\gamma}-1\right)\prod\nolimits_{\gamma=2}^{i_3}\epsilon^+_{\gamma},\,\,\,\,\left(\prod\nolimits_{\gamma=2}^{i_3}\epsilon^+_{\gamma}-1\right)\left(\prod\nolimits_{\gamma=i_2+1}^{i_4}\epsilon^+_{\gamma}-1\right)\prod\nolimits_{\gamma=2}^{i_2}\epsilon^+_{\gamma},\right.\\
&\,\,\,\,\,\,\,\,\,\,\,\,\,\,\,\,\,\,\,\,\,\,\,\,\left.\left(\prod\nolimits_{\gamma=2}^{i_4}\epsilon^+_{\gamma}-1\right)\left(\prod\nolimits_{\gamma=i_2+1}^{i_3}\epsilon^+_{\gamma}-1\right)\prod\nolimits_{\gamma=2}^{i_2}\epsilon^+_{\gamma}\right\}.\\
\end{split}
\end{equation*}
Since $\prod\nolimits_{\gamma=2}^{i_\alpha}\epsilon^+_{\gamma}-1\neq0$ in $C^2_N$ for $\alpha=2,3,4$, we can derive that $\left(\pi\circ j\right)^{-1}\mathscr I^{\underline s}_{3i_2i_3i_4}\cdot\mathcal O_{C^2_N}$ is a product of $\mathcal I_{1i_2i_3i_4}$ and the invertible ideal sheaf generated by $\prod\nolimits_{\gamma=2}^{i_2}\epsilon^+_{\gamma}$. 
Then the restriction of $\check {\mathcal R}_{2,2}$ to
$\left(\check{\mathcal R}_{2,1}\right)^{-1}\left(C^2_N\right)$ is an isomorphism.  Hence to prove  (\ref{p11}) for $k=2$, it suffices to show that $\left(\check{\mathcal R}_{2,1}\right)^{-1}\left(\mathbb A_+\right)$ is smooth.

Define a morphism  $\widetilde L:{\rm Spec}\,\mathbb Z\left[\epsilon^+_3,\epsilon^+_4,\cdots,\epsilon^+_{N-2}\right]\rightarrow \mathcal{G}(2,n-2)$ by
\begin{equation}\label{j3}
\left(\begin{matrix}
1&1&1&1&\cdots&1\\
0&1&\epsilon^+_3&\epsilon^+_3\epsilon^+_4&\cdots&\epsilon^+_3\epsilon^+_4\cdots\epsilon^+_{N-2}
\end{matrix}\right).  
\end{equation}
One can show that
\begin{equation*}
\left(\check{\mathcal R}_{2,1}\right)^{-1}\left(\mathbb A_+\right)\cong{\rm Spec}\,\mathbb Z\left[\epsilon^+_2\right]\times\mathbb M^{\underline {\widehat s}}_{\widetilde L},  
\end{equation*} where $\underline {\widehat s}=(1,1,\cdots,1)\in\mathbb Z^{N-2}$. By the assumption that Proposition \ref{ms} holds for all $N^{\prime}$ such that $4\leq N^{\prime}<N$, we can conclude
by (\ref{prod31}), (\ref{j1}), (\ref{prod32}) that $\mathbb M^{\underline {\widehat s}}_{\widetilde L}$ is smooth over ${\rm Spec}\,\mathbb Z$. 

The proof is complete. 

\smallskip

{\noindent\bf Case III ($3\leq k\leq N-2$).} Partition set (\ref{dudu1}) as follows. 
\begin{enumerate}[label={(A\arabic*)},ref=A\arabic*]
\item $\{(1,i_2,i_3,i_4)|\,k+2\leq i_2<i_3<i_4\leq n\}$.

\item $\{(i_1,i_2,i_3,i_4)|\,k+2\leq i_1<i_2<i_3<i_4\leq n\}$.

\end{enumerate}
\begin{enumerate}[label={(B\arabic*)},ref=B\arabic*]
\item $\{(1,3,i_3,i_4)|\,4\leq i_3<i_4\leq k+1\}$.

\item $\{(1,i_2,i_3,i_4)|\,4\leq i_2<i_3<i_4\leq k+1\}$.

\item $\{(3,i_2,i_3,i_4)|\,4\leq i_2<i_3<i_4\leq k+1\}$.

\item $\{(i_1,i_2,i_3,i_4)|\,4\leq i_1<i_2<i_3<i_4\leq k+1\}$.

\end{enumerate}
\begin{enumerate}[label={(C\arabic*)},ref=C\arabic*]
\item $\{(1,3,i_3,i_4)|\,k+2\leq i_3<i_4\leq n\}$.

\item $\{(3,i_2,i_3,i_4)|\,k+2\leq i_2<i_3<i_4\leq n\}$.

\item $\{(1,3,i_3,i_4)|\,4\leq i_3\leq k+1,\,\,k+2\leq i_4\leq n\}$.

\item $\{(1,i_2,i_3,i_4)|\,4\leq i_2<i_3\leq k+1,\,\,k+2\leq i_4\leq n\}$.

\item $\{(1,i_2,i_3,i_4)|\,4\leq i_2\leq k+1,\,\,k+2\leq i_3<i_4\leq n\}$.

\item $\{(3,i_2,i_3,i_4)|\,4\leq i_2<i_3\leq k+1,\,\,k+2\leq i_4\leq n\}$.

\item $\{(3,i_2,i_3,i_4)|\,4\leq i_2\leq k+1,\,\,k+2\leq i_3<i_4\leq n\}$.

\item $\{(i_1,i_2,i_3,i_4)|\,4\leq i_1<i_2<i_3\leq k+1,\,\,k+2\leq i_4\leq n\}$.

\item $\{(i_1,i_2,i_3,i_4)|\,4\leq i_1<i_2\leq k+1,\,\,k+2\leq i_3<i_4\leq n\}$.

\item $ \{(i_1,i_2,i_3,i_4)|\,4\leq i_1\leq k+1,\,\,k+2\leq i_2<i_3<i_4\leq n\}$.

\end{enumerate}Define a total order $\sigma$ so that the $4$-tuples in an upper subset are strictly smaller than the ones in a lower subset. 
Define
\begin{equation*}
\begin{split}
&\check {\mathcal R}_{3,1}:=\check\psi^{\sigma}_{1}\circ\check\psi^{\sigma}_{2}\circ\cdots\circ\check\psi^{\sigma}_{N_1},\\
&\check {\mathcal R}_{3,2}:=\check\psi^{\sigma}_{N_1+1}\circ\check\psi^{\sigma}_{N_1+2}\circ\cdots\circ\check\psi^{\sigma}_{N_1+N_2},\\
&\check{\mathcal R}_{3,3}:=\check\psi^{\sigma}_{N_1+N_2+1}\circ\check\psi^{\sigma}_{N_1+N_2+2}\circ\cdots\circ\check\psi^{\sigma}_{\frac{(n-2)!}{24\cdot(n-6)!}+\frac{(n-2)!}{6\cdot(n-5)!}},\\
\end{split}
\end{equation*}   
where $N_1$ (resp. $N_2$) is the sum of the cardinalities of (A1)-(A2) (resp. (B1)-(B4)). Equivalently, $\check{\mathcal R}_{3,1}$ ($\check{\mathcal R}_{3,2}$) is the sequence of blow-ups with respect to the (pull-backs of) ideal sheaves $\mathscr I^{\underline s}_{i_1i_2i_3i_4}$  with indices in (A1)-(A2) (resp. (B1)-(B4)).

\begin{claim}\label{3gi}
The restriction of  $\check{\mathcal R}_{3,3}$ to
$\left(\check{\mathcal R}_{3,1}\circ\check{\mathcal R}_{3,2}\right)^{-1}\left(C^k_N\right)$ is an isomorphism.
\end{claim}

{\bf\noindent Proof of Claim \ref{3gi}.} The proof is similar to that for Claim \ref{gi}. It suffices to show that for any $(i_1,i_2,i_3,i_4)$ in (C1)-(C10), the restriction  to 
$\left(\check{\mathcal R}_{3,1}\circ\check{\mathcal R}_{3,2}\right)^{-1}\left(C^k_N\right)$ of the ideal sheaf
\begin{equation*}
\left(\check{\mathcal R}_{3,1}\circ\check{\mathcal R}_{3,2}\right)^{-1}\left(\left(\pi\circ j\right)^{-1}\mathscr I^{\underline s}_{i_1i_2i_3i_4}\cdot\mathcal O_{\mathbb A_+}\right)\cdot\mathcal O_{\left(\check{\mathcal R}_{3,1}\circ\check{\mathcal R}_{3,2}\right)^{-1}\left(C^k_N\right)}=:\mathscr I_{i_1i_2i_3i_4}    
\end{equation*} is invertible.
The argument is on a case by case basis.

In case (C1), $\left(\pi\circ j\right)^{-1}\mathscr I^{\underline s}_{13i_3i_4}\cdot\mathcal O_{\mathbb A_+}$, $k+2\leq i_3<i_4\leq n$, is generated by
\begin{equation*}
\left\{\left(\prod\nolimits_{\gamma=i_3+1}^{i_4}\epsilon^+_{\gamma}-1\right)\prod\nolimits_{\gamma=2}^{i_3}\epsilon^+_{\gamma},\,\,\,\,\left(\prod\nolimits_{\gamma=2}^{i_4}\epsilon^+_{\gamma}-1\right)\prod\nolimits_{\gamma=2}^{i_3}\epsilon^+_{\gamma},\,\,\,\,\left(\prod\nolimits_{\gamma=2}^{i_3}\epsilon^+_{\gamma}-1\right)\prod\nolimits_{\gamma=2}^{i_4}\epsilon^+_{\gamma}\right\}.
\end{equation*}
Since $\prod\nolimits_{\gamma=2}^{i_3}\epsilon^+_{\gamma}-1\neq0$ and  $\prod\nolimits_{\gamma=2}^{i_4}\epsilon^+_{\gamma}-1\neq0$ in $C^k_N$, $\mathscr I_{13i_3i_4}$ is generated by $\prod\nolimits_{\gamma=2}^{i_3}\epsilon^+_{\gamma}$. Similarly, we can show that $\mathscr I_{i_1i_2i_3i_4}$ are invertible for indices in (C3), (C4), (C5), (C7), (C9). 

In case (C2), $\left(\pi\circ j\right)^{-1}\mathscr I^{\underline s}_{3i_2i_3i_4}\cdot\mathcal O_{\mathbb A_+}$, $k+2\leq i_2<i_3<i_4\leq n$, is generated by 
\begin{equation*}
\begin{aligned}
&\left\{\left(\prod\nolimits_{\gamma=2}^{i_2}\epsilon^+_{\gamma}-1\right)\left(\prod\nolimits_{\gamma=i_3+1}^{i_4}\epsilon^+_{\gamma}-1\right)\prod\nolimits_{\gamma=2}^{i_3}\epsilon^+_{\gamma},\,\,\,\,\left(\prod\nolimits_{\gamma=2}^{i_3}\epsilon^+_{\gamma}-1\right)\left(\prod\nolimits_{\gamma=i_2+1}^{i_4}\epsilon^+_{\gamma}-1\right)\prod\nolimits_{\gamma=2}^{i_2}\epsilon^+_{\gamma}\right.\\
&\left.\,\,\,\,\,\,\,\,\,\,\,\,\,\,\,\,\,\,\,\,\,\,\,\,\,\,\,\left(\prod\nolimits_{\gamma=2}^{i_4}\epsilon^+_{\gamma}-1\right)\left(\prod\nolimits_{\gamma=i_2+1}^{i_3}\epsilon^+_{\gamma}-1\right)\prod\nolimits_{\gamma=2}^{i_2}\epsilon^+_{\gamma}\right\}.\\ 
\end{aligned}
\end{equation*}
Computation yields that $\left(\pi\circ j\right)^{-1}\mathscr I^{\underline s}_{1i_2i_3i_4}\cdot\mathcal O_{\mathbb A_+}$, $4\leq i_2<i_3<i_4\leq k+1$, is a product of the invertible ideal sheaf generated by $\prod\nolimits_{\gamma=2}^{i_2}\epsilon^+_{\gamma}\cdot\prod\nolimits_{\gamma=2}^{i_3}\epsilon^+_{\gamma}$ and the ideal sheaf $\mathcal I_{i_2i_3i_4}$ generated by 
\begin{equation*}
\left\{\prod\nolimits_{\gamma=i_3+1}^{i_4}\epsilon^+_{\gamma}-1,\,\,\,\,\prod\nolimits_{\gamma=i_2+1}^{i_3}\epsilon^+_{\gamma}-1\right\}.
\end{equation*} 
Since $\prod\nolimits_{\gamma=2}^{i_2}\epsilon^+_{\gamma}-1\neq0$, $\prod\nolimits_{\gamma=2}^{i_3}\epsilon^+_{\gamma}-1\neq0$, $\prod\nolimits_{\gamma=2}^{i_4}\epsilon^+_{\gamma}-1\neq0$ in $C^k_N$,  we can derive that $\left(\pi\circ j\right)^{-1}\mathscr I^{\underline s}_{3i_2i_3i_4}\cdot\mathcal O_{\mathbb A_+}$ is a product of the invertible ideal sheaf generated by $\prod\nolimits_{\gamma=2}^{i_2}\epsilon^+_{\gamma}$ and the ideal sheaf $\mathcal I_{i_2i_3i_4}$, and hence $\mathscr I_{3i_2i_3i_4}$ is invertible. Similarly, we can show that $\mathscr I_{i_1i_2i_3i_4}$ are invertible for indices in (C6), (C8), (C10).

We complete the proof of Claim \ref{3gi}.\,\,\,\,$\endpf$

\smallskip

Define an open subscheme of $\mathbb A_+$ by
$D^{k}_{N}:=\left\{\mathfrak p\in \mathbb A_+\big|\,\epsilon^+_{\gamma}\notin\mathfrak p,\,\,\forall\, 2\leq \gamma\leq k-1 \right\}$.
Similarly to (\ref{j3}), we define a morphism $\widetilde L:{\rm Spec}\,\mathbb Z\left[\epsilon^+_{2},\epsilon^+_{3},\cdots,\epsilon^+_{k-1}\right]\longrightarrow \mathcal{G}(2,k)$ by
\begin{equation*}
\left(\begin{matrix}
1&1&1&1&\cdots&1\\
0&1&\epsilon^+_2&\epsilon^+_2\epsilon^+_3&\cdots&\epsilon^+_2\epsilon^+_3\cdots\epsilon^+_{k-1}
\end{matrix}\right),  
\end{equation*}
and a morphism $\widehat L:{\rm Spec}\,\mathbb Z\left[\epsilon^+_{k+1},\epsilon^+_{k+2},\cdots,\epsilon^+_{N-2}\right]\longrightarrow \mathcal{G}(2,n-k)$ by
\begin{equation*}
\left(\begin{matrix}
1&1&1&1&\cdots&1\\
0&1&\epsilon^+_{k+1}&\epsilon^+_{k+1}\epsilon^+_{k+2}&\cdots&\epsilon^+_{k+1}\epsilon^+_{k+2}\cdots\epsilon^+_{N-2}
\end{matrix}\right).   
\end{equation*}
Now by Claim \ref{3gi}, we can show that 
$\left(\check{\mathcal R}_{3,1}\circ\check{\mathcal R}_{3,2}\right)^{-1}\left(D^k_N\right)$ is isomorphic to an open subscheme of ${\rm Spec}\,\mathbb Z\left[\epsilon^+_{k}\right]\times\mathbb M^{\underline {\widetilde  s}}_{\widetilde L}\times\mathbb M^{\underline {\widehat s}}_{\widehat L}$, 
where $\underline {\widetilde s}=(1,1,\cdots,1)\in\mathbb Z^{k}$,  $\underline {\widehat s}=(1,1,\cdots,1)\in\mathbb Z^{N-k}$. By the same argument as that in {case II}, we can derive that $\mathbb M^{\underline {\widetilde s}}_{\widetilde L}$, $\mathbb M^{\underline {\widehat s}}_{\widehat L}$ are smooth over ${\rm Spec}\,\mathbb Z$. We can thus conclude the proof of (\ref{p11}) for case III.
\medskip

The proof of Lemma \ref{ud} is complete.\,\,\,\,$\endpf$

\section{The Case of Sub-torus Actions}\label{sta}

In this section, we are concerned with 
\begin{proposition}\label{3ms} Let $N\geq 2$ be an integer. Let $\underline s=(s_1,s_2,\cdots,s_{N})$ be a size vector. Denote $n:=s_1+s_2+\cdots s_N$.  Then, $\mathcal T_{n}^{\underline s}$ and $\mathcal M_{n}^{\underline s}$ are smooth over ${\rm Spec}\,\mathbb Z$ of relative dimensions $2n-4$ and $2n-N-3$, respectively.
Moreover, the morphism $\mathcal P^{\underline s}_{n}:\mathcal T_{n}^{\underline s}\rightarrow\mathcal M_{n}^{\underline s}$ is flat. \end{proposition}
We note that the structure of the section is in parallel with that of \S \ref{mta} to closely following the ideas therein, while the proofs will be concise to avoid repetition. 

\subsection{Parametrizations for \texorpdfstring{$\mathcal Q_{n}^{\underline s}$}{dd}}\label{3blf}

To better locate relative positions of columns, we introduce
\begin{definition}
Fix a size vector $\underline s=(s_1,s_2,\cdots,s_N)$. We associate to each integer $1\leq j\leq n$ a triple of positive integers $\left(j,t,p\right)$ such that
\begin{enumerate}
\item $1\leq t\leq N$,
\item $1\leq p \leq s_{t}$,
\item $j=s_1+s_2+\cdots+s_{t-1}+p$.
\end{enumerate}
Here we adopt the convention that $s_0=0$. We call $\left(j,t,p\right)$ the $\underline s$-triple of $j$.  We also write $\left(j,t,p\right)$ as $j$ for simplicity of notation when there is no ambiguity.
\end{definition}

We now define
\begin{definition}Let $\mathbb J^{\underline s}$ be the index set consisting of all tuples of $\underline s$-triples
\begin{equation}\label{ind}
\small
\left(\left(j_1,t_1,p_1\right),\left(j_2,t_2,p_2\right),\left(\left(j^+_1,t^+_1,p^+_1\right)\cdots,\left(j^+_l,t^+_l,p^+_l\right)\right),\left(\left(j^-_1,t^-_1,p^-_1\right),\cdots,\left(j^-_m,t^-_m,p^-_m\right)\right)\right)    
\end{equation}with the following properties.
\begin{enumerate}[label=(\arabic*)]

\item $l,m$ are non-negative integers such that $l+m=N-2$ or $N-1$. Moreover, if $l+m=N-2$, then $t_1<t_2$; if $l+m=N-1$, then $t_1=t_2$ and $p_1<p_2$.

\item  $t^+_1<t^+_2<\cdots<t^+_l$,  $t^-_1<t^-_2<\cdots<t^-_m$, and $\{t_1,t_2,t^+_1,\cdots,t^+_l,t^-_1,\cdots,t^-_{m}\}=\{1,2,\cdots,N\}$\,. 

\end{enumerate}
We write
(\ref{ind}) as $\left(j_1,j_2,\left(j^+_1,\cdots,j^+_l\right),\left( j^-_1,\cdots, j^-_{m}\right)\right)$ for simplicity of notation. We call an index $\tau=\left(j_1, j_2,\left(j^+_1,\cdots,j^+_l\right),\left(j^-_1,\cdots,j^-_{m}\right)\right)\in \mathbb J^{\underline s}$ a Class I index if $l+m=N-2$, and a Class II index if $l+m=N-1$.
\end{definition}

We associate to each Class I index $\tau=\left( j_1, j_2,\left(j^+_1,\cdots,j^+_l\right),\left(j^-_1,\cdots,j^-_{m}\right)\right)\in \mathbb J^{\underline s}$ an affine space ${\rm Spec}\,\mathbb Z\left[\overrightarrow A,\overrightarrow Y,\overrightarrow Z,\overrightarrow U,\overrightarrow V,\overrightarrow H^1,\overrightarrow H^2,\cdots,\overrightarrow H^l,\overrightarrow \Xi^1,\overrightarrow \Xi^2,\cdots,\overrightarrow \Xi^m\right]$. Here
\begin{equation}\label{t1ab}
\begin{split}
&\,\overrightarrow  A:=\left(a_{j^+_1},a_{j^+_2},\cdots,a_{j^+_l},a_{j^-_1},a_{j^-_2},\cdots,a_{j^-_m}\right);\\
\end{split}  
\end{equation}
\begin{equation}\label{317}
\begin{split}
&\overrightarrow Y:=\left(y_{1},y_{2},\cdots,y_{p_1-1},y_{p_1+1},\cdots,y_{s_{t_1}}\right),\,\,\overrightarrow U:=\left(u_{1},u_{2},\cdots,u_{p_1-1},u_{p_1+1},\cdots,u_{s_{t_1}}\right),\\
&\overrightarrow Z:=\left(z_{1},z_{2},\cdots,z_{p_2-1},z_{p_2+1},\cdots,z_{s_{t_2}}\right),\,\,\overrightarrow V:=\left(v_{1},v_{2},\cdots,v_{p_2-1},v_{p_2+1},\cdots,v_{s_{t_2}}\right);\\
\end{split}
\end{equation}
for $1\leq\alpha\leq l$,
\begin{equation}\label{ha}
\overrightarrow  H^{\alpha}:=\left(\eta^{\alpha}_{11},\eta^{\alpha}_{12},\cdots,\eta^{\alpha}_{1\left(p_{\alpha}^+-1\right)},\eta^{\alpha}_{1\left(p_{\alpha}^++1\right)},\cdots,\eta^{\alpha}_{1s_{t_{\alpha}^+}}, \eta^{\alpha}_{21},\eta^{\alpha}_{22},\cdots,\eta^{\alpha}_{2s_{t_{\alpha}^+}}\right); 
\end{equation}
for $1\leq\beta\leq m$,
\begin{equation}\label{xb}
\overrightarrow  \Xi^{\beta}:=\left(\xi^{\beta}_{11},\xi^{\beta}_{12},\cdots,\xi^{\beta}_{1s_{t_{\beta}^-}},\xi^{\beta}_{21},\xi^{\beta}_{22},\cdots,\xi^{\beta}_{2\left(p_{\beta}^--1\right)},\xi^{\beta}_{2\left(p_{\beta}^-+1\right)},\cdots,\xi^{\beta}_{2s_{t_{\beta}^-}}\right). 
\end{equation}
\begin{remark}
We adopt the convention that there is no such variable when the subscript is out of range. For instance,  there is no $y_{p_1-1}$, $u_{p_1-1}$  when $p_1=1$, 
and no $y_{p_1+1}$, $u_{p_1+1}$ when $p_1=s_{t_1}$.
\end{remark}
We next define a morphism $\Gamma^{\tau}$ from \begin{equation}\label{spab}
{\rm Spec}\,\mathbb Z\left[\overrightarrow A,\overrightarrow Y,\overrightarrow Z,\overrightarrow U,\overrightarrow V,\overrightarrow H^1,\overrightarrow H^2,\cdots,\overrightarrow H^l,\overrightarrow \Xi^1,\overrightarrow \Xi^2,\cdots,\overrightarrow \Xi^m\right]   
\end{equation} to $U_{j_1j_2}\subset G(2,n)$ as follows. For $j_1-p_1+1\leq i\leq j_1-p_1+s_{t_1}$, and $i\neq j_1$, set
\begin{equation}\label{dang1}
x_{1i}\mapsto y_{\left(i-j_1+p_1\right)},\,\,\,x_{2i}\mapsto u_{\left(i-j_1+p_1\right)};
\end{equation}
for $j_2-p_2+1\leq i\leq j_2-p_2+s_{t_2}$, and $i\neq j_2$, set
\begin{equation}\label{dang2}
x_{1i}\mapsto v_{\left(i-j_2+p_2\right)},\,\,\,x_{2i}\mapsto z_{\left(i-j_2+p_2\right)};
\end{equation}
for $1\leq\alpha\leq l$, set
\begin{equation}\label{ga}
\left\{
\begin{aligned}
&x_{1j_{\alpha}^+}\mapsto a_{j_{\alpha}^+},\\
&x_{1i}\mapsto a_{j_{\alpha}^+}\cdot\eta^{\alpha}_{1\left(i-j_{\alpha}^++p_{\alpha}^+\right)},&j_{\alpha}^+-p_{\alpha}^++1\leq i\leq j_{\alpha}^+-p_{\alpha}^++s_{t_{\alpha}^+}\,\,{\rm and}\,\, i\neq j_{\alpha}^+ \\
&x_{2i}\mapsto a_{j_{\alpha}^+}\cdot\eta^{\alpha}_{2\left(i-j_{\alpha}^++p_{\alpha}^+\right)},&j_{\alpha}^+-p_{\alpha}^++1\leq i\leq j_{\alpha}^+-p_{\alpha}^++s_{t_{\alpha}^+}\,\,\,\,\,\,\,\,\,\,\,\,\,\,\,\,\,\,\,\,\,\,\,\,\,\,\,\,\,\,\\
\end{aligned}\right.\,\,;
\end{equation}
for  $1\leq\beta\leq m$, set
\begin{equation}\label{gb}
\left\{
\begin{aligned}
&x_{1i}\mapsto a_{j_{\beta}^-}\cdot\xi^{\beta}_{1\left(i-j_{\beta}^-+p_{\beta}^-\right)},&j_{\beta}^--p_{\beta}^-+1\leq i\leq j_{\beta}^--p_{\beta}^-+s_{t_{\beta}^-}\,\,\,\,\,\,\,\,\,\,\,\,\,\,\,\,\,\,\,\,\,\,\,\,\,\,\,\,\,\,\\
&x_{2j_{\beta}^-}\mapsto a_{j_{\beta}^-},\\
&x_{2i}\mapsto a_{j_{\beta}^-}\cdot\xi^{\beta}_{2\left(i-j_{\beta}^-+p_{\beta}^-\right)},&j_{\beta}^--p_{\beta}^-+1\leq i\leq j_{\beta}^--p_{\beta}^-+s_{t_{\beta}^-}\,\,{\rm and}\,\, i\neq j_{\beta}^- \\
\end{aligned}\right.\,\,.
\end{equation}

Similarly, we associate to each Class II index $\tau=\left( j_1, j_2,\left(j^+_1,\cdots,j^+_l\right),\left(j^-_1,\cdots,j^-_{m}\right)\right)\in \mathbb J^{\underline s}$ an affine space ${\rm Spec}\,\mathbb Z\left[\overrightarrow A,\overrightarrow W,\overrightarrow H^1,\cdots,\overrightarrow H^l,\overrightarrow \Xi^1,\cdots,\overrightarrow \Xi^m\right]$. Here $\overrightarrow  A$ are defined by (\ref{t1ab}); $\overrightarrow  H^{\alpha}$ are defined by (\ref{ha}) for $1\leq\alpha\leq l$;  $\overrightarrow  \Xi^{\beta}$ are defined by (\ref{xb}) for $1\leq\beta\leq m$; 
\begin{equation}\label{t2w}
\begin{split}
&\overrightarrow W:=\left(w_{11},w_{12},\cdots,w_{1(p_1-1)},w_{1(p_1+1)},\cdots,w_{1(p_2-1)},w_{1(p_2+1)},\cdots,w_{s_{t_1}},\right.\\
&\,\,\,\,\,\,\,\,\,\,\,\,\,\,\,\,\,\,\,\,\,\,\left.w_{21},w_{22},\cdots,w_{2(p_1-1)},w_{2(p_1+1)},\cdots,w_{2(p_2-1)},w_{2(p_2+1)},\cdots,w_{2s_{t_1}}\right).\\
\end{split}   
\end{equation}
Define a morphism $\Gamma^{\tau}:{\rm Spec}\,\mathbb Z\left[\overrightarrow A,\overrightarrow W,\overrightarrow H^1,\cdots,\overrightarrow H^l,\overrightarrow \Xi^1,\cdots,\overrightarrow \Xi^m\right]\longrightarrow U_{j_1j_2}$  as follows. For $j_{\alpha}^+-p_{\alpha}^++1\leq i\leq j_{\alpha}^+-p_{\alpha}^++s_{t_{\alpha}^+}$ where $1\leq\alpha\leq l$, define the images of $x_{1i}$, $x_{2i}$ by (\ref{ga}); for $j_{\beta}^--p_{\beta}^-+1\leq i\leq j_{\beta}^--p_{\beta}^-+s_{t_{\beta}^-}$ where  $1\leq\beta\leq m$, define the images of $x_{1i}$, $x_{2i}$ by (\ref{gb}).
For $j_1-p_1+1\leq i\leq j_1-p_1+s_{t_1}$, and $i\neq j_1,j_2$, define
\begin{equation}\label{dang3}
x_{1i}\mapsto w_{1\left(i-j_1+p_1\right)},\,\,\,x_{2i}\mapsto w_{2\left(i-j_1+p_1\right)}.
\end{equation}

Similarly to Lemma \ref{em}, we can define for each Class I index $\tau\in\mathbb J^{\underline s}$ a locally closed embedding  \begin{equation*}
J^{\tau}:{\rm Spec}\,\mathbb Z\left[\overrightarrow A,\overrightarrow Y,\overrightarrow Z,\overrightarrow U,\overrightarrow V,\overrightarrow H^1,\cdots,\overrightarrow H^l,\overrightarrow \Xi^1,\cdots,\overrightarrow \Xi^m\right]\longrightarrow\mathbb P^{N_{2,n}}\times\prod\nolimits_{t=1}^N\mathbb P^{N^{\underline s}_{t}}    
\end{equation*} 
and for each Class II index $\tau\in\mathbb J^{\underline s}$ a locally closed embedding   \begin{equation*}
J^{\tau}:{\rm Spec}\,\mathbb Z\left[\overrightarrow A,\overrightarrow W,\overrightarrow H^1,\cdots,\overrightarrow H^l,\overrightarrow \Xi^1,\cdots,\overrightarrow \Xi^m\right]\longrightarrow\mathbb P^{N_{2,n}}\times\prod\nolimits_{t=1}^N\mathbb P^{N^{\underline s}_{t}}    
\end{equation*}
by extend the rational maps $\widetilde{\mathcal K}^{\underline s}\circ \Gamma^{\tau}=\left(e\circ\Gamma^{\tau},F^{\underline s}_{1}\circ e\circ \Gamma^{\tau},F^{\underline s}_{2}\circ e\circ \Gamma^{\tau},\cdots,F^{\underline s}_{N}\circ e\circ \Gamma^{\tau}\right)$. Denote by $A^{\tau}$ the corresponding image under  $J^{\tau}$ for each $\tau\in\mathbb J^{\underline s}$, which is equipped with the reduced scheme structure as a locally closed subscheme of $\mathcal Q_{n}^{\underline s}\subset\mathbb P^{N_{2,n}}\times\prod\nolimits_{t=1}^N\mathbb {P}^{N^{\underline s}_{t}}$. Moreover, 
\begin{proposition}\label{3coor2} $\mathcal Q_{n}^{\underline s}$ is a union of $A^{\tau}$ for all $\tau\in\mathbb J^{\underline s}$. $\mathcal Q_{n}^{\underline s}$ is smooth over ${\rm Spec}\,\mathbb Z$. 
\end{proposition}

{\bf\noindent Proof of Proposition \ref{3coor2}.} The proof is similar to that of Lemma \ref{coor2} and Corollary \ref{coor3}. We omit the details for brevity. 
\medskip

\subsection{Factorization}\label{3FAC} Fix a class I index $\tau=\left( j_1, j_2,\left(j^+_1,\cdots,j^+_l\right),\left(j^-_1,\cdots,j^-_{m}\right)\right)\in \mathbb J^{\underline s}$. In what follows, we shall write   ${\underline w}\in C^{\underline s}$ as \begin{equation}\label{wc}
\underline w=(w_{t_1},w_{t_2},w_{t^+_{1}},\cdots,w_{t^+_{\alpha}},\cdots,w_{t^+_{l}},w_{t^-_{1}},\cdots,w_{t^-_{\beta}},\cdots,w_{t^-_{m}}),    
\end{equation}
by rearranging its components corresponding to $\tau$. 

Let (${\rm AB}{\tau}$) be the subset of $C^{\underline s}$ consisting of the following elements.
\begin{enumerate}[label={(A\arabic*)},ref=A\arabic*]
\setcounter{enumi}{-1}

\item $w_{t_1}=2$, and the other components are $0$. \label{fa0}

\item $w_{t_1}=w_{t^+_{\alpha}}=1$ for a certain $1\leq \alpha\leq l$, and the other components are $0$. 
\label{fa1}

\item $w_{t_1}=w_{t_2}=w_{t^+_{\alpha_1}}=w_{t^+_{\alpha_2}}=1$ for certain $1\leq \alpha_1<\alpha_2\leq l$.
\label{fa2}
\end{enumerate}
\begin{enumerate}[label={(B\arabic*)},ref=B\arabic*]
\setcounter{enumi}{-1}

\item $w_{t_2}=2$, and the other components are $0$. \label{fb0}

\item $w_{t_2}=w_{t^-_{\beta}}=1$ for a certain $1\leq \beta\leq m$, and the other components are $0$. 
\label{fb1}

\item $w_{t_1}=w_{t_2}=w_{t^-_{\beta_1}}=w_{t^-_{\beta_2}}=1$ for certain $1\leq \beta_1<\beta_2\leq m$.
\label{fb2}

\end{enumerate}
Denote by $N_{\tau}$ the cardinality of subset (${\rm AB}{\tau}$). Take a total order $\sigma$ on $C^{\underline s}$ so that $\sigma(\gamma)$ are in subset ${\rm (AB}{\tau}{\rm )}$ for $1\leq\gamma\leq N_{\tau}$. 
As in \S \ref{pdm},
we define 
\begin{equation}\label{312a1}
\begin{split}
&\check {\mathcal R}_N:=\psi^{\sigma}_{1}\circ\psi^{\sigma}_{2}\circ\cdots\circ\psi^{\sigma}_{N_{\tau}},\\
&\ddot{\mathcal R}_N:=\psi^{\sigma}_{N_{\tau}+1}\circ\psi^{\sigma}_{N_{\tau}+2}\circ\cdots\circ\psi^{\sigma}_{N_{\underline s}}.\\
\end{split}
\end{equation}

To describe an open cover of $\left(\check{\mathcal R}_N\right)^{-1}\left(A^{\tau}\right)$,  we define $\Lambda^{\tau}$ to be the index set consisting of indices $\lambda:=\left(q^+,q^-,\lambda^+,\lambda^-,\Delta^+,\Delta^-\right)$ with the following properties.
\begin{enumerate}[label={(\arabic*)},ref=\arabic*]
\item $q^+$ and $q^-$ are integers such that 
\begin{enumerate}[label=$\bullet$]

\item $0\leq q^+\leq l$ and  $0\leq q^-\leq m$,

\item if $s_{t_1}=1$ and $s_{t_2}=1$, then $q^+=l$ and  $q^-=m$, respectively.

\end{enumerate}

\item $\lambda^+=(\lambda^+_1,\lambda^+_2,\cdots,\lambda^+_{l})$ is a permutation of $\{1,2,\cdots,l\}$,
and $\lambda^-=(\lambda^-_1,\lambda^-_2,\cdots,\lambda^-_{m})$ is a permutation of  $\{1,2,\cdots,m\}$.
\label{plus1}

\item  $\Delta^+:=\left(\delta^+_1,\cdots,\delta^+_{\rho},\cdots,\delta^+_{q^+},\delta^+_{q^++1}\right)$ is a $\left(q^++1\right)$-tuple of integers such that 
\begin{enumerate}[label=$\bullet$]

\item $1\leq\delta^+_{\rho}\leq s_{t^+_{\lambda^+_{\rho}}}$, for $1\leq\rho\leq q^+$,

\item $1\leq\delta^+_{q^++1}\leq s_{t_1}$ and $\delta^+_{q^++1}\neq p_1$.
\end{enumerate}
$\Delta^-:=\left(\delta_1^-,\cdots,\delta_{\rho}^-,\cdots,\delta_{q^-}^-,\delta_{q^-+1}^-\right)$ is a $(q^-+1)$-tuple of integers such that 
\begin{enumerate}[label=$\bullet$]
\item $1\leq\delta^-_{\rho}\leq s_{t^-_{\lambda^-_{\rho}}}$,  for $1\leq\rho\leq q^-$,

\item $1\leq\delta_{q^-+1}^-\leq s_{t_2}$ and $\delta_{q^-+1}^-\neq p_2$.
\end{enumerate}
\label{plus2}

\end{enumerate}

\begin{remark}
We also adopt the following conventions. \begin{enumerate}[label=$\bullet$]

\item
When $s_{t_1}=1$ (resp. $s_{t_2}=1$), there is no $\delta^+_{q^++1}$ (resp. $\delta^-_{q^-+1}$).

\item If $l=0$ (resp. $m=0$), then $\lambda^+$ and $\Delta^+$ (resp. $\lambda^-$ and $\Delta^-$) are null, and we write $\lambda=\left(q^+,q^-,\lambda^-,\Delta^-\right)$ (resp. $\lambda=\left(q^+,q^-,\lambda^+,\Delta^+\right)$). 

\end{enumerate}
\end{remark}

We associate to each  $\lambda=\left(q^+,q^-,\lambda^+,\lambda^-,\Delta^+,\Delta^-\right)\in \Lambda^{\tau}$ an affine space
\begin{equation}\label{bigsp}
\small
{\rm Spec}\,\mathbb Z\left[\overrightarrow A,\overrightarrow Y,\overrightarrow Z,\overrightarrow E_+,\overrightarrow E_-,\overrightarrow {\mathfrak U},\overrightarrow {\mathfrak V},\overrightarrow {\mathscr H}^1,\cdots,\overrightarrow {\mathscr H}^l,\overrightarrow {\mathfrak H}^1,\cdots,\overrightarrow {\mathfrak H}^l,\overrightarrow {\mathscr X}^1,\cdots,\overrightarrow {\mathscr X}^m,\overrightarrow {\mathfrak X}^1,\cdots,\overrightarrow {\mathfrak X}^m\right].   
\end{equation}
Here
\begin{equation}\label{348}
\begin{split}
&\overrightarrow  A:=\left(a_{j^+_1},a_{j^+_2},\cdots,a_{j^+_l},a_{j^-_1},a_{j^-_2},\cdots,a_{j^-_m}\right),\\
&\overrightarrow Y:=\left(y_{1},y_{2},\cdots,y_{p_1-1},y_{p_1+1},\cdots,y_{s_{t_1}}\right),\,\,\overrightarrow Z:=\left(z_{1},z_{2},\cdots,z_{p_2-1},z_{p_2+1},\cdots,z_{s_{t_2}}\right),\\ &\overrightarrow E_+:=\left(\epsilon_1^+,\epsilon_2^+,\cdots,\epsilon_{q^+}^+,\epsilon_{q^++1}^+\right),\,\,\,\,\,\,\,\,\,\,\,\,\,\,\,\,\,\,\,\,\,\,\,\overrightarrow E_-:=\left(\epsilon_1^-,\epsilon_2^-,\cdots,\epsilon_{q^-}^-,\epsilon^-_{q^-+1}\right),\\
\end{split}  
\end{equation}
\begin{equation*}
\begin{split}
&\overrightarrow {\mathfrak U}:=\left(\mathfrak u_{1},\mathfrak u_{2},\cdots,\mathfrak u_{p_1-1},\mathfrak u_{p_1+1},\cdots,\mathfrak u_{\delta^+_{q^++1}-1},\mathfrak u_{\delta^+_{q^++1}+1},\cdots,\mathfrak u_{s_1}\right),\\
&\overrightarrow {\mathfrak V}:=\left(\mathfrak v_{1},\mathfrak v_{2},\cdots,\mathfrak v_{p_2-1},\mathfrak v_{p_2+1},\cdots,\mathfrak v_{\delta^-_{q^-+1}-1},\mathfrak v_{\delta^-_{q^-+1}+1},\cdots,\mathfrak v_{s_2}\right);\,\,\,\,\,\,\,\,\,\,\,\,\,\,\,\,\,\,\,\,\,\,\,\,\,\,\,\,\,\,\,\,\,\,\,\,\,\,\,\,\,\,\,\,\\
\end{split}  
\end{equation*}
for $1\leq\alpha\leq l$,
\begin{equation}\label{3ha}
\begin{split}
&\overrightarrow  {\mathscr H}^{\alpha}:=\left(\eta^{\alpha}_{11},\eta^{\alpha}_{12},\cdots,\eta^{\alpha}_{1\left(p_{\alpha}^+-1\right)},\eta^{\alpha}_{1\left(p_{\alpha}^++1\right)},\cdots,\eta^{\alpha}_{1s_{t_{\alpha}^+}}\right);\\
\end{split}
\end{equation}
for $\alpha=\lambda^+_{\rho}$ where $1\leq\rho\leq q^+$,
\begin{equation}
\overrightarrow{
{\mathfrak H}}^{\alpha}:=\left({\mathfrak h}^{\alpha}_{21},{\mathfrak h}^{\alpha}_{22},\cdots,{\mathfrak h}^{\alpha}_{2\left(\delta^+_{\rho}-1\right)},{\mathfrak h}^{\alpha}_{2\left(\delta^+_{\rho}+1\right)},\cdots,{\mathfrak h}^{\alpha}_{2s_{t_{\alpha}^+}}\right),     
\end{equation}
and for $\alpha=\lambda^+_{\rho}$ where $q^++1\leq\rho\leq l$,
\begin{equation}
\overrightarrow{
{\mathfrak H}}^{\alpha}:=\left({\mathfrak h}^{\alpha}_{21},{\mathfrak h}^{\alpha}_{22},\cdots,{\mathfrak h}^{\alpha}_{2s_{t_{\alpha}^+}}\right);  
\end{equation}
for $1\leq\beta\leq m$,
\begin{equation}\overrightarrow  {\mathscr X}^{\beta}:=\left(\xi^{\beta}_{21},\xi^{\beta}_{22},\cdots,\xi^{\beta}_{2\left(p_{\beta}^--1\right)},\xi^{\beta}_{2\left(p_{\beta}^-+1\right)},\cdots,\xi^{\alpha}_{2s_{t_{\beta}^-}}\right);
\end{equation}
for $\beta=\lambda^-_{\rho}$ where $1\leq\rho\leq q^-$,
\begin{equation}
\overrightarrow{\mathfrak X}^{\beta}:=\left({\mathfrak x}^{\beta}_{11},{\mathfrak x}^{\beta}_{12},\cdots,{\mathfrak x}^{\beta}_{1\left(\delta^-_{\rho}-1\right)},{\mathfrak x}^{\beta}_{1\left(\delta^-_{\rho}+1\right)},\cdots,{\mathfrak x}^{\beta}_{1s_{t_{\beta}^-}}\right), 
\end{equation}
and for $\beta=\lambda^-_{\rho}$ where $q^-+1\leq\rho\leq m$,
\begin{equation}\label{3xb}
\overrightarrow{\mathfrak X}^{\beta}:=\left({\mathfrak x}^{\beta}_{11},{\mathfrak x}^{\beta}_{12},\cdots,{\mathfrak x}^{\beta}_{1s_{t_{\beta}^-}}\right).
\end{equation}
\begin{remark}
We adopt the following convention.

\begin{enumerate}[label=$\bullet$]

\item When $s_{t_1}=1$ (resp. $s_{t_1}=1$), there are no $\overrightarrow{\mathfrak U}$, $\epsilon^+_{q^++1}$ (resp. $\overrightarrow{\mathfrak V}$, $\epsilon^-_{q^-+1}$) variables.
    
\item When  $l=0$ (resp. $m=0$), there are no $a_{j^+_1},a_{j^+_2},\cdots,a_{j^+_l}$, $\overrightarrow E_+$, $\overrightarrow{\mathscr H}^{\alpha}$, $\overrightarrow{\mathfrak H}^{\alpha}$ (resp. $a_{j^-_1},a_{j^-_2},\cdots,a_{j^-_m}$, $\overrightarrow E_-$, $\overrightarrow{\mathscr X}^{\beta}$, $\overrightarrow{\mathfrak X}^{\beta}$) variables.
\end{enumerate}
\end{remark}

We define a morphism $\Sigma^{\lambda}$ from (\ref{bigsp}) to 
\begin{equation}\label{smallsp}
A^{\tau}\cong{\rm Spec}\,\mathbb Z\left[\overrightarrow A,\overrightarrow Y,\overrightarrow Z,\overrightarrow U,\overrightarrow V,\overrightarrow H^1,\cdots,\overrightarrow H^l,\overrightarrow \Xi^1,\cdots,\overrightarrow \Xi^m\right]    
\end{equation} by the following ring homomorphism.  $\overrightarrow A$, $\overrightarrow Y$, $\overrightarrow Z$, $\left(\eta^{\alpha}_{12},\eta^{\alpha}_{13},\cdots,\eta^{\alpha}_{1s_{t_{\alpha}^+}}\right)$ for  $1\leq\alpha\leq l$, and $\left(\xi^{\beta}_{22},\xi^{\beta}_{23},\cdots,\xi^{\beta}_{2s_{t_{\beta}^-}}\right)$ for $1\leq\beta\leq m$ are mapped to the same variables of (\ref{smallsp}).
When $s_{t_1}\geq 2$, 
\begin{equation*}
\footnotesize
\begin{split}
&\left(u_1,\cdots,u_{p_1-1},u_{p_1+1},\cdots,u_{\delta^+_{q^++1}-1},u_{\delta^+_{q^++1}},u_{\delta^+_{q^++1}+1},\cdots,u_{s_{t_1}}\right)\mapsto\left(\mathfrak u_{1}\prod\limits_{\gamma=1}^{q^++1}\epsilon^+_{\gamma},\cdots,\right.\\
&\left.\mathfrak u_{p_1-1}\prod\limits_{\gamma=1}^{q^++1}\epsilon^+_{\gamma},\,\,\mathfrak u_{p_1+1}\prod\limits_{\gamma=1}^{q^++1}\epsilon^+_{\gamma},\cdots,\mathfrak u_{\delta^+_{q^++1}-1}\prod\limits_{\gamma=1}^{q^++1}\epsilon^+_{\gamma},\,\,\prod\limits_{\gamma=1}^{q^++1}\epsilon^+_{\gamma},\,\,\mathfrak u_{\delta^+_{q^++1}+1}\prod\limits_{\gamma=1}^{q^++1}\epsilon^+_{\gamma},\cdots,\mathfrak u_{s_{t_1}}\prod\limits_{\gamma=1}^{q^++1}\epsilon^+_{\gamma}\right).\\  
\end{split}   
\end{equation*}
When $s_{t_2}\geq 2$, 
\begin{equation*}
\footnotesize
\begin{split}
&\left(v_1,\cdots,v_{p_2-1}, v_{p_2+1}, \cdots,v_{\delta^-_{q^-+1}-1},v_{\delta^-_{q^-+1}},v_{\delta^-_{q^-+1}+1},\cdots,v_{s_{t_2}}\right)\mapsto\left(\mathfrak v_{1}\prod\limits_{\gamma=1}^{q^-+1}\epsilon^-_{\gamma}, \cdots,\right.\\
&\left.\mathfrak v_{p_2-1}\prod\limits_{\gamma=1}^{q^-+1}\epsilon^-_{\gamma},\,\,\mathfrak v_{p_2+1}\prod\limits_{\gamma=1}^{q^-+1}\epsilon^-_{\gamma},\cdots,\mathfrak v_{\delta^-_{q^-+1}-1}\prod\limits_{\gamma=0}^{q^-}\epsilon^-_{\gamma},\,\,\prod\limits_{\gamma=1}^{q^-+1}\epsilon^-_{\gamma},\,\,\mathfrak v_{\delta^-_{q^-+1}+1}\prod\limits_{\gamma=0}^{q^-}\epsilon^-_{\gamma},\cdots,\mathfrak v_{s_{t_2}}\prod\limits_{\gamma=1}^{q^-+1}\epsilon^-_{\gamma}\right).\\  
\end{split}   
\end{equation*}
For $\alpha=\lambda^+_{\rho}$ where $1\leq\rho\leq q^+$, 
\begin{equation}\label{t1h1}\begin{split}
&\,\,\,\,\,\,\,\,\,\,\,\,\,\,\,\,\,\,\,\,\,\,\,\,\,\,\,\,\,\,\,\,\,\,\left(\eta^{\alpha}_{21},\eta^{\alpha}_{22},\cdots,\eta^{\alpha}_{2\left(\delta^+_{\rho}-1\right)},\eta^{\alpha}_{2\delta^+_{\rho}},\eta^{\alpha}_{2\left(\delta^+_{\rho}+1\right)},\cdots,\eta^{\alpha}_{2s_{t_{\alpha}^+}}\right)\mapsto\\
&\left({\mathfrak h}^{\alpha}_{21}\prod\limits_{\gamma=1}^{\rho}\epsilon^+_{\gamma},\,\,{\mathfrak h}^{\alpha}_{22}\prod\limits_{\gamma=1}^{\rho}\epsilon^+_{\gamma},\cdots,{\mathfrak h}^{\alpha}_{2\left(\delta^+_{\rho}-1\right)}\prod\limits_{\gamma=1}^{\rho}\epsilon^+_{\gamma},\,\,\prod\limits_{\gamma=1}^{\rho}\epsilon^+_{\gamma},\,\,{\mathfrak h}^{\alpha}_{2\left(\delta^+_{\rho}+1\right)}\prod\limits_{\gamma=1}^{\rho}\epsilon^+_{\gamma},\cdots,{\mathfrak h}^{\alpha}_{2s_{t_{\alpha}^+}}\prod\limits_{\gamma=1}^{\rho}\epsilon^+_{\gamma}\right),
\end{split}
\end{equation}
and for $\alpha=\lambda^+_{\rho}$ where $q^++1\leq\rho\leq l$,
\begin{equation}\label{t1h2}\begin{split}
&\left(\eta^{\alpha}_{21},\eta^{\alpha}_{22},\cdots,\eta^{\alpha}_{2s_{t_{\alpha}^+}}\right)\mapsto\left({\mathfrak h}^{\alpha}_{21}\prod\limits_{\gamma=1}^{q^++1}\epsilon^+_{\gamma},\,\,{\mathfrak h}^{\alpha}_{22}\prod\limits_{\gamma=1}^{q^++1}\epsilon^+_{\gamma},\cdots,{\mathfrak h}^{\alpha}_{2s_{t_{\alpha}^+}}\prod\limits_{\gamma=1}^{q^++1}\epsilon^+_{\gamma}\right).
\end{split}
\end{equation}
For $\beta=\lambda^-_{\rho}$ where $1\leq\rho\leq q^-$,
\begin{equation}
\begin{split}
&\,\,\,\,\,\,\,\,\,\,\,\,\,\,\,\,\,\,\,\,\,\,\,\,\,\,\,\,\,\,\,\,\,\,\left(\xi^{\beta}_{11},\xi^{\beta}_{12},\cdots,\xi^{\beta}_{1\left(\delta^-_{\rho}-1\right)},\xi^{\beta}_{1\delta^-_{\rho}},\xi^{\beta}_{1\left(\delta^-_{\rho}+1\right)},\cdots,\xi^{\beta}_{1s_{t_{\beta}^-}}\right)\mapsto\\
&\left({\mathfrak x}^{\beta}_{11}\prod\limits_{\gamma=1}^{\rho}\epsilon^-_{\gamma},\,\,{\mathfrak x}^{\beta}_{12}\prod\limits_{\gamma=1}^{\rho}\epsilon^-_{\gamma},\cdots,{\mathfrak x}^{\beta}_{1\left(\delta^-_{\rho}-1\right)}\prod\limits_{\gamma=1}^{\rho}\epsilon^-_{\gamma},\,\,\prod\limits_{\gamma=1}^{\rho}\epsilon^-_{\gamma},\,\,{\mathfrak x}^{\beta}_{1\left(\delta^-_{\rho}+1\right)}\prod\limits_{\gamma=1}^{\rho}\epsilon^-_{\gamma},\cdots,{\mathfrak x}^{\beta}_{1s_{t_{\beta}^-}}\prod\limits_{\gamma=1}^{\rho}\epsilon^-_{\gamma}\right),
\end{split}
\end{equation}
and for $\beta=\lambda^-_{\rho}$ where $q^-+1\leq\rho\leq m$,
\begin{equation*}
\begin{split}
&\left(\xi^{\beta}_{11},\xi^{\beta}_{12},\cdots,\xi^{\beta}_{1s_{t_{\beta}^-}}\right)\mapsto\left({\mathfrak x}^{\alpha}_{11}\prod\limits_{\gamma=1}^{q^-+1}\epsilon^-_{\gamma},\,\,{\mathfrak x}^{\alpha}_{12}\prod\limits_{\gamma=1}^{q^-+1}\epsilon^-_{\gamma},\cdots,{\mathfrak x}^{\beta}_{1s_{t_{\beta}^-}}\prod\limits_{\gamma=1}^{q^-+1}\epsilon^-_{\gamma}\right).
\end{split}
\end{equation*}

We can define a locally closed embedding $\Omega^{\lambda}$ from (\ref{bigsp}) to $\mathbb P^{N_{2,n}}\times\prod\nolimits_{t=1}^N\mathbb P^{N^{\underline s}_{t}}\times\prod\nolimits_{\substack{\underline w\,\,{\rm in\,\,(AB}{\tau}{\rm )}}}\mathbb P^{N^{\underline s}_{\underline w}}
$ by extending the rational map $\left(\widetilde{\mathcal K}^{\underline s}\circ \Gamma^{\tau}\circ\Sigma^{\lambda},\left(\cdots,F^{\underline s}_{\underline w}\circ e\circ\Gamma^{\tau}\circ \Sigma^{\lambda},\cdots\right)_{\underline w\,\,{\rm in\,\,(AB}{\tau}{\rm )}}\right)$ in the same way as Lemma \ref{em}.
Denote by $B^{\lambda}$ the corresponding image under $\Omega^{\lambda}$, which has the reduced scheme structure as a locally closed subscheme of $\left(\check{\mathcal R}_N\right)^{-1}\left(A^{\tau}\right)\subset\mathbb P^{N_{2,n}}\times\prod\nolimits_{t=1}^N\mathbb {P}^{N^{\underline s}_{t}}\times\prod\nolimits_{\substack{\underline w\,\,{\rm in\,\,(AB}{\tau}{\rm )}}}\mathbb P^{N^{\underline s}_{\underline w}}$.

\begin{lemma}\label{3fac}
For each $\lambda\in\Lambda^{\tau}$, $B^{\lambda}$ is an open subscheme of $\left(\check{\mathcal R}_N\right)^{-1}\left(A^{\tau}\right)$, and
\begin{equation*}
\bigcup\nolimits_{{\rm all\,\,indices\,\,}\lambda\in\Lambda^{\tau}}B^{\lambda}=\left(\check{\mathcal R}_N\right)^{-1}\left(A^{\tau}\right).     
\end{equation*}
In particular, $\left(\check{\mathcal R}_N\right)^{-1}\left(A^{\tau}\right)$ is smooth over ${\rm Spec}\,\mathbb Z$.
\end{lemma}
{\bf\noindent Proof of Lemma \ref{3fac}.} The idea is to apply the bubble sort algorithm in the same way as Lemma \ref{fac}. We omit the details here for brevity.\,\,\,\,\,\,$\endpf$

\subsection{Reduction to a local form and induction for Types (\ref{31mix}) and (\ref{32mix})}\label{UMX}
We modify the open cover $\left\{A^{\tau}\right\}_{\tau\in\mathbb J^{\underline s}}$ as follows. For each  $\tau=\left( j_1, j_2,\left(j^+_1,\cdots,j^+_l\right),\left(j^-_1,\cdots,j^-_{m}\right)\right)\in \mathbb J^{\underline s}$, let  $\mathring A^{\tau}$  be the open subscheme of $A^{\tau}$ defined by 
\begin{equation*}
\left\{\,\mathfrak p\in A^{\tau}\left|\,1-\eta^{\alpha}_{2p^+_{\alpha}}\cdot\xi^{\beta}_{1p^-_{\beta}}\notin\mathfrak p,\,\,\forall\,1\leq\alpha\leq l\,\,{\rm and}\,\,1\leq\beta\leq m\right.\right\}.  
\end{equation*}

By the same token, we can derive
\begin{definitionlemma}\label{3trun}
$\left\{\mathring A^{\tau}\right\}_{\tau\in\mathbb J^{\underline s}}$ is a finite open cover of  $\mathcal Q_{n}^{\underline s}$. We call $\mathring A^{\tau}$ the (truncated) coordinate charts of  $\mathcal Q_{n}^{\underline s}$.
\end{definitionlemma}

By analogy with Lemma \ref{red}, we can show that Proposition \ref{3ms} holds for a certain $N\geq 2$ if  the following lemma holds for the same $N$. Note that the lemma itself is an easy consequence of Proposition \ref{3ms}.

\begin{lemma}\label{3loc}  Let $N\geq 2$ be an integer. Let $\underline s\in\mathbb Z^N$ be a size vector.
For each (truncated) coordinate chart ${\mathring A^{\tau}}\subset\mathcal Q_{n}^{\underline s}$, there exists a finite open cover $\{T^{\tau}_{\alpha}\}_{\alpha}$ of  $\left(\overline{\mathcal R}^{\underline s}_n\right)^{-1}\left({\mathring A^{\tau}}\right)$ such that the following holds.

\begin{enumerate}[label=(\arabic*),ref=\arabic*]

\item Each $T^{\tau}_{\alpha}$ is an open subscheme of $\left(\overline{\mathcal R}^{\underline s}_n\right)^{-1}\left({\mathring A^{\tau}}\right)$, and is smooth over ${\rm Spec}\,\mathbb Z$.
\label{3loc1}   

\item  There is a locally closed subscheme $M^{\tau}_{\alpha}$ of $\prod\nolimits_{\underline w\in C^{\underline s}}\mathbb {P}^{N^{\underline s}_{\underline w}}$ for each $\alpha$, such that as a topological space $M^{\tau}_{\alpha}$ is the image of  $T^{\tau}_{\alpha}$ under the projection from $\mathbb {P}^{N_{2,n}}\times\prod\nolimits_{t=1}^N\mathbb {P}^{N^{\underline s}_{t}}\times\prod\nolimits_{\underline w\in C^{\underline s}}\mathbb {P}^{N^{\underline s}_{\underline w}}$ to  $\prod\nolimits_{\underline w\in C^{\underline s}}\mathbb {P}^{N^{\underline s}_{\underline w}}$. Moreover, $M^{\tau}_{\alpha}$ is smooth over ${\rm Spec}\,\mathbb Z$ of relative dimension  $2n-N-3$.
\label{3loc2}  

\item  The projection from $T^{\tau}_{\alpha}$ to $M^{\tau}_{\alpha}$ is a flat morphism.
\label{3loc3}  

\end{enumerate}
\end{lemma}

We define the type of  $\mathring A^{\tau}$ as follows.

\begin{enumerate}[label={\bf Type\,(I\arabic*)},ref=I\arabic*]

\item $\tau\in\mathbb J^{\underline s}$ is a Class I index such that $l,m\geq1$.
\label{31mix}
\end{enumerate}

\begin{enumerate}[label={\bf Type\,(I\arabic*)},ref=I\arabic*]
\setcounter{enumi}{1}

\item  $\tau\in\mathbb J^{\underline s}$ is a Class I index such that $m=0$ or $l=0$. 
\label{31up}
\end{enumerate}

\begin{enumerate}[label={\bf Type\,(II\arabic*)},ref=II\arabic*]

\item $\tau\in\mathbb J^{\underline s}$ is a Class II index such that $l,m\geq1$. 
\label{32mix}
\end{enumerate}

\begin{enumerate}[label={\bf Type\,(II\arabic*)},ref=II\arabic*]
\setcounter{enumi}{1}

\item  $\tau\in\mathbb J^{\underline s}$ is a Class II index such that $m=0$ or $l=0$.
\label{32up}
\end{enumerate}

\begin{lemma}\label{3uv}
Suppose that  Proposition \ref{3ms} holds for all integers $N^{\prime}$ such that $2\leq N^{\prime}<N$. Then Lemma \ref{3loc} holds for all (truncated) coordinate charts of {\bf Type (\ref{31mix})}.
\end{lemma}
{\bf\noindent Proof of Lemma \ref{3uv}.} Let $\tau=\left(j_1,j_2,(j^+_1,\cdots,j^+_l),(j^-_1,\cdots,j^-_{m})\right)\in\mathbb J^{\underline s}$ be an arbitrary Class I index such that $l,m\geq1$. We write vectors ${\underline w}\in C^{\underline s}$ in form (\ref{wc}).  

To define a total order $\sigma$ on $C^{\underline s}$, we partition $C^{\underline s}$ as follows.

\begin{enumerate}[label={(AB$\tau$)},ref={AB$\tau$}]
\item The union of subsets (\ref{fa0}), (\ref{fa1}), (\ref{fa2}), (\ref{fb0}), (\ref{fb1}), (\ref{fb2}) defined in \S  \ref{3FAC}.
\end{enumerate}
\begin{enumerate}[label={(AB\arabic*)},ref={AB\arabic*}]
\setcounter{enumi}{-1}
\item $w_{t_1}=w_{t_2}=1$, and the other components are $0$.
\label{ab0}
\end{enumerate}
\begin{enumerate}[label={(A\arabic*)},ref={A\arabic*}]
\setcounter{enumi}{2}
\item $w_{t^+_{\alpha}}=2$, for a certain $1\leq \alpha\leq l$, and other components are $0$.
\label{ab1}

\item $w_{t^+_{\alpha_1}}=w_{t^+_{\alpha_2}}=1$, for certain $1\leq \alpha_1<\alpha_2\leq l$, and other components are $0$.
\label{ab2}

\item $w_{t_1}=w_{t^+_{\alpha_1}}=w_{t^+_{\alpha_2}}=w_{t^+_{\alpha_3}}=1$, for certain $1\leq \alpha_1<\alpha_2<\alpha_3\leq l$.
\label{ab3}

\item $w_{t^+_{\alpha_1}}=w_{t^+_{\alpha_2}}=w_{t^+_{\alpha_3}}=w_{t^+_{\alpha_4}}=1$, for certain $1\leq \alpha_1<\alpha_2<\alpha_3<\alpha_4\leq l$.
\label{ab4}

\end{enumerate}
\begin{enumerate}[label={(B\arabic*)},ref={B\arabic*}]
\setcounter{enumi}{2} 
\item $w_{t^-_{\beta}}=2$, for a certain $1\leq \beta\leq m$, and other components are $0$. 

\item $w_{t^-_{\beta_1}}=w_{t^-_{\beta_2}}=1$, for certain $1\leq \beta_1<\beta_2\leq m$, and other components are $0$.

\item $w_{t_2}=w_{t^-_{\beta_1}}=w_{t^-_{\beta_2}}=w_{t^-_{\beta_3}}=1$, for certain $1\leq \beta_1<\beta_2<\beta_3\leq m$.

\item $w_{t^-_{\beta_1}}=w_{t^-_{\beta_2}}=w_{t^-_{\beta_3}}=w_{t^-_{\beta_4}}=1$, for certain $1\leq \beta_1<\beta_2<\beta_3<\beta_4\leq m$.
\end{enumerate}
\begin{enumerate}[label={(C\arabic*)},ref={C\arabic*}]
\setcounter{enumi}{0}
\item $w_{t_1}=w_{t^-_{\beta}}=1$, for a  certain $1\leq \beta\leq m$, and other components are $0$.

\item $w_{t_2}=w_{t^+_{\alpha}}=1$, for a certain $1\leq \alpha\leq l$, and other components are $0$.
\label{ab5}

\item $w_{t^+_{\alpha}}=w_{t^-_{\beta}}=1$, for certain $1\leq \alpha\leq l$ and $1\leq \beta\leq m$, and other components are $0$.

\item $w_{t_1}=w_{t^-_{\beta_1}}=w_{t^-_{\beta_2}}=w_{t^-_{\beta_3}}=1$, for  certain $1\leq \beta_1<\beta_2<\beta_3\leq m$.

\item $w_{t_2}=w_{t^+_{\alpha_1}}=w_{t^+_{\alpha_2}}=w_{t^+_{\alpha_3}}=1$, for certain $1\leq \alpha_1<\alpha_2<\alpha_3\leq l$.
\label{ab6}

\item $w_{t_1}=w_{t_2}=w_{t^+_{\alpha}}=w_{t^-_{\beta}}=1$, for  certain $1\leq \alpha\leq l$ and  $1\leq \beta\leq m$.

\item $w_{t_1}=w_{t^+_{\alpha_1}}=w_{t^+_{\alpha_2}}=w_{t^-_{\beta}}=1$, for certain $1\leq \alpha_1<\alpha_2\leq l$ and $1\leq \beta\leq m$.

\item $w_{t_1}=w_{t^+_{\alpha}}=w_{t^-_{\beta_1}}=w_{t^-_{\beta_2}}=1$, for  certain $1\leq \alpha\leq l$ and  $1\leq \beta_1<\beta_2\leq m$.

\item $w_{t_2}=w_{t^+_{\alpha_1}}=w_{t^+_{\alpha_2}}=w_{t^-_{\beta}}=1$, for certain $1\leq \alpha_1<\alpha_2\leq l$ and $1\leq \beta\leq m$.

\item $w_{t_2}=w_{t^+_{\alpha}}=w_{t^-_{\beta_1}}=w_{t^-_{\beta_2}}=1$, for certain $1\leq \alpha\leq l$ and  $1\leq \beta_1<\beta_2\leq m$.

\item $w_{t^+_{\alpha_1}}=w_{t^+_{\alpha_2}}=w_{t^-_{\beta_1}}=w_{t^-_{\beta_2}}=1$, for  certain $1\leq \alpha_1<\alpha_2\leq l$ and $1\leq \beta_1<\beta_2\leq m$.

\item $w_{t^+_{\alpha_1}}=w_{t^+_{\alpha_2}}=w_{t^+_{\alpha_3}}=w_{t^-_{\beta}}=1$, for certain $1\leq \alpha_1<\alpha_2<\alpha_3\leq l$ and  $1\leq \beta\leq m$.

\item $w_{t^+_{\alpha}}=w_{t^-_{\beta_1}}=w_{t^-_{\beta_2}}=w_{t^-_{\beta_3}}=1$, for certain $1\leq \alpha\leq l$ and  $1\leq \beta_1<\beta_2<\beta_3\leq m$.

\end{enumerate}
Define $\sigma$ so that the vectors in an upper subset are strictly smaller than the ones in a lower subset. Equivalently, we first blow up with respect to the (pull-backs of) ideal sheaves $\mathscr I^{\underline s}_{\underline w}$ with $\underline w$ in (AB$\tau$), (AB0), then (A3)-(A6),  (B3)-(B6), and finally (C1)-(C13).

It is clear that the blow-up with respect to the (pull-backs of) ideal sheaves $\mathscr I^{\underline s}_{\underline w}$ with $\underline w$ in (AB$\tau$), (AB0) yields $\check{\mathcal R}_N:W^{\sigma}_{N_{\tau}}\rightarrow \mathcal Q_{n}^{\underline s}$. Denote by $\mathcal R_{AB}$ the sequence of blow-ups with respect to the (pull-backs of) ideal sheaves $\mathscr I^{\underline s}_{\underline w}$ with $\underline w$ in (A3)-(A6) and (B3)-(B6); denote by $\mathcal R_{ABC}$ the sequence of blow-ups with respect to the (pull-backs of) ideal sheaves $\mathscr I^{\underline s}_{\underline w}$ with $\underline w$ in (C1)-(C13).

For each index  $\lambda\in\Lambda^{\tau}$, let $B^{\lambda}$ be the open subscheme of $\left(\check{\mathcal R}_N\right)^{-1}\left(A^{\tau}\right)$ defined in Lemma \ref{3fac}.

\begin{claim}\label{13gi}
The restriction of  $\mathcal R_{ABC}$ to
$\left(\mathcal R_{AB}\circ\mathcal R_{ABC}\right)^{-1}\left(B^{\lambda}\cap \left(\check{\mathcal R}_N\right)^{-1}\left({\mathring A^{\tau}}\right)\right)$ is an isomorphism.
\end{claim}

{\bf\noindent Proof of Claim \ref{13gi}.}  It suffices to show that, for any  ${\underline w}$ in (C1)-(C13), the restriction to 
$\left(\mathcal R_{AB}\right)^{-1}\left(B^{\lambda}\cap \left(\check{\mathcal R}_N\right)^{-1}\left({\mathring A^{\tau}}\right)\right)$ of the ideal sheaf $\left(\widetilde{\mathcal R}^{\underline s}\circ\check{\mathcal R}_N\circ\mathcal R_{AB}\right)^{-1}\mathscr I^{\underline s}_{\underline w}\cdot\mathcal O_{\left(\mathcal R_{AB}\right)^{-1}\left(B^{\lambda}\right)}$ is invertible.

By (\ref{348})-(\ref{3xb}), we can derive that  $\left(\widetilde{\mathcal R}^{\underline s}\circ\check{\mathcal R}_N\right)^{-1}\mathscr I^{\underline s}_{\underline w}\cdot\mathcal O_{B^{\lambda}\cap \left(\check{\mathcal R}_N\right)^{-1}\left({\mathring A^{\tau}}\right)}$ for ${\underline w}$ in (C1), (C2), (C3), (C6), (C8), (C9), (C11) are the invertible ideal sheaves generated by $a_{j^-_{\beta}}$, $a_{j^+_{\alpha}}$, $a_{j^+_{\alpha}}a_{j^-_{\beta}}$, $a_{j^+_{\alpha}}a_{j^-_{\beta}}$, $a_{j^+_{\alpha}}a_{j^-_{\beta_1}}a_{j^-_{\beta_2}}$, $a_{j^+_{\alpha}}a_{j^+_{\alpha_2}}a_{j^-_{\beta_1}}$, $a_{j^+_{\alpha_1}}a_{j^+_{\alpha_2}}a_{j^-_{\beta_1}}a_{j^-_{\beta_2}}$, respectively. For ${\underline w}$ in (C7), $\left(\widetilde{\mathcal R}^{\underline s}\circ\check{\mathcal R}_N\right)^{-1}\mathscr I^{\underline s}_{\underline w}\cdot\mathcal O_{B^{\lambda}}$ is the invertible ideal sheaf generated by
\begin{equation*}
a_{j^+_{\alpha_1}}a_{j^+_{\alpha_2}}a_{j^-_{\beta}}\prod\nolimits_{\gamma=1}^{\min\left\{q^++1,\,\,\rho_1\,\,{\rm such\,\,that\,\,}\lambda^+_{\rho_1}=a_1,\,\,\rho_2\,\,{\rm such\,\,that\,\,}\lambda^+_{\rho_2}=a_2\right\}}\epsilon^+_{\gamma}.   
\end{equation*}
Similar results hold for ${\underline w}$ in (C10) as well. 

For ${\underline w}$ in (C12),  we adopt the convention that $\eta_{1p^+_{\alpha}}^{\alpha}=\mathfrak h^{\lambda^+_{\rho}}_{2\delta^+_{\rho}}=1$ for $1\leq\alpha\leq l$ and $1\leq\rho\leq q^+$.  
Write $\alpha_1=\lambda^+_{\rho_1}$, $\alpha_2=\lambda^+_{\rho_2}$, $\alpha_3=\lambda^+_{\rho_3}$, and, without loss of generality,  assume that $\rho_1<\rho_2<\rho_3$.  It is clear that \begin{equation*}
\left(\widetilde{\mathcal R}^{\underline s}\circ\check{\mathcal R}_N\right)^{-1}\mathscr I^{\underline s}_{\underline w}\cdot\mathcal O_{B^{\lambda}}=\left(\prod\nolimits_{\gamma=1}^{\min\left\{q^++1,\,\,\rho_1\right\}}\epsilon^+_{\gamma}\right)\cdot\left(a_{j^+_{\alpha_1}}a_{j^+_{\alpha_2}}a_{j^+_{\alpha_3}}a_{j^-_{\beta}}\right)\cdot\mathfrak I,
\end{equation*} where $\mathfrak I$ is the ideal sheaf generated by  the following sub-determinants for all $1\leq q_1\leq s_{t_{\alpha_1}^+}$, $1\leq q_2\leq s_{t_{\alpha_2}^+}$, $1\leq q_3\leq s_{t_{\alpha_3}^+}$.
\begin{equation}\label{3il5}
\left|\begin{matrix}
\eta^{\alpha_1}_{1q_1}&\eta^{\alpha_2}_{1q_2}\\
{\mathfrak h}^{\alpha_1}_{2q_1}&\left(\prod\limits_{\gamma=\min\left\{q^++1,\,\,\rho_1\right\}+1}^{\min\left\{q^++1,\,\,\rho_2\right\}}\epsilon^+_{\gamma}\right)\cdot{\mathfrak h}^{\alpha_2}_{2q_2}\\
\end{matrix}\right|=:j^{12}_{q_1q_2}, \,\,\,\,\,\,\,\left|\begin{matrix}
\eta^{\alpha_1}_{1q_1}&\eta^{\alpha_3}_{1q_3}\\
{\mathfrak h}^{\alpha_1}_{2q_1}&\left(\prod\limits_{\gamma=\min\left\{q^++1,\,\,\rho_1\right\}+1}^{\min\left\{q^++1,\,\,\rho_3\right\}}\epsilon^+_{\gamma}\right)\cdot{\mathfrak h}^{\alpha_3}_{2q_3}\\
\end{matrix}\right|=:j^{13}_{q_1q_3}, 
\end{equation}
and
\begin{equation}\label{4il5}
\begin{split}
&\,\,\,\,\,\,\,\,\,\,\,\,\,\,\,\,\,\,\,\,\,\,\,\,\,\,\,\,\,\,\,\,\,\,\,\,\,\,\,\,\left|\begin{matrix}
\eta^{\alpha_2}_{1q_2}&\eta^{\alpha_3}_{1q_3}\\
\left(\prod\limits_{\gamma=\min\left\{q^++1,\,\,\rho_1\right\}+1}^{\min\left\{q^++1,\,\,\rho_2\right\}}\epsilon^+_{\gamma}\right)\cdot{\mathfrak h}^{\alpha_2}_{2q_2}&\left(\prod\limits_{\gamma=\min\left\{q^++1,\,\,\rho_1\right\}+1}^{\min\left\{q^++1,\,\,\rho_3\right\}}\epsilon^+_{\gamma}\right)\cdot{\mathfrak h}^{\alpha_3}_{2q_3}\\
\end{matrix}\right|=\\
&\left(\prod\limits_{\gamma=\min\left\{q^++1,\,\,\rho_1\right\}+1}^{\min\left\{q^++1,\,\,\rho_2\right\}}\epsilon^+_{\gamma}\right)\cdot\left|\begin{matrix}
\eta^{\alpha_2}_{1q_2}&\eta^{\alpha_3}_{1q_3}\\
{\mathfrak h}^{\alpha_2}_{2q_2}&\left(\prod\limits_{\gamma=\min\left\{q^++1,\,\,\rho_2\right\}+1}^{\min\left\{q^++1,\,\,\rho_3\right\}}\epsilon^+_{\gamma}\right)\cdot{\mathfrak h}^{\alpha_3}_{2q_3}\\
\end{matrix}\right|=:\left(\prod\limits_{\gamma=\min\left\{q^++1,\,\,\rho_1\right\}+1}^{\min\left\{q^++1,\,\,\rho_2\right\}}\epsilon^+_{\gamma}\right)\cdot j^{23}_{q_2q_3}.\\
\end{split} 
\end{equation}
Notice that by setting $q_1=\delta_{\rho_1}^+$ and $q_2=p^+_{\alpha_2}$ in $j^{12}_{q_1q_2}$, we can derive
\begin{equation*}
\begin{split}
&j^{23}_{q_2q_3}=-j^{12}_{\delta_{\rho_1}^+p^+_{\alpha_2}}\cdot j^{23}_{q_2q_3}-\eta^{\alpha_1}_{1\delta_{\rho_1}^+}\cdot{\mathfrak h}^{\alpha_2}_{2p^+_{\alpha_2}}\left(\left(\prod\nolimits_{\gamma=\min\left\{q^++1,\,\,\rho_1\right\}+1}^{\min\left\{q^++1,\,\,\rho_2\right\}}\epsilon^+_{\gamma}\right)\cdot j^{23}_{q_2q_3}\right).\\
\end{split}    
\end{equation*}
We can thus conclude that $\mathfrak I$ is generated by $j^{13}_{q_1q_3}$, $j^{23}_{q_2q_3}$, $j^{12}_{q_1q_2}$, where $q_1,q_2,q_3$ run over all $1\leq q_1\leq s_{t_{\alpha_1}^+}$, $1\leq q_2\leq s_{t_{\alpha_2}^+}$, $1\leq q_3\leq s_{t_{\alpha_3}^+}$.

Take  ${\underline w}^*\in C^{\underline s}$ such that $w^*_{t_1}=w^*_{t^+_{\alpha_1}}=w^*_{t^+_{\alpha_2}}=w^*_{t^+_{\alpha_3}}=1$.   Then,
\begin{equation*}
\left(\widetilde{\mathcal R}^{\underline s}\circ\check{\mathcal R}_N\right)^{-1}\mathscr I^{\underline s}_{\underline w^*}\cdot\mathcal O_{B^{\lambda}}=\left(\prod\limits_{\gamma=1}^{\min\left\{q^++1,\,\,\rho_1\right\}}\epsilon^+_{\gamma}\right)\cdot\left(\prod\limits_{\gamma=1}^{\min\left\{q^++1,\,\,\rho_2\right\}}\epsilon^+_{\gamma}\right)\cdot\left(a_{j^+_{\alpha_1}}a_{j^+_{\alpha_2}}a_{j^+_{\alpha_3}}\right)\cdot\mathfrak I^*,
\end{equation*} where $\mathfrak I^*$ is the ideal sheaf generated by $j^{13}_{q_1q_3}$, $j^{23}_{q_2q_3}$, and $\left(\prod\nolimits_{\gamma=\min\left\{q^++1,\,\,\rho_2\right\}+1}^{\min\left\{q^++1,\,\,\rho_3\right\}}\epsilon^+_{\gamma}\right)\cdot j^{12}_{q_1q_2}$, for all $1\leq q_1\leq s_{t_{\alpha_1}^+}$, $1\leq q_2\leq s_{t_{\alpha_2}^+}$, $1\leq q_3\leq s_{t_{\alpha_3}^+}$ (see (\ref{3il5} and (\ref{4il5})). Notice that by setting $q_2=\delta_{\rho_2}^+$ and $q_3=p^+_{\alpha_3}$ in $j^{23}_{q_2q_3}$, we can derive
\begin{equation*}
\begin{split}
&j^{12}_{q_1q_2}=j^{23}_{\delta_{\rho_2}^+p^+_{\alpha_3}}\cdot j^{12}_{q_1q_2}-\eta^{\alpha_1}_{1q_1}\cdot{\mathfrak h}^{\alpha_1}_{2q_1}\left(\left(\prod\nolimits_{\gamma=\min\left\{q^++1,\,\,\rho_2\right\}+1}^{\min\left\{q^++1,\,\,\rho_3\right\}}\epsilon^+_{\gamma}\right)\cdot j^{12}_{q_1q_2}\right).\\
\end{split}    
\end{equation*}
We can thus conclude that $\mathfrak I^*$ is generated by $j^{13}_{q_1q_3}$, $j^{23}_{q_2q_3}$, $j^{12}_{q_1q_2}$, $1\leq q_1\leq s_{t_{\alpha_1}^+}$, $1\leq q_2\leq s_{t_{\alpha_2}^+}$, $1\leq q_3\leq s_{t_{\alpha_3}^+}$, and hence $\mathfrak I=\mathfrak I^*$. Since the pull-back of $\left(\widetilde{\mathcal R}^{\underline s}\circ\check{\mathcal R}_N\right)^{-1}\mathscr I^{\underline s}_{\underline w^*}\cdot\mathcal O_{B^{\lambda}}$ under the blow-up $\mathcal R_{AB}$ is invertible, the ideal sheaf $\left(\widetilde{\mathcal R}^{\underline s}\circ\check{\mathcal R}_N\circ\mathcal R_{AB}\right)^{-1}\mathscr I^{\underline s}_{\underline w}\cdot\mathcal O_{\left(\mathcal R_{AB}\right)^{-1}\left(B^{\lambda}\cap \left(\check{\mathcal R}_N\right)^{-1}\left({\mathring A^{\tau}}\right)\right)}$ is invertible. Similar results hold for ${\underline w}$ in (C4), (C5), (C13) as well.

We complete the proof of Claim \ref{13gi}.\,\,\,\,$\endpf$
\smallskip

Then to prove Lemma \ref{3uv}, it suffices to show that for each $\lambda\in\Lambda^{\tau}$ the following holds. 

\begin{enumerate}[label={\bf(P\arabic*)},ref={\bf P\arabic*}]
\item $\left(\mathcal R_{AB}\right)^{-1}\left(B^{\lambda}\cap \left(\check{\mathcal R}_N\right)^{-1}\left({\mathring A^{\tau}}\right)\right)$
is  smooth over ${\rm Spec}\,\mathbb Z$.
\label{3one}

\item Under the projection from $\mathbb {P}^{N_{2,n}}\times\prod\nolimits_{t=1}^N\mathbb {P}^{N^{\underline s}_{t}}\times\prod\nolimits_{\underline w\in C^{\underline s}}\mathbb {P}^{N^{\underline s}_{\underline w}}$  to $\prod\nolimits_{\underline w\in C^{\underline s}}\mathbb {P}^{N^{\underline s}_{\underline w}}$ the image of \begin{equation}\label{3rabc}
\left(\mathcal R_{AB}\circ\mathcal R_{ABC}\right)^{-1}\left(B^{\lambda}\cap \left(\check{\mathcal R}_N\right)^{-1}\left({\mathring A^{\tau}}\right)\right)
\end{equation} has a reduced scheme structure as a locally closed subscheme of $\prod\nolimits_{\underline w\in C^{\underline s}}\mathbb {P}^{N^{\underline s}_{\underline w}}$; the subscheme is smooth over ${\rm Spec}\,\mathbb Z$ of relative dimension  $2n-N-3$; the projection from (\ref{3rabc}) to its image is a flat morphism.

\end{enumerate}
\medskip

Notice that the constructions of $\mathcal Q_{n}^{\underline s}$, $\mathcal T_{n}^{\underline s}$, $\mathcal M_{n}^{\underline s}$ are invariant under the automorphism induced by the permutations of blocks and that of columns within a block. Therefore, we may assume that $t_1=1$, $t_2=2$, and $p_1=p_2=p^+_{1}=\cdots=p^+_{l}=p^-_{1}=\cdots=p^-_{m}=1$ in the following.

Define indices $\widetilde{\underline s}:=\left(\widetilde s_1,\widetilde s_2,\cdots,\widetilde s_{l+2}\right)$ and  $\widehat{\underline s}:=\left(\widehat s_1,\widehat s_2,\cdots,\widehat s_{m+2}\right)$ respectively by 
\begin{equation*}
\widetilde s_1:=s_1,\,\,\,\,\,\,\widetilde s_2:=1,\,\,\,\,\,\,\widetilde s_{\alpha+2}:=s_{t^+_{\alpha}},\,\,\,{\rm for}\,\,1\leq\alpha\leq l,
\end{equation*}
and
\begin{equation*}
\widehat s_1:=1,\,\,\,\,\,\,\widehat s_2:=s_2,\,\,\,\,\,\,\widehat s_{\beta+2}:=s_{t^-_{\beta}},\,\,\,{\rm for}\,\,1\leq\beta\leq m.
\end{equation*}
Define $\widetilde\tau:=\left(\widetilde j_1,\widetilde j_2,\left(\widetilde j^+_{1},\cdots,\widetilde j^+_{l}\right)\right)\in\mathbb J^{\widetilde{\underline s}}$ and  $\widehat\tau:=\left(\widehat j_1,\widehat j_2,\left(\widehat j^-_{1},\cdots,\widehat j^-_{m}\right)\right)\in\mathbb J^{\widehat{\underline s}}$ respectively by
\begin{equation*}
\begin{aligned}
&\widetilde j_1:=j_1,\,\,\,\,\,\,\widetilde j_2:=s_1+1,\,\,\,\,\,\,\widetilde j_{\alpha}^+:=s_1+1+s_{t^+_1}+s_{t^+_2}+\cdots+s_{t^+_{\alpha-1}}+p_{\alpha},\,\,\,{\rm for}\,\,1\leq\alpha\leq l\\
\end{aligned},    
\end{equation*}
and
\begin{equation*}
\begin{aligned}
&\widehat j_1:=1,\,\,\,\,\,\,\widehat j_2:=1+p_2,\,\,\,\,\,\,\widehat j_{\beta}^-:=1+s_2+s_{t^-_1}+s_{t^-_2}+\cdots+s_{t^-_{\beta-1}}+p_{\beta},\,\,\,{\rm for}\,\,1\leq\beta\leq m\\
\end{aligned}.
\end{equation*}
Let $A^{\widetilde\tau}$ and $A^{\widehat\tau}$ respectively be the open subschemes of $\mathcal Q^{\widetilde{\underline s}}_{\widetilde n}$ and $\mathcal Q^{\widehat{\underline s}}_{\widehat n}$, where  $\widetilde n:=\widetilde s_1+\widetilde s_2+\cdots+\widetilde s_{l+2}$ and $\widehat n:=\widehat s_1+\widehat s_2+\cdots+\widehat s_{m+2}$.
Let $\check{\mathcal R}_{l+2}$, $\ddot{\mathcal R}_{l+2}$ be the blow-ups defined by setting $N=l+2$ in (\ref{312a1}); similarly, we have $\check{\mathcal R}_{m+2}$, $\ddot{\mathcal R}_{m+2}$. Define indices  $\widetilde\lambda:=\left(\widetilde q^+,\widetilde q^-,\widetilde\lambda^+,\widetilde\Delta^+\right)\in\Lambda^{\widetilde\tau}$ and $\widehat\lambda:=\left(\widehat q^+,\widehat q^-,\widehat\lambda^-,\widehat\Delta^-\right)\in\Lambda^{\widehat\tau}$ respectively by
\begin{equation*}
\widetilde q^+:=q^+,\,\,\,\,\,\,\widetilde q^-:=0,\,\,\,\,\,\,\widetilde \lambda^+:=\lambda^+,\,\,\,\,\,\,\widetilde \Delta^+:=\Delta^+,
\end{equation*}
and
\begin{equation*}
\widehat q^+:=0,\,\,\,\,\,\,\widehat q^-:=q^-,\,\,\,\,\,\,\widehat\lambda^-:=\lambda^-,\,\,\,\,\,\,\widehat\Delta^-:=\Delta^-.
\end{equation*}
Let $B^{\widetilde\lambda}$ and  $B^{\widehat\lambda}$ respectively be the open subschemes of $\left(\check{\mathcal R}_{l+2}\right)^{-1}\left(A^{\widetilde\tau}\right)$ and $\left(\check{\mathcal R}_{m+2}\right)^{-1}\left(A^{\widehat\tau}\right)$ defined by Lemma \ref{3fac}.

Notice that the ideal sheaf $\left(\widetilde{\mathcal R}^{\underline s}\circ\check{\mathcal R}_N\right)^{-1}\mathscr I^{\underline s}_{\underline w}\cdot\mathcal O_{B^{\lambda}}$ associated to vectors $\underline w$ in (A0)-(A6) can be generated by polynomials  independent of  $\epsilon_1^+,a_{j^-_1},\cdots,a_{j^-_m}$, $\overrightarrow E_-,\overrightarrow Z,\overrightarrow {\mathfrak V},\overrightarrow {\mathscr X}^1,\cdots,\overrightarrow {\mathscr X}^m,\overrightarrow {\mathfrak X}^1,$ $\cdots,\overrightarrow {\mathfrak X}^m$, and the ones associated to vectors $\underline w$ in (B0)-(B6) can be generated by polynomials independent of  $\epsilon_1^-,a_{j^+_1},\cdots,a_{j^+_l}$, $\overrightarrow E_+,\overrightarrow Y,\overrightarrow {\mathfrak U},\overrightarrow {\mathscr H}^1,\cdots,\overrightarrow {\mathscr H}^l,\overrightarrow {\mathfrak H}^1,\cdots,\overrightarrow {\mathfrak H}^l$. 
Similarly to \ref{rab}, we can derive that 
\begin{equation}\label{3dec}
\left(\mathcal R_{AB}\right)^{-1}\left(B^{\lambda}\right)\cong\left(\ddot{\mathcal R}_{l+2}\right)^{-1}\left(B^{\widetilde\lambda}\right)\times\left(\ddot{\mathcal R}_{m+2}\right)^{-1}\left(B^{\widehat\lambda}\right).
\end{equation}
By the assumption that Proposition \ref{3ms} holds for all $2\leq N^{\prime}\leq N-1$, we have that
\begin{enumerate}[label={$\bullet$}]
\item $\mathcal T_{l+2}^{\widetilde{\underline s}}$, $\mathcal T_{m+2}^{\widehat{\underline s}}$, $\mathcal M_{l+2}^{\widetilde{\underline s}}$, and $\mathcal M_{m+2}^{\widehat{\underline s}}$ are smooth over ${\rm Spec}\,\mathbb Z$,
    
\item the projections $\mathcal P_{l+2}^{\widetilde{\underline s}}:\mathcal T_{l+2}^{\widetilde{\underline s}}\rightarrow \mathcal M_{l+2}^{\widetilde{\underline s}}$ and $\mathcal P_{m+2}^{\widehat{\underline s}}:\mathcal T_{m+2}^{\widehat{\underline s}}\rightarrow \mathcal M_{m+2}^{\widehat{\underline s}}$ are flat.
    
\end{enumerate}
In particular, (\ref{3one}) holds.
\smallskip 

To derive a similar result as (\ref{cf}), we define a vector  ${\underline w}^*\in C^{\underline s}$ as follows.
\begin{enumerate}[label={$\bullet$}]
\item When $q^+=q^-=0$,  $w^*_{1}=w^*_{2}=1$,  and other components are $0$.

\item When $q^+=0$ and $q^-\geq1$, $w^*_{1}=w^*_{t^-_{\lambda^-_{1}}}=1$,  and other components are $0$.

\item When $q^+\geq1$ and $q^-=0$, $w^*_{2}=w^*_{t^+_{\lambda^+_{1}}}=1$,  and other components are $0$.

\item When $q^+\geq1$ and $q^-\geq1$, $w^*_1=w^*_2=w^*w^*_{t^+_{\lambda^+_{1}}}=w^*w^*_{t^-_{\lambda^-_{1}}}=1$.

\end{enumerate}
Write the homogeneous coordinates for $\mathbb P^{N^{\underline s}_{\underline w^*}}$ as $[\cdots,u_{\gamma},\cdots]_{\gamma=0}^{N^{\underline s}_{\underline w^*}}$. Then the rational map 
\begin{equation*}
\left.\left(F^{\underline s}_{\underline w^*}\circ e\circ\widetilde  {\mathcal R}^{\underline s}\circ\check{\mathcal R}_N\right)\right|_{B^{\lambda}}: B^{\lambda}\dashrightarrow \mathbb P^{N^{\underline s}_{\underline w^*}}  
\end{equation*}
can be given by
\begin{equation*}
\begin{split}
&[\cdots,u_{\gamma},\cdots]_{\gamma=0}^{N^{\underline s}_{\underline w^*}}\mapsto\left[\cdots,1,\cdots,\pm\left(\epsilon^+_0\cdot\epsilon^-_0\right)+h_{\underline w^*}\left(\overrightarrow Y,\overrightarrow Z,\overrightarrow {\mathscr H}^1,\cdots,\overrightarrow {\mathscr H}^l,\overrightarrow {\mathscr X}^1,\cdots,\overrightarrow {\mathscr X}^m\right),\cdots\right],\\
\end{split}
\end{equation*}
where one  homogeneous coordinate takes value $1$, one takes value $\pm\left(\epsilon^+_0\cdot\epsilon^-_0\right)+h_{\underline w^*}$ for a certain polynomials $h_{\underline w^*}$ in variables $\overrightarrow Y$, $\overrightarrow Z$, $\overrightarrow {\mathscr H}^1,\cdots,\overrightarrow {\mathscr H}^l$, $\overrightarrow {\mathscr X}^1$, $\cdots$, $\overrightarrow {\mathscr X}^m$.

We can then complete the proof of Lemma \ref{3uv} in the same way as Lemma \ref{uv}. \,\,\,\,\,\,$\endpf$

\begin{lemma}\label{3uv32}Suppose that  Proposition \ref{3ms} holds for all integers $N^{\prime}$ such that $2\leq N^{\prime}<N$. Then Lemma \ref{3loc} holds for all (truncated) coordinate charts of {\bf Type (\ref{32mix})}.
\end{lemma}
{\bf\noindent Proof of Lemma \ref{3uv32}.}  Fix a Class II index $\tau=\left(j_1,j_2,(j^+_1,\cdots,j^+_l),(j^-_1,\cdots,j^-_{m})\right)\in\mathbb J^{\underline s}$. Without loss of generality, we may assume that  $t_1=t_2=1$, $p_1=1$, $p_2=2$, $p^+_{1}=\cdots=p^+_{l}=p^-_{1}=\cdots=p^-_{m}=1$. Define indices $\widetilde{\underline s}:=\left(\widetilde s_1,\widetilde s_2,\cdots,\widetilde s_{l+1}\right)$ and  $\widehat{\underline s}:=\left(\widehat s_1,\widehat s_2,\cdots,\widehat s_{m+1}\right)$ respectively by 
\begin{equation*}\widetilde s_1:=2,\,\,\,\,\,\,\widetilde s_{\alpha+1}:=s_{t^+_{\alpha}},\,\,\,{\rm for}\,\,1\leq\alpha\leq l,
\end{equation*} and
\begin{equation*}
\widehat s_1=2,\,\,\,\,\,\,\widehat s_{\beta+1}:=s_{t^-_{\beta}},\,\,\,{\rm for}\,\,1\leq\beta\leq m.
\end{equation*}
Let $\widetilde n:=\widetilde s_1+\widetilde s_2+\cdots+\widetilde s_{l+1}$ and $\widehat n:=\widehat s_1+\widehat s_2+\cdots+\widehat s_{m+1}$. Define $\widetilde\tau:=\left(\widetilde j_1,\widetilde j_2,\left(\widetilde j^+_{1},\cdots,\widetilde j^+_{l}\right)\right)\in\mathbb J^{\widetilde{\underline s}}$ and  $\widehat\tau:=\left(\widehat j_1,\widehat j_2,\left(\widehat j^-_{1},\cdots,\widehat j^-_{m}\right)\right)\in\mathbb J^{\widehat{\underline s}}$ respectively by
\begin{equation*}
\begin{aligned}
&\widetilde j_1:=1,\,\,\,\,\,\,\widetilde j_2:=2,\,\,\,\,\,\,\widetilde j_{\alpha}^+:=s_1+s_{t^+_1}+s_{t^+_2}+\cdots+s_{t^+_{\alpha-1}}+p_{\alpha},\,\,\,{\rm for}\,\,1\leq\alpha\leq l\\
\end{aligned},  
\end{equation*}
and
\begin{equation*}
\begin{aligned}
&\widehat j_1:=1,\,\,\,\,\,\,\widehat j_2:=2,\,\,\,\,\,\,\widehat j_{\beta}^-:=s_1+s_{t^-_1}+s_{t^-_2}+\cdots+s_{t^-_{\beta-1}}+p_{\beta},\,\,\,{\rm for}\,\,1\leq\beta\leq m\\
\end{aligned}.
\end{equation*}
Let $A^{\widetilde\tau}$ and $A^{\widehat\tau}$ respectively be the open subschemes of $\mathcal Q^{\widetilde{\underline s}}_{\widetilde n}$ and $\mathcal Q^{\widehat{\underline s}}_{\widehat n}$.

Notice that the restriction of the ideal sheaf $\left(\widetilde{\mathcal R}^{\underline s}\circ\overline{\mathcal R}^{\underline s}_n\right)^{-1}\mathscr I^{\underline s}_{\underline w}\cdot\mathcal O_{\left(\overline{\mathcal R}^{\underline s}\right)^{-1}\left(\mathring A^{\tau}\right)}$  can be generated by polynomials independent of the variables $\overrightarrow W$. Similarly to Claim \ref{13gi}, we can show that $\left(\overline{\mathcal R}^{\underline s}_n\right)^{-1}\left({\mathring A^{\tau}}\right)$ is isomorphic to an open subscheme of
\begin{equation*}
{\rm Spec}\,\mathbb Z\left[\overrightarrow W\right]\times\left(\overline{\mathcal R}_{\widetilde n}^{{\widetilde{\underline s}}}\right)^{-1}\left({ A^{\widetilde\tau}}\right)\times \left(\overline{\mathcal R}_{\widehat n}^{{\widehat {\underline s}}}\right)^{-1}\left({ A^{\widehat \tau}}\right).  
\end{equation*}
Hence $\left(\overline{\mathcal R}^{\underline s}_n\right)^{-1}\left({\mathring A^{\tau}}\right)$ is smooth over ${\rm Spec}\,\mathbb Z$ by the assumption. By analogy with the proof of Lemma \ref{uv}, we can complete the proof of Lemma \ref{3uv32}.
\,\,\,\,$\endpf$ 

\subsection{Induction for Types (\ref{31up}), (\ref{32up})}
We define a morphism \begin{equation}\label{spnab}
j^{\tau}:{\rm Spec}\,\mathbb Z\left[\overrightarrow U,\overrightarrow V,\overrightarrow H^1,\cdots,\overrightarrow H^l,\overrightarrow \Xi^1,\cdots,\overrightarrow \Xi^m\right]\longrightarrow \mathcal G(2,n)   
\end{equation}
for any Class I index $\tau=\left(j_1, j_2,\left(j^+_1,\cdots,j^+_l\right),\left(j^-_1,\cdots,j^-_{m}\right)\right)\in \mathbb J^{\underline s}$ by 
\begin{equation}\label{4frakn5}
\begin{split}
&\left(\begin{matrix}
0&0&\cdots&0&1&0&\cdots\\
u_1&u_2&\cdots&u_{p_1-1}&0&u_{p_1+1}&\cdots\\
\end{matrix}\hspace{-0.12in}\begin{matrix} &\hfill\tikzmark{a2}\\
\\&\hfill\tikzmark{b2}
\end{matrix}\,\,\,\begin{matrix}
v_1&v_2&\cdots&v_{p_2-1}&0&v_{p_2+1}&\cdots\\
0&0&\cdots&0&1&0&\cdots\\
\end{matrix}\hspace{-0.12in}\begin{matrix} &\hfill\tikzmark{a12}\\
\\&\hfill\tikzmark{b12}
\end{matrix}\right.\\
&\,\,\,\,\,\,\,\,\,\,\,\,\begin{matrix} &\hfill\tikzmark{c12}\\
\\&\hfill\tikzmark{d12}
\end{matrix}\,\,\,\begin{matrix}\eta^1_{11}&\eta^1_{12}&\eta^1_{13}&\cdots\\
\eta^1_{21}&\eta^1_{22}&\eta^1_{23}&\cdots\end{matrix}\hspace{-0.12in}\begin{matrix} &\hfill\tikzmark{g2}\\
\\&\hfill\tikzmark{h2}\end{matrix}\,\,\,\begin{matrix}\eta^2_{11}&\eta^2_{12}&\eta^2_{13}&\cdots\\
\eta^2_{21}&\eta^2_{22}&\eta^2_{23}&\cdots
\end{matrix}\hspace{-0.12in}\begin{matrix} &\hfill\tikzmark{e2}\\
\\&\hfill\tikzmark{f2}\end{matrix}\,\,\,\begin{matrix}\cdots\cdots\\\cdots\cdots\\
\end{matrix}\hspace{-0.12in}\begin{matrix} &\hfill\tikzmark{c2}\\
\\&\hfill\tikzmark{d2}\end{matrix}\,\,\,\begin{matrix}\eta^l_{11}&\eta^l_{12}&\eta^l_{13}&\cdots\\
\eta^l_{21}&\eta^l_{22}&\eta^l_{23}&\cdots
\end{matrix}\hspace{-0.12in}\begin{matrix} &\hfill\tikzmark{e22}\\
\\&\hfill\tikzmark{f22}
\end{matrix}\\
&\left.\,\,\,\,\,\,\,\,\,\,\,\,\begin{matrix} &\hfill\tikzmark{e12}\\
\\&\hfill\tikzmark{f12}
\end{matrix}\,\,\,\begin{matrix}\xi^1_{11}&\xi^1_{12}&\xi^1_{13}&\cdots\\
\xi^1_{21}&\xi^1_{22}&\xi^1_{23}&\cdots\end{matrix}\hspace{-0.12in}\begin{matrix} &\hfill\tikzmark{g12}\\
\\&\hfill\tikzmark{h12}\end{matrix}\,\,\,\begin{matrix}\xi^2_{11}&\xi^2_{12}&\xi^2_{13}&\cdots\\
\xi^2_{21}&\xi^2_{22}&\xi^2_{23}&\cdots
\end{matrix}\hspace{-0.12in}\begin{matrix} &\hfill\tikzmark{a22}\\
\\&\hfill\tikzmark{b22}\end{matrix}\,\,\,\begin{matrix}\cdots\cdots\\\cdots\cdots\\
\end{matrix}\hspace{-0.12in}\begin{matrix} &\hfill\tikzmark{c22}\\
\\&\hfill\tikzmark{d22}\end{matrix}\,\,\,\begin{matrix}\xi^m_{11}&\xi^m_{12}&\xi^m_{13}&\cdots\\
\xi^m_{21}&\xi^m_{22}&\xi^m_{23}&\cdots
\end{matrix}\right),
\tikz[remember picture,overlay]   \draw[dashed,dash pattern={on 4pt off 2pt}] ([xshift=0.5\tabcolsep,yshift=7pt]a2.north) -- ([xshift=0.5\tabcolsep,yshift=-2pt]b2.south);\tikz[remember picture,overlay]   \draw[dashed,dash pattern={on 4pt off 2pt}] ([xshift=0.5\tabcolsep,yshift=7pt]c2.north) -- ([xshift=0.5\tabcolsep,yshift=-2pt]d2.south);\tikz[remember picture,overlay]   \draw[dashed,dash pattern={on 4pt off 2pt}] ([xshift=0.5\tabcolsep,yshift=7pt]e2.north) -- ([xshift=0.5\tabcolsep,yshift=-2pt]f2.south);\tikz[remember picture,overlay]   \draw[dashed,dash pattern={on 4pt off 2pt}] ([xshift=0.5\tabcolsep,yshift=7pt]g2.north) -- ([xshift=0.5\tabcolsep,yshift=-2pt]h2.south);\tikz[remember picture,overlay]   \draw[dashed,dash pattern={on 4pt off 2pt}] ([xshift=0.5\tabcolsep,yshift=7pt]a12.north) -- ([xshift=0.5\tabcolsep,yshift=-2pt]b12.south);\tikz[remember picture,overlay]   \draw[dashed,dash pattern={on 4pt off 2pt}] ([xshift=0.5\tabcolsep,yshift=7pt]c12.north) -- ([xshift=0.5\tabcolsep,yshift=-2pt]d12.south);\tikz[remember picture,overlay]   \draw[dashed,dash pattern={on 4pt off 2pt}] ([xshift=0.5\tabcolsep,yshift=7pt]e12.north) -- ([xshift=0.5\tabcolsep,yshift=-2pt]f12.south);\tikz[remember picture,overlay]   \draw[dashed,dash pattern={on 4pt off 2pt}] ([xshift=0.5\tabcolsep,yshift=7pt]g12.north) -- ([xshift=0.5\tabcolsep,yshift=-2pt]h12.south);\tikz[remember picture,overlay]   \draw[dashed,dash pattern={on 4pt off 2pt}] ([xshift=0.5\tabcolsep,yshift=7pt]a22.north) -- ([xshift=0.5\tabcolsep,yshift=-2pt]b22.south);\tikz[remember picture,overlay]   \draw[dashed,dash pattern={on 4pt off 2pt}] ([xshift=0.5\tabcolsep,yshift=7pt]c22.north) -- ([xshift=0.5\tabcolsep,yshift=-2pt]d22.south);\tikz[remember picture,overlay]   \draw[dashed,dash pattern={on 4pt off 2pt}] ([xshift=0.5\tabcolsep,yshift=7pt]e22.north) -- ([xshift=0.5\tabcolsep,yshift=-2pt]f22.south);\\
\end{split} \end{equation}
and define a morphism
$j^{\tau}:{\rm Spec}\,\mathbb Z\left[\overrightarrow H^1,\cdots,\overrightarrow H^l,\overrightarrow \Xi^1,\cdots,\overrightarrow \Xi^m\right]\longrightarrow \mathcal G(2,n)$
for any Class II index $\tau=\left(j_1, j_2,\left(j^+_1,\cdots,j^+_l\right),\left(j^-_1,\cdots,j^-_{m}\right)\right)\in \mathbb J^{\underline s}$ by
\begin{equation}\label{4frakn52}
\begin{split}
&\left(\begin{matrix}
1&0&0&\cdots&0\\
0&1&0&\cdots&0\\
\end{matrix}\hspace{-0.12in}\begin{matrix} &\hfill\tikzmark{c12}\\
\\&\hfill\tikzmark{d12}
\end{matrix}\,\,\,\begin{matrix}\eta^1_{11}&\eta^1_{12}&\eta^1_{13}&\cdots\\
\eta^1_{21}&\eta^1_{22}&\eta^1_{23}&\cdots\end{matrix}\hspace{-0.12in}\begin{matrix} &\hfill\tikzmark{g2}\\
\\&\hfill\tikzmark{h2}\end{matrix}\,\,\,\begin{matrix}\eta^2_{11}&\eta^2_{12}&\eta^2_{13}&\cdots\\
\eta^2_{21}&\eta^2_{22}&\eta^2_{23}&\cdots
\end{matrix}\hspace{-0.12in}\begin{matrix} &\hfill\tikzmark{e2}\\
\\&\hfill\tikzmark{f2}\end{matrix}\,\,\,\begin{matrix}\cdots\cdots\\\cdots\cdots\\
\end{matrix}\hspace{-0.12in}\begin{matrix} &\hfill\tikzmark{c2}\\
\\&\hfill\tikzmark{d2}\end{matrix}\,\,\,\begin{matrix}\eta^l_{11}&\eta^l_{12}&\eta^l_{13}&\cdots\\
\eta^l_{21}&\eta^l_{22}&\eta^l_{23}&\cdots
\end{matrix}\hspace{-0.12in}\begin{matrix} &\hfill\tikzmark{e22}\\
\\&\hfill\tikzmark{f22}
\end{matrix}\right.\\
&\left.\,\,\,\,\,\,\,\,\,\,\,\,\begin{matrix} &\hfill\tikzmark{e12}\\
\\&\hfill\tikzmark{f12}
\end{matrix}\,\,\,\begin{matrix}\xi^1_{11}&\xi^1_{12}&\xi^1_{13}&\cdots\\
\xi^1_{21}&\xi^1_{22}&\xi^1_{23}&\cdots\end{matrix}\hspace{-0.12in}\begin{matrix} &\hfill\tikzmark{g12}\\
\\&\hfill\tikzmark{h12}\end{matrix}\,\,\,\begin{matrix}\xi^2_{11}&\xi^2_{12}&\xi^2_{13}&\cdots\\
\xi^2_{21}&\xi^2_{22}&\xi^2_{23}&\cdots
\end{matrix}\hspace{-0.12in}\begin{matrix} &\hfill\tikzmark{a22}\\
\\&\hfill\tikzmark{b22}\end{matrix}\,\,\,\begin{matrix}\cdots\cdots\\\cdots\cdots\\
\end{matrix}\hspace{-0.12in}\begin{matrix} &\hfill\tikzmark{c22}\\
\\&\hfill\tikzmark{d22}\end{matrix}\,\,\,\begin{matrix}\xi^m_{11}&\xi^m_{12}&\xi^m_{13}&\cdots\\
\xi^m_{21}&\xi^m_{22}&\xi^m_{23}&\cdots
\end{matrix}\right).
\tikz[remember picture,overlay]   \draw[dashed,dash pattern={on 4pt off 2pt}] ([xshift=0.5\tabcolsep,yshift=7pt]c2.north) -- ([xshift=0.5\tabcolsep,yshift=-2pt]d2.south);\tikz[remember picture,overlay]   \draw[dashed,dash pattern={on 4pt off 2pt}] ([xshift=0.5\tabcolsep,yshift=7pt]e2.north) -- ([xshift=0.5\tabcolsep,yshift=-2pt]f2.south);\tikz[remember picture,overlay]   \draw[dashed,dash pattern={on 4pt off 2pt}] ([xshift=0.5\tabcolsep,yshift=7pt]g2.north) -- ([xshift=0.5\tabcolsep,yshift=-2pt]h2.south);\tikz[remember picture,overlay]   \draw[dashed,dash pattern={on 4pt off 2pt}] ([xshift=0.5\tabcolsep,yshift=7pt]c12.north) -- ([xshift=0.5\tabcolsep,yshift=-2pt]d12.south);\tikz[remember picture,overlay]   \draw[dashed,dash pattern={on 4pt off 2pt}] ([xshift=0.5\tabcolsep,yshift=7pt]e12.north) -- ([xshift=0.5\tabcolsep,yshift=-2pt]f12.south);\tikz[remember picture,overlay]   \draw[dashed,dash pattern={on 4pt off 2pt}] ([xshift=0.5\tabcolsep,yshift=7pt]g12.north) -- ([xshift=0.5\tabcolsep,yshift=-2pt]h12.south);\tikz[remember picture,overlay]   \draw[dashed,dash pattern={on 4pt off 2pt}] ([xshift=0.5\tabcolsep,yshift=7pt]a22.north) -- ([xshift=0.5\tabcolsep,yshift=-2pt]b22.south);\tikz[remember picture,overlay]   \draw[dashed,dash pattern={on 4pt off 2pt}] ([xshift=0.5\tabcolsep,yshift=7pt]c22.north) -- ([xshift=0.5\tabcolsep,yshift=-2pt]d22.south);\tikz[remember picture,overlay]   \draw[dashed,dash pattern={on 4pt off 2pt}] ([xshift=0.5\tabcolsep,yshift=7pt]e22.north) -- ([xshift=0.5\tabcolsep,yshift=-2pt]f22.south);\\
\end{split} \end{equation}
Here we adopt the convention that
$\eta_{1p^+_{\alpha}}^{\alpha}=\xi^{\beta}_{2p^-_{\beta}}=1$ 
for $1\leq\alpha\leq l$, $1\leq\beta\leq m$. Let $\mathbb M^{\underline s}_{j^{\tau}}$ be given by Definition \ref{mat}. Then, computation yields that for $\tau\in\mathbb J^{\underline s}$ of Class I (resp. Class II),
\begin{equation}\label{kabm}
\left(\overline{\mathcal R}^{\underline s}_n\right)^{-1}\left({A^{\tau}}\right)\cong{\rm Spec}\,\mathbb Z\left[\overrightarrow A,\overrightarrow Y,\overrightarrow Z\right]\times\mathbb M^{\underline s}_{j^{\tau}}\,\left({\rm resp.}\, \left(\overline{\mathcal R}^{\underline s}_n\right)^{-1}\left({A^{\tau}}\right)\cong{\rm Spec}\,\mathbb Z\left[\overrightarrow A,\overrightarrow W\right]\times\mathbb M^{\underline s}_{j^{\tau}}\right).
\end{equation}

We can derive that
\begin{lemma}\label{induc}
Suppose that  Proposition \ref{3ms} holds for all integers $N^{\prime}$ such that $2\leq N^{\prime}<N$. Let $\widehat{\underline s}\in\left(\mathbb Z\right)^{\widehat N}$ be any size vector such that $2\leq\widehat N<N$. Then $\mathbb M_{j^{\tau}}^{\widehat{\underline s}}$ constructed above is smooth over ${\rm Spec}\,\mathbb Z$ for any index $\widehat\tau\in\mathbb J^{\widehat{\underline s}}$.

\end{lemma}

In this subsection, we are concerned with
\begin{lemma}\label{3ud}Suppose that  Proposition \ref{3ms} holds for all integers $N^{\prime}$ such that $2\leq N^{\prime}<N$. Then Lemma \ref{3loc} holds for all (truncated) coordinate charts of {\bf Type (\ref{31up})}.
\end{lemma}

{\bf\noindent Proof of Lemma \ref{3ud}.} Fix $\tau=\left(j_1,j_2,(j_1^+,\cdots,j^+_{l})\right)$ $\in\mathbb J^{\underline s}$ where $l=N-2$. Without loss of generality, we can assume that  $t_1=1$, $t_2=2$, and $p_1=p_2=p^+_{1}=\cdots=p^+_{l}=1$. In what follows, we write vectors ${\underline w}\in C^{\underline s}$ in form (\ref{wc}), and, similarly, write $\underline s=(s_1,s_2,s_{t_1^+},\cdots,s_{t_l^+})$ by rearranging its components corresponding to $\tau$. 

We first partition $C^{\underline s}$ as follows. (It is derive from the partition given in the proof of Lemma \ref{3uv} by omitting the subsets where the vectors have non-trivial $w_{t^-_{\beta}}$ components.)
\begin{enumerate}[label={(A\arabic*)},ref=A\arabic*]
\setcounter{enumi}{-1}

\item $w_{t_1}=2$, and the other components are $0$. 

\item $w_{t_1}=1$, $w_{t^+_{\alpha}}=1$ for a certain $1\leq \alpha\leq l$, and the other components are $0$.

\item $w_{t_1}=w_{t_2}=1$, $w_{t^+_{\alpha_1}}=w_{t^+_{\alpha_2}}=1$ for certain $1\leq \alpha_1<\alpha_2\leq l$, and the other components are $0$.
\end{enumerate}
\begin{enumerate}[label={(B\arabic*)},ref=B\arabic*]
\setcounter{enumi}{-1}

\item $w_{t_2}=2$, and the other components are $0$.

\end{enumerate}

\begin{enumerate}[label={(AB\arabic*)},ref={AB\arabic*}]
\setcounter{enumi}{-1}
\item $w_{t_1}=w_{t_2}=1$, and the other components are $0$.
\end{enumerate}

\begin{enumerate}[label={(A\arabic*)},ref={A\arabic*}]
\setcounter{enumi}{2}
\item $w_{t^+_{\alpha}}=2$, for a certain $1\leq \alpha\leq l$, and other components are $0$.

\item $w_{t^+_{\alpha_1}}=w_{t^+_{\alpha_2}}=1$, for certain $1\leq \alpha_1<\alpha_2\leq l$, and other components are $0$.

\item $w_{t_1}=w_{t^+_{\alpha_1}}=w_{t^+_{\alpha_2}}=w_{t^+_{\alpha_3}}=1$, for certain $1\leq \alpha_1<\alpha_2<\alpha_3\leq l$.

\item $w_{t^+_{\alpha_1}}=w_{t^+_{\alpha_2}}=w_{t^+_{\alpha_3}}=w_{t^+_{\alpha_4}}=1$, for certain $1\leq \alpha_1<\alpha_2<\alpha_3<\alpha_4\leq l$.

\end{enumerate}

\begin{enumerate}[label={(C\arabic*)},ref={C\arabic*}]
\setcounter{enumi}{1}

\item $w_{t_2}=w_{t^+_{\alpha}}=1$, for a certain $1\leq \alpha\leq l$, and other components are $0$.

\end{enumerate}

\begin{enumerate}[label={(C\arabic*)},ref={C\arabic*}]
\setcounter{enumi}{4}

\item $w_{t_2}=w_{t^+_{\alpha_1}}=w_{t^+_{\alpha_2}}=w_{t^+_{\alpha_3}}=1$, for certain $1\leq \alpha_1<\alpha_2<\alpha_3\leq l$.

\end{enumerate}
It is clear that the blow-up with respect to the (pull-backs of) ideal sheaves $\mathscr I^{\underline s}_{\underline w}$ with $\underline w$ in (A0), (A1), (A2), (B0), (AB0) yields $\check{\mathcal R}_N:W^{\sigma}_{N_{\tau}}\rightarrow \mathcal Q_{n}^{\underline s}$. Denote by $\mathcal R_{AB}$ the sequence of blow-ups with respect to the (pull-backs of) ideal sheaves $\mathscr I^{\underline s}_{\underline w}$ with $\underline w$ in (A3)-(A6); denote by $\mathcal R_{ABC}$ the sequence of blow-ups with respect to the (pull-backs of) ideal sheaves $\mathscr I^{\underline s}_{\underline w}$ with $\underline w$ in (C2), (C5).

For each index  $\lambda:=\left(q^+,q^-,\lambda^+,\Delta^+\right)\in\Lambda^{\tau}$, let $B^{\lambda}$ be the subscheme of $\left(\check{\mathcal R}_N\right)^{-1}\left(A^{\tau}\right)$ defined in Lemma \ref{3fac}. Similarly to Claim \ref{13gi}, we can prove that

\begin{claim}\label{2iso}
$\mathcal R_{ABC}$ is an isomorphism.
\end{claim}

Define a new size vector  $\underline s^*:=(1,s_{t_1^+},\cdots,s_{t_l^+})$ when $s_1=1$, and $\underline s^*:=(2,s_{t_1^+},\cdots,s_{t_l^+})$
when $s_1\geq 2$. When $s_1=1$, we define a morphism 
\begin{equation*}
L_{\lambda}: {\rm Spec}\,\mathbb Z\left[\left(\epsilon_2^+,\epsilon_3^+,\cdots,\epsilon_{l}^+\right),\overrightarrow {\mathscr H}^1,\cdots,\overrightarrow {\mathscr H}^l,\overrightarrow {\mathfrak H}^1,\cdots,\overrightarrow {\mathfrak H}^l\right]\rightarrow  \mathcal G(2,1+s_{t_1^+}+\cdots+s_{t_l^+}) 
\end{equation*}  by
\begin{equation}\label{13frakn7}
\begin{split}
&\left(\begin{matrix}
1\\
0\\
\end{matrix}\hspace{-0.12in}\begin{matrix} &\hfill\tikzmark{a2}\\
\\&\hfill\tikzmark{b2}
\end{matrix}\,\,\,\begin{matrix}1&\eta^1_{12}&\eta^1_{13}&\cdots\\
\frac{\eta^1_{21}}{\epsilon^+_1}&\frac{\eta^1_{22}}{\epsilon^+_1}&\frac{\eta^1_{23}}{\epsilon^+_1}&\cdots\end{matrix}\hspace{-0.12in}\begin{matrix} &\hfill\tikzmark{g2}\\
\\&\hfill\tikzmark{h2}\end{matrix}\,\,\,\begin{matrix}1&\eta^2_{12}&\eta^2_{13}&\cdots\\
\frac{\eta^2_{21}}{\epsilon^+_1}&\frac{\eta^2_{22}}{\epsilon^+_1}&\frac{\eta^2_{23}}{\epsilon^+_1}&\cdots
\end{matrix}\hspace{-0.12in}\begin{matrix} &\hfill\tikzmark{e2}\\
\\&\hfill\tikzmark{f2}\end{matrix}\,\,\,\begin{matrix}\cdots\cdots\\\cdots\cdots\\
\end{matrix}\hspace{-0.12in}\begin{matrix} &\hfill\tikzmark{c2}\\
\\&\hfill\tikzmark{d2}\end{matrix}\,\,\,\begin{matrix}1&\eta^l_{12}&\eta^l_{13}&\cdots\\
\frac{\eta^l_{21}}{\epsilon^+_1}&\frac{\eta^l_{22}}{\epsilon^+_1}&\frac{\eta^l_{23}}{\epsilon^+_1}&\cdots
\end{matrix}\right);
\tikz[remember picture,overlay]   \draw[dashed,dash pattern={on 4pt off 2pt}] ([xshift=0.5\tabcolsep,yshift=7pt]a2.north) -- ([xshift=0.5\tabcolsep,yshift=-2pt]b2.south);\tikz[remember picture,overlay]   \draw[dashed,dash pattern={on 4pt off 2pt}] ([xshift=0.5\tabcolsep,yshift=7pt]c2.north) -- ([xshift=0.5\tabcolsep,yshift=-2pt]d2.south);\tikz[remember picture,overlay]   \draw[dashed,dash pattern={on 4pt off 2pt}] ([xshift=0.5\tabcolsep,yshift=7pt]e2.north) -- ([xshift=0.5\tabcolsep,yshift=-2pt]f2.south);\tikz[remember picture,overlay]   \draw[dashed,dash pattern={on 4pt off 2pt}] ([xshift=0.5\tabcolsep,yshift=7pt]g2.north) -- ([xshift=0.5\tabcolsep,yshift=-2pt]h2.south);\\    
\end{split}  
\end{equation}
when $s_1\geq2$, we define a morphism
\begin{equation*}
L_{\lambda}:{\rm Spec}\,\mathbb Z\left[\left(\epsilon_2^+,\epsilon_3^+,\cdots, \epsilon_{q^++1}^+\right),\overrightarrow {\mathscr H}^1,\cdots,\overrightarrow {\mathscr H}^l,\overrightarrow {\mathfrak H}^1,\cdots,\overrightarrow {\mathfrak H}^l\right]\rightarrow \mathcal {G}(2,2+s_{t_1^+}+\cdots+s_{t_l^+})   
\end{equation*} by
\begin{equation}\label{13frakn5}
\begin{split}
&\left(\begin{matrix}
1&0\\
0&\prod\limits_{\gamma=2}^{q^++1}\epsilon^+_{\gamma}\\
\end{matrix}\hspace{-0.12in}\begin{matrix} &\hfill\tikzmark{a2}\\\\
\\&\hfill\tikzmark{b2}
\end{matrix}\,\,\,\begin{matrix}1&\eta^1_{12}&\eta^1_{13}&\cdots\\\\
\frac{\eta^1_{21}}{\epsilon^+_1}&\frac{\eta^1_{22}}{\epsilon^+_1}&\frac{\eta^1_{23}}{\epsilon^+_1}&\cdots\end{matrix}\hspace{-0.12in}\begin{matrix} &\hfill\tikzmark{g2}\\\\
\\&\hfill\tikzmark{h2}\end{matrix}\,\,\,\begin{matrix}1&\eta^2_{12}&\eta^2_{13}&\cdots\\\\
\frac{\eta^2_{21}}{\epsilon^+_1}&\frac{\eta^2_{22}}{\epsilon^+_1}&\frac{\eta^2_{23}}{\epsilon^+_1}&\cdots
\end{matrix}\hspace{-0.12in}\begin{matrix} &\hfill\tikzmark{e2}\\\\
\\&\hfill\tikzmark{f2}\end{matrix}\,\,\,\begin{matrix}\cdots\cdots\\\\\cdots\cdots\\
\end{matrix}\hspace{-0.12in}\begin{matrix} &\hfill\tikzmark{c2}\\\\
\\&\hfill\tikzmark{d2}\end{matrix}\,\,\,\begin{matrix}1&\eta^l_{12}&\eta^l_{13}&\cdots\\\\
\frac{\eta^l_{21}}{\epsilon^+_1}&\frac{\eta^l_{22}}{\epsilon^+_1}&\frac{\eta^l_{23}}{\epsilon^+_1}&\cdots
\end{matrix}\right).
\tikz[remember picture,overlay]   \draw[dashed,dash pattern={on 4pt off 2pt}] ([xshift=0.5\tabcolsep,yshift=7pt]a2.north) -- ([xshift=0.5\tabcolsep,yshift=-2pt]b2.south);\tikz[remember picture,overlay]   \draw[dashed,dash pattern={on 4pt off 2pt}] ([xshift=0.5\tabcolsep,yshift=7pt]c2.north) -- ([xshift=0.5\tabcolsep,yshift=-2pt]d2.south);\tikz[remember picture,overlay]   \draw[dashed,dash pattern={on 4pt off 2pt}] ([xshift=0.5\tabcolsep,yshift=7pt]e2.north) -- ([xshift=0.5\tabcolsep,yshift=-2pt]f2.south);\tikz[remember picture,overlay]   \draw[dashed,dash pattern={on 4pt off 2pt}] ([xshift=0.5\tabcolsep,yshift=7pt]g2.north) -- ([xshift=0.5\tabcolsep,yshift=-2pt]h2.south);\\    
\end{split}  
\end{equation}
Here, by a slight abuse of notation, we denote by $\eta^{\alpha}_{21},\eta^{\alpha}_{22},\cdots,\eta^{\alpha}_{2s_{t_{\alpha}^+}}$, $1\leq\alpha\leq l$, their images under $\Sigma^{\lambda}$ (see (\ref{t1h1}), (\ref{t1h2})), which thus have a common factor $\epsilon^+_1$.

We can show that the ideal sheaf $\left(\widetilde{\mathcal R}^{\underline s}\circ\check{\mathcal R}_N\right)^{-1}\mathscr I^{\underline s}_{\underline w}\cdot\mathcal O_{B^{\lambda}}$ is a product of an invertible ideal sheaf and an ideal sheaf defined by polynomials in $\epsilon_2^+,\epsilon_3^+,\cdots,\epsilon_{q^++1}^+$, $\overrightarrow {\mathscr H}^1,\cdots,\overrightarrow {\mathscr H}^l,\overrightarrow {\mathfrak H}^1,\cdots,\overrightarrow {\mathfrak H}^l$. Then, we can derive that
\begin{equation}\label{subsch}
\begin{split}
&{\mathcal R}_{AB}^{-1}\left(B^{\lambda}\right)\cong{\rm Spec}\,\mathbb Z\left[\epsilon^+_1,\epsilon^-_1,\overrightarrow A,\overrightarrow Y,\overrightarrow Z,\overrightarrow{\mathfrak U},\overrightarrow{\mathfrak V}\right]\times\mathbb M^{\underline s^*}_{L_{\lambda}},\\
\end{split}  
\end{equation}
where $\mathbb M^{\underline s^*}_{L_{\lambda}}$ is defined by Definition \ref{mat}.

By analogy with Lemma \ref{ud}, we can reduce the proof of Lemma \ref{3ud} to 
\begin{claim}\label{strip} $\mathbb M^{\underline s^*}_{L_{\lambda}}$ is smooth  over ${\rm Spec}\,\mathbb Z$ for any $\lambda\in\Lambda^{\tau}$.
\end{claim}

We first prove  that

\begin{claim}\label{special}
When $s_1=1$, $\mathbb M^{\underline s^*}_{L_{\lambda}}$ is smooth  over ${\rm Spec}\,\mathbb Z$   for any $\lambda\in\Lambda^{\tau}$.
\end{claim}
{\bf\noindent Proof of Claim \ref{special}.} Without loss of generality, we assume that $\lambda^+=(1,2,\cdots,l)$ where $l=N-2$.  We adopt the convention that $\eta_{11}^{\alpha}=\mathfrak h^{\alpha}_{2\delta^+_{\alpha}}=1$ 
for $1\leq\alpha\leq l$.  

Define
an open cover $\left\{C^k_N\right\}_{k=2}^{l+1}$ of $\mathbb A_+:={\rm Spec}\,\mathbb Z\left[\epsilon_2^+,\cdots,\epsilon_{l}^+,\overrightarrow {\mathscr H}^1,\cdots,\overrightarrow {\mathscr H}^l,\overrightarrow {\mathfrak H}^1,\cdots,\overrightarrow {\mathfrak H}^l\right]$ by
\begin{equation*}
\begin{split}
C^k_N:=&\left\{\mathfrak p\in \mathbb A_+\left|\begin{matrix}
&\epsilon^+_{\gamma}\notin\mathfrak p,\,\forall\,2\leq\gamma\leq k-1,\,{\rm and}\\
&1-\eta_{1\delta^+_{\alpha}}^{\alpha}\mathfrak h^{\beta}_{21}\prod\limits_{\gamma=\alpha+1}^{\beta}\epsilon^+_{\gamma}\notin\mathfrak p,\,\forall\,1\leq\alpha\leq k-1\,\,{\rm and}\,\,k\leq \beta\leq  l\\
\end{matrix}\right.\right\}.\\
\end{split}
\end{equation*}
Let $\mathcal R^{\underline s^*}_{L_{\lambda}}:\mathbb M^{\underline s^*}_{L_{\lambda}}\rightarrow \mathbb A_+$ be the blow-up given by Definition \ref{mat}. Then, it suffices to show that
\begin{enumerate}[label={({\bf SQ1})},ref={\bf SQ1}]
\item for each $2\leq k\leq l+1=N-1$,  $\left(\mathcal R^{\underline s^*}_{L_{\lambda}}\right)^{-1}\left(C^k_N\right)$ is smooth over ${\rm Spec}\,\mathbb Z$.  \label{3dudu}
\end{enumerate}

By analogy with  {\bf (\ref{p11})}, 
we proceed to prove  on a case by case basis.

\medskip

{\bf\noindent Case I ($k=l+1=N-1$).} We first study the case of $\delta^+_1\geq 2$. We cover $C^{l+1}_N$ by two open subschemes $D_1:=\left\{\mathfrak p\in C^{l+1}_N\left|\,\mathfrak h^1_{21}\notin\mathfrak p\right.\right\}$ and $D_2:=\left\{\mathfrak p\in C^{l+1}_N\left|\,1-\mathfrak h^1_{21}\eta^{1}_{1\delta^+_{1}}\notin\mathfrak p\right.\right\}$.

Left multiplying (\ref{13frakn7}) by $\left(\begin{matrix}
1&-\eta^1_{1\delta^+_1}\\
0&1\\
\end{matrix}\right)$ and substituting in  (\ref{t1h1}), (\ref{t1h2}), we get
\begin{equation}\small
\begin{split}
&\left(\begin{matrix}
1\\\\0\\\end{matrix}\hspace{-0.12in}\begin{matrix} &\hfill\tikzmark{a12}\\
\\\\&\hfill\tikzmark{b12}
\end{matrix}\,\,\,\begin{matrix}\eta^1_{11}-\mathfrak h^1_{21}\cdot\eta^1_{1\delta^+_1}&\cdots&\eta^1_{1\left(\delta^+_1-1\right)}-\mathfrak h^1_{2\left(\delta^+_1-1\right)}\cdot\eta^1_{1\delta^+_1}&0&\eta^1_{1\left(\delta^+_1+1\right)}-\mathfrak h^1_{2\left(\delta^+_1+1\right)}\cdot\eta^1_{1\delta^+_1}&\cdots\\\\
\mathfrak h^1_{21}&\cdots&\mathfrak h^1_{2\left(\delta^+_1-1\right)}&1&\mathfrak h^1_{2\left(\delta^+_1+1\right)}&\cdots\end{matrix}\hspace{-0.12in}\begin{matrix} &\hfill\tikzmark{g12}\\\\
\\&\hfill\tikzmark{h12}\end{matrix}\right.\\
&\,\,\,\,\,\,\,\,\,\,\,\,\,\,\,\,\,\,\,\,\begin{matrix} &\hfill\tikzmark{x12}\\\\
\\&\hfill\tikzmark{y12}\end{matrix}\,\,\,\begin{matrix}\eta^2_{11}-\left(\mathfrak h^2_{21}\cdot\epsilon^+_2\right)\eta^1_{1\delta^+_1}&\cdots&\eta^2_{1\delta^+_2}-\left(\epsilon^+_2\right)\eta^1_{1\delta^+_1}&\cdots\\\\\mathfrak h^2_{21}\cdot\epsilon^+_2&\cdots&\epsilon^+_2&\cdots\\
\end{matrix}\hspace{-0.12in}\begin{matrix} &\hfill\tikzmark{e12}\\\\
\\&\hfill\tikzmark{f12}\end{matrix}\,\,\,\begin{matrix}\cdots\cdots\\\\\cdots\cdots\\
\end{matrix}\\
&\left.\,\,\,\,\,\,\,\,\,\,\,\,\,\,\,\,\,\,\,\begin{matrix} &\hfill\tikzmark{c12}\\\\
\\&\hfill\tikzmark{d12}\end{matrix}\,\,\,\begin{matrix}\eta^l_{11}-\left(\mathfrak h^l_{21}\cdot\prod\nolimits_{\gamma=2}^{l}{\epsilon^+_{\gamma}}\right)\eta^1_{1\delta^+_1}&\cdots&\eta^l_{1\delta^+_l}-\left(\prod\nolimits_{\gamma=2}^{l}{\epsilon^+_{\gamma}}\right)\eta^1_{1\delta^+_1}&\cdots\\
\mathfrak h^l_{21}\cdot\prod\nolimits_{\gamma=2}^{l}{\epsilon^+_{\gamma}}&\cdots&\prod\nolimits_{\gamma=2}^{l}{\epsilon^+_{\gamma}}&\cdots
\end{matrix}\right)\\
&=\left(\begin{matrix}
1\\\\0\\\end{matrix}\hspace{-0.12in}\begin{matrix} &\hfill\tikzmark{a22}\\
\\\\&\hfill\tikzmark{b22}
\end{matrix}\,\,\,\begin{matrix}1-\mathfrak h^1_{21}\cdot\eta^1_{1\delta^+_1}&\cdots&\eta^1_{1\left(\delta^+_1-1\right)}-\mathfrak h^1_{2\left(\delta^+_1-1\right)}\cdot\eta^1_{1\delta^+_1}&0&\eta^1_{1\left(\delta^+_1+1\right)}-\mathfrak h^1_{2\left(\delta^+_1+1\right)}\cdot\eta^1_{1\delta^+_1}&\cdots\\\\
\mathfrak h^1_{21}&\cdots&\mathfrak h^1_{2\left(\delta^+_1-1\right)}&1&\mathfrak h^1_{2\left(\delta^+_1+1\right)}&\cdots\end{matrix}\hspace{-0.12in}\begin{matrix} &\hfill\tikzmark{g22}\\\\
\\&\hfill\tikzmark{h22}\end{matrix}\right.\\
&\left.\begin{matrix} &\hfill\tikzmark{x22}\\\\
\\&\hfill\tikzmark{y22}\end{matrix}\,\,\,\begin{matrix}\frac{1}{\epsilon^+_2}-\mathfrak h^2_{21}\eta^1_{1\delta^+_1}&\cdots&\frac{\eta^2_{1\delta^+_2}}{\epsilon^+_2}-\eta^1_{1\delta^+_1}&\cdots\\\\\mathfrak h^2_{21}&\cdots&1&\cdots\\
\end{matrix}\hspace{-0.12in}\begin{matrix} &\hfill\tikzmark{e22}\\\\
\\&\hfill\tikzmark{f22}\end{matrix}\,\,\,\begin{matrix}\cdots\\\\\cdots\\
\end{matrix}\hspace{-0.12in}\begin{matrix} &\hfill\tikzmark{c22}\\\\
\\&\hfill\tikzmark{d22}\end{matrix}\,\,\,\begin{matrix}\frac{1}{\prod\nolimits_{\gamma=2}^{l}{\epsilon^+_{\gamma}}}-\mathfrak h^l_{21}\eta^1_{1\delta^+_1}&\cdots&\frac{\eta^l_{1\delta^+_1}}{\prod\nolimits_{\gamma=2}^{l}{\epsilon^+_{\gamma}}}-\eta^1_{1\delta^+_1}&\cdots\\
\mathfrak h^l_{21}&\cdots&1&\cdots
\end{matrix}\right)\\
&\,\,\,\,\,\,\,\,\,\,\cdot{\rm diag}\left(1,\,\,\,I_{s_{t_1^+}\times s_{t_1^+}},\,\,\,\epsilon^+_2\cdot I_{s_{t_2^+}\times s_{t_2^+}},\,\,\,\cdots,\,\,\,\left(\prod\nolimits_{\gamma=2}^{l}{\epsilon^+_{\gamma}}\right)\cdot I_{s_{t_l^+}\times s_{t_l^+}}\right).\\
\tikz[remember picture,overlay]   \draw[dashed,dash pattern={on 4pt off 2pt}] ([xshift=0.5\tabcolsep,yshift=7pt]g12.north) -- ([xshift=0.5\tabcolsep,yshift=-2pt]h12.south);\tikz[remember picture,overlay]   \draw[dashed,dash pattern={on 4pt off 2pt}] ([xshift=0.5\tabcolsep,yshift=7pt]a22.north) -- ([xshift=0.5\tabcolsep,yshift=-2pt]b22.south);\tikz[remember picture,overlay]   \draw[dashed,dash pattern={on 4pt off 2pt}] ([xshift=0.5\tabcolsep,yshift=7pt]c22.north) -- ([xshift=0.5\tabcolsep,yshift=-2pt]d22.south);\tikz[remember picture,overlay]   \draw[dashed,dash pattern={on 4pt off 2pt}] ([xshift=0.5\tabcolsep,yshift=7pt]e22.north) -- ([xshift=0.5\tabcolsep,yshift=-2pt]f22.south);\tikz[remember picture,overlay]   \draw[dashed,dash pattern={on 4pt off 2pt}] ([xshift=0.5\tabcolsep,yshift=7pt]g22.north) -- ([xshift=0.5\tabcolsep,yshift=-2pt]h22.south);\tikz[remember picture,overlay]   \draw[dashed,dash pattern={on 4pt off 2pt}] ([xshift=0.5\tabcolsep,yshift=7pt]x22.north) -- ([xshift=0.5\tabcolsep,yshift=-2pt]y22.south);\tikz[remember picture,overlay]   \draw[dashed,dash pattern={on 4pt off 2pt}] ([xshift=0.5\tabcolsep,yshift=7pt]c12.north) -- ([xshift=0.5\tabcolsep,yshift=-2pt]d12.south);\tikz[remember picture,overlay]   \draw[dashed,dash pattern={on 4pt off 2pt}] ([xshift=0.5\tabcolsep,yshift=7pt]e12.north) -- ([xshift=0.5\tabcolsep,yshift=-2pt]f12.south);\tikz[remember picture,overlay]   \draw[dashed,dash pattern={on 4pt off 2pt}] ([xshift=0.5\tabcolsep,yshift=7pt]a12.north) -- ([xshift=0.5\tabcolsep,yshift=-2pt]b12.south);\tikz[remember picture,overlay]   \draw[dashed,dash pattern={on 4pt off 2pt}] ([xshift=0.5\tabcolsep,yshift=7pt]x12.north) -- ([xshift=0.5\tabcolsep,yshift=-2pt]y12.south);
\end{split}   
\end{equation}

Define a Class I index $\widehat\tau=\left(\widehat j_1, \widehat j_2,\left(\widehat j^-_1,\cdots,\widehat j^-_{m}\right)\right)\in \mathbb J^{\underline s^*}$ where $m=l-1$ by
\begin{equation*}
\widehat j_1:=1,\,\,\,\,\,\,\widehat j_2:=1+\delta_1^+,\,\,\,\,\,\,\widehat j_{\beta}^-:=1+s_{t^+_1}+s_{t^+_2}+\cdots+s_{t^+_{\beta}}+\delta^+_{\beta+1},\,\,\,{\rm for}\,\,1\leq\beta\leq m. 
\end{equation*}
We perform the following change of coordinates in $D_1$ such that $D_1$ is isomorphic to an open subscheme of ${\rm Spec}\,\mathbb Z\left[\overrightarrow Z,\overrightarrow V,\overrightarrow \Xi^1,\cdots,\overrightarrow \Xi^m\right]$. 
For all $1\leq\gamma\leq s_{t_1^+}$ such that $\gamma\neq\widehat p_2=\delta^+_1$, set
\begin{equation*}
\begin{split}
&z_{\gamma}:=\mathfrak h^1_{2\gamma},\,\,\,\,\,\,\,\,v_{\gamma}:=\eta^1_{1\gamma}-\mathfrak h^1_{2\gamma}\cdot\eta^1_{1\delta^+_1},\\
\end{split}
\end{equation*}
and for $1\leq\beta\leq m$, set
\begin{equation*}
\begin{split}
&\xi^{\beta}_{1\gamma}:=\eta^{\beta+1}_{1\gamma}\left(\prod\nolimits_{\gamma=2}^{\beta+1}{\epsilon^+_{\gamma}}\right)^{-1}-\mathfrak h^{\beta+1}_{2\gamma}\cdot\eta^1_{1\delta^+_1},\,\,\,\,\,\,\,\,\,\,\,\,\,\forall\,\,1\leq\gamma\leq s_{t_{\beta+1}^-},\\
&\xi^{\beta}_{2\rho}:=\mathfrak h^{\beta+1}_{2\rho},\,\,\,\,\,\,\,\,\,\,\,\,\,\,\,\,\,\,\,\,\,\,\,\,\,\,\,\,\,\,\,\,\,\,\,\forall\,\,1\leq\rho\leq s_{t_{\beta+1}^-}\,\,{\rm such\,\,that\,\,}\rho\neq\widehat p_{\beta}^-=\delta^+_{\beta+1}.\\    
\end{split} 
\end{equation*}
Define a morphism
$j^{\widehat\tau}:{\rm Spec}\,\mathbb Z\left[\overrightarrow V,\overrightarrow \Xi^1,\cdots,\overrightarrow \Xi^m\right]\longrightarrow \mathcal G(2,1+s_{t_1^+}+\cdots+s_{t_l^+})$
as (\ref{spnab}), (\ref{4frakn5}) by
\begin{equation*}
\begin{split}
&\left(\begin{matrix}
1\\
0\\
\end{matrix}\hspace{-0.12in}\begin{matrix} &\hfill\tikzmark{a2}\\
\\&\hfill\tikzmark{b2}
\end{matrix}\,\,\,\begin{matrix}v_1&v_2&\cdots&v_{p_2-1}&0&v_{p_2+1}&\cdots\\
0&0&\cdots&0&1&0&\cdots\\\end{matrix}\hspace{-0.12in}\begin{matrix} &\hfill\tikzmark{g2}\\
\\&\hfill\tikzmark{h2}\end{matrix}\,\,\,\begin{matrix}\xi^1_{11}&\xi^1_{12}&\cdots\\
\xi^1_{21}&\xi^1_{22}&\cdots
\end{matrix}\hspace{-0.12in}\begin{matrix} &\hfill\tikzmark{e2}\\
\\&\hfill\tikzmark{f2}\end{matrix}\,\,\,\begin{matrix}\cdots\cdots\\\cdots\cdots\\
\end{matrix}\hspace{-0.12in}\begin{matrix} &\hfill\tikzmark{c2}\\
\\&\hfill\tikzmark{d2}\end{matrix}\,\,\,\begin{matrix}\xi^m_{11}&\xi^m_{12}&\cdots\\
\xi^m_{21}&\xi^m_{22}&\cdots
\end{matrix}\right),
\tikz[remember picture,overlay]   \draw[dashed,dash pattern={on 4pt off 2pt}] ([xshift=0.5\tabcolsep,yshift=7pt]a2.north) -- ([xshift=0.5\tabcolsep,yshift=-2pt]b2.south);\tikz[remember picture,overlay]   \draw[dashed,dash pattern={on 4pt off 2pt}] ([xshift=0.5\tabcolsep,yshift=7pt]c2.north) -- ([xshift=0.5\tabcolsep,yshift=-2pt]d2.south);\tikz[remember picture,overlay]   \draw[dashed,dash pattern={on 4pt off 2pt}] ([xshift=0.5\tabcolsep,yshift=7pt]e2.north) -- ([xshift=0.5\tabcolsep,yshift=-2pt]f2.south);\tikz[remember picture,overlay]   \draw[dashed,dash pattern={on 4pt off 2pt}] ([xshift=0.5\tabcolsep,yshift=7pt]g2.north) -- ([xshift=0.5\tabcolsep,yshift=-2pt]h2.south);\\ \end{split} \end{equation*}
where
$\xi^{\beta}_{2p^-_{\beta}}=1$ 
for $1\leq\beta\leq m$ by convention.
Similarly to {case I} in Lemma \ref{ud},  we can show that $\left({\mathcal R}^{\underline s^*}_{L_{\lambda}}\right)^{-1}\left(D_1\right)$ is isomorphic to an open subscheme of ${\rm Spec}\,\mathbb Z\left[\overrightarrow Z\right]\times\mathbb M_{j^{\widehat\tau}}^{\underline s^*}$. Therefore, $\left({\mathcal R}^{\underline s^*}_{L_{\lambda}}\right)^{-1}\left(D_1\right)$ is smooth over ${\rm Spec}\,\mathbb Z$
by Lemma \ref{induc}.

Regarding $D_2$, we define a size vector $\widetilde{\underline s}:=\left(1,s_{t^+_2},s_{t^+_3},\cdots,s_{t^+_l}\right)$. 
We define a morphism 
\begin{equation*}
\widetilde{L}_{\lambda}: {\rm Spec}\,\mathbb Z\left[\left(\epsilon_3^+,\epsilon_4^+,\cdots,\epsilon_{l}^+\right),\overrightarrow {\mathscr H}^2,\cdots,\overrightarrow {\mathscr H}^l,\overrightarrow {\mathfrak H}^2,\cdots,\overrightarrow {\mathfrak H}^l\right]\rightarrow  \mathcal G(2,1+s_{t_2^+}+\cdots+s_{t_l^+}) 
\end{equation*}  by 
\begin{equation*}
\begin{split}
&\left(\begin{matrix}
1\\
0\\
\end{matrix}\hspace{-0.12in}\begin{matrix} &\hfill\tikzmark{a2}\\
\\&\hfill\tikzmark{b2}
\end{matrix}\,\,\,\begin{matrix}1&\eta^2_{12}&\eta^2_{13}&\cdots\\
\frac{\eta^2_{21}}{\epsilon^+_1\epsilon^+_2}&\frac{\eta^2_{22}}{\epsilon^+_1\epsilon^+_2}&\frac{\eta^2_{23}}{\epsilon^+_1\epsilon^+_2}&\cdots\end{matrix}\hspace{-0.12in}\begin{matrix} &\hfill\tikzmark{g2}\\
\\&\hfill\tikzmark{h2}\end{matrix}\,\,\,\begin{matrix}1&\eta^3_{12}&\eta^3_{13}&\cdots\\
\frac{\eta^3_{21}}{\epsilon^+_1\epsilon^+_2}&\frac{\eta^3_{22}}{\epsilon^+_1\epsilon^+_2}&\frac{\eta^3_{23}}{\epsilon^+_1\epsilon^+_2}&\cdots
\end{matrix}\hspace{-0.12in}\begin{matrix} &\hfill\tikzmark{e2}\\
\\&\hfill\tikzmark{f2}\end{matrix}\,\,\,\begin{matrix}\cdots\cdots\\\cdots\cdots\\
\end{matrix}\hspace{-0.12in}\begin{matrix} &\hfill\tikzmark{c2}\\
\\&\hfill\tikzmark{d2}\end{matrix}\,\,\,\begin{matrix}1&\eta^l_{12}&\eta^l_{13}&\cdots\\
\frac{\eta^l_{21}}{\epsilon^+_1\epsilon^+_2}&\frac{\eta^l_{22}}{\epsilon^+_1\epsilon^+_2}&\frac{\eta^l_{23}}{\epsilon^+_1\epsilon^+_2}&\cdots
\end{matrix}\right),
\tikz[remember picture,overlay]   \draw[dashed,dash pattern={on 4pt off 2pt}] ([xshift=0.5\tabcolsep,yshift=7pt]a2.north) -- ([xshift=0.5\tabcolsep,yshift=-2pt]b2.south);\tikz[remember picture,overlay]   \draw[dashed,dash pattern={on 4pt off 2pt}] ([xshift=0.5\tabcolsep,yshift=7pt]c2.north) -- ([xshift=0.5\tabcolsep,yshift=-2pt]d2.south);\tikz[remember picture,overlay]   \draw[dashed,dash pattern={on 4pt off 2pt}] ([xshift=0.5\tabcolsep,yshift=7pt]e2.north) -- ([xshift=0.5\tabcolsep,yshift=-2pt]f2.south);\tikz[remember picture,overlay]   \draw[dashed,dash pattern={on 4pt off 2pt}] ([xshift=0.5\tabcolsep,yshift=7pt]g2.north) -- ([xshift=0.5\tabcolsep,yshift=-2pt]h2.south);\\ \end{split} \end{equation*}
which is derived from (\ref{13frakn7}) by deleting the second block and dividing the second row by $\epsilon_2^+$. Similarly to Claim \ref{13gi}, we can show that $\left({\mathcal R}^{\underline s^*}_{L_{\lambda}}\right)^{-1}\left(D_2\right)$ is isomorphic to an open subscheme of ${\rm Spec}\,\mathbb Z\left[\epsilon_2^+, \overrightarrow {\mathscr H}^1,\overrightarrow {\mathfrak H}^1\right]\times \mathbb M^{\widetilde{\underline  s}}_{ \widetilde{L}_{\lambda}}$. Therefore,  $\left({\mathcal R}^{\underline s^*}_{L_{\lambda}}\right)^{-1}\left(D_2\right)$ is smooth over ${\rm Spec}\,\mathbb Z$ by (\ref{subsch}) and the assumption that Proposition \ref{3ms} holds for all $2\leq N^{\prime}<N$.

The proof for the case of $\delta_1^+=1$ is similar, and  we omit the details  here for brevity.
\smallskip

{\bf\noindent Case II ($2\leq k\leq l$).}  Define $\widetilde{\underline s}:=\left(1,s_{t_1^+},s_{t_2^+},\cdots,s_{t_{k-1}^+}\right)$, and $\widehat{\underline s}:=\left(1,s_{t_{k}^+},s_{t_{k+1}^+},\cdots,s_{t_l^+}\right)$. Define a morphism 
\begin{equation*}
\widetilde{L}_{\lambda}: {\rm Spec}\,\mathbb Z\left[\left(\epsilon_2^+,\epsilon_3^+,\cdots,\epsilon_{k-1}^+\right),\overrightarrow {\mathscr H}^1,\cdots,\overrightarrow {\mathscr H}^{k-1},\overrightarrow {\mathfrak H}^1,\cdots,\overrightarrow {\mathfrak H}^{k-1}\right]\rightarrow  \mathcal{G}(2,1+s_{t_1^+}+\cdots+s_{t_{k-1}^+}) 
\end{equation*}
by
\begin{equation}\label{13frakn9}
\begin{split}
&\left(\begin{matrix}
1\\
0\\
\end{matrix}\hspace{-0.12in}\begin{matrix} &\hfill\tikzmark{a2}\\
\\&\hfill\tikzmark{b2}
\end{matrix}\,\,\,\begin{matrix}1&\eta^1_{12}&\eta^1_{13}&\cdots\\
\frac{\eta^1_{21}}{\epsilon^+_1}&\frac{\eta^1_{22}}{\epsilon^+_1}&\frac{\eta^1_{23}}{\epsilon^+_1}&\cdots\end{matrix}\hspace{-0.12in}\begin{matrix} &\hfill\tikzmark{g2}\\
\\&\hfill\tikzmark{h2}\end{matrix}\,\,\,\begin{matrix}1&\eta^2_{12}&\eta^2_{13}&\cdots\\
\frac{\eta^2_{21}}{\epsilon^+_1}&\frac{\eta^2_{22}}{\epsilon^+_1}&\frac{\eta^2_{23}}{\epsilon^+_1}&\cdots
\end{matrix}\hspace{-0.12in}\begin{matrix} &\hfill\tikzmark{e2}\\
\\&\hfill\tikzmark{f2}\end{matrix}\,\,\,\begin{matrix}\cdots\cdots\\\cdots\cdots\\
\end{matrix}\hspace{-0.12in}\begin{matrix} &\hfill\tikzmark{c2}\\
\\&\hfill\tikzmark{d2}\end{matrix}\,\,\,\begin{matrix}1&\eta^{k-1}_{12}&\eta^{k-1}_{13}&\cdots\\
\frac{\eta^{k-1}_{21}}{\epsilon^+_1}&\frac{\eta^{k-1}_{22}}{\epsilon^+_1}&\frac{\eta^{k-1}_{23}}{\epsilon^+_1}&\cdots
\end{matrix}\right),
\tikz[remember picture,overlay]   \draw[dashed,dash pattern={on 4pt off 2pt}] ([xshift=0.5\tabcolsep,yshift=7pt]a2.north) -- ([xshift=0.5\tabcolsep,yshift=-2pt]b2.south);\tikz[remember picture,overlay]   \draw[dashed,dash pattern={on 4pt off 2pt}] ([xshift=0.5\tabcolsep,yshift=7pt]c2.north) -- ([xshift=0.5\tabcolsep,yshift=-2pt]d2.south);\tikz[remember picture,overlay]   \draw[dashed,dash pattern={on 4pt off 2pt}] ([xshift=0.5\tabcolsep,yshift=7pt]e2.north) -- ([xshift=0.5\tabcolsep,yshift=-2pt]f2.south);\tikz[remember picture,overlay]   \draw[dashed,dash pattern={on 4pt off 2pt}] ([xshift=0.5\tabcolsep,yshift=7pt]g2.north) -- ([xshift=0.5\tabcolsep,yshift=-2pt]h2.south);\\    
\end{split}  
\end{equation}
and a morphism 
\begin{equation*}
\widehat{L}_{\lambda}: {\rm Spec}\,\mathbb Z\left[\left(\epsilon_{k+1}^+,\epsilon_{k+2}^+,\cdots,\epsilon_{l}^+\right),\overrightarrow {\mathscr H}^k,\cdots,\overrightarrow {\mathscr H}^{l},\overrightarrow {\mathfrak H}^k,\cdots,\overrightarrow {\mathfrak H}^{l}\right]\rightarrow  \mathcal {G}(2,1+s_{t_k^+}+\cdots+s_{t_{l}^+}) 
\end{equation*}
 by
\begin{equation}\label{13frakn6}
\begin{split}
&\left(\begin{matrix}1\\\\0\\\end{matrix}\hspace{-0.12in}\begin{matrix} &\hfill\tikzmark{a2}\\\\
\\&\hfill\tikzmark{b2}
\end{matrix}\,\,\,\begin{matrix}1&\eta^k_{12}&\cdots\\
\frac{\eta^k_{21}}{\prod\limits_{\gamma=2}^{k}{\epsilon^+_{\gamma}}}&\frac{\eta^k_{22}}{\prod\limits_{\gamma=2}^{k}{\epsilon^+_{\gamma}}}&\cdots\end{matrix}\hspace{-0.12in}\begin{matrix} &\hfill\tikzmark{g2}\\\\
\\&\hfill\tikzmark{h2}\end{matrix}\,\,\,\begin{matrix}1&\eta^{k+1}_{12}&\cdots\\
\frac{\eta^{k+1}_{21}}{\prod\limits_{\gamma=2}^{k}{\epsilon^+_{\gamma}}}&\frac{\eta^{k+1}_{22}}{\prod\limits_{\gamma=2}^{k}{\epsilon^+_{\gamma}}}&\cdots
\end{matrix}\hspace{-0.12in}\begin{matrix} &\hfill\tikzmark{e2}\\\\
\\&\hfill\tikzmark{f2}\end{matrix}\,\,\,\begin{matrix}\cdots\cdots\\\\\cdots\cdots\\
\end{matrix}\hspace{-0.12in}\begin{matrix} &\hfill\tikzmark{c2}\\\\
\\&\hfill\tikzmark{d2}\end{matrix}\,\,\,\begin{matrix}1&\eta^{l}_{12}&\cdots\\
\frac{\eta^{l}_{21}}{\prod\limits_{\gamma=2}^{k}{\epsilon^+_{\gamma}}}&\frac{\eta^{l}_{22}}{\prod\limits_{\gamma=2}^{k}{\epsilon^+_{\gamma}}}&\cdots
\end{matrix}\right),
\tikz[remember picture,overlay]   \draw[dashed,dash pattern={on 4pt off 2pt}] ([xshift=0.5\tabcolsep,yshift=7pt]a2.north) -- ([xshift=0.5\tabcolsep,yshift=-2pt]b2.south);\tikz[remember picture,overlay]   \draw[dashed,dash pattern={on 4pt off 2pt}] ([xshift=0.5\tabcolsep,yshift=7pt]c2.north) -- ([xshift=0.5\tabcolsep,yshift=-2pt]d2.south);\tikz[remember picture,overlay]   \draw[dashed,dash pattern={on 4pt off 2pt}] ([xshift=0.5\tabcolsep,yshift=7pt]e2.north) -- ([xshift=0.5\tabcolsep,yshift=-2pt]f2.south);\tikz[remember picture,overlay]   \draw[dashed,dash pattern={on 4pt off 2pt}] ([xshift=0.5\tabcolsep,yshift=7pt]g2.north) -- ([xshift=0.5\tabcolsep,yshift=-2pt]h2.south);\\    
\end{split}  
\end{equation}
where $\eta^{\alpha}_{21},\eta^{\alpha}_{22},\cdots,\eta^{\alpha}_{2s_{t_{\alpha}^+}}$, $1\leq\alpha\leq l$, are viewed as their images under $\Sigma^{\lambda}$ (see (\ref{t1h1})).

Similarly to {case III} in Lemma \ref{ud}, we can show that $\left({\mathcal R}^{\underline s^*}_{L_{\lambda}}\right)^{-1}\left(C^k_N\right)$ is isomorphic to an open subscheme of ${\rm Spec}\,\mathbb Z\left[\epsilon_k^+\right]\times \mathbb M^{\widetilde{\underline  s}}_{\widetilde L_{\lambda}}\times \mathbb M^{\widehat{\underline  s}}_{\widehat L_{\lambda}}$.
By (\ref{subsch}) and the assumption that Proposition \ref{ms} holds for all $2\leq N^{\prime}<N$, we can conclude (\ref{3dudu}) for case II.
\smallskip

The proof of Claim \ref{special} is complete.\,\,\,\,$\endpf$
\medskip

{\bf\noindent Proof of Claim \ref{strip}.}  Similarly to Claim \ref{special}, we shall prove on a case by case basis. \smallskip

{\bf\noindent Case I ($q^+=0$).}  Computation yields that $L_{\lambda}$ takes the form 
\begin{equation*}
\begin{split}
&\left(\begin{matrix}1&0\\0&1\\\end{matrix}\hspace{-0.12in}\begin{matrix} &\hfill\tikzmark{a2}\\\\&\hfill\tikzmark{b2}\end{matrix}\,\,\,\begin{matrix}1&\eta^1_{12}&\eta^1_{13}&\cdots\\\mathfrak h^1_{21}&\mathfrak h^1_{22}&\mathfrak h^1_{23}&\cdots\end{matrix}\hspace{-0.12in}\begin{matrix} &\hfill\tikzmark{g2}\\\\&\hfill\tikzmark{h2}\end{matrix}\,\,\,\begin{matrix}1&\eta^2_{12}&\eta^2_{13}&\cdots\\\mathfrak h^2_{21}&\mathfrak h^2_{22}&\mathfrak h^2_{23}&\cdots\end{matrix}\hspace{-0.12in}\begin{matrix} &\hfill\tikzmark{e2}\\\\&\hfill\tikzmark{f2}\end{matrix}\,\,\,\begin{matrix}\cdots\cdots\\\cdots\cdots\\\end{matrix}\hspace{-0.12in}\begin{matrix} &\hfill\tikzmark{c2}\\\\&\hfill\tikzmark{d2}\end{matrix}\,\,\,\begin{matrix}1&\eta^l_{12}&\eta^l_{13}&\cdots\\\mathfrak h^l_{21}&\mathfrak h^l_{22}&\mathfrak h^l_{23}&\cdots\end{matrix}\right).\tikz[remember picture,overlay]   \draw[dashed,dash pattern={on 4pt off 2pt}] ([xshift=0.5\tabcolsep,yshift=7pt]a2.north) -- ([xshift=0.5\tabcolsep,yshift=-2pt]b2.south);\tikz[remember picture,overlay]   \draw[dashed,dash pattern={on 4pt off 2pt}] ([xshift=0.5\tabcolsep,yshift=7pt]c2.north) -- ([xshift=0.5\tabcolsep,yshift=-2pt]d2.south);\tikz[remember picture,overlay]   \draw[dashed,dash pattern={on 4pt off 2pt}] ([xshift=0.5\tabcolsep,yshift=7pt]e2.north) -- ([xshift=0.5\tabcolsep,yshift=-2pt]f2.south);\tikz[remember picture,overlay]   \draw[dashed,dash pattern={on 4pt off 2pt}] ([xshift=0.5\tabcolsep,yshift=7pt]g2.north) -- ([xshift=0.5\tabcolsep,yshift=-2pt]h2.south);\\
\end{split}  \end{equation*}
Then $\mathbb M^{\underline s^*}_{L_{\lambda}}$ is smooth over ${\rm Spec}\,\mathbb Z$ by (\ref{4frakn52}) and Lemma \ref{induc}.
\smallskip

{\bf\noindent Case II ($N=3$ and $q^+=1$).} Note that $L_{\lambda}$ takes the form
\begin{equation*}
\left(\begin{matrix}
1&0\\0&\epsilon^+_2\\\end{matrix}\hspace{-0.12in}\begin{matrix} &\hfill\tikzmark{a2}\\
\\&\hfill\tikzmark{b2}
\end{matrix}\,\,\,\begin{matrix}1&\eta^1_{12}&\cdots&\eta^1_{1\delta^+_1}&\cdots&\eta^1_{1s_{t_1^+}}\\
\mathfrak h^1_{21}&\mathfrak h^1_{22}&\cdots&\mathfrak h^1_{2\delta^+_1}&\cdots&\mathfrak h^1_{2s_{t_1^+}}\end{matrix}\right).
\tikz[remember picture,overlay]   \draw[dashed,dash pattern={on 4pt off 2pt}] ([xshift=0.5\tabcolsep,yshift=7pt]a2.north) -- ([xshift=0.5\tabcolsep,yshift=-2pt]b2.south);
\end{equation*}
Set $\widetilde s:=(1,s_{t_1^+})$ and define a morphism $\widetilde L_{\lambda}:{\rm Spec}\,\mathbb Z\left[\epsilon_2^+,\overrightarrow {\mathscr H}^1,\overrightarrow {\mathfrak H}^1\right]\rightarrow \mathcal{G}(2,1+s_{t_1^+})$ by
\begin{equation*}
\small
\begin{split}
&\left(\begin{matrix}
1\\0\end{matrix}\hspace{-0.12in}\begin{matrix} &\hfill\tikzmark{a12}\\
\\&\hfill\tikzmark{b12}
\end{matrix}\,\,\,\begin{matrix}1-\mathfrak h^1_{21}\eta^1_{1\delta^+_1}&\cdots&\eta^1_{1\left(\delta^+_1-1\right)}-\mathfrak h^1_{2\left(\delta^+_1-1\right)}\eta^1_{1\delta^+_1}&0&\eta^1_{1\left(\delta^+_1+1\right)}-\mathfrak h^1_{2\left(\delta^+_1+1\right)}\eta^1_{1\delta^+_1}&\cdots&\eta^1_{1s_{t_1^+}}-\mathfrak h^1_{2s_{t_1^+}}\eta^1_{1\delta^+_1}\\
\mathfrak h^1_{21}&\cdots&\mathfrak h^1_{2\left(\delta^+_1-1\right)}&1&\mathfrak h^1_{2\left(\delta^+_1+1\right)}&\cdots&\mathfrak h^1_{2s_{t_1^+}}\end{matrix}\right).\\
\tikz[remember picture,overlay]   \draw[dashed,dash pattern={on 4pt off 2pt}] ([xshift=0.5\tabcolsep,yshift=7pt]a12.north) -- ([xshift=0.5\tabcolsep,yshift=-2pt]b12.south);
\end{split}
\end{equation*}
Computation yields that $\mathbb M^{\underline s^*}_{L_{\lambda}}\cong{\rm Spec}\,\mathbb Z\left[\epsilon^+_2\right]\times\mathbb M^{\underline {\widetilde s}}_{\widetilde L_{\lambda}}$. Hence, $\mathbb M^{\underline s^*}_{L_{\lambda}}$ is smooth over ${\rm Spec}\,\mathbb Z$ by {case I} in Claim \ref{special}.
\smallskip

{\bf\noindent Case III ($q^+\geq 1$).} 
Without loss of generality, we assume that $\lambda^+$ is the identity permutation.  We adopt the convention that $\eta_{11}^{\alpha}=\mathfrak h^{\rho}_{2\delta^+_{\rho}}=1$ for $1\leq\alpha\leq l=N-2$ and $1\leq\rho\leq q^+$. 

Define
an open cover $\left\{C^k_N\right\}_{k=2}^{q^++2}$  of $\mathbb A_+:={\rm Spec}\,\mathbb Z\left[\epsilon_2^+,\cdots,\epsilon_{q^++2}^+,\overrightarrow {\mathscr H}^1,\cdots,\overrightarrow {\mathscr H}^l,\overrightarrow {\mathfrak H}^1,\cdots,\overrightarrow {\mathfrak H}^l\right]$ by $C^{q^++2}_N:=\left\{\mathfrak p\in \mathbb A_+\left|\,\epsilon^+_{\gamma}\notin\mathfrak p,\,\forall\,\,2\leq\gamma\leq q^++1\right.\right\}$, and, for $2\leq k\leq q^++1$,
\begin{equation*}
\begin{split}
C^k_N:=&\left\{\mathfrak p\in \mathbb A_+\left|\begin{matrix}
&\epsilon^+_{\gamma}\notin\mathfrak p,\,\forall\,2\leq\gamma\leq k-1,\,{\rm and}\\
&1-\eta_{1\delta^+_{\alpha}}^{\alpha}\mathfrak h^{\beta}_{21}\prod\limits_{\gamma=\alpha+1}^{\min\{q^++1,\,\beta\}}\epsilon^+_{\gamma}\notin\mathfrak p,\,\forall\,1\leq\alpha\leq k-1\,\,{\rm and}\,\,k\leq \beta\leq  l\\
\end{matrix}\right.\right\}.\\
\end{split}
\end{equation*}
Then it suffices to show that
\begin{enumerate}[label={({\bf SQ2})},ref={\bf SQ2}]
\item for each $2\leq k\leq q^++2$,  $\left(\mathcal R^{\underline s^*}_{L_{\lambda}}\right)^{-1}\left(C^k_N\right)$ is smooth over ${\rm Spec}\,\mathbb Z$. 
\label{4dudu}
\end{enumerate}

\smallskip

{\bf Subcase IIIA ($k=q^++2$).}  We first study the case of $\delta^+_1\geq 2$. For each subset $I$ of $\left\{q^++1,q^++2,\cdots,l\right\}$, we define an open subscheme $D^I_1$ of $C^{q^++2}_N$ by \begin{equation*}
D^I_1:=\left\{\mathfrak p\in C^{q^++2}_N\left|\,\,\begin{matrix}
\mathfrak h^1_{21}\notin\mathfrak p,\,\, \mathfrak h^j_{21}\notin\mathfrak p,\,\,\forall\,j\in \left\{2,3,\cdots,l\right\}-I\\
\,{\rm and}\,\,1-\mathfrak h^{i}_{21}\eta^1_{1\delta^+_1}\prod\nolimits_{\gamma=2}^{q^++1}{\epsilon^+_{\gamma}}\notin\mathfrak p,\,\forall\,\,i\in I
\end{matrix}\right.\right\};   
\end{equation*}
define $D_2:=\left\{\mathfrak p\in C^{q^++2}_N\left|\,1-\mathfrak h^1_{21}\eta^{1}_{1\delta^+_{1}}\notin\mathfrak p\right.\right\}$.
Note that $C^{q^++2}_N=D_2\bigcup_{I\subset\{q^++1,q^++2,\cdots,l\}}D^I_1$.

We can show that $\mathbb M^{\underline s^*}_{L_{\lambda}}\cong\mathbb M^{\underline s^*}_{\widetilde L_{\lambda}}$, where
\begin{equation*}
\widetilde L_{\lambda}:{\rm Spec}\,\mathbb Z\left[\left(\epsilon_2^+,\epsilon_3^+,\cdots, \epsilon_{q^++1}^+\right),\overrightarrow {\mathscr H}^1,\cdots,\overrightarrow {\mathscr H}^l,\overrightarrow {\mathfrak H}^1,\cdots,\overrightarrow {\mathfrak H}^l\right]\rightarrow \mathcal {G}(2,2+s_{t_1^+}+\cdots+s_{t_l^+})   
\end{equation*}
is defined by
\begin{equation*}
\small
\begin{split}
&\left(\begin{matrix}
1&0\\0&\prod\limits_{\gamma=2}^{q^++1}{\epsilon^+_{\gamma}}\end{matrix}\hspace{-0.12in}\begin{matrix} &\hfill\tikzmark{a12}\\
\\\\&\hfill\tikzmark{b12}
\end{matrix}\,\,\,\begin{matrix}\eta^1_{11}-\mathfrak h^1_{21}\cdot\eta^1_{1\delta^+_1}&\cdots&\eta^1_{1\left(\delta^+_1-1\right)}-\mathfrak h^1_{2\left(\delta^+_1-1\right)}\cdot\eta^1_{1\delta^+_1}&0&\eta^1_{1\left(\delta^+_1+1\right)}-\mathfrak h^1_{2\left(\delta^+_1+1\right)}\cdot\eta^1_{1\delta^+_1}&\cdots\\\\
\mathfrak h^1_{21}&\cdots&\mathfrak h^1_{2\left(\delta^+_1-1\right)}&1&\mathfrak h^1_{2\left(\delta^+_1+1\right)}&\cdots\end{matrix}\hspace{-0.12in}\begin{matrix} &\hfill\tikzmark{g12}\\\\
\\&\hfill\tikzmark{h12}\end{matrix}\right.\\
&\,\,\,\,\,\,\,\,\begin{matrix} &\hfill\tikzmark{x12}\\\\
\\&\hfill\tikzmark{y12}\end{matrix}\,\,\,\begin{matrix}\frac{1}{\epsilon^+_2}-\mathfrak h^2_{21}\cdot\eta^1_{1\delta^+_1}&\frac{\eta^2_{12}}{\epsilon^+_2}-\mathfrak h^2_{22}\cdot\eta^1_{1\delta^+_1}&\cdots\\\\\mathfrak h^2_{21}&\mathfrak h^2_{22}&\cdots\\
\end{matrix}\hspace{-0.12in}\begin{matrix} &\hfill\tikzmark{e12}\\\\
\\&\hfill\tikzmark{f12}\end{matrix}\,\,\,\begin{matrix}\cdots\cdots\\\\\cdots\cdots\\
\end{matrix}\left.\hspace{-0.18in}\begin{matrix} &\hfill\tikzmark{c12}\\\\
\\&\hfill\tikzmark{d12}\end{matrix}\,\,\,\begin{matrix}\frac{1}{\prod\limits_{\gamma=2}^{q^++1}{\epsilon^+_{\gamma}}}-\mathfrak h^l_{21}\cdot\eta^1_{1\delta^+_1}&\frac{\eta^l_{12}}{\prod\limits_{\gamma=2}^{q^++1}{\epsilon^+_{\gamma}}}-\mathfrak h^l_{22}\cdot\eta^1_{1\delta^+_1}&\cdots\\
\mathfrak h^l_{21}&\mathfrak h^l_{22}&\cdots
\end{matrix}\right).\\
\tikz[remember picture,overlay]   \draw[dashed,dash pattern={on 4pt off 2pt}] ([xshift=0.5\tabcolsep,yshift=7pt]a12.north) -- ([xshift=0.5\tabcolsep,yshift=-2pt]b12.south);\tikz[remember picture,overlay]   \draw[dashed,dash pattern={on 4pt off 2pt}] ([xshift=0.5\tabcolsep,yshift=7pt]c12.north) -- ([xshift=0.5\tabcolsep,yshift=-2pt]d12.south);\tikz[remember picture,overlay]   \draw[dashed,dash pattern={on 4pt off 2pt}] ([xshift=0.5\tabcolsep,yshift=7pt]e12.north) -- ([xshift=0.5\tabcolsep,yshift=-2pt]f12.south);\tikz[remember picture,overlay]   \draw[dashed,dash pattern={on 4pt off 2pt}] ([xshift=0.5\tabcolsep,yshift=7pt]g12.north) -- ([xshift=0.5\tabcolsep,yshift=-2pt]h12.south);\tikz[remember picture,overlay]   \draw[dashed,dash pattern={on 4pt off 2pt}] ([xshift=0.5\tabcolsep,yshift=7pt]x12.north) -- ([xshift=0.5\tabcolsep,yshift=-2pt]y12.south);
\end{split}
\end{equation*}

Define a Class I index $\widehat\tau=\left(\widehat j_1, \widehat j_2,\left(\widehat j^+_1,\cdots,\widehat j^+_{|I|}\right),\left(\widehat j^-_1,\cdots,\widehat j^-_{l-|I|-1}\right)\right)\in \mathbb J^{\underline s^*}$, for each subset $I$ of $\left\{q^++1,q^++2,\cdots,l\right\}$, as follows. Write $I=\left\{i_1<i_2<\cdots<i_{|I|}\right\}$, and $\left\{2,3,\cdots,l\right\}-I=\left\{2<3<\cdots<q^+<k_{q^+}<k_{q^++1}<\cdots<k_{l-|I|-1}\right\}$, where  $|I|$ denotes the cardinality of $I$. Set $\widehat j_1:=1$, $\widehat j_2:=2+\delta_1^+$, and
\begin{equation*}
\left\{\begin{aligned}
&\widehat j_{\beta}^-:=2+s_{t^+_1}+s_{t^+_2}+\cdots+s_{t^+_{\beta}}+\delta^+_{\beta+1},\,\,\,{\rm for}\,\,1\leq\beta\leq q^+-1\\
&\widehat j_{\beta}^-:=2+s_{t^+_1}+s_{t^+_2}+\cdots+s_{t^+_{k_{\beta}}}+1,\,\,\,\,{\rm for}\,\,q^+\leq\beta\leq l-|I|-1\\
&\widehat j_{\alpha}^+:=2+s_{t^+_1}+s_{t^+_2}+\cdots+s_{t^+_{i_{\alpha}-1}}+1,\,\,\,\,{\rm for}\,\,1\leq\alpha\leq |I|\\
\end{aligned}\right..    
\end{equation*}
In what follows, we perform a change of coordinates in $D^I_1$ such that $D^I_1$ is isomorphic to an open subscheme of ${\rm Spec}\,\mathbb Z\left[\overrightarrow U,\overrightarrow Z,\overrightarrow V,\overrightarrow H^1,\cdots,\overrightarrow H^{|I|},\overrightarrow \Xi^1,\cdots,\overrightarrow \Xi^{l-|I|-1}\right]$.  Set $\widehat u_{2}:=\prod\limits_{\gamma=2}^{q^++1}{\epsilon^+_{\gamma}}$; for all $1\leq\gamma\leq s_{t_1}^+$ such that $\gamma\neq\delta^+_1$, set
\begin{equation*}
\widehat z_{\gamma}:=\mathfrak h^1_{2\gamma},\,\,\,\,\,\,\,\,\,\,\,\,\,\,\,\,\,\,\,\,\,\,\,\,
\widehat v_{\gamma}:=\eta^1_{1\gamma}-\mathfrak h^1_{2\gamma}\cdot\eta^1_{1\delta^+_1};
\end{equation*}
for $1\leq\beta\leq q^+-1$, set
\begin{equation*}
\begin{split}
&\widehat\xi^{\beta}_{1\gamma}:=\eta^{\beta+1}_{1\gamma}\left({\prod\nolimits_{\gamma=2}^{\beta+1}{\epsilon^+_{\gamma}}}\right)^{-1}-\mathfrak h^{\beta+1}_{2\gamma}\cdot\eta^1_{1\delta^+_1},\,\,\,\,\,\,\,\,\,\,\,\,\,\,\,\,\,\forall \,\,2\leq\gamma\leq s_{t_{\beta+1}^+},\\
&\widehat\xi^{\beta}_{2\rho}:=\mathfrak h^{\beta+1}_{2\rho},\,\,\,\,\,\,\,\,\,\,\,\,\,\,\,\,\,\,\,\,\,\,\,\,\,\,\,\,\,\,\,\,\,\,\,\,\,\,\,\,\,\,\,\,\,\forall\,\, 1\leq\rho\leq s_{t_{\beta+1}^+}{\rm\,\,such\,\,that\,\,}\rho\neq\delta_{\beta+1}^+;\\ 
\end{split} 
\end{equation*}
for $q^+\leq\beta\leq l-|I|-1$, set
\begin{equation*}
\begin{split}
&\widehat\xi^{\beta}_{1\gamma}:=\eta^{k_{\beta}}_{1\gamma}\left(\mathfrak h^{k_{\beta}}_{21}\prod\nolimits_{\gamma=2}^{q^++1}{\epsilon^+_{\gamma}}\right)^{-1}-\mathfrak h^{k_{\beta}}_{2\gamma}\left(\mathfrak h^{k_{\beta}}_{21}\right)^{-1}\eta^1_{1\delta^+_1},\,\,\forall\,\,1\leq\gamma\leq s_{t_{k_{\beta}+1}^+},\\
&\widehat\xi^{\beta}_{2\rho}:=\mathfrak h^{k_{\beta}}_{2\rho}\left(\mathfrak h^{k_{\beta}}_{21}\right)^{-1},\,\,\,\,\,\,\,\,\,\,\,\,\forall\,\,2\leq\rho\leq s_{t_{k_{\beta}+1}^+};\\    
\end{split} 
\end{equation*}
for $1\leq\alpha\leq |I|$,  set
\begin{equation*}
\begin{split}
&\widehat\eta^{\alpha}_{1\gamma}:=\left(\eta^{i_{\alpha}}_{1\gamma}-\mathfrak h^{i_{\alpha}}_{2\gamma}\eta^1_{1\delta^+_1}\prod\nolimits_{\gamma=2}^{q^++1}{\epsilon^+_{\gamma}}\right)\left(1-\mathfrak h^{i_{\alpha}}_{21}\eta^1_{1\delta^+_{1}}\prod\nolimits_{\gamma=2}^{q^++1}{\epsilon^+_{\gamma}}\right)^{-1},\,\,\forall\,\, 2\leq\gamma\leq s_{t_{i_{\alpha}}^+},\\
&\widehat\eta^{\alpha}_{2\rho}:=\mathfrak h^{i_{\alpha}}_{21}\left(1-\mathfrak h^{i_{\alpha}}_{21}\cdot\eta^1_{1\delta^+_{1}}\prod\nolimits_{\gamma=2}^{q^++1}{\epsilon^+_{\gamma}}\right)^{-1}\prod\nolimits_{\gamma=2}^{q^++1}{\epsilon^+_{\gamma}},\,\,\,\,\,\,\,\,\,\,\,\,\,\,\,\,\,\,\,\,\,\,\,\,\forall\,\, 1\leq\rho\leq s_{t_{i_{\alpha}}^+}.\\
\end{split}
\end{equation*}
Similarly,  we can show that $\left({\mathcal R}^{\underline s^*}_{L_{\lambda}}\right)^{-1}\left(D^I_1\right)$ is isomorphic to an open subscheme of $\mathbb M_{j^{\widehat\tau}}^{\underline s^*}$, where $j^{\widehat\tau}$ is the morphism defined by (\ref{spnab}), (\ref{4frakn5}), 
hence it is smooth over ${\rm Spec}\,\mathbb Z$ by Lemma \ref{induc}.

Regarding $D_2$, define $\widetilde{\underline s}:=(2,s_{t_2^+},\cdots,s_{t_l^+})$. Define a morphism 
\begin{equation*}
\widetilde L_{\lambda}:{\rm Spec}\,\mathbb Z\left[\left(\epsilon_3^+,\epsilon_4^+,\cdots, \epsilon_{q^++1}^+\right),\overrightarrow {\mathscr H}^2,\cdots,\overrightarrow {\mathscr H}^l,\overrightarrow {\mathfrak H}^2,\cdots,\overrightarrow {\mathfrak H}^l\right]\rightarrow \mathcal {G}(2,2+s_{t_2^+}+\cdots+s_{t_l^+}) 
\end{equation*}
by 
\begin{equation*}
\begin{split}
&\left(\begin{matrix}
1&0\\
0&\prod\nolimits_{\gamma=3}^{q^++1}\epsilon^+_{\gamma}\\
\end{matrix}\hspace{-0.12in}\begin{matrix} &\hfill\tikzmark{a2}\\
\\&\hfill\tikzmark{b2}
\end{matrix}\,\,\,\begin{matrix}1&\eta^2_{12}&\cdots\\
\frac{\eta^2_{21}}{\epsilon^+_1\epsilon^+_2}&\frac{\eta^2_{22}}{\epsilon^+_1\epsilon^+_2}&\cdots\end{matrix}\hspace{-0.12in}\begin{matrix} &\hfill\tikzmark{g2}\\
\\&\hfill\tikzmark{h2}\end{matrix}\,\,\,\begin{matrix}1&\eta^3_{12}&\cdots\\
\frac{\eta^3_{21}}{\epsilon^+_1\epsilon^+_2}&\frac{\eta^3_{22}}{\epsilon^+_1\epsilon^+_2}&\cdots
\end{matrix}\hspace{-0.12in}\begin{matrix} &\hfill\tikzmark{e2}\\
\\&\hfill\tikzmark{f2}\end{matrix}\,\,\,\begin{matrix}\cdots\\\cdots\\
\end{matrix}\hspace{-0.12in}\begin{matrix} &\hfill\tikzmark{c2}\\
\\&\hfill\tikzmark{d2}\end{matrix}\,\,\,\begin{matrix}1&\eta^l_{12}&\cdots\\
\frac{\eta^l_{21}}{\epsilon^+_1\epsilon^+_2}&\frac{\eta^l_{22}}{\epsilon^+_1\epsilon^+_2}&\cdots
\end{matrix}\right),
\tikz[remember picture,overlay]   \draw[dashed,dash pattern={on 4pt off 2pt}] ([xshift=0.5\tabcolsep,yshift=7pt]a2.north) -- ([xshift=0.5\tabcolsep,yshift=-2pt]b2.south);\tikz[remember picture,overlay]   \draw[dashed,dash pattern={on 4pt off 2pt}] ([xshift=0.5\tabcolsep,yshift=7pt]c2.north) -- ([xshift=0.5\tabcolsep,yshift=-2pt]d2.south);\tikz[remember picture,overlay]   \draw[dashed,dash pattern={on 4pt off 2pt}] ([xshift=0.5\tabcolsep,yshift=7pt]e2.north) -- ([xshift=0.5\tabcolsep,yshift=-2pt]f2.south);\tikz[remember picture,overlay]   \draw[dashed,dash pattern={on 4pt off 2pt}] ([xshift=0.5\tabcolsep,yshift=7pt]g2.north) -- ([xshift=0.5\tabcolsep,yshift=-2pt]h2.south);\\ \end{split} \end{equation*}
which is derived from (\ref{13frakn5}) by deleting the second block and dividing the second row by $\epsilon_2^+$. Similarly to Claim \ref{13gi}, we can show that $\left({\mathcal R}^{\underline s^*}_{L_{\lambda}}\right)^{-1}\left(D_2\right)$ is isomorphic to an open subscheme of ${\rm Spec}\,\mathbb Z\left[\epsilon_2^+,\overrightarrow {\mathscr H}^1,\overrightarrow {\mathfrak H}^1\right]\times \mathbb M^{\widetilde{\underline  s}}_{\widetilde L_{\lambda}}$, hence $\left({\mathcal R}^{\underline s^*}_{L_{\lambda}}\right)^{-1}\left(D_2\right)$ is smooth over ${\rm Spec}\,\mathbb Z$ by the assumption that Proposition \ref{3ms} holds for all $2\leq N^{\prime}<N$.

Similarly, we can prove (\ref{4dudu}) for the case of $\delta_1^+=1$. We omit the details here for brevity.

\smallskip

{\bf Subcase IIIB ($2\leq k\leq q^++1$).} Define new size vectors $\widetilde{\underline s}:=\left(1,s_{t_1^+},s_{t_2^+},\cdots,s_{t_{k-1}^+}\right)$, and $\widehat{\underline s}:=\left(2,s_{t_{k}^+},s_{t_{k+1}^+},\cdots,s_{t_l^+}\right)$. Define a morphism 
\begin{equation*}
\widetilde{L}_{\lambda}: {\rm Spec}\,\mathbb Z\left[\left(\epsilon_2^+,\epsilon_3^+,\cdots,\epsilon_{k-1}^+\right),\overrightarrow {\mathscr H}^1,\cdots,\overrightarrow {\mathscr H}^{k-1},\overrightarrow {\mathfrak H}^1,\cdots,\overrightarrow {\mathfrak H}^{k-1}\right]\rightarrow  \mathcal{G}(2,1+s_{t_2^+}+\cdots+s_{t_{k-1}^+}) 
\end{equation*}
by (\ref{13frakn9}), and a morphism 
\begin{equation*}
\widehat{L}_{\lambda}: {\rm Spec}\,\mathbb Z\left[\left(\epsilon_{k+1}^+,\epsilon_{k+2}^+,\cdots,\epsilon_{l}^+\right),\overrightarrow {\mathscr H}^k,\cdots,\overrightarrow {\mathscr H}^{l},\overrightarrow {\mathfrak H}^k,\cdots,\overrightarrow {\mathfrak H}^{l}\right]\rightarrow  \mathcal {G}(2,2+s_{t_k^+}+\cdots+s_{t_{l}^+}) \end{equation*}
similarly to (\ref{13frakn6}) by
\begin{equation*}
\left(\begin{matrix}
1&0\\
0&\prod\limits_{\gamma=k+1}^{q^++1}{\epsilon^+_{\gamma}}\\
\end{matrix}\hspace{-0.12in}\begin{matrix} &\hfill\tikzmark{a2}\\\\
\\&\hfill\tikzmark{b2}
\end{matrix}\,\,\,\begin{matrix}1&\eta^k_{12}&\cdots\\
\frac{\eta^k_{21}}{\prod\limits_{\gamma=2}^{k}{\epsilon^+_{\gamma}}}&\frac{\eta^k_{22}}{\prod\limits_{\gamma=2}^{k}{\epsilon^+_{\gamma}}}&\cdots\end{matrix}\hspace{-0.12in}\begin{matrix} &\hfill\tikzmark{g2}\\\\
\\&\hfill\tikzmark{h2}\end{matrix}\,\,\,\begin{matrix}1&\eta^{k+1}_{12}&\cdots\\
\frac{\eta^{k+1}_{21}}{\prod\limits_{\gamma=2}^{k}{\epsilon^+_{\gamma}}}&\frac{\eta^{k+1}_{22}}{\prod\limits_{\gamma=2}^{k}{\epsilon^+_{\gamma}}}&\cdots
\end{matrix}\hspace{-0.12in}\begin{matrix} &\hfill\tikzmark{e2}\\\\
\\&\hfill\tikzmark{f2}\end{matrix}\,\,\,\begin{matrix}\cdots\cdots\\\\\cdots\cdots\\
\end{matrix}\hspace{-0.12in}\begin{matrix} &\hfill\tikzmark{c2}\\\\
\\&\hfill\tikzmark{d2}\end{matrix}\,\,\,\begin{matrix}1&\eta^{l}_{12}&\cdots\\
\frac{\eta^{l}_{21}}{\prod\limits_{\gamma=2}^{k}{\epsilon^+_{\gamma}}}&\frac{\eta^{l}_{22}}{\prod\limits_{\gamma=2}^{k}{\epsilon^+_{\gamma}}}&\cdots
\end{matrix}\right).
\tikz[remember picture,overlay]   \draw[dashed,dash pattern={on 4pt off 2pt}] ([xshift=0.5\tabcolsep,yshift=7pt]a2.north) -- ([xshift=0.5\tabcolsep,yshift=-2pt]b2.south);\tikz[remember picture,overlay]   \draw[dashed,dash pattern={on 4pt off 2pt}] ([xshift=0.5\tabcolsep,yshift=7pt]c2.north) -- ([xshift=0.5\tabcolsep,yshift=-2pt]d2.south);\tikz[remember picture,overlay]   \draw[dashed,dash pattern={on 4pt off 2pt}] ([xshift=0.5\tabcolsep,yshift=7pt]e2.north) -- ([xshift=0.5\tabcolsep,yshift=-2pt]f2.south);\tikz[remember picture,overlay]   \draw[dashed,dash pattern={on 4pt off 2pt}] ([xshift=0.5\tabcolsep,yshift=7pt]g2.north) -- ([xshift=0.5\tabcolsep,yshift=-2pt]h2.south);
\end{equation*}
where $\eta^{\alpha}_{21},\eta^{\alpha}_{22},\cdots,\eta^{\alpha}_{2s_{t_{\alpha}^+}}$, $k\leq\alpha\leq l$, are viewed as their images under $\Sigma^{\lambda}$ (see (\ref{t1h1})). 

Similarly to {case III} in Lemma \ref{ud}, we can show that $\left({\mathcal R}^{\underline s^*}_{L_{\lambda}}\right)^{-1}\left(C^k_N\right)$ is isomorphic to an open subscheme of ${\rm Spec}\,\mathbb Z\left[\epsilon_k^+\right]\times \mathbb M^{\widetilde{\underline  s}}_{\widetilde L_{\lambda}}\times \mathbb M^{\widehat{\underline  s}}_{L_{\widehat\lambda}}$, which is smooth by the assumption that Proposition \ref{3ms} holds for all $2\leq N^{\prime}\leq N-1$.
\smallskip

In conclusion, we complete the proof of Claim \ref{strip}.\,\,\,\,$\endpf$
\medskip

The proof of Lemma \ref{3ud} is complete.\,\,\,\,$\endpf$

\begin{lemma}\label{32ud}
Suppose that  Proposition \ref{3ms} holds for all integers $N^{\prime}$ such that $2\leq N^{\prime}<N$. Then Lemma \ref{3loc} holds for all (truncated) coordinate charts of  {\bf Type (\ref{32up})}.
\end{lemma}
{\bf\noindent Proof of Lemma \ref{32ud}.} It suffices to consider  $\tau=\left(1,2,(j^+_1,\cdots,j^+_l)\right)$ $\in\mathbb J^{\underline s}$ where $l=N-1$,  and $p^+_{1}=\cdots=p^+_{l}=1$. Define a new size vector $\underline {\widetilde s}=(2,s_{t_1}^+,s_{t_2}^+,\cdots,s_{t_l}^+)$. Let
$j^{\tau}:{\rm Spec}\,\mathbb Z\left[\overrightarrow H^1,\cdots,\overrightarrow H^l\right]\rightarrow \mathcal {G}(2,n)$ be given by (\ref{4frakn52}).
Define a morphism\begin{equation*}\mathring L^{\tau}:{\rm Spec}\,\mathbb Z\left[\overrightarrow H^1,\cdots,\overrightarrow H^l,a,a^{-1}\right]\rightarrow \mathcal {G}(2,n)   \end{equation*} by 
\begin{equation}\label{wma}
\left(\begin{matrix}
1&0\\
0&1\\
\end{matrix}\hspace{-0.12in}\begin{matrix} &\hfill\tikzmark{a2}\\
\\&\hfill\tikzmark{b2}
\end{matrix}\,\,\,\begin{matrix}a&a\cdot\eta^1_{12}&\cdots\\
a\cdot\eta^1_{21}&a\cdot\eta^1_{22}&\cdots\end{matrix}\hspace{-0.12in}\begin{matrix} &\hfill\tikzmark{g2}\\
\\&\hfill\tikzmark{h2}\end{matrix}\,\,\,\begin{matrix}a&a\cdot\eta^2_{12}&\cdots\\
a\cdot\eta^2_{21}&a\cdot\eta^2_{22}&\cdots
\end{matrix}\hspace{-0.12in}\begin{matrix} &\hfill\tikzmark{e2}\\
\\&\hfill\tikzmark{f2}\end{matrix}\,\,\,\begin{matrix}\cdots\cdots\\\cdots\cdots\\
\end{matrix}\hspace{-0.12in}\begin{matrix} &\hfill\tikzmark{c2}\\
\\&\hfill\tikzmark{d2}\end{matrix}\,\,\,\begin{matrix}a&a\cdot\eta^l_{12}&\cdots\\
a\cdot\eta^l_{21}&a\cdot\eta^l_{22}&\cdots
\end{matrix}\right).
\tikz[remember picture,overlay]   \draw[dashed,dash pattern={on 4pt off 2pt}] ([xshift=0.5\tabcolsep,yshift=7pt]a2.north) -- ([xshift=0.5\tabcolsep,yshift=-2pt]b2.south);\tikz[remember picture,overlay]   \draw[dashed,dash pattern={on 4pt off 2pt}] ([xshift=0.5\tabcolsep,yshift=7pt]c2.north) -- ([xshift=0.5\tabcolsep,yshift=-2pt]d2.south);\tikz[remember picture,overlay]   \draw[dashed,dash pattern={on 4pt off 2pt}] ([xshift=0.5\tabcolsep,yshift=7pt]e2.north) -- ([xshift=0.5\tabcolsep,yshift=-2pt]f2.south);\tikz[remember picture,overlay]   \draw[dashed,dash pattern={on 4pt off 2pt}] ([xshift=0.5\tabcolsep,yshift=7pt]g2.north) -- ([xshift=0.5\tabcolsep,yshift=-2pt]h2.south);
\end{equation}
It is clear that $\mathbb M^{\underline {\widetilde s}}_{\mathring L^{\tau}}\cong {\rm Spec}\,\mathbb Z\left[a,a^{-1}\right]\times\mathbb M^{\underline {\widetilde s}}_{j^{\tau}}$. Left multiplying (\ref{wma}) by $\left(\begin{matrix}
-\eta^1_{21}&1\\
a^{-1}&0\\
\end{matrix}\right)$,  we derive
\begin{equation*}
\begin{split}
&\left(\begin{matrix}
-\eta^1_{21}&1\\
a^{-1}&0\\
\end{matrix}\hspace{-0.12in}\begin{matrix} &\hfill\tikzmark{a12}\\
\\&\hfill\tikzmark{b12}
\end{matrix}\,\,\,\begin{matrix}0&a\left(\eta^1_{22}-\eta^1_{12}\cdot\eta^1_{21}\right)&\cdots\\1&\eta^1_{12}&\cdots\\
\end{matrix}\hspace{-0.12in}\begin{matrix} &\hfill\tikzmark{g12}\\
\\&\hfill\tikzmark{h12}\end{matrix}\,\,\,\begin{matrix}
a\left(\eta^2_{21}-\eta^1_{21}\right)&a\left(\eta^2_{22}-\eta^2_{12}\cdot\eta^1_{21}\right)&\cdots\\
1&\eta^2_{12}&\cdots\\
\end{matrix}\hspace{-0.12in}\begin{matrix} &\hfill\tikzmark{e12}\\
\\&\hfill\tikzmark{f12}\end{matrix}\right.\\
&\left.\,\,\,\,\,\,\,\,\,\,\,\,\,\,\,\,\,\,\,\,\,\,\,\,\,\,\,\,\begin{matrix} &\hfill\tikzmark{x12}\\
\\&\hfill\tikzmark{y12}\end{matrix}\,\,\,\begin{matrix}\cdots\cdots\\\cdots\cdots\\
\end{matrix}\hspace{-0.12in}\begin{matrix} &\hfill\tikzmark{c12}\\
\\&\hfill\tikzmark{d12}\end{matrix}\,\,\,\begin{matrix}
a\left(\eta^l_{21}-\eta^1_{21}\right)&a\left(\eta^l_{22}-\eta^l_{12}\cdot\eta^1_{21}\right)&\cdots\\
1&\eta^l_{12}&\cdots\\
\end{matrix}\right).\\
\tikz[remember picture,overlay]   \draw[dashed,dash pattern={on 4pt off 2pt}] ([xshift=0.5\tabcolsep,yshift=7pt]a12.north) -- ([xshift=0.5\tabcolsep,yshift=-2pt]b12.south);\tikz[remember picture,overlay]   \draw[dashed,dash pattern={on 4pt off 2pt}] ([xshift=0.5\tabcolsep,yshift=7pt]c12.north) -- ([xshift=0.5\tabcolsep,yshift=-2pt]d12.south);\tikz[remember picture,overlay]   \draw[dashed,dash pattern={on 4pt off 2pt}] ([xshift=0.5\tabcolsep,yshift=7pt]e12.north) -- ([xshift=0.5\tabcolsep,yshift=-2pt]f12.south);\tikz[remember picture,overlay]   \draw[dashed,dash pattern={on 4pt off 2pt}] ([xshift=0.5\tabcolsep,yshift=7pt]g12.north) -- ([xshift=0.5\tabcolsep,yshift=-2pt]h12.south);\tikz[remember picture,overlay]   \draw[dashed,dash pattern={on 4pt off 2pt}] ([xshift=0.5\tabcolsep,yshift=7pt]x12.north) -- ([xshift=0.5\tabcolsep,yshift=-2pt]y12.south);\\ 
\end{split}    
\end{equation*}
Define a Class I index $\widetilde\tau=\left(\widetilde j_1,\widetilde j_2,\left(\widetilde j^-_1,\cdots,\widetilde j^-_{m}\right)\right)$ $\in \mathbb J^{\underline {\widetilde s}}$ where $m=N-2$ by
\begin{equation*}
\widetilde j_1:=2,\,\,\,\,\,\,\widetilde j_2:=s_1+1,\,\,\,\,\,\,\widetilde j_{\beta}^-:=1+s_1+s_{t^+_1}+s_{t^+_2}+\cdots+s_{t^+_{\beta-1}},\,\,\,{\rm for}\,\,1\leq\beta\leq m=N-2.    
\end{equation*}
Similarly to (\ref{spnab}), (\ref{4frakn5}), we define a morphism $j^{\widetilde\tau}:{\rm Spec}\,\mathbb Z\left[\overrightarrow U,\overrightarrow V,\overrightarrow \Xi^1,\cdots,\overrightarrow \Xi^m\right]\longrightarrow \mathcal G(2,n)$
by 
\begin{equation*}
\begin{split}
&\left(\begin{matrix}
0&1\\
u_1&0\\
\end{matrix}\hspace{-0.12in}\begin{matrix} &\hfill\tikzmark{a2}\\
\\&\hfill\tikzmark{b2}
\end{matrix}\,\,\,\begin{matrix}
0&v_2&\cdots\\
1&0&\cdots\\
\end{matrix}\hspace{-0.12in}\begin{matrix} &\hfill\tikzmark{a12}\\
\\&\hfill\tikzmark{b12}
\end{matrix}\,\,\,\begin{matrix}\xi^1_{11}&\xi^1_{12}&\xi^1_{13}&\cdots\\
1&\xi^1_{22}&\xi^1_{23}&\cdots\end{matrix}\hspace{-0.12in}\begin{matrix} &\hfill\tikzmark{g12}\\
\\&\hfill\tikzmark{h12}\end{matrix}\,\,\,\begin{matrix}\xi^2_{11}&\xi^2_{12}&\xi^2_{13}&\cdots\\
1&\xi^2_{22}&\xi^2_{23}&\cdots
\end{matrix}\hspace{-0.12in}\begin{matrix} &\hfill\tikzmark{a22}\\
\\&\hfill\tikzmark{b22}\end{matrix}\,\,\,\begin{matrix}\cdots\cdots\\\cdots\cdots\\
\end{matrix}\hspace{-0.12in}\begin{matrix} &\hfill\tikzmark{c22}\\
\\&\hfill\tikzmark{d22}\end{matrix}\,\,\,\begin{matrix}\xi^m_{11}&\xi^m_{12}&\xi^m_{13}&\cdots\\
1&\xi^m_{22}&\xi^m_{23}&\cdots
\end{matrix}\right).
\tikz[remember picture,overlay]   \draw[dashed,dash pattern={on 4pt off 2pt}] ([xshift=0.5\tabcolsep,yshift=7pt]a2.north) -- ([xshift=0.5\tabcolsep,yshift=-2pt]b2.south);\tikz[remember picture,overlay]   \draw[dashed,dash pattern={on 4pt off 2pt}] ([xshift=0.5\tabcolsep,yshift=7pt]a12.north) -- ([xshift=0.5\tabcolsep,yshift=-2pt]b12.south);\tikz[remember picture,overlay]   \draw[dashed,dash pattern={on 4pt off 2pt}] ([xshift=0.5\tabcolsep,yshift=7pt]g12.north) -- ([xshift=0.5\tabcolsep,yshift=-2pt]h12.south);\tikz[remember picture,overlay]   \draw[dashed,dash pattern={on 4pt off 2pt}] ([xshift=0.5\tabcolsep,yshift=7pt]a22.north) -- ([xshift=0.5\tabcolsep,yshift=-2pt]b22.south);\tikz[remember picture,overlay]   \draw[dashed,dash pattern={on 4pt off 2pt}] ([xshift=0.5\tabcolsep,yshift=7pt]c22.north) -- ([xshift=0.5\tabcolsep,yshift=-2pt]d22.south);\\
\end{split} \end{equation*}
By a change of coordinates as  for {case I} in Claim \ref{special} and for {subcase IIIA} in Claim \ref{strip},
we can conclude that $\mathbb M^{\underline {\widetilde s}}_{\mathring L^{\tau}}$ is isomorphic to an open subscheme of ${\rm Spec}\,\mathbb Z\left[\eta^1_{12},\eta^1_{13},\cdots,\eta^1_{1s_{t_1^+}}\right]\times\mathbb M^{\underline {\widetilde s}}_{j_{\widetilde\tau}}$. 

We then conclude Lemma \ref{32ud}
by Lemma \ref{3ud}.\,\,\,\,$\endpf$

\section{Compactifications of Homogeneous Varieties}\label{ahom}

Similarly to Lemma 2.10 of \cite{Fan}, we can lift the  torus action $T_{\underline s}$ from $G(2,n)$ to $\mathcal T^{\underline s}_{n}$ as follows.
Define a map \begin{equation}\label{intau}
\tau:\{1,2,\cdots,n\}\longrightarrow\{1,2,\cdots,N\}
\end{equation}
such that $1+\sum_{t=1}^{\tau(i)-1}s_t\leq i\leq\sum_{t=1}^{\tau(i)}s_t$ for $1\leq i\leq n$.   Write $T_{\underline s}={\rm Spec}\,\mathbb Z\left[\delta_1,\cdots,\delta_N,\frac{1}{\delta_1},\cdots,\frac{1}{\delta_N}\right]$. Define a torus action $T_{\underline s}\times\mathbb P^{N_{2,n}}\rightarrow\mathbb P^{N_{2,n}}$ by the ring homomorphism
\begin{equation}\label{tpt}
z_{(i,j)}\mapsto\delta_{\tau(i)}\cdot\delta_{\tau(j)}\cdot z_{(i,j)},\,\,\,\forall\,\,(i,j)\in\mathbb I_{2,n}. 
\end{equation}
Similarly,  define torus actions $T_{\underline s}\times\mathbb P^{N_{t}^{\underline s}}\rightarrow\mathbb P^{N_{t}^{\underline s}}$, $1\leq t\leq N$, by
\begin{equation}\label{bfj}
z_{(i,j)}\mapsto\delta_{\tau(i)}\cdot\delta_{\tau(j)}\cdot z_{(i,j)},\,\,\,\forall\,\,(i,j)\in\mathbb I^{\underline s}_{t}.
\end{equation}
Combining (\ref{tpt}), (\ref{bfj}), and the trivial $T_{\underline s}$-action on $\prod\nolimits_{\underline w\in C^{\underline s}}\mathbb {P}^{N^{\underline s}_{\underline w}}$, we derive a $T_{\underline s}$-action on $\mathbb P^{N_{2,n}}\times\prod\nolimits_{t=1}^N\mathbb P^{N^{\underline s}_{t}}\times\prod\nolimits_{\underline w\in C^{\underline s}}\mathbb {P}^{N^{\underline s}_{\underline w}}$. It is easy to verify the $T_{\underline s}$-action descends to the subscheme $\mathcal T^{\underline s}_{n}$ of $\mathbb P^{N_{2,n}}\times\prod\nolimits_{t=1}^N\mathbb P^{N^{\underline s}_{t}}\times\prod\nolimits_{\underline w\in C^{\underline s}}\mathbb {P}^{N^{\underline s}_{\underline w}}$, which is compactible with the standard $T_{\underline s}$-action on $G(2,n)$.

We call an $S$-scheme $X$ a $G$-scheme, if there is an action of a group scheme $G$ on $X$ over $S$. We say that $G$ acts on $X$ trivially if the action is $p_2:G\times_S X\rightarrow X$ the projection on $X$. According to Theorem 2.3 of \cite{Fo}, we can derive that
\begin{definitionlemma}
There exists the largest closed subscheme of $\mathcal T^{\underline s}_n$ on which $T_{\underline s}$ acts trivially. We denote the fixed point scheme by $\left(\mathcal T^{\underline s}_n\right)^{T_{\underline s}}$.
\end{definitionlemma}

We then have

\begin{lemma}\label{homeward}
There is a closed subscheme $\mathfrak M^{\underline s}_n$ of $\mathcal T^{\underline s}_n$ such that: 
\begin{enumerate}

\item The restriction of $\mathcal P^{\underline s}_n:\mathcal T^{\underline s}_n\rightarrow\mathcal M^{\underline s}_n$ to $\mathfrak M^{\underline s}_n$ is an isomorphism from  $\mathfrak M^{\underline s}_n$ to $\mathcal M^{\underline s}_n$.
    
\item $\mathfrak M^{\underline s}_n$ is a connected component of the fixed point scheme $\left(\mathcal T^{\underline s}_n\right)^{T_{\underline s}}$. 
    
\end{enumerate}
 
\end{lemma}
{\noindent\bf Proof of Lemma \ref{homeward}.}
The proof is on a case by case basis in terms of the size vector.  We first study the case of $\underline s=(s_1,\cdots,s_N)$ such that there exists $s_t\geq 2$ for a certain $1\leq t\leq N$. Without loss of generality, we can assume that $s_1\geq 2$.

Define a sub-Grassmannian $\mathcal S$ by
\begin{equation}\label{cals1}
\mathcal S:=\bigcap\nolimits_{2\leq t\leq N} S_t,
\end{equation}
where $S_t$ are the subs-Grassmannians of $G(2,n)$ associated to the ideal sheaf $\mathscr I^{\underline s}_t$. Now define \begin{equation}\label{2s1s2}
\mathfrak M^{\underline s}_n:=\left(\widetilde{\mathcal R}^{\underline s}\circ\overline{\mathcal R}^{\underline s}_n\right)^{-1}\left(\mathcal S\right).  
\end{equation}

Notice that $\left(\widetilde{\mathcal R}^{\underline s}\right)^{-1}\left(\mathcal S\right)$ is covered by  $A^{\tau}$ where $\tau=\left( j_1, j_2,\left(j^+_1,\cdots,j^+_l\right),\left(j^-_1,\cdots,j^-_{m}\right)\right)\in \mathbb J^{\underline s}$ run over all Class II indices such that $1\leq j_1<j_2\leq s_1$.
Recall (\ref{t1ab}), (\ref{ha}), (\ref{xb}), (\ref{t2w}) that
\begin{equation*}
A^{\tau}\cong{\rm Spec}\,\mathbb Z\left[\overrightarrow A,\overrightarrow W,\overrightarrow H^1,\cdots,\overrightarrow H^l,\overrightarrow \Xi^1,\cdots,\overrightarrow \Xi^m\right]. 
\end{equation*}

Define a morphism 
$L^{\tau}$ from $U:={\rm Spec}\,\mathbb Z\left[\overrightarrow W,\overrightarrow H^1,\cdots,\overrightarrow H^l,\overrightarrow \Xi^1,\cdots,\overrightarrow \Xi^m\right]$ to $G(2,n)$ by $\Gamma^{\tau}\circ j$, where $j$ is the embedding from $U$ to ${\rm Spec}\,\mathbb Z\left[\overrightarrow A,\overrightarrow W,\overrightarrow H^1,\cdots,\overrightarrow H^l,\overrightarrow \Xi^1,\cdots,\overrightarrow \Xi^m\right]$ given by
\begin{equation*}
\overrightarrow  A=\left(a_{j^+_1},a_{j^+_2},\cdots,a_{j^+_l},a_{j^-_1},a_{j^-_2},\cdots,a_{j^-_m}\right)\mapsto(1,1,\cdots,1),
\end{equation*}
and $\Gamma^{\tau}$ is the morphism defined by (\ref{ga}), (\ref{gb}), (\ref{dang3}). As in Definition \ref{mat}, we can define a  blow-up of $U$ with respect to the product of ideal sheaves $\prod\nolimits_{\underline w\in C^{\underline s}}\left(L^{\tau}\right)^{-1}\mathscr I^{\underline s}_{\underline w}\cdot\mathcal O_U$; by a slight abuse of notation, we  denote the blow-up by $\mathcal R^{\underline s}_{L^{\tau}}:\mathbb M^{\underline s}_{L^{\tau}}\rightarrow U$. Similarly to (\ref{kabm}), we can conclude that
\begin{equation}\label{ppro}
\left(\overline{\mathcal R}^{\underline s}_n\right)^{-1}\left({A^{\tau}}\right)\cong{\rm Spec\,}\mathbb Z\left[\overrightarrow  A\right]\times\mathbb M^{\underline s}_{L^{\tau}}, 
\end{equation}
and, moreover, the embedding $\left(\overline{\mathcal R}^{\underline s}_n\right)^{-1}\left(A^{\tau}\bigcap\left(\widetilde{\mathcal R}^{\underline s}\right)^{-1}\left(\mathcal S\right)\right)\hookrightarrow\left(\overline{\mathcal R}^{\underline s}_n\right)^{-1}\left({A^{\tau}}\right)$ is induced by the ring homorphism
\begin{equation*}
\overrightarrow  A=\left(a_{j^+_1},a_{j^+_2},\cdots,a_{j^+_l},a_{j^-_1},a_{j^-_2},\cdots,a_{j^-_m}\right)\mapsto(0,0,\cdots,0).
\end{equation*}

One can then show that the restriction of $\mathcal P^{\underline s}_n$ to $\left(\overline{\mathcal R}^{\underline s}_n\right)^{-1}\left(A^{\tau}\bigcap\left(\widetilde{\mathcal R}^{\underline s}\right)^{-1}\left(\mathcal S\right)\right)$ is a locally closed embedding to $\mathcal M^{\underline s}_n$. Since $\tau$ is arbitrary, we can conclude that  the restriction of $\mathcal P^{\underline s}_n$ to $\mathfrak M^{\underline s}_n$ is an isomorphism from  $\mathfrak M^{\underline s}_n$ to $\mathcal M^{\underline s}_n$.

Computation yields that the $T_{\underline s}$-action on (\ref{ppro}) is a product of the trivial action on $\mathbb M^{\underline s}_{L^{\tau}}$ and the action on ${\rm Spec\,}\mathbb Z\left[\overrightarrow  A\right]$ defined by
\begin{equation*}\left(a_{j^+_1},\cdots,a_{j^+_l},a_{j^-_1},\cdots,a_{j^-_m}\right)\mapsto\left(\delta_{j^+_1}a_{j^+_1},\,\,\,\cdots,\,\,\,\delta_{j^+_l}a_{j^+_l},\,\,\,\delta_{j^-_1}a_{j^-_1},\,\,\,\cdots,\,\,\,\delta_{j^-_m}a_{j^-_m}\right).
\end{equation*}
Then, it is easy to verify that $\mathfrak M^{\underline s}_n$ is a connected component of the fixed point scheme of $\mathcal T^{\underline s}_n$ under the $T_{\underline s}$-action.

\medskip 

We next assume that $s_1=s_2=\cdots s_N=1$. Similarly to (\ref{cals1}), we define a sub-Grassmannian 
\begin{equation*}
\mathcal S_{12}:=\bigcap\nolimits_{3\leq i\leq N} S_t.
\end{equation*}
For $1\leq t\leq N$, we write as $\left[u_{t},v_{t}\right]$ the homogeneous coordinates for $\mathbb {P}^{N^{\underline s}_{t}}$. By viewing as closed subschemes of $P^{N_{2,n}}\times\prod\nolimits_{t=1}^N\mathbb {P}^{N^{\underline s}_{t}}$, we can define a closed subscheme $\check{\mathcal S}_{12}$ of $\left(\widetilde{\mathcal R}^{\underline s}\right)^{-1}\left(\mathcal S_{12}\right)$ by vanishing of the monomials $u_3=u_4=\cdots=u_N=0$. We now define
\begin{equation}\label{s1s21}
\mathfrak M^{\underline s}_n:=\left(\overline{\mathcal R}^{\underline s}_n\right)^{-1}\left(\check{\mathcal S}_{12}\right).  
\end{equation} 

Note that $\check {\mathcal S}_{12}$ is isomorphic to a point and is contained in $A^{\tau}$ where $\tau=\left(j_1, j_2,\left(j^+_1,\cdots,j^+_l\right)\right)=\left(1, 2,\left(3,4,\cdots,
N\right)\right)\in \mathbb J^{\underline s}$. Let $\check {\mathcal R}_n$, $\ddot{\mathcal R}_n$ be defined by (\ref{12a1}).
Recall Lemma \ref{fac} that 
$\left(\check{\mathcal R}_n\right)^{-1}\left(A^{\tau}\right)$ is covered by $B^{\lambda}$ where $\lambda$ run over all permutations of $\{1,2,\cdots,n-2\}$.

Define a morphism 
$L_{\lambda}$ from $U:={\rm Spec}\,\mathbb Z\left[\epsilon^+_{2},\epsilon^+_{3},\cdots,\epsilon^+_{N-2}\right]$  to $G(2,n)$ by $\Gamma^{\tau}\circ\check{\mathcal R}_n\circ j$, where $j$ is the embedding from $U$ to ${\rm Spec}\,\mathbb Z\left[\overrightarrow A,\overrightarrow E_+\right]$ given by
\begin{equation*}
\epsilon^+_1\mapsto1,\,\,{\rm and}\,\,\overrightarrow  A=\left(a_{3},a_{4},\cdots,a_{N}\right)\mapsto(1,1,\cdots,1),
\end{equation*}
and $\Gamma^{\tau}$ is the morphism defined by (\ref{ngamma}).
We can define a  blow-up of $U$ with respect to the product of ideal sheaves $\prod\nolimits_{\underline w\in C^{\underline s}}\left(L_{\lambda}\right)^{-1}\mathscr I^{\underline s}_{\underline w}\cdot\mathcal O_U$, which we denote $\mathcal R^{\underline s}_{L_{\lambda}}:\mathbb M^{\underline s}_{L_{\lambda}}\rightarrow U$. Similarly, we can conclude that
\begin{equation}\label{ppta}
\left(\ddot{\mathcal R}_N\right)^{-1}\left(B^{\lambda}\right)\cong{\rm Spec\,}\mathbb Z\left[\overrightarrow  A,\epsilon^+_1\right]\times\mathbb M^{\underline s}_{L_{\lambda}}, 
\end{equation}
and, moreover, the embedding $\left(\ddot{\mathcal R}_N\right)^{-1}\left(\left(\check {\mathcal R}_N\right)^{-1}\left(\check{\mathcal S}_{12}\right)\cap B^{\lambda}\right)\hookrightarrow\left(\ddot{\mathcal R}_N\right)^{-1}\left(B^{\lambda}\right)$ is induced by the ring homorphism
\begin{equation*}\epsilon^+_1\mapsto0,\,\,{\rm and}\,\,\overrightarrow  A=\left(a_{3},a_{4},\cdots,a_{N}\right)\mapsto(0,0,\cdots,0).
\end{equation*}
It is clear that the restriction of  $\mathcal P^{\underline s}_n$ to $\left(\ddot{\mathcal R}_N\right)^{-1}\left(\left(\check {\mathcal R}_N\right)^{-1}\left(\check{\mathcal S}_{12}\right)\cap B^{\lambda}\right)$ is a locally closed embedding to $\mathcal M^{\underline s}_n$, hence the restriction of $\mathcal P^{\underline s}_n$ to $\mathfrak M^{\underline s}_n$ is an isomorphism from  $\mathfrak M^{\underline s}_n$ to $\mathcal M^{\underline s}_n$.

Computation yields that the $T_{\underline s}$-action on (\ref{ppta}) is a product of the trivial action on $\mathbb M^{\underline s}_{L_{\lambda}}$ and the action on ${\rm Spec\,}\mathbb Z\left[\overrightarrow  A,\epsilon^+_1\right]$ defined by
\begin{equation*}
\begin{split}
&\epsilon^+_1\mapsto\delta_2^{-1}\epsilon^+_1,\,\,{\rm and}\,\,\left(a_{3},a_{4},\cdots,a_{N}\right)\mapsto\left(\delta_{3}a_{3},\,\,\,\delta_{4}a_{4},\,\,\,\cdots,\,\,\,\delta_{N}a_{N}\right).\\
\end{split}
\end{equation*}
Then, it is easy to verify that $\mathfrak M^{\underline s}_n$ is a connected component of $\left(\mathcal T^{\underline s}_n\right)^{T_{\underline s}}$.

We complete the proof of Lemma \ref{homeward}.\,\,\,\,\,\,$\endpf$
\medskip 


\begin{corollary}\label{dimen}
The connected components of the fixed point scheme $\left(\mathcal T^{\underline s}_n\right)^{T_{\underline s}}$ are smooth over ${\rm Spec}\,\mathbb Z$, and of relative dimensions $2n-N-3$ or $2n-N-4$. 
\end{corollary}

{\noindent\bf Proof of Corollary \ref{dimen}.} We first assume that $\underline s=(s_1,\cdots,s_t,\cdots,s_N)$ such that $s_t\geq 2$ for $1\leq t\leq N$. By the proof of Lemma \ref{homeward}, it suffices to consider the restriction of $\left(\mathcal T^{\underline s}_n\right)^{T_{\underline s}}$ to $\left(\overline{\mathcal R}^{\underline s}_n\right)^{-1}\left({A^{\tau}}\right)$, where $\tau=\left( j_1, j_2,\left(j^+_1,\cdots,j^+_l\right),\left(j^-_1,\cdots,j^-_{m}\right)\right)\in \mathbb J^{\underline s}$ is a Class I index. Without loss of generality, we may assume that $j_1=1$ and $j_2=s_1+1$.

Let $\check {\mathcal R}_N$ and $\ddot{\mathcal R}_N$ be defined by (\ref{312a1}).
Recall Lemma \ref{3fac} that 
$\left(\check{\mathcal R}_N\right)^{-1}\left(A^{\tau}\right)$ is covered by $B^{\lambda}$ where $\lambda$ run over all indices of  $\Lambda^{\tau}$. According to (\ref{3dec}) and (\ref{subsch}), we can conclude that
\begin{equation}\label{3ppta}
\left(\ddot{\mathcal R}_N\right)^{-1}\left(B^{\lambda}\right)\cong{\rm Spec\,}\mathbb Z\left[\overrightarrow  A,\epsilon^+_1,\epsilon^-_1\right]\times U, 
\end{equation}
where $U$ is a certain quasi-projective scheme over ${\rm Spec}\,\mathbb Z$, and, moreover, the $T_{\underline s}$-action on (\ref{3ppta}) is a product of the trivial action on $U$ and the action on ${\rm Spec\,}\mathbb Z\left[\overrightarrow  A,\epsilon^+_1,\epsilon^-_1\right]$ defined by
\begin{equation*}
\left\{\begin{aligned}
&\epsilon^+_1\mapsto\delta_1\delta_2^{-1}\epsilon^+_1,\,\,\,\,\,\,\,\,\,\,\,\,\,\,\epsilon^-_1\mapsto\delta_2\delta_1^{-1}\epsilon^-_1,\\
&\left(a_{j^+_1},a_{j^+_2},\cdots,a_{j^+_l}\right)\mapsto\left(\delta_{j^+_1}\delta_2^{-1}a_{j^+_1},\,\,\,\delta_{j^+_2}\delta_2^{-1}a_{j^+_2},\,\,\,\cdots,\delta_{j^+_l}\delta_2^{-1}a_{j^+_l}\right),\\
&\left(a_{j^-_1},a_{j^-_2},\cdots,a_{j^-_m}\right)\mapsto\left(\delta_{j^-_1}\delta_1^{-1}a_{j^-_1},\,\,\,\delta_{j^-_2}\delta_1^{-1}a_{j^-_2},\,\,\,\cdots,\delta_{j^-_m}\delta_1^{-1}a_{j^-_m}\right).\\
\end{aligned}\right.
\end{equation*}
Then, it is easy to verify that the restriction of $\left(\mathcal T^{\underline s}_n\right)^{T_{\underline s}}$ to $\left(\ddot{\mathcal R}_N\right)^{-1}\left(B^{\lambda}\right)$ is defined by
\begin{equation*}\epsilon^+_1\mapsto0,\,\,\epsilon^+_1\mapsto0,\,\,\overrightarrow  A=\left(a_{j^+_1},a_{j^+_2},\cdots,a_{j^+_l},a_{j^-_1},a_{j^-_2},\cdots,a_{j^-_m}\right)\mapsto(0,0,\cdots,0),
\end{equation*}
which is smooth over ${\rm Spec}\,\mathbb Z$, and of relative dimensions $2n-N-4$. 

The proof is the same for the remaining cases. For brevity we omit the details here. \,\,\,\,\,\,$\endpf$

\begin{remark}
It is interesting to investigate the relation between the family $\mathcal P^{\underline s}_n:\mathcal T^{\underline s}_n\rightarrow \mathcal M^{\underline s}_n$ and the Chow quotients of Grassmannians by sub-tori (see \cite{GiWu}). 
\end{remark}

\begin{remark} Denote by  $\mathfrak M_1,\cdots,\mathfrak M_{\rho},\cdots,\mathfrak M_p,\mathfrak m_1,\cdots,\mathfrak m_{\gamma},\cdots,\mathfrak m_q$ the connected components of $\left(\mathcal T^{\underline s}_n\right)^{T_{\underline s}}$, where $\mathfrak M_{\rho}$ are of relative dimension $2n-N-3$ and $\mathfrak m_{\gamma}$ are of relative dimension $2n-N-4$.   Then the restriction of $\mathcal P^{\underline s}_n$ to  $\mathfrak M_{\rho}$ is an isomorphism, and  to $\mathfrak m_{\gamma}$ is an embedding such that the image is contained in the boundary $\mathcal M^{\underline s}_n\mathbin{\scaleobj{1.5}{\backslash}}\overline {\rm Gr}^{2,E}_{\mathbb V^{\underline s}}$. Moreover, for each irreducible component $D$ of $\mathcal M^{\underline s}_n\mathbin{\scaleobj{1.5}{\backslash}}\overline {\rm Gr}^{2,E}_{\mathbb V^{\underline s}}$, there exists $\mathfrak m_{\rho}$ such that its image under $\mathcal P^{\underline s}_n$ is $D$.
\end{remark}

{\noindent\bf Proof of Theorem \ref{HCR}.} According to Proposition \ref{3ms} it is clear that  $\mathcal T^{\underline s}_n$ is a smooth, projective scheme over ${\rm Spec}\,\mathbb Z$ , and that the morphism $\mathcal P^{\underline s}_n:\mathcal T^{\underline s}_n\rightarrow\mathcal M^{\underline s}_n$ is flat and ${\rm Aut}(E_1)\times\cdots\times{\rm Aut}(E_N)$-equivariant. Property (C) holds by Lemma \ref{homeward}.

Notice that by our construction of $\widehat {\mathcal R}^{\underline s}_n:\mathcal T^{\underline s}_n\rightarrow G(2,n)$ locally as a certain sequence of blow-ups, the center of each  blow-up is always chosen to be smooth. We can thus prove that the boundary $\mathcal T^{\underline s}_n\mathbin{\scaleobj{1.5}{\backslash}}{\rm Gr}^{2,E}_{\mathbb V^{\underline s}}$ is a simple normal crossing divisor, by repeating the process in the proof of Proposition \ref{3ms}. We complete the proof.\,\,\,\,\,\,$\endpf$

\section{Proof of Theorem \ref{CR} }\label{1short}

In this section, we will identify $\mathcal M^{\underline s}_n$ with Lafforgue's spaces when $r=2$. In what follows, we fix a size vector $\underline s=(s_1,s_2,\cdots,s_N)$ and denote $n=\sum_{t=1}^Ns_{t}$.

\subsection{{Properties of Lafforgue's spaces}}\label{recall}

We first recall notations in \cite{L2}. 

For any family $(d_I)_{I\subset\{1,2,\cdots,N\}}$ of non-negative integers $d_I$ such that
\begin{enumerate}
\item $d_{\emptyset}=0$, $d_{\{1,2,\cdots,N\}}=2$,
\item $d_I+d_J\leq d_{I\cup J}+d_{I\cap J}$ for all $I,J\subset\{1,2,\cdots,N\}$,
\end{enumerate}
we define an entire convex  $S^{(d_I)}\subset\mathbb R^N$ by 
\begin{equation*}
\left\{(x_1,\cdots,x_N)\in\left(\mathbb Z_{\geq0}\right)^N\left|\,\sum\nolimits_{i=1}^Nx_i=2\,\,{\rm and}\,\,\sum\nolimits_{i\in I}x_i\geq d_I\,\,{\rm for\,\,all}\,\,I\subset\{1,2,\cdots,N\}\right.\right\}.
\end{equation*}
We denote by 
$S^{(d_I)}_{\mathbb R}$ the convex hull of $S^{(d_I)}$ in $\mathbb R^N$, and call $S^{(d_I)}$ an entire pave if $S^{(d_I)}_{\mathbb R}$ is of the maximal dimension $N-1$.
For each entire convex  $S:=S^{(d_I)}\subset\mathbb V^{\underline s}$, we define a locally closed subscheme of the Grassmannian ${\rm Gr}^{2,E}\cong G(2,n)$ by
\begin{equation*}
 {\rm Gr}^{2,E}_S:=\left\{F\in E\left|\,\dim(F\cap E_I)= d_{I},\,\,\forall I\subset\{1,2,\cdots,N\}\right.\right\}.   
\end{equation*} 
Since ${\rm Gr}^{2,E}_S$ is preserved by the action of the group ${\rm Aut}(E_1)\times\cdots{\rm Aut}(E_N)$, and in particular the action of its center $\mathbb G^N_m$, we can
define a quasi-projective scheme
$\overline {\rm Gr}^{2,E}_S:=  {\rm Gr}^{2,E}_S/\mathbb G^N_{m}$.

Let $\mathbb G_m$ be the multiplicative group scheme over ${\rm Spec}\,\mathbb Z$, and $\mathbb G^S_m$ be the $\mathbb Z$-torus consisting of  functions $S\rightarrow\mathbb G_m$, and $(\mathbb G^S_m)_{\emptyset}\subset\mathbb G^S_m$ be the sub-torus consisting of affine functions $S\rightarrow\mathbb G_m$. For each entire convex $S\subset \mathbb V^{\underline s}$, we denote by $\mathbb R^S$ the real vector space of functions $S\rightarrow \mathbb R$. Let $\mathcal C^{S}\subset\mathbb R^S$ be the
cone of functions $v:S\rightarrow \mathbb R$ such that for any affine function $l:S\rightarrow \mathbb R$ satisfying
$l\leq v$, the set $\left\{i\in S\left |\, l(i) = v(i)\right.\right\}$ is an entire convex if it is not empty. We call $\underline{S}$ an entire convex paving of $S$ if $\underline{S}$ is a set of entire convexes
$S^{(d_I)}$ such that the following holds.
\begin{enumerate}

\item The convex hulls $S^{(d_I)}_{\mathbb R}$ are of the same dimension as  $S_{\mathbb R}$.

\item The convex hulls $S^{(d_I)}_{\mathbb R}$ form a paving of $S_{\mathbb R}$.

\item $\mathcal C^S_{\underline S}$ is not empty. Here $\mathcal C^S_{\underline S}$ is the convex cone of functions $v:S\rightarrow\mathbb R$ such that for
any element $S^{\prime}$ of $\underline S$ there exists an affine map $l_{S^{\prime}}:S\rightarrow\mathbb R$ satisfying $l_{S^{\prime}}\leq v$ and
$S^{\prime}=\{i\in S |\, l_{S^{\prime}}(i) = v(i)\}$.   
\end{enumerate}

Note that for any entire convex paving $\underline S$ of $S$, 
$\overline{\mathcal C}^S_{\underline S}/\mathcal C_{\emptyset}^S $ is a rational polyhedral convex cone, where $\overline{\mathcal C}^S_{\underline S}$ is the closure of $\mathcal C^S_{\underline S}$ in $\mathbb R^S$, and, moreover, the family of rational polyhedral convex cones
$\overline{\mathcal C}^S_{\underline S}/\mathcal C_{\emptyset}^S $ is a fan in the quotient  $\mathbb R^S/\mathcal C_{\emptyset}^S $. Then, the general theory of toric varieties (see  \cite{De} and  \cite{KKMS})
associates to this fan a normal toric variety $\mathcal A^S$ of torus $\mathcal A^S_{\emptyset}= \mathbb G^S_m/(\mathbb G^S_m)_{\emptyset}$ over ${\rm Spec}\,\mathbb Z$. The orbits in $\mathcal A^S$ are locally closed subschemes  indexed by the entire convex pavings $\underline S$ of $S$, which we denote by $\mathcal A^S_{\underline S}$.

Lafforgue constructed a $\mathbb G^S_m/\mathbb G_m$-torsor $\Omega^{S,E}$  in terms of explicit equations obtained by twisting the Pl\"ucker relations via characters
of $\mathcal A^S$ as follows.
Let $\{P\}$ be a family of homogeneous polynomials defining $\cong{\rm Gr}^{2,E}$ as a closed subscheme of $\mathbb G_m\Big\backslash\left(\prod\nolimits_{\underline v\in \mathbb V^{\underline s}}\wedge^{\underline v}E_{\bullet}\right)-\{0\}$. It is clear that the sub-torus $(\mathbb G^S_m)_{\emptyset}$ acts on each polynomial
$P\in\{P\}$ by a character $\chi_P:(\mathbb G^S_m)_{\emptyset}\rightarrow\mathbb G_m$. For any element $S^{\prime}$ of a paving $\underline S$, choose in $S^{\prime}$ a family $e_{S^{\prime}}$ of $s+1$ points that generates the lattice of integer points in the affine subspace of
$\mathbb R^N$ generated by $S$. The homomorphism
$\mathbb G^S_m\rightarrow \mathbb G^{s+1}_m$ of restriction to points of $e_{S^{\prime}}$ induces an isomorphism $(\mathbb G^S_m)_{\emptyset}\rightarrow\mathbb G^{s+1}_m$
and defines a splitting $b_{e_{S^{\prime}}}: \mathbb G^S_m/(\mathbb G^S_m)_{\emptyset}\rightarrow \mathbb G^S_m$
of the exact sequence $1\rightarrow (\mathbb G^S_m)_{\emptyset} \rightarrow  \mathbb G^S_m\rightarrow \mathbb G^S_m/(\mathbb G^S_m)_{\emptyset} \rightarrow 1$. For any polynomial $P\in\{P\}$ and any element $S^{\prime}$ of a paving $\underline S$, we can thus define a polynomial $P_{e_{S^{\prime}}}$ on the product
$\mathcal A^S_{\emptyset}\times\prod_{\underline v\in S}\wedge^{\underline v}E_{\bullet}$ by 
\begin{equation}\label{twist}P_{e_{S^{\prime}}}\left(\lambda,(x_{\underline v})\right)= P\left(b_{e_{S^{\prime}}}(\lambda)\cdot(x_{\underline v})\right)
\end{equation}
for $\lambda\in\mathcal A^S_{\emptyset}$; $P_{e_{S^{\prime}}}$ extends over the affine open subset $\mathcal A^{\prime}$ of $\mathcal A^S$ associated with
the paving $\underline S$. 
We define the scheme above each affine open $\mathcal A^{\prime}$ of $\mathcal A^S$ associated with
an entire convex paving $\underline S$ as the closed subscheme of
$\mathcal A^{\prime}\times\mathbb G_m \prod_{\underline v\in S}\wedge^{\underline v}E_{\bullet}$
where all polynomials $P_{e_{S^{\prime}}}$ vanish.  They fit together to define $\Omega^{S,E}$ over the whole $\mathcal A^S$.

\begin{theorem}[Theorem 2.4 of \cite{L2}]
\begin{enumerate}[label={\rm (\arabic*)}]

\item $\Omega^{S,E}$ is a closed subscheme of the product scheme
\begin{equation*}
\mathcal A^S\times \mathbb G_m{\Big\backslash\prod\nolimits_{\underline v\in S}\left(\wedge^{\underline v}E_{\bullet}-\{0\}\right)}   
\end{equation*}
over ${\rm Spec}\,\mathbb Z$ such that:
\begin{itemize}
    \item $\Omega^{S,E}$ is stabilized by the double action of ${\rm Aut}\,(E_1)\times\cdots {\rm Aut}\,(E_N)$ and of torus $\mathbb G^S_m/\mathbb G_m$, and thus it comes with an equivariant morphism $\Omega^{S,E}\rightarrow \mathcal A^S$.
    
    \item  For any entire convex paving $\underline S$ of $S$, the fiber of $\Omega^{S,E}$ above the distinguished point $\alpha_{\underline S}$ of the orbit $\mathcal A^S_{\underline S}$ is  ${\rm Gr}^{2,E}_{\underline S}$.
\end{itemize}

\item The quotient $\overline{\Omega}^{S,E}$ of $\Omega^{S,E}$ by the free action of the torus $\mathbb G^S_m/\mathbb G_m$ is a projective
scheme.
\end{enumerate}
\end{theorem}
\begin{theorem}[Corollary 3.3 of \cite{L2}] For any $N\geq 2$ and any entire pave $S$ of $\mathbb V^{\underline s}$, the toric
variety $\mathcal A^S$ is smooth. The codimension of an orbit associated with a paving $\underline S$ is equal to $|\underline S|-1$.
\end{theorem}
\begin{theorem}[Theorem 3.6 of \cite{L2}]\label{3.6}
Let $S$ be an entire pave of
$\mathbb V^{\underline s}$. The equivariant structure morphism
$\Omega^{S,E}\rightarrow\mathcal A^S$ is smooth and surjective. By consequence, the equivariant compactification
$\Omega^{S,E}$  of $\overline {\rm Gr}^{2,E}_S$ is smooth over ${\rm Spec}\,\mathbb Z$ and its boundary is a normal crossing divisor.
\end{theorem}

\subsection{{Identification for \texorpdfstring{$N=4$}{dd}}}\label{n4}  
In this subsection, we shall prove that
\begin{proposition}\label{N=4} $\mathcal M^{\underline s}_n$ is isomorphic to $\overline{\Omega}^{\mathbb V^{\underline s},E}$ over ${\rm Spec}\,\mathbb Z$ when $N=4$.
\end{proposition}
{\bf\noindent Proof of Proposition \ref{N=4}.} Write $\underline s=(s_1,s_2,s_3,s_4)$ and $n=\sum_{t=1}^4s_{t}$. When $s_1=s_2=s_3=s_4=1$, $\mathcal M^{\underline s}_n$ and $\Omega^{\mathbb V^{\underline s},E}$ are both isomorphic to $\mathbb P^1$ over ${\rm Spec}\,\mathbb Z$ (see Example in \S 2.3 of \cite{L2}).

Next, we assume that $s_1=s_2=s_3=1$ and $s_4\geq2$. One can show that $\mathcal M^{\underline s}_n$ is covered by open subschemes $X_{1A}$, $X_{1B}$, $X_2$, $X_{3}$, $X_4$, $X_5$ given in Appendix \ref{x11}, up to permutations of blocks and that of columns within a block. Denote  $\mathcal P:=\{1A,1B,2,3,4,5\}$. 

For $i\in\mathcal P$, denote by $p_{j}$ the projection from the product scheme $X_{i}\times\mathbb G_m^{\mathbb V^{\underline s}}/\mathbb G_m$ to its first factor $X_{i}$, and define a torus action $\mathbb G^{\mathbb V^{\underline s}}_m/\mathbb G_m$ on $X_i\times\mathbb G_m^{\mathbb V^{\underline s}}/\mathbb G_m$ as a product of the standard group action on $\mathbb G_m^{\mathbb V^{\underline s}}/\mathbb G_m$ and the trivial action on $X_i$.

We reduce the proof to the construction of an isomorphism between torsors as follows.
  
\begin{claim}\label{torsor}
Assume that for all $i\in\mathcal P$ there exist locally closed embeddings \begin{equation}\label{f1a}
f_i:X_i\times\mathbb G_m^{\mathbb V^{\underline s}}/\mathbb G_m\longrightarrow \mathcal A^{ \mathbb V^{\underline s}}\times \mathbb G_m{\Big\backslash\prod\nolimits_{\underline i\in \mathbb V^{\underline s}}\left(\wedge^{\underline i}E_{\bullet}-\{0\}\right)}   
\end{equation}
over ${\rm Spec}\,\mathbb Z$ such that the following holds.
\begin{enumerate}[label={(\arabic*)}]
\item The image of $X_i\times\mathbb G_m^{\mathbb V^{\underline s}}/\mathbb G_m$ under $f_i$ is   contained in $\Omega^{\mathbb V^{\underline s},E}$ for a certain $i\in\mathcal P$.
    
\item The morphisms $f_i$ are $\mathbb G^{\mathbb V^{\underline s}}_m/\mathbb G_m$-equivariant.
    
\item For any $j_1,j_2\in\mathcal P$, the images $f_{j_1}\left((X_{j_1}\cap X_{j_2})\times\mathbb G_m^{\mathbb V^{\underline s}}/\mathbb G_m\right)$, $f_{j_2}\left((X_{j_1}\cap X_{j_2})\times\mathbb G_m^{\mathbb V^{\underline s}}/\mathbb G_m\right)$ coincide as locally closed subschemes.  Denote by $f^{-1}_{j_2}$ the inverse isomorphism from the image 
$f_{j_2}\left( X_{j_2}\times\mathbb G_m^{\mathbb V^{\underline s}}/\mathbb G_m\right)$ to $X_{j_2}\times\mathbb G_m^{\mathbb V^{\underline s}}/\mathbb G_m$. Then
$p_{j_2}\circ f_{j_2}^{-1}\circ f_{j_1}=p_{j_1}$ when restricted to $(X_{j_1}\cap X_{j_2})\times\mathbb G_m^{\mathbb V^{\underline s}}/\mathbb G_m$. 
\end{enumerate}
Then there exists an isomorphism from $\mathcal M^{\underline s}_n$ to $\overline{\Omega}^{\mathbb V^{\underline s},E}$ over ${\rm Spec}\,\mathbb Z$.
\end{claim}
{\bf\noindent Proof of Claim \ref{torsor}.} 
Notice that, according to the assumption, from each open subscheme in an open cover of $\mathcal M^{\underline s}_n$, we can construct a locally closed embedding to $\overline{\Omega}^{\mathbb V^{\underline s},E}$ by the standard projection $\Omega^{\mathbb V^{\underline s},E}\rightarrow \overline{\Omega}^{\mathbb V^{\underline s},E}$. Moreover, such locally closed embeddings glue up to a morphism $F$ from $\mathcal M^{\underline s}_n$ to $\overline{\Omega}^{\mathbb V^{\underline s},E}$  over ${\rm Spec}\,\mathbb Z$. 

It is clear that 
$\mathcal M^{\underline s}_n$ and $\overline{\Omega}^{\mathbb V^{\underline s},E}$  are irreducible, and of the same dimension. Then $F$ is a projective, surjective, birational, finite morphism,  for each $f_i$ is a locally closed embedding.
Recall Proposition \ref{3ms} and Theorem \ref{3.6} that  $\mathcal M^{\underline s}_n$ and $\overline{\Omega}^{\mathbb V^{\underline s},E}$  are smooth over ${\rm Spec}\,\mathbb Z$,  respectively. We can thus conclude that $F$ is an isomorphism.

The proof is complete.\,\,\,\,\,$\endpf$
\smallskip

In what follows, we shall construct such embeddings. 

We first compute the toric variety $\mathcal A^S$. Note that $\mathbb V^{\underline s}$ consists of points
\begin{equation*}
\begin{split}
&t_{14}:=(1,0,0,1),\,\, t_{24}:=(0,1,0,1),\,\, t_{34}:=(0,0,1,1),\,\,
t_{44}:=(0,0,0,2),\\
&t_{12}:=(1,1,0,0),\,\, t_{13}:=(1,0,1,0),\,\, t_{23}:=(0,1,1,0).\\
\end{split}
\end{equation*}
We have the following entire paves.
\begin{equation*}
\begin{split}
&S_{12}:=\{t_{12},t_{13},t_{14},t_{23},t_{24}\},\,\,S_{13}:=\{t_{12},t_{13},t_{14},t_{23},t_{34}\},\,\,S_{23}:=\{t_{12},t_{23},t_{24},t_{13},t_{34}\},\\
&S_{24}:=\{t_{12},t_{23},t_{24},t_{14},t_{34}\},\,\,S_{34}:=\{t_{13},t_{23},t_{34},t_{14},t_{24}\},\,\,S_{14}:=\{t_{12},t_{13},t_{14},t_{24},t_{34}\},\\
&S_{4}:=\{t_{14},t_{24},t_{34},t_{44}\},\,\,S_{144}:=\{t_{12},t_{13},t_{14},t_{24},t_{34},t_{44}\},\,\,S_{244}:=\{t_{12},t_{23},t_{24},t_{14},t_{34},t_{44}\},\\
&S_{344}:=\{t_{13},t_{23},t_{34},t_{14},t_{24},t_{44}\},\,\,S_{1234}:=\{t_{12},t_{13},t_{23},t_{24},t_{14},t_{34}\},\\
&S_{12344}:=\{t_{12},t_{13},t_{23},t_{24},t_{14},t_{34},t_{44}\}.\\
\end{split}    
\end{equation*}
Then $\mathbb V^{\underline s}$ has the following entire convex pavings.
\begin{equation*}
\begin{split}
&\underline S^{12}:=\left\{S_{12},S_{344}\right\},\,\,\underline S^{13}:=\left\{S_{13},S_{244}\right\},\,\,\underline S^{23}:=\left\{S_{23},S_{144}\right\},\,\,\underline S^{4}:=\left\{S_{4},S_{1234}\right\},\\
&\underline S^{1244}:=\left\{S_{12},S_{34},S_4\right\},\,\,\underline S^{1344}:=\left\{S_{13},S_{24},S_4\right\},\,\,\underline S^{2344}:=\left\{S_{23},S_{14},S_4\right\},\emptyset:=\{S_{12344}\}.\\
\end{split}    
\end{equation*}
It is clear that
$\mathcal C_{\emptyset}^{\mathbb V^{\underline s}}$ consists of all affine functions $l:\mathbb V^{\underline s}\rightarrow\mathbb R$ such that
\begin{equation*}
\begin{split}
&l(t_{12})=a,\,\, l(t_{13})=b,\,\, l(t_{23})=c,\,\, l(t_{34})=d,\,\, \\
&l(t_{14})=a+d-c,\,\, l(t_{24})=a+d-b,\,\, l(t_{44})=a+2d-b-c,\\
\end{split}
\end{equation*}
where $a,b,c,d\in\mathbb R$. Computation yields that $\mathcal C_{\underline S^{12}}^{\mathbb V^{\underline s}}$ is the convex cone generated by $\mathcal C_{\emptyset}^{\mathbb V^{\underline s}}$ and the function $l_{12}$ defined by
\begin{equation*}
\begin{split}
&l_{12}(t_{12})=l_{12}(t_{13})=l_{12}(t_{23})=l_{12}(t_{34})=0,\,\, \,l_{12}(t_{14})=l_{12}(t_{24})=l_{12}(t_{44})=-1;\\
\end{split}
\end{equation*}
$\mathcal C_{\underline S^{13}}^{\mathbb V^{\underline s}}$ is the convex cone generated by $\mathcal C_{\emptyset}^{\mathbb V^{\underline s}}$ and the function $l_{13}$ defined by
\begin{equation*}
\begin{split}
&l_{13}(t_{12})=l_{13}(t_{13})=l_{13}(t_{14})=l_{13}(t_{23})=l_{13}(t_{34})=0,\,\, \,l_{13}(t_{24})=l_{13}(t_{44})=1;\\
\end{split}
\end{equation*}
$\mathcal C_{\underline S^{23}}^{\mathbb V^{\underline s}}$ is the convex cone generated by $\mathcal C_{\emptyset}^{\mathbb V^{\underline s}}$ and the function $l_{23}$ defined by
\begin{equation*}
\begin{split}
&l_{23}(t_{12})=l_{23}(t_{13})=l_{23}(t_{23})=l_{23}(t_{24})=l_{23}(t_{34})=0,\,\, \,l_{23}(t_{14})=l_{23}(t_{44})=1;\\
\end{split}
\end{equation*}
$\mathcal C_{\underline S^{4}}^{\mathbb V^{\underline s}}$ is the convex cone generated by $\mathcal C_{\emptyset}^{\mathbb V^{\underline s}}$ and the function $l_{4}$ defined by
\begin{equation*}
\begin{split}
&l_{4}(t_{12})=l_{4}(t_{13})=l_{4}(t_{23})=l_{4}(t_{34})=l_{4}(t_{14})=l_{4}(t_{24})=0,\,\, \,l_{4}(t_{44})=1.\\
\end{split}
\end{equation*}
$\overline{\mathcal C}_{\underline S^{1244}}^{\mathbb V^{\underline s}}$ is the convex cone generated by $\mathcal C_{\emptyset}^{\mathbb V^{\underline s}}$ and $l_{12}$, $l_4$;  $\overline{\mathcal C}_{\underline S^{1344}}^{\mathbb V^{\underline s}}$ is the convex cone generated by $\mathcal C_{\emptyset}^{\mathbb V^{\underline s}}$ and $l_{13}$, $l_4$; $\overline{\mathcal C}_{\underline S^{2344}}^{\mathbb V^{\underline s}}$ is the convex cone generated by $\mathcal C_{\emptyset}^{\mathbb V^{\underline s}}$ and $l_{23}$, $l_4$. One can show that $\mathcal A^{\mathbb V^{\underline s}}$ is covered by affine open subschemes \begin{equation*}
\mathcal A_{12}:={\rm Spec}\,\mathbb Z\left[\left(\overline{\mathcal C}_{\underline S^{1244}}^{\mathbb V^{\underline s}}\right)^*\cap M\right],\,\mathcal A_{13}:={\rm Spec}\,\mathbb Z\left[\left(\overline{\mathcal C}_{\underline S^{1344}}^{\mathbb V^{\underline s}}\right)^*\cap M\right],\,\mathcal A_{23}:={\rm Spec}\,\mathbb Z\left[\left(\overline{\mathcal C}_{\underline S^{2344}}^{\mathbb V^{\underline s}}\right)^*\cap M\right].    
\end{equation*} 
Here $M$ is the abelian group of characters of  $\mathbb G^{\mathbb V^{\underline s}}_m$.

Write the coordinate ring of the torus $\mathbb G^{\mathbb V^{\underline s}}_m$ as
\begin{equation*}
\mathbb Z\left[\cdots,x_{\underline v},\frac{1}{x_{\underline v}},\cdots\right]_{\underline v\in\mathbb V^{\underline s}}:=\mathbb Z\left[x_{12},\frac{1}{x_{12}},x_{13},\frac{1}{x_{13}},x_{23},\frac{1}{x_{23}},x_{14},\frac{1}{x_{14}},x_{24},\frac{1}{x_{24}},x_{34},\frac{1}{x_{34}},x_{44},\frac{1}{x_{44}}\right].   
\end{equation*}
For $(i_1,i_2)=(1,2),(1,3),(2,3),(1,4),(2,4),(3,4),(4,4)$, define characters $\chi_{i_1i_2}$ of $\mathbb G^{\mathbb V^{\underline s}}_m$ by
\begin{equation*}
\begin{split}
&\chi_{i_1i_2}(x_{12},x_{13},x_{23},x_{14},x_{24},x_{34},x_{44})=x_{i_1i_2}.\\
\end{split}    
\end{equation*}
Then $M\subset(\mathbb R^{\mathbb V^{\underline s}})^*$ is  generated by  $\chi_{12}$, $\chi_{13}$, $\chi_{23}$, $\chi_{14}$, $\chi_{24}$, $\chi_{34}$, $\chi_{44}$. Computation yields that
\begin{equation*}
\begin{split}
&\left(\overline{\mathcal C}_{\underline S^{1244}}^{\mathbb V^{\underline s}}\right)^*\cap M=\left\{x\left(\chi_{12}+\chi_{34}-\chi_{23}-\chi_{14}\right)+y\left(\chi_{12}+\chi_{34}-\chi_{13}-\chi_{24}\right)\right.\\
&\,\,\,\,\,\,\,\,\,\,\,\,\left.+z\left(\chi_{12}-\chi_{13}-\chi_{23}+2\chi_{34}-\chi_{44}\right)\left|\,x+y+z\geq0,z\leq0,(x,y,z)\in\mathbb Z^3\right.\right\},\\
&\left(\overline{\mathcal C}_{\underline S^{1344}}^{\mathbb V^{\underline s}}\right)^*\cap M=\left\{x\left(\chi_{12}+\chi_{34}-\chi_{23}-\chi_{14}\right)+y\left(\chi_{12}+\chi_{34}-\chi_{13}-\chi_{24}\right)\right.\\
&\,\,\,\,\,\,\,\,\,\,\,\,\left.+z\left(\chi_{12}-\chi_{13}-\chi_{23}+2\chi_{34}-\chi_{44}\right)\left|\,y+z\leq0,z\leq0,(x,y,z)\in\mathbb Z^3\right.\right\},\\
&\left(\overline{\mathcal C}_{\underline S^{2344}}^{\mathbb V^{\underline s}}\right)^*\cap M=\left\{x\left(\chi_{12}+\chi_{34}-\chi_{23}-\chi_{14}\right)+y\left(\chi_{12}+\chi_{34}-\chi_{13}-\chi_{24}\right)\right.\\
&\,\,\,\,\,\,\,\,\,\,\,\,\left.+z\left(\chi_{12}-\chi_{13}-\chi_{23}+2\chi_{34}-\chi_{44}\right)\left|\,x+z\leq0,z\leq0,(x,y,z)\in\mathbb Z^3\right.\right\}.\\
\end{split}    
\end{equation*}
Therefore,
\begin{equation}\label{a12}
\begin{split}
&\mathcal A_{12}:={\rm Spec}\,\mathbb Z\left[\frac{\chi_{13}\chi_{24}}{\chi_{23}\chi_{14}},\frac{\chi_{23}\chi_{14}}{\chi_{13}\chi_{24}},\frac{\chi_{12}\chi_{34}}{\chi_{23}\chi_{14}},\frac{\chi_{13}\chi_{44}}{\chi_{14}\chi_{34}}\right],\\
&\mathcal A_{13}:={\rm Spec}\,\mathbb Z\left[\frac{\chi_{23}\chi_{14}}{\chi_{12}\chi_{34}},\frac{\chi_{12}\chi_{34}}{\chi_{23}\chi_{14}},\frac{\chi_{13}\chi_{24}}{\chi_{12}\chi_{34}},\frac{\chi_{23}\chi_{44}}{\chi_{24}\chi_{34}}\right],\\
&\mathcal A_{23}:={\rm Spec}\,\mathbb Z\left[\frac{\chi_{12}\chi_{34}}{\chi_{13}\chi_{24}},\frac{\chi_{13}\chi_{24}}{\chi_{12}\chi_{34}},\frac{\chi_{23}\chi_{14}}{\chi_{12}\chi_{34}},\frac{\chi_{13}\chi_{44}}{\chi_{14}\chi_{34}}\right].      
\end{split}  
\end{equation} 

Recall the Pl\"ucker coordinate functions in Definition \ref{pluc}. 
For convenience, we denote by $\mathcal P^{\Theta}_{(i_1,i_2)}$ the regular functions  $P_{(i_1,i_2)}\circ e_{12}\circ L_{1A}$  on $X_{1A}$ for all $(i_1,i_2)\in\mathbb I_{2,n}$. Here the embedding $e_{12}:U_{12}\hookrightarrow G(2,n)$ is given by (\ref{ej1j2}), $L_{1A}:X_{1A}\longrightarrow G(2,n)$ is the morphism given by (\ref{x1a}), $\Theta$ is the $2\times n$ matrix defining $L_{1A}$.   Notice that,  $P_{(i_1,i_2)}\circ e_{12}\circ L_{1A}$ can be given by the determinants of $2\times 2$ submatrices of $\Theta$ consisting of the $i_1^{\rm th}$ and $i_2^{\rm th}$ columns.

We view $\mathbb G_m^{\mathbb V^{\underline s}}/\mathbb G_m$ as a subscheme of ${\rm Proj}\,\mathbb Z\left[x_{12},x_{13},x_{23},x_{14},x_{24},x_{34},x_{44}\right]$.
It is clear that $\mathbb G_m\Big\backslash\prod\nolimits_{\underline i\in \mathbb V^{\underline s}}\left(\wedge^{\underline i}E_{\bullet}-\{0\}\right)$ is an open subscheme of  ${\rm Proj}\,\mathbb Z\left[\cdots,z_{(i_1,i_2)},\cdots\right]_{1\leq i_1<i_2\leq n}$. Now, by exploiting the coordinate ring  (\ref{a12}) and that for $X_{1A}$ defined in Appendix \ref{x11}, we can construct a locally closed embedding $f_{1A}$ from $X_{1A}\times\mathbb G_m^{\mathbb V^{\underline s}}/\mathbb G_m$ to $\mathcal A_{12}\times \mathbb G_m{\Big\backslash\prod\nolimits_{\underline i\in \mathbb V^{\underline s}}\left(\wedge^{\underline i}E_{\bullet}-\{0\}\right)}$ by the following ring homomorphism, which is clearly $\mathbb G^{\mathbb V^{\underline s}}_m/\mathbb G_m$-equivariant.
\begin{equation}\label{firsti}
\begin{split}
&\frac{\chi_{13}\chi_{24}}{\chi_{23}\chi_{14}}\mapsto\frac{\mathcal P^{\Theta}_{(1,3)}\mathcal P^{\Theta}_{(2,4)}}{\mathcal P^{\Theta}_{(2,3)}\mathcal P^{\Theta}_{(1,4)}}\frac{x_{23}x_{14}}{x_{13}x_{24}}=\epsilon^+_2\frac{x_{23}x_{14}}{x_{13}x_{24}},\\
&\frac{\chi_{23}\chi_{14}}{\chi_{13}\chi_{24}}\mapsto\frac{\mathcal P^{\Theta}_{(2,3)}\mathcal P^{\Theta}_{(1,4)}}{\mathcal P^{\Theta}_{(1,3)}\mathcal P^{\Theta}_{(2,4)}}\frac{x_{13}x_{24}}{x_{23}x_{14}}=\frac{1}{\epsilon^+_2}\frac{x_{13}x_{24}}{x_{23}x_{14}},\\
&\frac{\chi_{12}\chi_{34}}{\chi_{23}\chi_{14}}\mapsto\frac{\mathcal P^{\Theta}_{(1,2)}\mathcal P^{\Theta}_{(3,4)}}{\mathcal P^{\Theta}_{(2,3)}\mathcal P^{\Theta}_{(1,4)}}\frac{x_{23}x_{14}}{x_{12}x_{34}}=-(1-\epsilon^+_2)\frac{x_{23}x_{14}}{x_{12}x_{34}},\\
&\frac{\chi_{13}\chi_{44}}{\chi_{14}\chi_{34}}\mapsto\frac{\mathcal P^{\Theta}_{(1,3)}\mathcal P^{\Theta}_{(4,5)}}{\mathcal P^{\Theta}_{(1,4)}\mathcal P^{\Theta}_{(3,4)}}\frac{x_{14}x_{34}}{x_{13}x_{44}}=\epsilon^+_2\mathfrak x_{12}\frac{x_{14}x_{34}}{x_{13}x_{44}},\,\,\,\,\,\,\,\,\,\,\,\,\,\,\,\,\,\,\,\,\,\,\,\,\,\,\,\,\,\,\,\,\,\,\,\,\,\,\,\,\,\,\,\,\,\,\,\,\,\,\,\,\,\,\,\,\,\,\,\,\,\,\,\,\,\,\,\,\,\,\,\,\,\,\,\,\,\\
\end{split}
\end{equation}
\begin{equation*}
\begin{split}
&z_{(1,2)}\mapsto\frac{\mathcal P^{\Theta}_{(1,2)}}{\mathcal P^{\Theta}_{(1,2)}}x_{12}=x_{12},\,\,\,z_{(1,3)}\mapsto\frac{\mathcal P^{\Theta}_{(1,3)}}{\mathcal P^{\Theta}_{(1,3)}}x_{13}=x_{13},\,\,\,z_{(2,3)}\mapsto \frac{\mathcal P^{\Theta}_{(2,3)}}{\mathcal P^{\Theta}_{(2,3)}}x_{23}=x_{23},\\
&z_{(1,4)}\mapsto \frac{\mathcal P^{\Theta}_{(1,4)}}{\mathcal P^{\Theta}_{(1,4)}}x_{14}=x_{14},\,\,\,z_{(1,i)}\mapsto \frac{\mathcal P^{\Theta}_{(1,i)}}{\mathcal P^{\Theta}_{(1,4)}}x_{14}=\mathfrak{h}_{22}x_{14},\,\,{\rm for}\,\,i\geq5,\\
&z_{(2,4)}\mapsto\frac{\mathcal P^{\Theta}_{(2,4)}}{\mathcal P^{\Theta}_{(2,4)}}x_{24}=x_{24},\,\,\,z_{(2,5)}\mapsto\frac{\mathcal P^{\Theta}_{(2,5)}}{\mathcal P^{\Theta}_{(2,4)}}x_{24}=\mathfrak x_{12}(1-\epsilon^+_2)x_{24},\\
&z_{(2,i)}\mapsto \frac{\mathcal P^{\Theta}_{(2,i)}}{\mathcal P^{\Theta}_{(2,4)}}x_{24}=\mathfrak x_{12}\mathfrak x_{1(i-3)}(1-\epsilon^+_2)x_{24},\,\,{\rm for}\,\,i\geq6,\\
&z_{(3,4)}\mapsto \frac{\mathcal P^{\Theta}_{(3,4)}}{\mathcal P^{\Theta}_{(3,4)}}x_{34}=x_{34},\,\,\,z_{(3,5)}\mapsto \frac{\mathcal P^{\Theta}_{(3,5)}}{\mathcal P^{\Theta}_{(3,4)}}x_{34}=(\mathfrak h_{22}-\epsilon^+_2\mathfrak x_{12})x_{34},\\
&z_{(3,i)}\mapsto\frac{\mathcal P^{\Theta}_{(3,i)}}{\mathcal P^{\Theta}_{(3,4)}}x_{34}=(\mathfrak h_{2(i-3)}-\epsilon^+_2\mathfrak x_{12}\mathfrak x_{1(i-3)})x_{34},\,\,{\rm for}\,\,i\geq6,\\
&z_{(4,5)}\mapsto\frac{\mathcal P^{\Theta}_{(4,5)}}{\mathcal P^{\Theta}_{(4,5)}}x_{44}=x_{44},\,\,\,z_{(4,i)}\mapsto\frac{\mathcal P^{\Theta}_{(4,i)}}{\mathcal P^{\Theta}_{(4,5)}}x_{44}=\mathfrak x_{1(i-3)}x_{44},\,\,{\rm for}\,\,i\geq6,\\
&z_{(5,i)}\mapsto\frac{\mathcal P^{\Theta}_{(5,i)}}{\mathcal P^{\Theta}_{(4,5)}}x_{44}=(\mathfrak h_{22}\mathfrak x_{1(i-3)}-\mathfrak h_{2(i-3)})x_{44},\,\,{\rm for}\,\,i\geq6,\\
&z_{(i_1,i_2)}\mapsto\frac{\mathcal P^{\Theta}_{(i_1,i_2)}}{\mathcal P^{\Theta}_{(4,5)}}x_{44}=(\mathfrak x_{1(i_2-3)}\mathfrak h_{2(i_1-3)}-\mathfrak x_{1(i_1-3)}\mathfrak h_{2(i_2-3)})x_{44},\,\,{\rm for}\,\,6\leq i_1< i_2.\\
\end{split}    
\end{equation*}

Define a splitting  
$b_{e_{\mathbb V^{\underline s}}}$ from $\mathbb G^{\mathbb V^{\underline s}}_m/(\mathbb G^{\mathbb V^{\underline s}}_m)_{\emptyset}:={\rm Spec}\,\mathbb Z\left[\lambda_1,\frac{1}{\lambda_1},\lambda_2,\frac{1}{\lambda_2},\lambda_3,\frac{1}{\lambda_3}\right]$ to $\mathbb G^{\mathbb V^{\underline s}}_m$
by 
\begin{equation*}
\begin{split}
&x_{12}\mapsto1,\,\,\,x_{13}\mapsto1,\,\,\,x_{23}\mapsto1,\,\,\,x_{34}\mapsto1,\,\,\,x_{14}\mapsto\lambda_2^{-1},\,\,\,x_{24}\mapsto\lambda_1^{-1}\lambda_2^{-1},\,\,\,x_{44}\mapsto\lambda_3\lambda_2^{-1}.\\
\end{split}    
\end{equation*}
Then for any polynomial $P\in\{P\}$, the polynomial $P_{e_{\mathbb V^{\underline s}}}$ on 
$\mathcal A^{\mathbb V^{\underline s}}_{\emptyset}\times\prod_{\underline v\in \mathbb V^{\underline s}}\wedge^{\underline v}E_{\bullet}$ takes the form
\begin{equation*}
P_{e_{\mathbb V^{\underline s}}}\left((\lambda_1,\lambda_2,\lambda_3),(x_{\underline v})\right)= P\left(b_{e_{\mathbb V^{\underline s}}}(\lambda_1,\lambda_2,\lambda_3)\cdot(x_{\underline v})\right).    
\end{equation*}
One can verify that all $P_{e_{\mathbb V^{\underline s}}}$ vanish on $f\left(X_{1A}\times\mathbb G_m^{\mathbb V^{\underline s}}/\mathbb G_m\right)\bigcap\left( \mathcal A^{\mathbb V^{\underline s}}_{\emptyset}\times \mathbb G_m{\Big\backslash\prod\nolimits_{\underline i\in \mathbb V^{\underline s}}\left(\wedge^{\underline i}E_{\bullet}-\{0\}\right)}\right)$. Therefore,  $f_{1A}\left(X_{1A}\times\mathbb G_m^{\mathbb V^{\underline s}}/\mathbb G_m\right)$ is contained in $\Omega^{\mathbb V^{\underline s},E}$.

In the following,  we shall construct a $\mathbb G^{\mathbb V^{\underline s}}_m/\mathbb G_m$-equivariant locally closed embedding $f_2$ from $X_2\times\mathbb G_m^{\mathbb V^{\underline s}}/\mathbb G_m$ to $\mathcal A_{12}\times \mathbb G_m{\Big\backslash\prod\nolimits_{\underline i\in \mathbb V^{\underline s}}\left(\wedge^{\underline i}E_{\bullet}-\{0\}\right)}$.  For convenience, we denote by $\mathcal P^{\Theta}_{(i_1,i_2)}$ the regular functions  $P_{(i_1,i_2)}\circ e_{12}\circ L_{2}$  on $X_{2}$ for all $(i_1,i_2)\in\mathbb I_{2,n}$. Here $L_{2}:X_{2}\longrightarrow G(2,n)$ is the morphism given by (\ref{x2}), $\Theta$ is the $2\times n$ matrix defining $L_{2}$. Then, the morphism $f_2$ is induced by the following homomorphism between coordinate rings.

\begin{equation}
\begin{split}
&\frac{\chi_{13}\chi_{24}}{\chi_{23}\chi_{14}}\mapsto\frac{\mathcal P^{\Theta}_{(1,3)}\mathcal P^{\Theta}_{(2,4)}}{\mathcal P^{\Theta}_{(2,3)}\mathcal P^{\Theta}_{(1,4)}}\frac{x_{23}x_{14}}{x_{13}x_{24}}=\left(1+\mathfrak x_{11}\left(\mathfrak h_{22}-\eta_{12}\right)\right)\frac{x_{23}x_{14}}{x_{13}x_{24}},\\
&\frac{\chi_{23}\chi_{14}}{\chi_{13}\chi_{24}}\mapsto\frac{\mathcal P^{\Theta}_{(2,3)}\mathcal P^{\Theta}_{(1,4)}}{\mathcal P^{\Theta}_{(1,3)}\mathcal P^{\Theta}_{(2,4)}}\frac{x_{13}x_{24}}{x_{23}x_{14}}=\frac{1}{1+\mathfrak x_{11}\left(\mathfrak h_{22}-\eta_{12}\right)}\frac{x_{13}x_{24}}{x_{23}x_{14}},\\
&\frac{\chi_{12}\chi_{34}}{\chi_{23}\chi_{14}}\mapsto\frac{\mathcal P^{\Theta}_{(1,2)}\mathcal P^{\Theta}_{(3,5)}}{\mathcal P^{\Theta}_{(2,3)}\mathcal P^{\Theta}_{(1,4)}}\frac{x_{23}x_{14}}{x_{12}x_{34}}=-\left(1-\mathfrak x_{11}\eta_{12}\right)\left(\mathfrak h_{22}-\eta_{12}\right)\frac{x_{23}x_{14}}{x_{12}x_{34}},\\
&\frac{\chi_{13}\chi_{44}}{\chi_{14}\chi_{34}}\mapsto\frac{\mathcal P^{\Theta}_{(1,3)}\mathcal P^{\Theta}_{(4,5)}}{\mathcal P^{\Theta}_{(1,4)}\mathcal P^{\Theta}_{(3,5)}}\frac{x_{14}x_{34}}{x_{13}x_{44}}=\frac{1+\mathfrak x_{11}\left(\mathfrak h_{22}-\eta_{12}\right)}{1-\mathfrak x_{11}\eta_{12}}\frac{x_{14}x_{34}}{x_{13}x_{44}},\,\,\,\,\,\,\,\,\,\,\,\,\,\,\,\,\,\,\,\,\,\,\,\,\,\,\,\,\,\,\,\,\,\,\,\,\,\,\,\,\,\,\,\,\,\,\,\,\,\,\,\,\,\,\,\,\,\,\,\,\,\,\\
&z_{(1,2)}\mapsto\frac{\mathcal P^{\Theta}_{(1,2)}}{\mathcal P^{\Theta}_{(1,2)}}x_{12}=x_{12},\,\,\,z_{(1,3)}\mapsto\frac{\mathcal P^{\Theta}_{(1,3)}}{\mathcal P^{\Theta}_{(1,3)}}x_{13}=x_{13},\,\,\,z_{(2,3)}\mapsto \frac{\mathcal P^{\Theta}_{(2,3)}}{\mathcal P^{\Theta}_{(2,3)}}x_{23}=x_{23},\\
&z_{(1,4)}\mapsto \frac{\mathcal P^{\Theta}_{(1,4)}}{\mathcal P^{\Theta}_{(1,4)}}x_{14}=x_{14},\,\,\,z_{(1,i)}\mapsto \frac{\mathcal P^{\Theta}_{(1,i)}}{\mathcal P^{\Theta}_{(1,4)}}x_{14}=\mathfrak{h}_{2(i-3)}x_{14},\,\,{\rm for}\,\,i\geq5,\\
&z_{(2,4)}\mapsto\frac{\mathcal P^{\Theta}_{(2,4)}}{\mathcal P^{\Theta}_{(2,4)}}x_{24}=x_{24},\,\,\,z_{(2,5)}\mapsto \frac{\mathcal P^{\Theta}_{(2,5)}}{\mathcal P^{\Theta}_{(2,4)}}x_{24}=\eta_{12}x_{24},\\
&z_{(2,i)}\mapsto\frac{\mathcal P^{\Theta}_{(2,i)}}{\mathcal P^{\Theta}_{(2,4)}}x_{24}=\mathfrak h_{2(i-3)}+\mathfrak{x}_{1(i-3)}(\mathfrak{h}_{22}-\eta_{12}),\,\,{\rm for}\,\,i\geq6,\,\,\,\\
&z_{(3,5)}\mapsto \frac{\mathcal P^{\Theta}_{(3,5)}}{\mathcal P^{\Theta}_{(3,5)}}x_{34}=x_{34},\,\,\,z_{(3,4)}\mapsto \frac{\mathcal P^{\Theta}_{(3,4)}}{\mathcal P^{\Theta}_{(3,5)}}x_{34}=-\frac{\mathfrak x_{11}}{1-\mathfrak x_{11}\eta_{12}}x_{34},\\
&z_{(3,i)}\mapsto\frac{\mathcal P^{\Theta}_{(3,i)}}{\mathcal P^{\Theta}_{(3,5)}}x_{34}=-\frac{\mathfrak x_{1(i-3)}+\mathfrak x_{11}\mathfrak h_{2(i-3)}+\mathfrak x_{11}\mathfrak x_{1(i-3)}(\mathfrak h_{22}-\eta_{12})}{1-\mathfrak x_{11}\eta_{12}}x_{34},\,\,{\rm for}\,\,i\geq6,\\
&z_{(4,5)}\mapsto\frac{\mathcal P^{\Theta}_{(4,5)}}{\mathcal P^{\Theta}_{(4,5)}}x_{44}=x_{44},\,\,\,z_{(4,i)}\mapsto\frac{\mathcal P^{\Theta}_{(4,i)}}{\mathcal P^{\Theta}_{(4,5)}}x_{44}=\mathfrak x_{1(i-3)}x_{44},\,\,{\rm for}\,\,i\geq6,\\
&z_{(5,i)}\mapsto\frac{\mathcal P^{\Theta}_{(5,i)}}{\mathcal P^{\Theta}_{(4,5)}}x_{44}=\left(\mathfrak x_{1(i-3)}\mathfrak h_{22}+\mathfrak h_{2(i-3)}\right)x_{44},\,\,{\rm for}\,\,i\geq6,\\
&z_{(i_1,i_2)}\mapsto\frac{\mathcal P^{\Theta}_{(i_1,i_2)}}{\mathcal P^{\Theta}_{(4,5)}}x_{44}=\left(\mathfrak x_{1(i_1-3)}\mathfrak h_{2(i_2-3)}-\mathfrak x_{1(i_2-3)}\mathfrak h_{2(i_1-3)}\right)x_{44},\,\,{\rm for}\,\,6\leq i_1< i_2.\\
\end{split}    
\end{equation}
Notice that when restricted to
$X_{1A}\cap X_2$, we have 
\begin{equation*}
\left(P_{(i_1,i_2)}\circ e_{12}\circ L_{1A}\right)\left(P_{(i_3,i_4)}\circ e_{12}\circ L_{1A}\right)^{-1}=\left(P_{(i_1,i_2)}\circ e_{12}\circ L_{2}\right)\left(P_{(i_3,i_4)}\circ e_{12}\circ L_{2}\right)^{-1}.
\end{equation*}
Then  $f_{1A}\left((X_{1A}\cap X_{2})\times\mathbb G_m^{\mathbb V^{\underline s}}/\mathbb G_m\right)$ and $f_{2}\left((X_{1A}\cap X_{2})\times\mathbb G_m^{\mathbb V^{\underline s}}/\mathbb G_m\right)$ coincide as locally closed subschemes, and $p_{2}\circ f_{2}^{-1}\circ f_{1A}=p_{1A}$ when restricted to $(X_{1A}\cap X_{2})\times\mathbb G_m^{\mathbb V^{\underline s}}/\mathbb G_m$.

The construction of the remaining $f_i$ are the same, and we leave the details to Appendix \ref{lx1}. The proof of Proposition \ref{N=4} when $s_1=s_2=s_3=1$ and $s_4\geq2$ is now complete.
\smallskip

Notice that, up to permutations of blocks and that of columns within a block, when $s_1=s_2=1$ and $s_3,s_4\geq2$, $\mathcal M^{\underline s}_n$ is covered by open subschemes $Y_{1A}$, $Y_{1B}$, $\cdots$, $Y_6$ given in Appendix \ref{y1}; when $s_1=1$ and $s_2,s_3,s_4\geq2$, $\mathcal M^{\underline s}_n$ is covered by open subschemes $Z_{1A}$, $Z_{1B}$, $\cdots$, $Z_9$ given in Appendix \ref{z1}; when $s_1,s_2,s_3,s_4\geq2$, $\mathcal M^{\underline s}_n$ is covered by open subschemes $W_{1A}$, $W_{1B}$, $\cdots$, $W_{11D}$ given in Appendix \ref{w1}. Similarly, we can reduce the proof of Proposition \ref{N=4} for the remaining cases to the construction of the local embeddings which satisfy the properties as that in Claim \ref{torsor}. For brevity, we leave the details to Appendices \ref{ly1}, \ref{lz1}, \ref{lw1} correspondingly.

We complete the proof of Proposition \ref{N=4}. 
\,\,\,\,\,$\endpf$
\medskip

Similarly, we can show that
\begin{corollary}\label{N=23} Let $\underline s=(s_1,\cdots,s_N)$ be size vector with $N=2,3$. Then $\mathcal M^{\underline s}_n$ is isomorphic to $\overline{\Omega}^{\mathbb V^{\underline s},E}$ over ${\rm Spec}\,\mathbb Z$, where $n=\sum_{t=1}^Ns_{t}$. 
\end{corollary}
{\bf\noindent Proof of Corollary \ref{N=23}.} The proof is similar, and hence we omit the details for brevity.\,\,\,\,\,$\endpf$
\medskip

Denote the morphism constructed in Proposition \ref{N=4} by $\widehat F^{\underline s}_n:\mathcal M^{\underline s}_n\longrightarrow\overline{\Omega}^{\mathbb V^{\underline s},E}$. Let $\alpha_{\emptyset}$ be the distinguished point corresponding to  $\overline{\mathcal C}^{\mathbb V^{\underline s}}_{\emptyset}/\mathcal C_{\emptyset}^{\mathbb V^{\underline s}}$. Denote by ${\rm Gr}^{2,E}_{\emptyset}$
the fiber  of $\Omega^{\mathbb V^{\underline s},E}$ above $\alpha_{\emptyset}$, and by $\overline {\rm Gr}^{2,E}_{\emptyset}$ the quotient $ {\rm Gr}^{2,E}_{\emptyset}/\mathbb G^N_{m}$. Define a morphism 
\begin{equation*}
G^{\underline s}_n:{\rm Gr}^{2,E}_{\emptyset}\longrightarrow\mathcal M^{\underline s}_n\subset\prod\nolimits_{\underline w\in C^{\underline s}}\mathbb {P}^{N^{\underline s}_{\underline w}}
\end{equation*}
by $\left(\cdots,F^{\underline s}_{\underline w}\circ e ,\cdots\right)_{\underline w\in C^{\underline s}}$ (see (\ref{Fw})). It is easy to verify that $G^{\underline s}_n$ descends to 
\begin{equation*}
{\overline G}^{\underline s}_n:{\overline{\rm Gr}}^{2,E}_{\emptyset}\longrightarrow\mathcal M^{\underline s}_n\subset\prod\nolimits_{\underline w\in C^{\underline s}}\mathbb {P}^{N^{\underline s}_{\underline w}}.
\end{equation*}
\begin{lemma}\label{inverse}
Let $N=2,3,4$. There exists an open subscheme $\check{\mathcal M}$ of $\mathcal M^{\underline s}_n$ such that the restriction to $\check{\mathcal M}$ of the morphism $\left({\overline G}^{\underline s}_n\circ\widehat F^{\underline s}_n\right)$
coincides with the identity morphism on $\check{\mathcal M}$.
\end{lemma}
{\bf\noindent Proof of Lemma \ref{inverse}.}
We first assume that $N=4$, and $s_1=s_2=s_3=1$, $s_4\geq2$. 

Let $\check{\mathcal M}$ be the open subscheme of $X_{1A}\subset\mathcal M^{\underline s}_n$ (see Appendix \ref{x11} for $X_{1A}$) defined by
\begin{equation*}
\check{\mathcal M}:={\rm Spec}\,\mathbb Z\left[\epsilon^+_2,\eta_{12},\eta_{13},\cdots,\eta_{1s_4},\mathfrak h_{22},\mathfrak h_{23},\cdots,\mathfrak h_{2s_4},\frac{1}{\epsilon^+_2},\frac{1}{1-\epsilon^+_2},\frac{1}{\mathfrak h_{22}}\right].    
\end{equation*}
Take a locally closed subscheme of $\check{\mathcal M}\times\mathbb G_m^{\mathbb V^{\underline s}}/\mathbb G_m$ defined by
\begin{equation*}
\begin{split}
&\frac{P_{(1,3)}P_{(2,4)}}{P_{(2,3)}P_{(1,4)}}\frac{x_{23}x_{14}}{x_{13}x_{24}}=\frac{P_{(1,2)}P_{(3,4)}}{P_{(2,3)}P_{(1,4)}}\frac{x_{23}x_{14}}{x_{12}x_{34}}=\frac{P_{(1,3)}P_{(4,5)}}{P_{(1,4)}P_{(3,4)}}\frac{x_{14}x_{34}}{x_{13}x_{44}}=1,\\
\end{split}    
\end{equation*}
where we denote by $P_{(i_1,i_2)}$ the regular functions $P_{(i_1,i_2)}\circ L_{1A}$. It is easy to verify that this subscheme is isomorphic to $\check{\mathcal M}\times\mathbb G^N_m/\mathbb G_m$, and thus we denote it by $\check{\mathcal M}\times\mathbb G^N_m/\mathbb G_m$ for convenience.

One can verify by (\ref{firsti}) that the restriction of the morphism $f_{1A}$ defined in Proposition \ref{N=4} induces an open immersion from $\check{\mathcal M}\times\mathbb G^N_m/\mathbb G_m$ to the fiber ${\rm Gr}^{2,E}_{\emptyset}$ of $\Omega^{\mathbb V^{\underline s},E}$ over the distinguished point $\alpha_{\emptyset}$. Moreover,  the morphism 
$\left.\left(G^{\underline s}_{n}\circ f_{1A}\right)\right|_{\left(\check{\mathcal M}\times\mathbb G^N_m/\mathbb G_m\right)}$ descends to $\check{\mathcal M}$. Comparing with (\ref{x1a}), we can conclude that $\left.\left({\overline G}^{\underline s}_n\circ\widehat F^{\underline s}_n\right)\right|_{\check{\mathcal M}}$ is the identity morphism.

The other cases can be proved in the same way. We omit the details for brevity. \,\,\,\,\,$\endpf$

\subsection{Identification for \texorpdfstring{$N\geq5$}{dd}}\label{n5}
Define lower-dimensional entire convexes $\mathbb V^{\underline s}_j$, $1\leq j\leq N$, by
\begin{equation*}
\mathbb V^{\underline s}_j:=\left\{(x_1,\cdots,x_N)\in \mathbb V^{\underline s}\left|\,x_j=0\right.\right\}.    
\end{equation*} 
Then for any subdivision of $\mathbb V^{\underline s}$ into entire convexes $\mathscr{P}:=\cup P_{\alpha}$, there is a natural subdivision of $\mathbb V^{\underline s}_j$ into entire convexes given by  
\begin{equation*}
\mathscr{P}|_{\mathbb V^{\underline s}_j}:=\bigcup\left( P_{\alpha}\cap \mathbb V^{\underline s}_j\right).    
\end{equation*}

\begin{lemma}
\label{rest}
Assume that $N\geq 5$. Let $\mathscr{P},\mathscr{P}'$ be two subdivisions of $\mathbb V^{\underline s}$ into entire convexes such that $\mathscr{P}'$ is coarser than $\mathscr{P}$. Then there exists an integer $1\leq j\leq N$ such that $\mathscr{P}'|_{ \mathbb V^{\underline s}_j}$ is coarser than $\mathscr{P}|_{\mathbb V^{\underline s}_j}$.
\end{lemma}

{\bf\noindent Proof of Lemma \ref{rest}}
We prove by contradiction. Assume that for each $1\leq j\leq N$, $\mathscr{P}'|_{\mathbb V^{\underline s}_j}$ coincides with $\mathscr{P}|_{\mathbb V^{\underline s}_j}$. Since $\mathscr{P}'$ is coarser than $\mathscr{P}$, we can find entire convexes $P\in\mathscr{P}$ and $P'\in\mathscr{P}'$, and a hyperplane $H$ in the linear span of $\mathbb V^{\underline s}$ such that the following holds.
\begin{enumerate}
    \item  The intersection of $H$ and the convex hull $P_{\mathbb R}$ of $P$ is a facet of $P_{\mathbb R}$;

    \item $H$ contains an interior point of the convex hull $P^{\prime}_{\mathbb R}$ of $P^{\prime}$.
    
\end{enumerate}
After a certain permutation of coordinates $(x_1,\cdots,x_N)$ of $\mathbb R^N$, we can conclude by Lemmas 1.7 and 3.1 of \cite{L2} that $H$ is defined for a certain $1\leq l\leq N-3$ by \begin{equation*}
H:=\left\{(x_1,\cdots,x_N)\in \mathbb R^N\left|\,\sum\nolimits_{i=1}^lx_i=1\,\,{\rm and}\,\,\sum\nolimits_{i=l+1}^Nx_i=1\right.\right\}.
\end{equation*}

If $2\leq l\leq N-3$, $H\cap \mathbb V^{\underline s}_{l+1}$ contains an interior point $\mathbb V^{\underline s}_{l+1}$ . If $l=1$, $\mathbb V^{\underline s}$ must contain the point $(2,0,0,\cdots,0)$, for otherwise $H$ does not contain interior points of $P^{\prime}$. Then $H\cap \mathbb V^{\underline s}_{l+1}$ contain an interior point of $\mathbb V^{\underline s}_{l+1}$
as well.  Hence, there exists a hyperplane $H_{l+1}$ such that the following holds.
\begin{enumerate}
    \item  There exists an entire convex $P^{l+1}\in\mathscr P^{\prime}$ such that the intersection of $H_{l+1}$ and the convex hull $P^{l+1}_{\mathbb R}$ of $P^{l+1}$ is a facet of $P^{l+1}_{\mathbb R}$;

    \item $H_{l+1}\cap \mathbb V^{\underline s}_{l+1}=H\cap \mathbb V^{\underline s}_{l+1}$.
    
\end{enumerate}
By Lemmas 1.7 and 3.1 of \cite{L2} again, we can conclude that 
\begin{equation*}
H_{l+1} =\left\{(x_1,\cdots,x_N)\in \mathbb R^N\left|\,\sum\nolimits_{i=1}^{l+1}x_i=1\,\,{\rm and}\,\,\sum\nolimits_{i=l+2}^Nx_i=1\right.\right\}.
\end{equation*}

Similarly, by considering $\mathbb V^{\underline s}_{l+2}$ we can conclude that there exists an entire convex $P^{l+2}\in\mathscr P^{\prime}$ such that the intersection of the hyperplane
\begin{equation*}
H_{l+2}:=\left\{(x_1,\cdots,x_N)\in\mathbb R^N\left|\,x_{l+2}+\sum\nolimits_{i=1}^{l}x_i=1\,\,{\rm and}\,\,x_{l+1}+\sum\nolimits_{i=l+3}^Nx_i=1\right.\right\},
\end{equation*}
and the convex hull $P^{l+2}_{\mathbb R}$ of $P^{l+2}$ is a facet of $P^{l+2}_{\mathbb R}$.

Now it is easy to show that $H_{l+1}\cap H_{l+2}$ contains an interior point of the convex hull of $\mathbb V^{\underline s}$, which contradicts the fact that $\mathcal P^{\prime}$
is a subdivison. We complete the proof of Lemma \ref{rest}.\,\,\,\,$\endpf$
\medskip

For an entire convex $S\subset\mathbb V^{\underline s}$ and $S^{\prime}$ a face of S,
forgetting the coordinates outside of $S^{\prime}$
\begin{equation*}
\begin{split}
\mathbb G_m{\Big\backslash\prod\nolimits_{\underline v\in S}\left(\wedge^{\underline v}E_{\bullet}-\{0\}\right)}&\longrightarrow \mathbb G_m{\Big\backslash\prod\nolimits_{\underline v\in S^{\prime}}\left(\wedge^{\underline v}E_{\bullet}-\{0\}\right)}\\
(x_v)_{v\in S}&\mapsto(x_v)_{v\in S^{\prime}}\\ 
\end{split}
\end{equation*}
defines a morphism $f_{S,S^{\prime}}:{\rm Gr}^{2,E}_S\rightarrow {\rm Gr}^{2,E}_{S^{\prime}}$. Moreover, the restriction $\mathbb G^S_m\rightarrow\mathbb G^S_{m^{\prime}}$ induces a homomorphism $\mathcal A^S_{\emptyset}\rightarrow \mathcal A^{S^{\prime}}_{\emptyset}$, 
which extends into an equivariant morphism of toric varieties
\begin{equation}\label{toricmorphism}
\mathcal A^S\rightarrow \mathcal A^{S^{\prime}}.
\end{equation}
Lafforgue extended $f_{S,S^{\prime}}$ to a morphism from $\Omega^{S,E}\rightarrow \Omega^{S^{\prime},E}$ above $\mathcal A^S\rightarrow A^{S^{\prime}}$ (see Proposition 2.5 of \cite{L2}). It preserves the action of the torus
$\mathbb G^S_m/\mathbb G_m$, $\mathbb G^{S^{\prime}}_m/\mathbb G_m$ related by the restriction $\mathbb G^S_m\rightarrow\mathbb G^{S^{\prime}}_m$. Moreover, the morphism $f_{S,S^{\prime}}$ sends the fiber ${\rm Gr}^{2,E}_{\emptyset}$ of $\Omega^{S,E}$ above $\alpha_{\emptyset}$ to the fiber ${\rm Gr}^{2,E}_{\emptyset^{\prime}}$ of $\Omega^{S^{\prime},E}$ above $\alpha_{\emptyset^{\prime}}$, where $\emptyset$ and $\emptyset^{\prime}$
are trivial pavings of $S$ and $S^{\prime}$, respectively.

\begin{lemma}
\label{finiteness}
Assume that $N\geq 5$. Then the fibers of the product of face maps
\begin{equation}
\mathcal{A}^{\mathbb V^{\underline s}}\longrightarrow\prod\nolimits_{j=1}^N\mathcal{A}^{\mathbb V^{\underline s}_j}
\end{equation}
do not contain any projective subscheme of positive dimension.
\end{lemma}

{\bf\noindent Proof of Lemma \ref{finiteness}.} Note that the proof of Lemma 3.3 of \cite{ST} works for schemes over ${\rm Spec}\,\mathbb Z$ as well. We repeat the proof therein as follows. 

Let $C\subseteq\mathcal{A}^{\mathbb V^{\underline s}}$ be an irreducible, projective subscheme of positive dimension. Let $S'\subseteq\mathcal{A}^{\mathbb V^{\underline s}}$ be the minimal closed toric stratum containing $C$. More explicitly, we can write $S'=\overline{T'}$ where $\overline{T'}\subseteq\mathcal{A}^{\mathbb V^{\underline s}}$ is the closure of a torus $T'$ of the appropriate dimension. Since the torus $T'$ is affine, $C$ is projective, and $C\cap T'\neq\emptyset$ by the minimality of $S'$, we cannot have that $C\subseteq T'$, hence $C$ intersects the closure of another torus orbit $\overline{T}\subsetneq\overline{T'}$. Then $S=\overline{T}$ is a toric stratum properly contained in $S'$ which intersects $C$ nontrivially. Let $p\in S\cap C$ and $p'\in(S'\setminus S)\cap C$.

The torus $T'$ (resp. $T$) corresponds to a subdivision $\mathscr{P}'$ (resp. $\mathscr{P}$) of~$Q^{\underline s}$ into entire convexes. Since $\overline T\subsetneq \overline {T^{\prime}}$,  $\mathscr{P}'$ is coarser than $\mathscr{P}$. By Lemma \ref{rest} there exists a $1\leq j\leq N$ such that $\mathscr{P}'|_{\mathbb V^{\underline s}_j}$ is coarser than $\mathscr{P}|_{\mathbb V^{\underline s}_j}$. Then $\mathscr{P}'|_{\mathbb V^{\underline s}_j}$ and $\mathscr{P}|_{\mathbb V^{\underline s}_j}$ correspond to two toric strata $\Sigma\subsetneq\Sigma'$ of $\mathcal{A}^{\mathbb V^{\underline s}_j}$. Let $f_j$ be the face map $\mathcal{A}^{\mathbb V^{\underline s}}\rightarrow\mathcal{A}^{\mathbb V^{\underline s}_j}$. Then $f_j(p')\in\Sigma'\setminus\Sigma$ and $f_j(p)\in\Sigma$, which is a contradiction. 

The proof is complete.\,\,\,\,$\endpf$.\medskip

For any $1\leq i_1<i_2<i_3<i_4\leq N$, define lower-dimensional entire convexes $\mathbb V^{\underline s}_{i_1i_2i_3i_4}$ by
\begin{equation*}
\mathbb V^{\underline s}_{i_1i_2i_3i_4}:=\left\{(x_1,\cdots,x_N)\in \mathbb V^{\underline s}\left|\,x_j=0,\,\,\forall\,1\leq j\leq N\,\,{\rm such\,\,that\,\,}j\neq i_1,i_2,i_3,i_4\right.\right\}.    
\end{equation*}
Applying Lemma \ref{finiteness} repeatedly, we can conclude that
\begin{corollary}\label{4finite}
Assume that $N\geq 5$. Then the fibers of the product of face maps
\begin{equation*}
\mathcal{A}^{\mathbb V^{\underline s}}\longrightarrow\prod\nolimits_{1\leq i_1<i_2<i_3<i_4\leq N}\mathcal{A}^{\mathbb V^{\underline s}_{i_1i_2i_3i_4}}
\end{equation*}
do not contain any projective subscheme of positive dimension.
\end{corollary}

{\bf\noindent Proof of Theorem \ref{CR}.}
By Proposition 2.5 of \cite{L2}, we have a morphism 
\begin{equation*}
\small
f:\mathcal{A}^{\mathbb V^{\underline s}}\times\mathbb G_m{\Big\backslash\prod\nolimits_{\underline v\in \mathbb V^{\underline s}}\left(\wedge^{\underline v}E_{\bullet}-\{0\}\right)}\longrightarrow\prod\limits_{1\leq i_1<i_2<i_3<i_4\leq N}\left(\mathcal{A}^{\mathbb V^{\underline s}_{i_1i_2i_3i_4}}\times\mathbb G_m{\Big\backslash\prod\nolimits_{\underline v\in \mathbb V^{\underline s}_{i_1i_2i_3i_4}}\left(\wedge^{\underline v}E_{\bullet}-\{0\}\right)}\right).
\end{equation*}
It descends to a morphism $\mathcal F:\overline{\Omega}^{\mathbb V^{\underline s},E}\rightarrow\prod\nolimits_{1\leq i_1<i_2<i_3<i_4\leq N}\overline{\Omega}^{\mathbb V^{\underline s}_{i_1i_2i_3i_4},E_{i_1i_2i_3i_4}}$, where $E_{i_1i_2i_3i_4}:=E_{i_1}\oplus E_{i_2}\oplus E_{i_3}\oplus E_{i_4}$.
According to Corollary \ref{4finite}, we can conclude that $\mathcal F$ is a finite morphism.  

By Proposition \ref{N=4}, we have the isomorphisms
\begin{equation*}
f_{i_1i_2i_3i_4}:\overline{\Omega}^{\mathbb V^{\underline s}_{i_1i_2i_3i_4},E_{i_1i_2i_3i_4}}\longrightarrow \mathcal M^{(s_{i_1},s_{i_2},s_{i_3},s_{i_4})}_{s_{i_1}+s_{i_2}+s_{i_3}+s_{i_4}}.
\end{equation*}
Define a morphism
\begin{equation*}
\mathcal G:\overline{\Omega}^{\mathbb V^{\underline s},E}\longrightarrow\prod\nolimits_{1\leq i_1<i_2<i_3<i_4\leq N}\mathcal M^{(s_{i_1},s_{i_2},s_{i_3},s_{i_4})}_{s_{i_1}+s_{i_2}+s_{i_3}+s_{i_4}}
\end{equation*}
by $\mathcal G:=\left(\cdots,f_{i_1i_2i_3i_4},\cdots\right)_{1\leq i_1<i_2<i_3<i_4\leq N}\circ\mathcal F$. Notice that the morphism $f$ sends the fiber ${\rm Gr}^{2,E}_{\emptyset}$ of $\Omega^{\mathbb V^{\underline s},E}$ above $\alpha_{\emptyset}$ to the product of fibers ${\rm Gr}^{2,E}_{\emptyset\mathbb V^{\underline s}_{i_1i_2i_3i_4}}$ of $\Omega^{\mathbb V^{\underline s}_{i_1i_2i_3i_4},E_{i_1i_2i_3i_4}}$ above $\alpha_{\emptyset_{i_1i_2i_3i_4}}$, where $\emptyset$ and $\emptyset_{i_1i_2i_3i_4}$
are trivial pavings of $\mathbb V^{\underline s}$ and $\mathbb V^{\underline s}_{i_1i_2i_3i_4}$, respectively. Then by Lemma \ref{inverse} and Remark \ref{last}, we can conclude that the image of $\overline{\Omega}^{\mathbb V^{\underline s},E}$ under $\mathcal G$ in $\prod\nolimits_{1\leq i_1<i_2<i_3<i_4\leq N}\mathcal M^{(s_{i_1},s_{i_2},s_{i_3},s_{i_4})}_{s_{i_1}+s_{i_2}+s_{i_3}+s_{i_4}}$
is isomorphic to $\mathcal M^{\underline s}_n$. 

Since $\mathcal M^{\underline s}$ is smooth and $\mathcal G$ is a finite morphism, we can conclude that $\mathcal G$ is an isomorphism. We complete the proof of Theorem \ref{CR}.\,\,\,\,$\endpf$.

\appendix
\section{Coordinate charts for \texorpdfstring{$\mathcal M^{\underline s}_n$}{dd} with {\texorpdfstring{$N=4$}{dd}}}\label{ccm4}

In this section, we will provide an explicit open cover of $\mathcal M^{\underline s}_n$ with $N=4$, up to permutations of blocks and that of columns within a block, on a case by case basis in terms of $\underline s=(s_1,s_2,s_3,s_4)$. Note that $\mathcal M^{\underline s}_n\cong\mathbb P^1$ when $n=4$.

Notice that thanks to Proposition \ref{homeward} that $\mathcal M^{\underline s}_n$ can be embedded into $\mathcal T^{\underline s}_n$ as a closed subscheme, we can exploit the explicit coordinate charts of $\mathcal T^{\underline s}_n$ to to give local coordinates for  $\mathcal M^{\underline s}_n$.

\subsection{{\texorpdfstring{$s_1=s_2=s_3=1$ and $s_4\geq2$}{dd}}}\label{x11}

Define an affine scheme $X_{1A}$ by
\begin{equation*}
X_{1A}:={\rm Spec}\,\mathbb Z\left[\epsilon^+_2,\mathfrak x_{12},\mathfrak x_{13},\cdots,\mathfrak x_{1s_4},\mathfrak h_{22},\mathfrak h_{23},\cdots,\mathfrak h_{2s_4},(\epsilon^+_2)^{-1}\right],  
\end{equation*}
and an associated morphism $L:X_{1A}\longrightarrow G(2,n)$ by 
\begin{equation}\label{x1a}
\begin{split}
&\left(\begin{matrix}
1\\
0\\
\end{matrix}\hspace{-0.12in}\begin{matrix} &\hfill\tikzmark{a2}\\
&\hfill\tikzmark{b2}\\
\end{matrix}\,\,\,\begin{matrix}
0\\
1\\
\end{matrix}\hspace{-0.12in}\begin{matrix} &\hfill\tikzmark{c2}\\
&\hfill\tikzmark{d2}\\
\end{matrix}\,\,\,\begin{matrix}
1\\
\epsilon_2^+\\
\end{matrix}\hspace{-0.12in}\begin{matrix} &\hfill\tikzmark{e2}\\
&\hfill\tikzmark{f2}\\
\end{matrix}\,\,\,\begin{matrix}
1&\mathfrak x_{12}(1-\epsilon^+_2)+\mathfrak h_{22}&\mathfrak x_{12}\mathfrak x_{13}(1-\epsilon^+_2)+\mathfrak h_{23}&\cdots\\
1&\mathfrak{h}_{22}&\mathfrak{h}_{23}&\cdots\\
\end{matrix}\right)=:\Theta.
\tikz[remember picture,overlay]   
\draw[dashed,dash pattern={on 4pt off 2pt}] ([xshift=0.5\tabcolsep,yshift=7pt]a2.north) -- ([xshift=0.5\tabcolsep,yshift=-2pt]b2.south);
\tikz[remember picture,overlay]   
\draw[dashed,dash pattern={on 4pt off 2pt}] ([xshift=0.5\tabcolsep,yshift=7pt]c2.north) -- ([xshift=0.5\tabcolsep,yshift=-2pt]d2.south);
\tikz[remember picture,overlay]   
\draw[dashed,dash pattern={on 4pt off 2pt}] ([xshift=0.5\tabcolsep,yshift=7pt]e2.north) -- ([xshift=0.5\tabcolsep,yshift=-2pt]f2.south);    
\end{split}  
\end{equation}
We can define a morphism \begin{equation}\label{klx}
\mathfrak K^{\underline s}_{L}:X_{1A}\longrightarrow\prod\nolimits_{\underline w\in C^{\underline s}}\mathbb {P}^{N^{\underline s}_{\underline w}}    
\end{equation} by extending the rational map $\left(\cdots,F^{\underline s}_{\underline w}\circ e\circ {L},\cdots\right)_{\underline w\in C^{\underline s}}$ (see Definition \ref{mat} for  essentially the same construction). One can show that $\mathfrak K^{\underline s}_{L}$ is a locally closed embedding from $X_{1A}$ to $\prod\nolimits_{\underline w\in C^{\underline s}}\mathbb {P}^{N^{\underline s}_{\underline w}}$. By viewing $\mathcal M^{\underline s}_n$ as a subscheme of $\prod\nolimits_{\underline w\in C^{\underline s}}\mathbb {P}^{N^{\underline s}_{\underline w}}$ as in Definition \ref{msss}, we can show that the image of $X_{1A}$ under the morphism $\mathfrak K^{\underline s}_{L}$ is an open subscheme of $\mathcal M^{\underline s}_n$. By a slight abuse of notation, we denote the corresponding open subscheme  by $X_{1A}$ as well. Notice that the parametrization of $X_{1A}$ here can be derived from the embedding given by (\ref{s1s21}).

In what follows, we describe certain affine schemes and their associated morphisms to $G(2,n)$, such that the corresponding images give an open cover of $\mathcal M^{\underline s}_n$, up to permutations of blocks and that of columns within a block.
For brevity, we denote such affine schemes and their images in $\mathcal M^{\underline s}_n$ by the same notation.

\begin{enumerate}
\item[($X_{1B}$)] $X_{1B}:={\rm Spec}\,\mathbb Z\left[\epsilon^+_2,\mathfrak x_{12},\mathfrak x_{13},\cdots,\mathfrak x_{1s_4},\mathfrak h_{22},\mathfrak h_{23},\cdots,\mathfrak h_{2s_4},\frac{1}{1-\epsilon^+_2}\right]$. The associated morphism $L:X_{1B}\longrightarrow G(2,n)$ is defined by (\ref{x1a}).
\smallskip

\item[($X_2$)]  $X_{2}:={\rm Spec}\,\mathbb Z\left[\mathfrak x_{11},\eta_{12},\mathfrak x_{13},\mathfrak x_{14},\cdots,\mathfrak x_{1s_4},\mathfrak h_{22},\mathfrak h_{23},\cdots,\mathfrak h_{2s_4},\frac{1}{1+\mathfrak{x}_{11}(\mathfrak h_{22}-\eta_{12})},\frac{1}{1-\mathfrak{x}_{11}\eta_{12}}\right]$. The associated morphism $L:X_{2}\longrightarrow G(2,n)$ is defined by 
\begin{equation}\label{x2}
\begin{split}
&\left(\begin{matrix}
1\\
0\\
\end{matrix}\hspace{-0.12in}\begin{matrix} &\hfill\tikzmark{a2}\\
&\hfill\tikzmark{b2}\\
\end{matrix}\,\,\,\begin{matrix}
0\\
1\\
\end{matrix}\hspace{-0.12in}\begin{matrix} &\hfill\tikzmark{c2}\\
&\hfill\tikzmark{d2}\\
\end{matrix}\,\,\,\begin{matrix}
1\\
1+\mathfrak{x}_{11}(\mathfrak h_{22}-\eta_{12})\\
\end{matrix}\hspace{-0.12in}\begin{matrix} &\hfill\tikzmark{e2}\\
&\hfill\tikzmark{f2}\\
\end{matrix}\,\,\,\begin{matrix}
1&\eta_{12}&\mathfrak h_{23}+\mathfrak{x}_{13}(\mathfrak{h}_{22}-\eta_{12})&\cdots\\
1&\mathfrak{h}_{22}&\mathfrak{h}_{23}&\cdots\\
\end{matrix}\right)=:\Theta.
\tikz[remember picture,overlay]   
\draw[dashed,dash pattern={on 4pt off 2pt}] ([xshift=0.5\tabcolsep,yshift=7pt]a2.north) -- ([xshift=0.5\tabcolsep,yshift=-2pt]b2.south);
\tikz[remember picture,overlay]   
\draw[dashed,dash pattern={on 4pt off 2pt}] ([xshift=0.5\tabcolsep,yshift=7pt]c2.north) -- ([xshift=0.5\tabcolsep,yshift=-2pt]d2.south);
\tikz[remember picture,overlay]   
\draw[dashed,dash pattern={on 4pt off 2pt}] ([xshift=0.5\tabcolsep,yshift=7pt]e2.north) -- ([xshift=0.5\tabcolsep,yshift=-2pt]f2.south);    
\end{split} 
\end{equation}
\label{x13}


\item[($X_{3}$)] $X_{3}:={\rm Spec}\,\mathbb Z\left[\epsilon^+_2,\eta_{12},\eta_{13},\cdots,\eta_{1s_4},\mathfrak h_{21},\mathfrak h_{23},\mathfrak h_{24},\cdots,\mathfrak h_{2s_4},\frac{1}{1-\eta_{12}\mathfrak h_{21}},\frac{1}{1-\epsilon_2^+\eta_{12}}\right]$. The associated morphism $L:X_{3}\longrightarrow G(2,n)$ is defined by  \begin{equation}\label{x3a}
\begin{split}
&\left(\begin{matrix}
1\\
0\\
\end{matrix}\hspace{-0.12in}\begin{matrix} &\hfill\tikzmark{a2}\\
&\hfill\tikzmark{b2}\\
\end{matrix}\,\,\,\begin{matrix}
0\\
1\\
\end{matrix}\hspace{-0.12in}\begin{matrix} &\hfill\tikzmark{c2}\\
&\hfill\tikzmark{d2}\\
\end{matrix}\,\,\,\begin{matrix}
1\\
\epsilon_2^+\\
\end{matrix}\hspace{-0.12in}\begin{matrix} &\hfill\tikzmark{e2}\\
&\hfill\tikzmark{f2}\\
\end{matrix}\,\,\,\begin{matrix}
1&\eta_{12}&\eta_{13}&\cdots\\
\mathfrak{h}_{21}&1&\mathfrak{h}_{23}&\cdots\\
\end{matrix}\right)=:\Theta.
\tikz[remember picture,overlay]   
\draw[dashed,dash pattern={on 4pt off 2pt}] ([xshift=0.5\tabcolsep,yshift=7pt]a2.north) -- ([xshift=0.5\tabcolsep,yshift=-2pt]b2.south);
\tikz[remember picture,overlay]   
\draw[dashed,dash pattern={on 4pt off 2pt}] ([xshift=0.5\tabcolsep,yshift=7pt]c2.north) -- ([xshift=0.5\tabcolsep,yshift=-2pt]d2.south);
\tikz[remember picture,overlay]   
\draw[dashed,dash pattern={on 4pt off 2pt}] ([xshift=0.5\tabcolsep,yshift=7pt]e2.north) -- ([xshift=0.5\tabcolsep,yshift=-2pt]f2.south);    
\end{split}  
\end{equation}
\label{x14}

\item[($X_4$)] 
$X_{4}:={\rm Spec}\,\mathbb Z\left[\epsilon^+_2,\eta_{12},\eta_{13},\cdots,\eta_{1s_4},\mathfrak h_{22},\mathfrak h_{23},\cdots,\mathfrak h_{2s_4},\frac{1}{1-\epsilon_2^+}\right]$. The associated morphism $L:X_{4}\longrightarrow G(2,n)$ is defined by
\begin{equation}\label{x4}
\begin{split}
&\left(\begin{matrix}
1\\
0\\
\end{matrix}\hspace{-0.12in}\begin{matrix} &\hfill\tikzmark{a2}\\
&\hfill\tikzmark{b2}\\
\end{matrix}\,\,\,\begin{matrix}
0\\
1\\
\end{matrix}\hspace{-0.12in}\begin{matrix} &\hfill\tikzmark{c2}\\
&\hfill\tikzmark{d2}\\
\end{matrix}\,\,\,\begin{matrix}
1\\
1\\
\end{matrix}\hspace{-0.12in}\begin{matrix} &\hfill\tikzmark{e2}\\
&\hfill\tikzmark{f2}\\
\end{matrix}\,\,\,\begin{matrix}
1&\eta_{12}&\eta_{13}&\cdots\\
\epsilon_2^+&\epsilon_2^+\left(\mathfrak{h}_{22}+\eta_{12}\right)&\epsilon_2^+\left(\mathfrak{h}_{23}\mathfrak h_{22}+\eta_{13}\right)&\cdots\\
\end{matrix}\right)=:\Theta.
\tikz[remember picture,overlay]   
\draw[dashed,dash pattern={on 4pt off 2pt}] ([xshift=0.5\tabcolsep,yshift=7pt]a2.north) -- ([xshift=0.5\tabcolsep,yshift=-2pt]b2.south);
\tikz[remember picture,overlay]   
\draw[dashed,dash pattern={on 4pt off 2pt}] ([xshift=0.5\tabcolsep,yshift=7pt]c2.north) -- ([xshift=0.5\tabcolsep,yshift=-2pt]d2.south);
\tikz[remember picture,overlay]   
\draw[dashed,dash pattern={on 4pt off 2pt}] ([xshift=0.5\tabcolsep,yshift=7pt]e2.north) -- ([xshift=0.5\tabcolsep,yshift=-2pt]f2.south);    
\end{split}  
\end{equation}
\label{x15}

\item[($X_5$)] 
$X_{5}:={\rm Spec}\,\mathbb Z\left[\epsilon^+_2,\eta_{12},\eta_{13},\cdots,\eta_{1s_4},\mathfrak h_{21},\mathfrak h_{23},\mathfrak h_{24},\cdots,\mathfrak h_{2s_4},\frac{1}{1-\eta_{12}\mathfrak h_{21}},\frac{1}{1-\epsilon^+_2\mathfrak h_{21}}\right]$. The associated morphism $L:X_{5}\longrightarrow G(2,n)$ is defined by
\begin{equation}\label{x5}
\begin{split}
&\left(\begin{matrix}
1\\
0\\
\end{matrix}\hspace{-0.12in}\begin{matrix} &\hfill\tikzmark{a2}\\
&\hfill\tikzmark{b2}\\
\end{matrix}\,\,\,\begin{matrix}
0\\
1\\
\end{matrix}\hspace{-0.12in}\begin{matrix} &\hfill\tikzmark{c2}\\
&\hfill\tikzmark{d2}\\
\end{matrix}\,\,\,\begin{matrix}
1\\
1\\
\end{matrix}\hspace{-0.12in}\begin{matrix} &\hfill\tikzmark{e2}\\
&\hfill\tikzmark{f2}\\
\end{matrix}\,\,\,\begin{matrix}
1&\eta_{12}&\eta_{13}&\cdots\\
\epsilon_2^+\mathfrak{h}_{21}&\epsilon_2^+&\epsilon_2^+\mathfrak{h}_{23}&\cdots\\
\end{matrix}\right)=:\Theta.
\tikz[remember picture,overlay]   
\draw[dashed,dash pattern={on 4pt off 2pt}] ([xshift=0.5\tabcolsep,yshift=7pt]a2.north) -- ([xshift=0.5\tabcolsep,yshift=-2pt]b2.south);
\tikz[remember picture,overlay]   
\draw[dashed,dash pattern={on 4pt off 2pt}] ([xshift=0.5\tabcolsep,yshift=7pt]c2.north) -- ([xshift=0.5\tabcolsep,yshift=-2pt]d2.south);
\tikz[remember picture,overlay]   
\draw[dashed,dash pattern={on 4pt off 2pt}] ([xshift=0.5\tabcolsep,yshift=7pt]e2.north) -- ([xshift=0.5\tabcolsep,yshift=-2pt]f2.south);    
\end{split}  
\end{equation}

\end{enumerate}

\subsection{{\texorpdfstring{$s_1=s_2=1$  and $s_3, s_4\geq2$}{dd}}}\label{y1} Similarly,  we will describe certain affine schemes 
and the associated morphisms such that  the images give an open cover of $\mathcal M^{\underline s}_n$ up to permutations.

\begin{enumerate}
\item[($Y_{1A}$)] 
$Y_{1A}:={\rm Spec}\,\mathbb Z\left[\epsilon^+_2,\mathfrak x^1_{12},\cdots,\mathfrak x^1_{1s_3},\mathfrak h^1_{22},\cdots,\mathfrak h^1_{2s_3},\mathfrak x^2_{12},\cdots,\mathfrak x^2_{1s_4},\mathfrak h^2_{22},\cdots,\mathfrak h^2_{2s_4},\frac{1}{1-\epsilon^+_2}\right]$. The associated morphism $L:Y_{1}\longrightarrow G(2,n)$ is defined by
\begin{equation}\label{didi}
\begin{split}
&\left(\begin{matrix}1\\
0\\
\end{matrix}\hspace{-0.12in}\begin{matrix} &\hfill\tikzmark{a12}
\\&\hfill\tikzmark{b12}
\end{matrix}\,\,\,\begin{matrix}0\\
1\\
\end{matrix}\hspace{-0.12in}\begin{matrix} &\hfill\tikzmark{g12}
\\&\hfill\tikzmark{h12}\end{matrix}\,\,\,\begin{matrix}1&\mathfrak h^1_{22}+\mathfrak x_{12}^1(1-\epsilon^+_2)&\mathfrak h^1_{23}+\mathfrak x_{12}^1\mathfrak x_{13}^1(1-\epsilon^+_2)&\cdots\\
1&\mathfrak h^1_{22}&\mathfrak h^1_{23}&\cdots\\
\end{matrix}\hspace{-0.12in}\begin{matrix} &\hfill\tikzmark{e12}
\\&\hfill\tikzmark{f12}\end{matrix}\right.\\
&\left.\,\,\,\,\,\,\,\,\,\,\,\,\,\,\,\,\,\,\,\,\begin{matrix} &\hfill\tikzmark{c12}
\\&\hfill\tikzmark{d12}\end{matrix}\,\,\,\begin{matrix}1&\mathfrak x_{12}^2(1-\epsilon^+_2)+\mathfrak h^2_{22}&\mathfrak x_{12}^2\mathfrak x_{13}^2(1-\epsilon^+_2)+\mathfrak h^2_{23}&\cdots\\
\epsilon_2^+&\epsilon_2^+\mathfrak h_{22}^2&\epsilon_2^+\mathfrak h_{23}^2&\cdots\\
\end{matrix}\right)=:\Theta.
\tikz[remember picture,overlay]   \draw[dashed,dash pattern={on 4pt off 2pt}] ([xshift=0.5\tabcolsep,yshift=7pt]a12.north) -- ([xshift=0.5\tabcolsep,yshift=-2pt]b12.south);\tikz[remember picture,overlay]   \draw[dashed,dash pattern={on 4pt off 2pt}] ([xshift=0.5\tabcolsep,yshift=7pt]c12.north) -- ([xshift=0.5\tabcolsep,yshift=-2pt]d12.south);\tikz[remember picture,overlay]   \draw[dashed,dash pattern={on 4pt off 2pt}] ([xshift=0.5\tabcolsep,yshift=7pt]e12.north) -- ([xshift=0.5\tabcolsep,yshift=-2pt]f12.south);\tikz[remember picture,overlay]   \draw[dashed,dash pattern={on 4pt off 2pt}] ([xshift=0.5\tabcolsep,yshift=7pt]g12.north) -- ([xshift=0.5\tabcolsep,yshift=-2pt]h12.south); 
\end{split}    
\end{equation}
\item[($Y_{1B}$)] $Y_{1B}:={\rm Spec}\,\mathbb Z\left[\epsilon^+_2,\mathfrak x^1_{12},\cdots,\mathfrak x^1_{1s_3},\mathfrak h^1_{22},\cdots,\mathfrak h^1_{2s_3},\mathfrak x^2_{12},\cdots,\mathfrak x^2_{1s_4},\mathfrak h^2_{22},\cdots,\mathfrak h^2_{2s_4},\frac{1}{\epsilon^+_2}\right]$. The associated morphism $L:Y_{1B}\longrightarrow G(2,n)$ is defined by (\ref{didi}).
\smallskip 

\item[($Y_{2A}$)] $Y_{2A}:={\rm Spec}\,\mathbb Z\left[\eta^1_{12},\mathfrak x^1_{13},\cdots,\mathfrak x^1_{1s_3},\mathfrak h^1_{22},\cdots,\mathfrak h^1_{2s_3},\mathfrak x^2_{11},\cdots,\mathfrak x^2_{1s_4},\mathfrak h^2_{22},\cdots,\mathfrak h^2_{2s_4},\frac{1}{1+\mathfrak x_{11}^2\left(\eta_{12}^1-\mathfrak h^1_{22}\right)},\frac{1}{\mathfrak x^2_{11}}\right]$. The associated morphism $L:Y_{2A}\longrightarrow G(2,n)$ is defined by
\begin{equation}\label{y3a}
\begin{split}
&\left(\begin{matrix}1\\
0\\
\end{matrix}\hspace{-0.12in}\begin{matrix} &\hfill\tikzmark{a12}\\&\hfill\tikzmark{b12}
\end{matrix}\,\,\,\begin{matrix}0\\
1\\
\end{matrix}\hspace{-0.12in}\begin{matrix} &\hfill\tikzmark{g12}\\&\hfill\tikzmark{h12}\end{matrix}\,\,\,\begin{matrix}1&\eta_{12}^1&\mathfrak h^1_{23}+\mathfrak x_{13}^1\left(\eta_{12}^1-\mathfrak h_{22}^1\right)&\cdots\\
1&\mathfrak h^1_{22}&\mathfrak h^1_{23}&\cdots\\
\end{matrix}\hspace{-0.12in}\begin{matrix} &\hfill\tikzmark{e12}\\&\hfill\tikzmark{f12}\end{matrix}\,\,\,\begin{matrix}
  1\\1+\mathfrak x_{11}^2(\eta_{12}^1-\mathfrak h^1_{22})
\end{matrix}\right.\\
&\left.\,\,\,\,\,\,\,\,\,\begin{matrix}\mathfrak h_{22}^2+(\mathfrak x_{12}^2+\mathfrak x_{11}^2\mathfrak h^2_{22})(\eta^1_{12}-\mathfrak h^1_{22})&\mathfrak h_{23}^2+(\mathfrak x_{12}^2\mathfrak x_{13}^2+\mathfrak x_{11}^2\mathfrak x_{13}^2\mathfrak h^2_{23})(\eta^1_{12}-\mathfrak h^1_{22})&\cdots\\
(1+\mathfrak x_{11}^2(\eta_{12}^1-\mathfrak h^1_{22}))\mathfrak h_{22}^2&(1+\mathfrak x_{11}^2(\eta_{12}^1-\mathfrak h^1_{22}))\mathfrak h_{23}^2&\cdots\\
\end{matrix}\right)=:\Theta.
\tikz[remember picture,overlay]   \draw[dashed,dash pattern={on 4pt off 2pt}] ([xshift=0.5\tabcolsep,yshift=7pt]e12.north) -- ([xshift=0.5\tabcolsep,yshift=-2pt]f12.south);\tikz[remember picture,overlay]   \draw[dashed,dash pattern={on 4pt off 2pt}] ([xshift=0.5\tabcolsep,yshift=7pt]g12.north) -- ([xshift=0.5\tabcolsep,yshift=-2pt]h12.south); \tikz[remember picture,overlay]   \draw[dashed,dash pattern={on 4pt off 2pt}] ([xshift=0.5\tabcolsep,yshift=7pt]a12.north) -- ([xshift=0.5\tabcolsep,yshift=-2pt]b12.south);
\end{split}    
\end{equation}

\item[($Y_{2B}$)] $Y_{2B}:={\rm Spec}\,\mathbb Z\left[\eta^1_{12},\mathfrak x^1_{13},\cdots,\mathfrak x^1_{1s_3},\mathfrak h^1_{22},\cdots,\mathfrak h^1_{2s_3},\mathfrak x^2_{11},\cdots,\mathfrak x^2_{1s_4},\mathfrak h^2_{22},\cdots,\mathfrak h^2_{2s_4},\frac{1}{1+\mathfrak x_{11}^2\left(\eta_{12}^1-\mathfrak h^1_{22}\right)},\frac{1}{1+\mathfrak x_{11}^2\eta^1_{12}}\right]$. The associated morphism $L:Y_{2B}\longrightarrow G(2,n)$ is defined  by (\ref{y3a}).
\smallskip

\item[($Y_3$)] $Y_{3}:={\rm Spec}\,\mathbb Z\left[\mathfrak x^1_{12},\cdots,\mathfrak x^1_{1s_3},\mathfrak h^1_{22},\cdots,\mathfrak h^1_{2s_3},\mathfrak x^2_{11},\eta^2_{12},\mathfrak x^2_{13},\mathfrak x^2_{14},\cdots,\mathfrak x^2_{1s_4},\mathfrak h^2_{22},\cdots,\mathfrak h^2_{2s_4},\frac{1}{1-\mathfrak x^2_{11}\eta^2_{12}},\right.$
$\left.\frac{1}{1-\mathfrak x^2_{11}\mathfrak h^2_{22}}\right]$. The associated morphism $L:Y_{3}\longrightarrow G(2,n)$ is defined  by
\begin{equation}\label{y3}
\begin{split}
&\left(\begin{matrix}1\\
0\\
\end{matrix}\hspace{-0.12in}\begin{matrix} &\hfill\tikzmark{a12}\\
\\&\hfill\tikzmark{b12}
\end{matrix}\,\,\,\begin{matrix}0\\
1\\
\end{matrix}\hspace{-0.12in}\begin{matrix} &\hfill\tikzmark{g12}\\
\\&\hfill\tikzmark{h12}\end{matrix}\,\,\,\begin{matrix}1&\mathfrak h^1_{22}+\mathfrak x_{12}^1\frac{\eta_{12}^2-\mathfrak h^2_{22}}{1-\mathfrak x^2_{11}\mathfrak h^2_{22}}&\mathfrak h^1_{23}+\mathfrak x_{12}^1\mathfrak x_{13}^1\frac{\eta_{12}^2-\mathfrak h^2_{22}}{1-\mathfrak x^2_{11}\mathfrak h^2_{22}}&\cdots\\
1&\mathfrak h^1_{22}&\mathfrak h^1_{23}&\cdots\\
\end{matrix}\hspace{-0.12in}\begin{matrix} &\hfill\tikzmark{e12}\\
\\&\hfill\tikzmark{f12}\end{matrix}\right.\\
&\left.\,\,\,\,\,\,\,\,\,\,\,\,\,\,\,\,\,\,\,\,\,\,\,\,\begin{matrix} &\hfill\tikzmark{c12}\\
\\&\hfill\tikzmark{d12}\end{matrix}\,\,\,\begin{matrix}1&\eta_{12}^2&\frac{\eta^2_{12}-\mathfrak h^2_{22}}{1-\mathfrak x_{11}^2\mathfrak h^2_{22}}\mathfrak x_{13}^2+\frac{1-\mathfrak x_{11}^2\eta^2_{12}}{1-\mathfrak x_{11}^2\mathfrak h^2_{22}}\mathfrak h^2_{23}&\cdots\\
\frac{1-\mathfrak x_{11}^2\eta^2_{12}}{1-\mathfrak x_{11}^2\mathfrak h^2_{22}}&\frac{1-\mathfrak x_{11}^2\eta^2_{12}}{1-\mathfrak x_{11}^2\mathfrak h^2_{22}}\mathfrak h_{22}^2&\frac{1-\mathfrak x_{11}^2\eta^2_{12}}{1-\mathfrak x_{11}^2\mathfrak h^2_{22}}\mathfrak h_{23}^2&\cdots\\
\end{matrix}\right)=:\Theta.
\tikz[remember picture,overlay]   \draw[dashed,dash pattern={on 4pt off 2pt}] ([xshift=0.5\tabcolsep,yshift=7pt]a12.north) -- ([xshift=0.5\tabcolsep,yshift=-2pt]b12.south);\tikz[remember picture,overlay]   \draw[dashed,dash pattern={on 4pt off 2pt}] ([xshift=0.5\tabcolsep,yshift=7pt]c12.north) -- ([xshift=0.5\tabcolsep,yshift=-2pt]d12.south);\tikz[remember picture,overlay]   \draw[dashed,dash pattern={on 4pt off 2pt}] ([xshift=0.5\tabcolsep,yshift=7pt]e12.north) -- ([xshift=0.5\tabcolsep,yshift=-2pt]f12.south);\tikz[remember picture,overlay]   \draw[dashed,dash pattern={on 4pt off 2pt}] ([xshift=0.5\tabcolsep,yshift=7pt]g12.north) -- ([xshift=0.5\tabcolsep,yshift=-2pt]h12.south); 
\end{split}    
\end{equation}

\item[($Y_4$)] $Y_{4}:={\rm Spec}\,\mathbb Z\left[\epsilon^+_2,\eta^1_{12},\cdots,\eta^1_{1s_3},\mathfrak h^1_{21},\mathfrak h^1_{23},\mathfrak h^1_{24},\cdots,\mathfrak h^1_{2s_3},\eta^2_{12},\cdots,\eta^2_{1s_4},\mathfrak h^2_{22},\cdots,\mathfrak h^2_{2s_4},\frac{1}{1-\eta_{12}^1\epsilon^+_2},\right.$ $\left.\frac{1}{1-\eta_{12}^1\mathfrak h_{21}^1}\right]$. The associated morphism $L:Y_{4}\longrightarrow G(2,n)$ is defined by
\begin{equation}\label{y4}
\begin{split}
&\left(\begin{matrix}
1\\
0\\
\end{matrix}\hspace{-0.12in}\begin{matrix} &\hfill\tikzmark{a2}\\
&\hfill\tikzmark{b2}\\
\end{matrix}\,\,\,\begin{matrix}
0\\
1\\
\end{matrix}\hspace{-0.12in}\begin{matrix} &\hfill\tikzmark{c2}\\
&\hfill\tikzmark{d2}\\
\end{matrix}\,\,\,\begin{matrix}
1&\eta_{12}^1&\eta_{13}^1&\cdots\\
\mathfrak h^1_{21}&1&\mathfrak h^1_{23}&\cdots\\
\end{matrix}\hspace{-0.12in}\begin{matrix} &\hfill\tikzmark{e2}\\
&\hfill\tikzmark{f2}\\
\end{matrix}\,\,\,\begin{matrix}
1&\eta_{12}^2&\eta_{13}^2&\cdots\\
\epsilon_2^+&\epsilon_2^+\left(\mathfrak{h}^2_{22}+\eta^2_{12}\right)&\epsilon_2^+\left(\mathfrak{h}^2_{23}\mathfrak h^2_{22}+\eta^2_{13}\right)&\cdots\\
\end{matrix}\right)=:\Theta.
\tikz[remember picture,overlay]   
\draw[dashed,dash pattern={on 4pt off 2pt}] ([xshift=0.5\tabcolsep,yshift=7pt]a2.north) -- ([xshift=0.5\tabcolsep,yshift=-2pt]b2.south);
\tikz[remember picture,overlay]   
\draw[dashed,dash pattern={on 4pt off 2pt}] ([xshift=0.5\tabcolsep,yshift=7pt]c2.north) -- ([xshift=0.5\tabcolsep,yshift=-2pt]d2.south);
\tikz[remember picture,overlay]   
\draw[dashed,dash pattern={on 4pt off 2pt}] ([xshift=0.5\tabcolsep,yshift=7pt]e2.north) -- ([xshift=0.5\tabcolsep,yshift=-2pt]f2.south);    
\end{split}  
\end{equation}

\item[($Y_5$)] $Y_{5}:={\rm Spec}\,\mathbb Z\left[\epsilon^+_2,\eta^1_{12},\cdots,\eta^1_{1s_3},\mathfrak h^1_{22},\cdots,\mathfrak h^1_{2s_3},\eta^2_{12},\cdots,\eta^2_{1s_4},\mathfrak h^2_{21},\mathfrak h^2_{23},\mathfrak h^2_{24},\cdots,\mathfrak h^2_{2s_4},\frac{1}{1-\epsilon^+_2\mathfrak h_{21}^2},\right.$ $\left.\frac{1}{1-\eta_{12}^2\mathfrak h_{21}^2}\right]$. The associated morphism $L:Y_{5}\longrightarrow G(2,n)$ is defined by
\begin{equation}\label{y5}
\begin{split}
&\left(\begin{matrix}
1\\
0\\
\end{matrix}\hspace{-0.12in}\begin{matrix} &\hfill\tikzmark{a2}\\
&\hfill\tikzmark{b2}\\
\end{matrix}\,\,\,\begin{matrix}
0\\
1\\
\end{matrix}\hspace{-0.12in}\begin{matrix} &\hfill\tikzmark{c2}\\
&\hfill\tikzmark{d2}\\
\end{matrix}\,\,\,\begin{matrix}
1&\eta_{12}^1&\eta_{13}^1&\cdots\\
1&\mathfrak{h}^1_{22}+\eta^1_{12}&\mathfrak{h}^1_{23}\mathfrak h^1_{22}+\eta^1_{13}&\cdots\\
\end{matrix}\hspace{-0.12in}\begin{matrix} &\hfill\tikzmark{e2}\\
&\hfill\tikzmark{f2}\\
\end{matrix}\,\,\,\begin{matrix}
1&\eta_{12}^2&\eta_{13}^2&\cdots\\
\epsilon_2^+\mathfrak{h}^2_{21}&\epsilon_2^+&\epsilon_2^+\mathfrak{h}^2_{23}&\cdots\\
\end{matrix}\right)=:\Theta.
\tikz[remember picture,overlay]   
\draw[dashed,dash pattern={on 4pt off 2pt}] ([xshift=0.5\tabcolsep,yshift=7pt]a2.north) -- ([xshift=0.5\tabcolsep,yshift=-2pt]b2.south);
\tikz[remember picture,overlay]   
\draw[dashed,dash pattern={on 4pt off 2pt}] ([xshift=0.5\tabcolsep,yshift=7pt]c2.north) -- ([xshift=0.5\tabcolsep,yshift=-2pt]d2.south);
\tikz[remember picture,overlay]   
\draw[dashed,dash pattern={on 4pt off 2pt}] ([xshift=0.5\tabcolsep,yshift=7pt]e2.north) -- ([xshift=0.5\tabcolsep,yshift=-2pt]f2.south);    
\end{split}  
\end{equation}

\item[($Y_6$)] $Y_{6}:={\rm Spec}\,\mathbb Z\Bigl[\epsilon^+_2,\eta^1_{12},\cdots,\eta^1_{1s_3},\mathfrak h^1_{21},\mathfrak h^1_{23},\mathfrak h^1_{24},\cdots,\mathfrak h^1_{2s_3},\eta^2_{12},\cdots,\eta^2_{1s_4},\mathfrak h^2_{21},\mathfrak h^2_{23},\mathfrak h^2_{24},\cdots,\mathfrak h^2_{2s_4},$ $\left.\frac{1}{1-\eta_{12}^1\mathfrak h_{21}^1},\frac{1}{1-\eta_{12}^2\mathfrak h_{21}^2},\frac{1}{1-\eta_{12}^1\epsilon^+_2\mathfrak h_{21}^2}\right]$. The associated morphism $L:Y_{6}\longrightarrow G(2,n)$ is defined by
\begin{equation}\label{y6}
\begin{split}
&\left(\begin{matrix}
1\\
0\\
\end{matrix}\hspace{-0.12in}\begin{matrix} &\hfill\tikzmark{a2}\\
&\hfill\tikzmark{b2}\\
\end{matrix}\,\,\,\begin{matrix}
0\\
1\\
\end{matrix}\hspace{-0.12in}\begin{matrix} &\hfill\tikzmark{c2}\\
&\hfill\tikzmark{d2}\\
\end{matrix}\,\,\,\begin{matrix}
1&\eta_{12}^1&\eta_{13}^1&\cdots\\
\mathfrak h^1_{21}&1&\mathfrak{h}^1_{23}&\cdots\\
\end{matrix}\hspace{-0.12in}\begin{matrix} &\hfill\tikzmark{e2}\\
&\hfill\tikzmark{f2}\\
\end{matrix}\,\,\,\begin{matrix}
1&\eta_{12}^2&\eta_{13}^2&\cdots\\
\epsilon_2^+\mathfrak{h}^2_{21}&\epsilon_2^+&\epsilon_2^+\mathfrak{h}^2_{23}&\cdots\\
\end{matrix}\right)=:\Theta.
\tikz[remember picture,overlay]   
\draw[dashed,dash pattern={on 4pt off 2pt}] ([xshift=0.5\tabcolsep,yshift=7pt]a2.north) -- ([xshift=0.5\tabcolsep,yshift=-2pt]b2.south);
\tikz[remember picture,overlay]   
\draw[dashed,dash pattern={on 4pt off 2pt}] ([xshift=0.5\tabcolsep,yshift=7pt]c2.north) -- ([xshift=0.5\tabcolsep,yshift=-2pt]d2.south);
\tikz[remember picture,overlay]   
\draw[dashed,dash pattern={on 4pt off 2pt}] ([xshift=0.5\tabcolsep,yshift=7pt]e2.north) -- ([xshift=0.5\tabcolsep,yshift=-2pt]f2.south);    
\end{split}  
\end{equation}

\end{enumerate}

\subsection{{\texorpdfstring{$s_1=1$   and $s_2,s_3, s_4\geq2$}{dd}}}\label{z1} In the following, we shall describe certain affine schemes 
and the associated morphisms such that the images give an open cover of $\mathcal M^{\underline s}_n$ up to permutations. Notice that the parametrizations here can be derived as follows. We first embed $\mathcal M^{\underline s}_n$ into $\mathcal T^{\underline s}_n$, and then transform from $A^{\tau}$ with a Class II index to $A^{\tau^{\prime}}$ with a Class I index as in Lemma \ref{32ud}. 

\begin{enumerate}
\item[($Z_{1A}$)] $Z_{1A}:={\rm Spec}\,\mathbb Z\Bigl[\epsilon^+_2,\mathfrak x^1_{12},\cdots,\mathfrak x^1_{1s_3},\mathfrak h^1_{22},\cdots,\mathfrak h^1_{2s_3},\mathfrak x^2_{12},\cdots,\mathfrak x^2_{1s_4},\mathfrak h^2_{22},\cdots,\mathfrak h^2_{2s_4},\epsilon^-_1,z_2,\cdots,z_{s_2},\mathfrak v_3,$  $\left.\cdots,\mathfrak v_{s_2},\frac{1}{1-\epsilon^+_2}\right]$. The associated morphism $L:Z_{1A}\longrightarrow G(2,n)$ is defined by

\begin{equation}\label{didi2}
\begin{split}
&\left(\begin{matrix}1\\
0\\
\end{matrix}\hspace{-0.12in}\begin{matrix} &\hfill\tikzmark{a12}\\&\hfill\tikzmark{b12}
\end{matrix}\,\,\,\begin{matrix}0&\epsilon^-_1&\epsilon^-_1\mathfrak v_3&\cdots\\
1&z_2&z_3&\cdots\\
\end{matrix}\hspace{-0.12in}\begin{matrix} &\hfill\tikzmark{g12}\\&\hfill\tikzmark{h12}\end{matrix}\,\,\,\begin{matrix}1&\mathfrak h^1_{22}+\mathfrak x_{12}^1(1-\epsilon_2^+)&\mathfrak h^1_{23}+\mathfrak x_{12}^1\mathfrak x_{13}^1(1-\epsilon_2^+)&\cdots\\
1&\mathfrak h^1_{22}&\mathfrak h^1_{23}&\cdots\\
\end{matrix}\hspace{-0.12in}\begin{matrix} &\hfill\tikzmark{e12}\\&\hfill\tikzmark{f12}\end{matrix}\right.\\
&\left.\,\,\,\,\,\,\,\,\,\,\,\,\,\,\,\,\,\,\,\,\,\,\,\,\,\begin{matrix} &\hfill\tikzmark{c12}\\&\hfill\tikzmark{d12}\end{matrix}\,\,\,\begin{matrix}1&\mathfrak x_{12}^2(1-\epsilon^+_2)+\mathfrak h^2_{22}&\mathfrak x_{12}^2\mathfrak x_{13}^2(1-\epsilon^+_2)+\mathfrak h^2_{23}&\cdots\\
\epsilon_2^+&\epsilon_2^+\mathfrak h_{22}^2&\epsilon_2^+\mathfrak h_{23}^2&\cdots\\
\end{matrix}\right)=:\Theta.
\tikz[remember picture,overlay]   \draw[dashed,dash pattern={on 4pt off 2pt}] ([xshift=0.5\tabcolsep,yshift=7pt]a12.north) -- ([xshift=0.5\tabcolsep,yshift=-2pt]b12.south);\tikz[remember picture,overlay]   \draw[dashed,dash pattern={on 4pt off 2pt}] ([xshift=0.5\tabcolsep,yshift=7pt]c12.north) -- ([xshift=0.5\tabcolsep,yshift=-2pt]d12.south);\tikz[remember picture,overlay]   \draw[dashed,dash pattern={on 4pt off 2pt}] ([xshift=0.5\tabcolsep,yshift=7pt]e12.north) -- ([xshift=0.5\tabcolsep,yshift=-2pt]f12.south);\tikz[remember picture,overlay]   \draw[dashed,dash pattern={on 4pt off 2pt}] ([xshift=0.5\tabcolsep,yshift=7pt]g12.north) -- ([xshift=0.5\tabcolsep,yshift=-2pt]h12.south); 
\end{split}    
\end{equation}

\item[($Z_{1B}$)] $Z_{1B}:={\rm Spec}\,\mathbb Z\Bigl[\epsilon^+_2,\mathfrak x^1_{12},\cdots,\mathfrak x^1_{1s_3},\mathfrak h^1_{22},\cdots,\mathfrak h^1_{2s_3},\mathfrak x^2_{12},\cdots,\mathfrak x^2_{1s_4},\mathfrak h^2_{22},\cdots,\mathfrak h^2_{2s_4},\epsilon^-_1,z_2,\cdots,z_{s_2},\mathfrak v_3,$  $\left.\cdots,\mathfrak v_{s_2},\frac{1}{\epsilon^+_2}\right]$. The associated morphism $L:Z_{1B}\longrightarrow G(2,n)$ is defined by (\ref{didi2}).

\item[($Z_{2A}$)] $Z_{2A}:={\rm Spec}\,\mathbb Z\Bigl[\eta^1_{12},\mathfrak x^1_{13},\cdots,\mathfrak x^1_{1s_3},\mathfrak h^1_{22},\cdots,\mathfrak h^1_{2s_3},\mathfrak x^2_{11},\cdots,\mathfrak x^2_{1s_4},\mathfrak h^2_{22},\cdots,\mathfrak h^2_{2s_4},\epsilon^-_1,z_2,\cdots,z_{s_2},\mathfrak v_3,$ 
$\left.\cdots,\mathfrak v_{s_2},\frac{1}{1+\mathfrak x_{11}^2\left(\eta_{12}^1-\mathfrak h^1_{22}\right)},\frac{1}{\mathfrak x^2_{11}}\right]$. The associated morphism $L:Z_{2A}\longrightarrow G(2,n)$ is defined by
\begin{equation}\label{z3a}
\begin{split}
&\left(\begin{matrix}1\\
0\\
\end{matrix}\hspace{-0.12in}\begin{matrix} &\hfill\tikzmark{a12}\\&\hfill\tikzmark{b12}
\end{matrix}\,\,\,\begin{matrix}0&\epsilon^-_1&\epsilon^-_1\mathfrak v_3&\cdots\\
1&z_2&z_3&\cdots\\
\end{matrix}\hspace{-0.12in}\begin{matrix} &\hfill\tikzmark{g12}\\&\hfill\tikzmark{h12}\end{matrix}\,\,\,\begin{matrix}1&\eta_{12}^1&\mathfrak h^1_{23}+\mathfrak x_{13}^1\left(\eta_{12}^1-\mathfrak h_{22}^1\right)&\cdots\\
1&\mathfrak h^1_{22}&\mathfrak h^1_{23}&\cdots\\
\end{matrix}\hspace{-0.12in}\begin{matrix} &\hfill\tikzmark{e12}\\&\hfill\tikzmark{f12}\end{matrix}\,\,\,\begin{matrix}1\\1+\mathfrak x_{11}^2(\eta_{12}^1-\mathfrak h^1_{22})
\end{matrix}\right.\\
&\left.\,\,\,\,\,\,\begin{matrix}\mathfrak h_{22}^2+(\mathfrak x_{12}^2+\mathfrak x_{11}^2\mathfrak h^2_{22})(\eta^1_{12}-\mathfrak h^1_{22})&\mathfrak h_{23}^2+(\mathfrak x_{12}^2\mathfrak x_{13}^2+\mathfrak x_{11}^2\mathfrak x_{13}^2\mathfrak h^2_{23})(\eta^1_{12}-\mathfrak h^1_{22})&\cdots\\
(1+\mathfrak x_{11}^2(\eta_{12}^1-\mathfrak h^1_{22}))\mathfrak h_{22}^2&(1+\mathfrak x_{11}^2(\eta_{12}^1-\mathfrak h^1_{22}))\mathfrak h_{23}^2&\cdots\\
\end{matrix}\right)=:\Theta.
\tikz[remember picture,overlay]   \draw[dashed,dash pattern={on 4pt off 2pt}] ([xshift=0.5\tabcolsep,yshift=7pt]a12.north) -- ([xshift=0.5\tabcolsep,yshift=-2pt]b12.south);\tikz[remember picture,overlay]   \draw[dashed,dash pattern={on 4pt off 2pt}] ([xshift=0.5\tabcolsep,yshift=7pt]c12.north) -- ([xshift=0.5\tabcolsep,yshift=-2pt]d12.south);\tikz[remember picture,overlay]   \draw[dashed,dash pattern={on 4pt off 2pt}] ([xshift=0.5\tabcolsep,yshift=7pt]e12.north) -- ([xshift=0.5\tabcolsep,yshift=-2pt]f12.south);\tikz[remember picture,overlay]   \draw[dashed,dash pattern={on 4pt off 2pt}] ([xshift=0.5\tabcolsep,yshift=7pt]g12.north) -- ([xshift=0.5\tabcolsep,yshift=-2pt]h12.south); 
\end{split}    
\end{equation}

\item[($Z_{2B}$)] $Z_{2B}:={\rm Spec}\,\mathbb Z\Bigl[\eta^1_{12},\mathfrak x^1_{13},\cdots,\mathfrak x^1_{1s_3},\mathfrak h^1_{22},\cdots,\mathfrak h^1_{2s_3},\mathfrak x^2_{11},\cdots,\mathfrak x^2_{1s_4},\mathfrak h^2_{22},\cdots,\mathfrak h^2_{2s_4},\epsilon^-_1,z_2,\cdots,z_{s_2},$  $\left.\mathfrak v_3,\cdots,\mathfrak v_{s_2},\frac{1}{1+\mathfrak x_{11}^2\left(\eta_{12}^1-\mathfrak h^1_{22}\right)},\frac{1}{1+\mathfrak x_{11}^2\eta^1_{12}}\right]$. The associated morphism $L:Z_{2B}\longrightarrow G(2,n)$ is defined  by (\ref{z3a}).
\smallskip

\item[($Z_3$)] 
$Z_{3}:={\rm Spec}\,\mathbb Z\Bigl[\mathfrak x^1_{12},\cdots,\mathfrak x^1_{1s_3},\mathfrak h^1_{22},\cdots,\mathfrak h^1_{2s_3},\mathfrak x^2_{11},\eta^2_{12},\mathfrak x^2_{13},\mathfrak x^2_{14},\cdots,\mathfrak x^2_{1s_4},\mathfrak h^2_{22},\cdots,\mathfrak h^2_{2s_4},\epsilon^-_1,z_2,\cdots,$  $\left.z_{s_2},\mathfrak v_3,\cdots,\mathfrak v_{s_2},\frac{1}{1-\mathfrak x^2_{11}\eta^2_{12}},\frac{1}{1-\mathfrak x^2_{11}\mathfrak h^2_{22}}\right]$. The associated morphism $L:Z_{3}\longrightarrow G(2,n)$ is defined  by

\begin{equation}\label{z3}
\begin{split}
&\left(\begin{matrix}1\\
0\\
\end{matrix}\hspace{-0.12in}\begin{matrix} &\hfill\tikzmark{a12}\\
\\&\hfill\tikzmark{b12}
\end{matrix}\,\,\,\begin{matrix}0&\epsilon^-_1&\epsilon^-_1\mathfrak v_3&\cdots\\
1&z_2&z_3&\cdots\\
\end{matrix}\hspace{-0.12in}\begin{matrix} &\hfill\tikzmark{g12}\\
\\&\hfill\tikzmark{h12}\end{matrix}\,\,\,\begin{matrix}1&\mathfrak h^1_{22}+\mathfrak x_{12}^1\frac{\eta_{12}^2-\mathfrak h^2_{22}}{1-\mathfrak x^2_{11}\mathfrak h^2_{12}}&\mathfrak h^1_{23}+\mathfrak x_{12}^1\mathfrak x_{13}^1\frac{\eta_{12}^2-\mathfrak h^2_{22}}{1-\mathfrak x^2_{11}\mathfrak h^2_{12}}&\cdots\\
1&\mathfrak h^1_{22}&\mathfrak h^1_{23}&\cdots\\
\end{matrix}\hspace{-0.12in}\begin{matrix} &\hfill\tikzmark{e12}\\
\\&\hfill\tikzmark{f12}\end{matrix}\right.\\
&\left.\,\,\,\,\,\,\,\,\,\,\,\,\,\,\,\,\,\,\,\,\,\,\,\,\,\,\,\,\,\,\begin{matrix} &\hfill\tikzmark{c12}\\
\\&\hfill\tikzmark{d12}\end{matrix}\,\,\,\begin{matrix}1&\eta_{12}^2&\frac{\eta^2_{12}-\mathfrak h^2_{22}}{1-\mathfrak x_{11}^2\mathfrak h^2_{22}}\mathfrak x_{13}^2+\frac{1-\mathfrak x_{11}^2\eta^2_{12}}{1-\mathfrak x_{11}^2\mathfrak h^2_{22}}\mathfrak h^2_{23}&\cdots\\
\frac{1-\mathfrak x_{11}^2\eta^2_{12}}{1-\mathfrak x_{11}^2\mathfrak h^2_{22}}&\frac{1-\mathfrak x_{11}^2\eta^2_{12}}{1-\mathfrak x_{11}^2\mathfrak h^2_{22}}\mathfrak h_{22}^2&\frac{1-\mathfrak x_{11}^2\eta^2_{12}}{1-\mathfrak x_{11}^2\mathfrak h^2_{22}}\mathfrak h_{23}^2&\cdots\\
\end{matrix}\right)=:\Theta.
\tikz[remember picture,overlay]   \draw[dashed,dash pattern={on 4pt off 2pt}] ([xshift=0.5\tabcolsep,yshift=7pt]a12.north) -- ([xshift=0.5\tabcolsep,yshift=-2pt]b12.south);\tikz[remember picture,overlay]   \draw[dashed,dash pattern={on 4pt off 2pt}] ([xshift=0.5\tabcolsep,yshift=7pt]c12.north) -- ([xshift=0.5\tabcolsep,yshift=-2pt]d12.south);\tikz[remember picture,overlay]   \draw[dashed,dash pattern={on 4pt off 2pt}] ([xshift=0.5\tabcolsep,yshift=7pt]e12.north) -- ([xshift=0.5\tabcolsep,yshift=-2pt]f12.south);\tikz[remember picture,overlay]   \draw[dashed,dash pattern={on 4pt off 2pt}] ([xshift=0.5\tabcolsep,yshift=7pt]g12.north) -- ([xshift=0.5\tabcolsep,yshift=-2pt]h12.south); 
\end{split}    
\end{equation}

\item[($Z_4$)] 
$Z_{4}:={\rm Spec}\,\mathbb Z\Bigl[\epsilon^+_2,\eta^1_{12},\cdots,\eta^1_{1s_3},\mathfrak h^1_{21},\mathfrak h^1_{23},\mathfrak h^1_{24},\cdots,\mathfrak h^1_{2s_3},\eta^2_{12},\cdots,\eta^2_{1s_4},\mathfrak h^2_{22},\cdots,\mathfrak h^2_{2s_4},\epsilon^-_1,z_2,\cdots,$  $\left.z_{s_2},\mathfrak v_3,\cdots,\mathfrak v_{s_2},\frac{1}{1-\eta_{12}^1\epsilon^+_2},\frac{1}{1-\eta_{12}^1\mathfrak h_{21}^1}\right]$. The associated morphism $L:Z_{4}\longrightarrow G(2,n)$ is defined by
\begin{equation}\label{z4}
\begin{split}
&\left(\begin{matrix}
1\\
0\\
\end{matrix}\hspace{-0.12in}\begin{matrix} &\hfill\tikzmark{a2}\\
&\hfill\tikzmark{b2}\\
\end{matrix}\,\,\,\begin{matrix}
0&\epsilon^-_1&\epsilon^-_1\mathfrak v_3&\cdots\\
1&z_2&z_3&\cdots\\
\end{matrix}\hspace{-0.12in}\begin{matrix} &\hfill\tikzmark{c2}\\
&\hfill\tikzmark{d2}\\
\end{matrix}\,\,\,\begin{matrix}
1&\eta_{12}^1&\eta_{13}^1&\cdots\\
\mathfrak h^1_{21}&1&\mathfrak h^1_{23}&\cdots\\
\end{matrix}\hspace{-0.12in}\begin{matrix} &\hfill\tikzmark{e2}\\
&\hfill\tikzmark{f2}\\
\end{matrix}\,\,\,\begin{matrix}
1&\eta_{12}^2&\eta_{13}^2&\cdots\\
\epsilon_2^+&\epsilon_2^+\left(\mathfrak{h}^2_{22}+\eta^2_{12}\right)&\epsilon_2^+\left(\mathfrak{h}^2_{23}\mathfrak h^2_{22}+\eta^2_{13}\right)&\cdots\\
\end{matrix}\right)=:\Theta.
\tikz[remember picture,overlay]   
\draw[dashed,dash pattern={on 4pt off 2pt}] ([xshift=0.5\tabcolsep,yshift=7pt]a2.north) -- ([xshift=0.5\tabcolsep,yshift=-2pt]b2.south);
\tikz[remember picture,overlay]   
\draw[dashed,dash pattern={on 4pt off 2pt}] ([xshift=0.5\tabcolsep,yshift=7pt]c2.north) -- ([xshift=0.5\tabcolsep,yshift=-2pt]d2.south);
\tikz[remember picture,overlay]   
\draw[dashed,dash pattern={on 4pt off 2pt}] ([xshift=0.5\tabcolsep,yshift=7pt]e2.north) -- ([xshift=0.5\tabcolsep,yshift=-2pt]f2.south);    
\end{split}  
\end{equation}

\item[($Z_5$)] 

$Z_{5}:={\rm Spec}\,\mathbb Z\Bigl[\epsilon^+_2,\eta^1_{12},\cdots,\eta^1_{1s_3},\mathfrak h^1_{22},\cdots,\mathfrak h^1_{2s_3},\eta^2_{12},\cdots,\eta^2_{1s_4},\mathfrak h^2_{21},\mathfrak h^2_{23},\mathfrak h^2_{24},\cdots,\mathfrak h^2_{2s_4},\epsilon^-_1,z_2,\cdots,$  $\left.z_{s_2},\mathfrak v_3,\cdots,\mathfrak v_{s_2},\frac{1}{1-\epsilon^+_2\mathfrak h_{21}^2},\frac{1}{1-\eta_{12}^2\mathfrak h_{21}^2}\right]$. The associated morphism $L:Z_{5}\longrightarrow G(2,n)$ is defined by
\begin{equation}\label{z5}
\begin{split}
&\left(\begin{matrix}
1\\
0\\
\end{matrix}\hspace{-0.12in}\begin{matrix} &\hfill\tikzmark{a2}\\
&\hfill\tikzmark{b2}\\
\end{matrix}\,\,\,\begin{matrix}
0&\epsilon^-_1&\epsilon^-_1\mathfrak v_3&\cdots\\
1&z_2&z_3&\cdots\\
\end{matrix}\hspace{-0.12in}\begin{matrix} &\hfill\tikzmark{c2}\\
&\hfill\tikzmark{d2}\\
\end{matrix}\,\,\,\begin{matrix}
1&\eta_{12}^1&\eta_{13}^1&\cdots\\
1&\mathfrak{h}^1_{22}+\eta^1_{12}&\mathfrak{h}^1_{23}\mathfrak h^1_{22}+\eta^1_{13}&\cdots\\
\end{matrix}\hspace{-0.12in}\begin{matrix} &\hfill\tikzmark{e2}\\
&\hfill\tikzmark{f2}\\
\end{matrix}\,\,\,\begin{matrix}
1&\eta_{12}^2&\eta_{13}^2&\cdots\\
\epsilon_2^+\mathfrak{h}^2_{21}&\epsilon_2^+&\epsilon_2^+\mathfrak{h}^2_{23}&\cdots\\
\end{matrix}\right)=:\Theta.
\tikz[remember picture,overlay]   
\draw[dashed,dash pattern={on 4pt off 2pt}] ([xshift=0.5\tabcolsep,yshift=7pt]a2.north) -- ([xshift=0.5\tabcolsep,yshift=-2pt]b2.south);
\tikz[remember picture,overlay]   
\draw[dashed,dash pattern={on 4pt off 2pt}] ([xshift=0.5\tabcolsep,yshift=7pt]c2.north) -- ([xshift=0.5\tabcolsep,yshift=-2pt]d2.south);
\tikz[remember picture,overlay]   
\draw[dashed,dash pattern={on 4pt off 2pt}] ([xshift=0.5\tabcolsep,yshift=7pt]e2.north) -- ([xshift=0.5\tabcolsep,yshift=-2pt]f2.south);    
\end{split}  
\end{equation}

\item[($Z_6$)] 
$Z_{6}:={\rm Spec}\,\mathbb Z\Bigl[\epsilon^+_2,\eta^1_{12},\cdots,\eta^1_{1s_3},\mathfrak h^1_{21},\mathfrak h^1_{23},\mathfrak h^1_{24},\cdots,\mathfrak h^1_{2s_3},\eta^2_{12},\cdots,\eta^2_{1s_4},\mathfrak h^2_{21},\mathfrak h^2_{23},\mathfrak h^2_{24},\cdots,\mathfrak h^2_{2s_4},\epsilon^-_1,$  $\left.z_2,\cdots,z_{s_2},\mathfrak v_3,\cdots,\mathfrak v_{s_2},\frac{1}{1-\eta_{12}^1\mathfrak h_{21}^1},\frac{1}{1-\eta_{12}^2\mathfrak h_{21}^2},\frac{1}{1-\eta_{12}^1\epsilon^+_2\mathfrak h_{21}^2}\right]$. The associated morphism $L:Z_{6}\longrightarrow G(2,n)$ is defined by
\begin{equation}\label{z6}
\begin{split}
&\left(\begin{matrix}
1\\
0\\
\end{matrix}\hspace{-0.12in}\begin{matrix} &\hfill\tikzmark{a2}\\
&\hfill\tikzmark{b2}\\
\end{matrix}\,\,\,\begin{matrix}
0&\epsilon^-_1&\epsilon^-_1\mathfrak v_3&\cdots\\
1&z_2&z_3&\cdots\\
\end{matrix}\hspace{-0.12in}\begin{matrix} &\hfill\tikzmark{c2}\\
&\hfill\tikzmark{d2}\\
\end{matrix}\,\,\,\begin{matrix}
1&\eta_{12}^1&\eta_{13}^1&\cdots\\
\mathfrak h^1_{21}&1&\mathfrak{h}^1_{23}&\cdots\\
\end{matrix}\hspace{-0.12in}\begin{matrix} &\hfill\tikzmark{e2}\\
&\hfill\tikzmark{f2}\\
\end{matrix}\,\,\,\begin{matrix}
1&\eta_{12}^2&\eta_{13}^2&\cdots\\
\epsilon_2^+\mathfrak{h}^2_{21}&\epsilon_2^+&\epsilon_2^+\mathfrak{h}^2_{23}&\cdots\\
\end{matrix}\right)=:\Theta.
\tikz[remember picture,overlay]   
\draw[dashed,dash pattern={on 4pt off 2pt}] ([xshift=0.5\tabcolsep,yshift=7pt]a2.north) -- ([xshift=0.5\tabcolsep,yshift=-2pt]b2.south);
\tikz[remember picture,overlay]   
\draw[dashed,dash pattern={on 4pt off 2pt}] ([xshift=0.5\tabcolsep,yshift=7pt]c2.north) -- ([xshift=0.5\tabcolsep,yshift=-2pt]d2.south);
\tikz[remember picture,overlay]   
\draw[dashed,dash pattern={on 4pt off 2pt}] ([xshift=0.5\tabcolsep,yshift=7pt]e2.north) -- ([xshift=0.5\tabcolsep,yshift=-2pt]f2.south);    
\end{split}  
\end{equation}

\item[($Z_7$)] 
$Z_{7}:={\rm Spec}\,\mathbb Z\Bigl[\epsilon^+_1,\eta^1_{12},\cdots,\eta^1_{1s_4},\mathfrak h^1_{21},\mathfrak h^1_{23},\mathfrak h^1_{24},\cdots,\mathfrak h^1_{2s_4},\mathfrak x^1_{12},\cdots,\mathfrak x^1_{1s_3},\xi^1_{22},\cdots,\xi^1_{2s_3},\epsilon^-_2,z_2,\cdots,$  $\left.z_{s_2},\mathfrak v_3,\cdots,\mathfrak v_{s_2},\frac{1}{1-\eta_{12}^1\mathfrak h_{21}^1},\frac{1}{1-\epsilon^+_1\mathfrak h_{21}^1}\right]$. The associated morphism $L:Z_{7}\longrightarrow G(2,n)$ is defined by
\begin{equation}\label{z7}
\begin{split}
&\left(\begin{matrix}
1\\
0\\
\end{matrix}\hspace{-0.12in}\begin{matrix} &\hfill\tikzmark{a2}\\
&\hfill\tikzmark{b2}\\
\end{matrix}\,\,\,\begin{matrix}
0&\epsilon^-_2&\epsilon^-_2\mathfrak v_3&\cdots\\
1&z_2&z_3&\cdots\\
\end{matrix}\hspace{-0.12in}\begin{matrix} &\hfill\tikzmark{c2}\\
&\hfill\tikzmark{d2}\\
\end{matrix}\,\,\,\begin{matrix}
1&\mathfrak x_{12}^1+\xi_{22}^1&\mathfrak x_{12}^1\mathfrak x_{13}^1+\xi_{23}^1&\cdots\\
1&\xi^1_{22}&\xi^1_{23}&\cdots\\
\end{matrix}\hspace{-0.12in}\begin{matrix} &\hfill\tikzmark{e2}\\
&\hfill\tikzmark{f2}\\
\end{matrix}\,\,\,\begin{matrix}
1&\eta_{12}^1&\eta_{13}^1&\cdots\\
\epsilon_1^+\mathfrak{h}^1_{21}&\epsilon_1^+&\epsilon_1^+\mathfrak{h}^1_{23}&\cdots\\
\end{matrix}\right)=:\Theta.
\tikz[remember picture,overlay]   
\draw[dashed,dash pattern={on 4pt off 2pt}] ([xshift=0.5\tabcolsep,yshift=7pt]a2.north) -- ([xshift=0.5\tabcolsep,yshift=-2pt]b2.south);
\tikz[remember picture,overlay]   
\draw[dashed,dash pattern={on 4pt off 2pt}] ([xshift=0.5\tabcolsep,yshift=7pt]c2.north) -- ([xshift=0.5\tabcolsep,yshift=-2pt]d2.south);
\tikz[remember picture,overlay]   
\draw[dashed,dash pattern={on 4pt off 2pt}] ([xshift=0.5\tabcolsep,yshift=7pt]e2.north) -- ([xshift=0.5\tabcolsep,yshift=-2pt]f2.south);    
\end{split}  
\end{equation}

\item[($Z_8$)] 
$Z_{8}:={\rm Spec}\,\mathbb Z\Bigl[\epsilon^+_1,\eta^1_{12},\cdots,\eta^1_{1s_4},\mathfrak h^1_{21},\mathfrak h^1_{23},\mathfrak h^1_{24},\cdots,\mathfrak h^1_{2s_4},\mathfrak x^1_{11},\mathfrak x^1_{13},\mathfrak x^1_{14},\cdots,\mathfrak x^1_{1s_3},\xi^1_{22},\cdots,\xi^1_{2s_3},\epsilon^-_2,$  $\left.z_2,\cdots,z_{s_2},\mathfrak v_3,\cdots,\mathfrak v_{s_2},\frac{1}{1-\eta_{12}^1\mathfrak h_{21}^1},\frac{1}{1-\mathfrak x^1_{11}\xi^1_{22}},\frac{1}{1-\mathfrak x^1_{11}\epsilon^+_1\mathfrak h^1_{21}}\right]$. The associated morphism $L:Z_{8}\longrightarrow G(2,n)$ is defined by
\begin{equation}\label{z8}
\begin{split}
&\left(\begin{matrix}
1\\
0\\
\end{matrix}\hspace{-0.12in}\begin{matrix} &\hfill\tikzmark{a2}\\
&\hfill\tikzmark{b2}\\
\end{matrix}\,\,\,\begin{matrix}
0&\epsilon^-_2&\epsilon^-_2\mathfrak v_3&\cdots\\
1&z_2&z_3&\cdots\\
\end{matrix}\hspace{-0.12in}\begin{matrix} &\hfill\tikzmark{c2}\\
&\hfill\tikzmark{d2}\\
\end{matrix}\,\,\,\begin{matrix}
\mathfrak x_{11}^1&1&\mathfrak x_{13}^1&\cdots\\
1&\xi^1_{22}&\xi^1_{23}&\cdots\\
\end{matrix}\hspace{-0.12in}\begin{matrix} &\hfill\tikzmark{e2}\\
&\hfill\tikzmark{f2}\\
\end{matrix}\,\,\,\begin{matrix}
1&\eta_{12}^1&\eta_{13}^1&\cdots\\
\epsilon_1^+\mathfrak{h}^1_{21}&\epsilon_1^+&\epsilon_1^+\mathfrak{h}^1_{23}&\cdots\\
\end{matrix}\right)=:\Theta.
\tikz[remember picture,overlay]   
\draw[dashed,dash pattern={on 4pt off 2pt}] ([xshift=0.5\tabcolsep,yshift=7pt]a2.north) -- ([xshift=0.5\tabcolsep,yshift=-2pt]b2.south);
\tikz[remember picture,overlay]   
\draw[dashed,dash pattern={on 4pt off 2pt}] ([xshift=0.5\tabcolsep,yshift=7pt]c2.north) -- ([xshift=0.5\tabcolsep,yshift=-2pt]d2.south);
\tikz[remember picture,overlay]   
\draw[dashed,dash pattern={on 4pt off 2pt}] ([xshift=0.5\tabcolsep,yshift=7pt]e2.north) -- ([xshift=0.5\tabcolsep,yshift=-2pt]f2.south);    
\end{split}  
\end{equation}

\item[($Z_{9}$)] 
$Z_{9}:={\rm Spec}\,\mathbb Z\Bigl[\epsilon^+_1,\eta^1_{12},\cdots,\eta^1_{1s_4},\mathfrak h^1_{21},\mathfrak h^1_{23},\mathfrak h^1_{24},\cdots,\mathfrak h^1_{2s_4},\mathfrak x^1_{11},\cdots,\mathfrak x^1_{1s_3},\xi^1_{22},\cdots,\xi^1_{2s_3},z_2,\cdots,z_{s_2},$  $\left.\mathfrak v_3,\cdots,\mathfrak v_{s_2},\frac{1}{1-\eta_{12}^1\mathfrak h_{21}^1},\frac{1}{1-\mathfrak x^1_{11}z_2},\frac{1}{1-\mathfrak x^1_{11}\epsilon^+_1\mathfrak h_{21}^1}\right]$. The associated morphism $L:Z_{9}\longrightarrow G(2,n)$ is defined by
\begin{equation}\label{z9}
\begin{split}
&\left(\begin{matrix}
1\\
0\\
\end{matrix}\hspace{-0.12in}\begin{matrix} &\hfill\tikzmark{a2}\\
&\hfill\tikzmark{b2}\\
\end{matrix}\,\,\,\begin{matrix}
0&1&\mathfrak v_3&\cdots\\
1&z_2&z_3&\cdots\\
\end{matrix}\hspace{-0.12in}\begin{matrix} &\hfill\tikzmark{c2}\\
&\hfill\tikzmark{d2}\\
\end{matrix}\,\,\,\begin{matrix}
\mathfrak x_{11}^1&\mathfrak x_{12}^1+\mathfrak x_{11}^1\xi_{22}^1&\mathfrak x_{12}^1\mathfrak x_{13}^1+\mathfrak x_{11}^1\xi_{23}^1&\cdots\\
1&\xi^1_{22}&\xi^1_{23}&\cdots\\
\end{matrix}\hspace{-0.12in}\begin{matrix} &\hfill\tikzmark{e2}\\
&\hfill\tikzmark{f2}\\
\end{matrix}\,\,\,\begin{matrix}
1&\eta_{12}^1&\eta_{13}^1&\cdots\\
\epsilon_1^+\mathfrak{h}^1_{21}&\epsilon_1^+&\epsilon_1^+\mathfrak{h}^1_{23}&\cdots\\
\end{matrix}\right)=:\Theta.
\tikz[remember picture,overlay]   
\draw[dashed,dash pattern={on 4pt off 2pt}] ([xshift=0.5\tabcolsep,yshift=7pt]a2.north) -- ([xshift=0.5\tabcolsep,yshift=-2pt]b2.south);
\tikz[remember picture,overlay]   
\draw[dashed,dash pattern={on 4pt off 2pt}] ([xshift=0.5\tabcolsep,yshift=7pt]c2.north) -- ([xshift=0.5\tabcolsep,yshift=-2pt]d2.south);
\tikz[remember picture,overlay]   
\draw[dashed,dash pattern={on 4pt off 2pt}] ([xshift=0.5\tabcolsep,yshift=7pt]e2.north) -- ([xshift=0.5\tabcolsep,yshift=-2pt]f2.south);    
\end{split}  
\end{equation}
\end{enumerate}

\subsection{{\texorpdfstring{$s_1,s_2,s_3,s_4\geq2$}{dd}}}\label{w1}
Similarly to Appendix \ref{z1}, we shall describe certain affine schemes 
and the associated morphisms such that the images give an open cover of $\mathcal M^{\underline s}_n$ up to permutations.
\begin{enumerate}
\item[($W_{1A}$)] $W_{1A}:={\rm Spec}\,\mathbb Z\Bigl[\epsilon^+_2,\mathfrak x^1_{12},\cdots,\mathfrak x^1_{1s_3},\mathfrak h^1_{22},\cdots,\mathfrak h^1_{2s_3},\mathfrak x^2_{12},\cdots,\mathfrak x^2_{1s_4},\mathfrak h^2_{22},\cdots,\mathfrak h^2_{2s_4},\epsilon^-_1,z_2,\cdots,z_{s_2},\mathfrak v_3,$  $\left.\cdots,\mathfrak v_{s_2},\epsilon_3^+,y_2,\cdots,y_{s_1},\mathfrak u_3,\cdots,\mathfrak u_{s_1},\frac{1}{1-\epsilon^+_2}\right]$. The associated morphism $L:W_{1A}\longrightarrow G(2,n)$ is defined by

\begin{equation}\label{didi3}
\begin{split}
&\left(\begin{matrix}1&y_2&y_3&\cdots\\
0&\epsilon_2^+\epsilon^+_3&\epsilon_2^+\epsilon^+_3\mathfrak u_3&\cdots\\
\end{matrix}\hspace{-0.12in}\begin{matrix} &\hfill\tikzmark{a12}\\&\hfill\tikzmark{b12}
\end{matrix}\,\,\,\begin{matrix}0&\epsilon^-_1&\epsilon^-_1\mathfrak v_3&\cdots\\
1&z_2&z_3&\cdots\\
\end{matrix}\hspace{-0.12in}\begin{matrix} &\hfill\tikzmark{g12}\\&\hfill\tikzmark{h12}\end{matrix}\,\,\,\begin{matrix}1&\mathfrak h^1_{22}+\mathfrak x_{12}^1(1-\epsilon_2^+)&\mathfrak h^1_{23}+\mathfrak x_{12}^1\mathfrak x_{13}^1(1-\epsilon_2^+)\\
1&\mathfrak h^1_{22}&\mathfrak h^1_{23}\\
\end{matrix}\right.\\
&\left.\,\,\,\,\,\begin{matrix}
  \cdots\\\cdots
\end{matrix}\hspace{-0.12in}\begin{matrix} &\hfill\tikzmark{c12}\\&\hfill\tikzmark{d12}\end{matrix}\,\,\,\begin{matrix}1&\mathfrak x_{12}^2(1-\epsilon^+_2)+\mathfrak h^2_{22}&\mathfrak x_{12}^2\mathfrak x_{13}^2(1-\epsilon^+_2)+\mathfrak h^2_{23}&\cdots\\
\epsilon_2^+&\epsilon_2^+\mathfrak h_{22}^2&\epsilon_2^+\mathfrak h_{23}^2&\cdots\\
\end{matrix}\right)=:\Theta.
\tikz[remember picture,overlay]   \draw[dashed,dash pattern={on 4pt off 2pt}] ([xshift=0.5\tabcolsep,yshift=7pt]a12.north) -- ([xshift=0.5\tabcolsep,yshift=-2pt]b12.south);\tikz[remember picture,overlay]   \draw[dashed,dash pattern={on 4pt off 2pt}] ([xshift=0.5\tabcolsep,yshift=7pt]c12.north) -- ([xshift=0.5\tabcolsep,yshift=-2pt]d12.south);\tikz[remember picture,overlay]   \draw[dashed,dash pattern={on 4pt off 2pt}] ([xshift=0.5\tabcolsep,yshift=7pt]g12.north) -- ([xshift=0.5\tabcolsep,yshift=-2pt]h12.south); 
\end{split}    
\end{equation}

\item[($W_{1B}$)] $W_{1B}:={\rm Spec}\,\mathbb Z\Bigl[\epsilon^+_2,\mathfrak x^1_{12},\cdots,\mathfrak x^1_{1s_3},\mathfrak h^1_{22},\cdots,\mathfrak h^1_{2s_3},\mathfrak x^2_{12},\cdots,\mathfrak x^2_{1s_4},\mathfrak h^2_{22},\cdots,\mathfrak h^2_{2s_4},\epsilon^-_1,z_2,\cdots,z_{s_2},\mathfrak v_3,$  $\left.\cdots,\mathfrak v_{s_2},\epsilon_3^+,y_2,\cdots,y_{s_1},\mathfrak u_3,\cdots,\mathfrak u_{s_1},\frac{1}{\epsilon^+_2}\right]$. The associated morphism $L:W_{2}\longrightarrow G(2,n)$ is defined by (\ref{didi3}).

\item[($W_{2A}$)] $W_{2A}:={\rm Spec}\,\mathbb Z\Bigl[\eta^1_{12},\mathfrak x^1_{13},\cdots,\mathfrak x^1_{1s_3},\mathfrak h^1_{22},\cdots,\mathfrak h^1_{2s_3},\mathfrak x^2_{11},\cdots,\mathfrak x^2_{1s_4},\mathfrak h^2_{22},\cdots,\mathfrak h^2_{2s_4},\epsilon^-_1,z_2,\cdots,z_{s_2},$  $\left.\mathfrak v_3,\cdots,\mathfrak v_{s_2},\epsilon_3^+,y_2,\cdots,y_{s_1},\mathfrak u_3,\cdots,\mathfrak u_{s_1},\frac{1}{1+\mathfrak x_{11}^2\left(\eta_{12}^1-\mathfrak h^1_{22}\right)},\frac{1}{\mathfrak x^2_{11}}\right]$. The associated morphism $L:W_{2A}\longrightarrow G(2,n)$ is defined by
\begin{equation}\label{w3a}
\begin{split}
&\left(\begin{matrix}1&y_2&y_3&\cdots\\
0&(1+\mathfrak x_{11}^2(\eta_{12}^1-\mathfrak h^1_{22}))\epsilon^+_3&(1+\mathfrak x_{11}^2(\eta_{12}^1-\mathfrak h^1_{22}))\epsilon^+_3\mathfrak u_3&\cdots\\
\end{matrix}\hspace{-0.12in}\begin{matrix} &\hfill\tikzmark{a12}\\&\hfill\tikzmark{b12}
\end{matrix}\,\,\,\begin{matrix}0&\epsilon^-_1&\epsilon^-_1\mathfrak v_3&\cdots\\
1&z_2&z_3&\cdots\\
\end{matrix}\hspace{-0.12in}\begin{matrix} &\hfill\tikzmark{g12}\\&\hfill\tikzmark{h12}\end{matrix}\right.\\
&\,\,\,\,\,\,\,\,\,\,\,\,\begin{matrix} &\hfill\tikzmark{c12}\\&\hfill\tikzmark{d12}\end{matrix}\,\,\,\begin{matrix}1&\eta_{12}^1&\mathfrak h^1_{23}+\mathfrak x_{13}^1(\eta_{12}^1-\mathfrak h_{22}^1)&\cdots\\
1&\mathfrak h^1_{22}&\mathfrak h^1_{23}&\cdots\\
\end{matrix}\hspace{-0.12in}\begin{matrix} &\hfill\tikzmark{e12}\\&\hfill\tikzmark{f12}\end{matrix}\,\,\,\begin{matrix}1\\1+\mathfrak x_{11}^2(\eta_{12}^1-\mathfrak h^1_{22})
\end{matrix}\\
&\left.\,\,\,\,\,\begin{matrix}\mathfrak h_{22}^2+(\mathfrak x_{12}^2+\mathfrak x_{11}^2\mathfrak h^2_{22})(\eta^1_{12}-\mathfrak h^1_{22})&\mathfrak h_{23}^2+(\mathfrak x_{12}^2\mathfrak x_{13}^2+\mathfrak x_{11}^2\mathfrak x_{13}^2\mathfrak h^2_{23})(\eta^1_{12}-\mathfrak h^1_{22})&\cdots\\
(1+\mathfrak x_{11}^2(\eta_{12}^1-\mathfrak h^1_{22}))\mathfrak h_{22}^2&(1+\mathfrak x_{11}^2(\eta_{12}^1-\mathfrak h^1_{22}))\mathfrak h^2_{23}&\cdots\\
\end{matrix}\right)=:\Theta.
\tikz[remember picture,overlay]   \draw[dashed,dash pattern={on 4pt off 2pt}] ([xshift=0.5\tabcolsep,yshift=7pt]a12.north) -- ([xshift=0.5\tabcolsep,yshift=-2pt]b12.south);\tikz[remember picture,overlay]   \draw[dashed,dash pattern={on 4pt off 2pt}] ([xshift=0.5\tabcolsep,yshift=7pt]c12.north) -- ([xshift=0.5\tabcolsep,yshift=-2pt]d12.south);\tikz[remember picture,overlay]   \draw[dashed,dash pattern={on 4pt off 2pt}] ([xshift=0.5\tabcolsep,yshift=7pt]e12.north) -- ([xshift=0.5\tabcolsep,yshift=-2pt]f12.south);\tikz[remember picture,overlay]   \draw[dashed,dash pattern={on 4pt off 2pt}] ([xshift=0.5\tabcolsep,yshift=7pt]g12.north) -- ([xshift=0.5\tabcolsep,yshift=-2pt]h12.south); 
\end{split}    
\end{equation}

\item[($W_{2B}$)] $W_{2B}:={\rm Spec}\,\mathbb Z\Bigl[\eta^1_{12},\mathfrak x^1_{13},\cdots,\mathfrak x^1_{1s_3},\mathfrak h^1_{22},\cdots,\mathfrak h^1_{2s_3},\mathfrak x^2_{11},\cdots,\mathfrak x^2_{1s_4},\mathfrak h^2_{22},\cdots,\mathfrak h^2_{2s_4},\epsilon^-_1,z_2,\cdots,z_{s_2},$  $\left.\mathfrak v_3,\cdots,\mathfrak v_{s_2},\epsilon_3^+,y_2,\cdots,y_{s_1},\mathfrak u_3,\cdots,\mathfrak u_{s_1},\frac{1}{1+\mathfrak x_{11}^2\left(\eta_{12}^1-\mathfrak h^1_{22}\right)},\frac{1}{1+\mathfrak x_{11}^2\eta^1_{12}}\right]$. The associated morphism $L:W_{2B}\longrightarrow G(2,n)$ is defined by  (\ref{w3a}).
\smallskip

\item[($W_3$)] 
$W_{3}:={\rm Spec}\,\mathbb Z\Bigl[\mathfrak x^1_{12},\cdots,\mathfrak x^1_{1s_3},\mathfrak h^1_{22},\cdots,\mathfrak h^1_{2s_3},\mathfrak x^2_{11},\eta^2_{12},\mathfrak x^2_{13},\mathfrak x^2_{14},\cdots,\mathfrak x^2_{1s_4},\mathfrak h^2_{22},\cdots,\mathfrak h^2_{2s_4},\epsilon^-_1,z_2,\cdots,$  $\left.z_{s_2},\mathfrak v_3,\cdots,\mathfrak v_{s_2},\epsilon_3^+,y_2,\cdots,y_{s_1},\mathfrak u_3,\cdots,\mathfrak u_{s_1},\frac{1}{1-\mathfrak x^2_{11}\eta^2_{12}},\frac{1}{1-\mathfrak x^2_{11}\mathfrak h^2_{22}}\right]$. The associated morphism $L:W_{3}\longrightarrow G(2,n)$ is defined  by
\begin{equation}\label{w3}
\begin{split}
&\left(\begin{matrix}1&y_2&y_3&\cdots\\
0&\frac{1-\mathfrak x_{11}^2\eta^2_{12}}{1-\mathfrak x_{11}^2\mathfrak h^2_{22}}\epsilon^+_3&\frac{1-\mathfrak x_{11}^2\eta^2_{12}}{1-\mathfrak x_{11}^2\mathfrak h^2_{22}}\epsilon^+_3\mathfrak u_3&\cdots\\
\end{matrix}\hspace{-0.12in}\begin{matrix} &\hfill\tikzmark{a12}\\
\\&\hfill\tikzmark{b12}
\end{matrix}\,\,\,\begin{matrix}0&\epsilon^-_1&\epsilon^-_1\mathfrak v_3&\cdots\\
1&z_2&z_3&\cdots\\
\end{matrix}\hspace{-0.12in}\begin{matrix} &\hfill\tikzmark{g12}\\
\\&\hfill\tikzmark{h12}\end{matrix}\,\,\,\begin{matrix}1&\mathfrak h^1_{22}+\mathfrak x_{12}^1\frac{\eta_{12}^2-\mathfrak h^2_{22}}{1-\mathfrak x^2_{11}\mathfrak h^2_{22}}&\mathfrak h^1_{23}+\mathfrak x_{12}^1\mathfrak x_{13}^1\frac{\eta_{12}^2-\mathfrak h^2_{22}}{1-\mathfrak x^2_{11}\mathfrak h^2_{22}}\\
1&\mathfrak h^1_{22}&\mathfrak h^1_{23}\\
\end{matrix}\right.\\
&\left.\,\,\,\,\,\,\,\,\,\,\,\,\,\,\,\begin{matrix}\cdots\\
\cdots\\
\end{matrix}\hspace{-0.12in}\begin{matrix} &\hfill\tikzmark{c12}\\
\\&\hfill\tikzmark{d12}\end{matrix}\,\,\,\begin{matrix}1&\eta_{12}^2&\frac{\eta^2_{12}-\mathfrak h^2_{22}}{1-\mathfrak x_{11}^2\mathfrak h^2_{22}}\mathfrak x_{13}^2+\frac{1-\mathfrak x_{11}^2\eta^2_{12}}{1-\mathfrak x_{11}^2\mathfrak h^2_{22}}\mathfrak h^2_{23}&\cdots\\
\frac{1-\mathfrak x_{11}^2\eta^2_{12}}{1-\mathfrak x_{11}^2\mathfrak h^2_{22}}&\frac{1-\mathfrak x_{11}^2\eta^2_{12}}{1-\mathfrak x_{11}^2\mathfrak h^2_{22}}\mathfrak h_{22}^2&\frac{1-\mathfrak x_{11}^2\eta^2_{12}}{1-\mathfrak x_{11}^2\mathfrak h^2_{22}}\mathfrak h_{23}^2&\cdots\\
\end{matrix}\right)=:\Theta.
\tikz[remember picture,overlay]   \draw[dashed,dash pattern={on 4pt off 2pt}] ([xshift=0.5\tabcolsep,yshift=7pt]a12.north) -- ([xshift=0.5\tabcolsep,yshift=-2pt]b12.south);\tikz[remember picture,overlay]   \draw[dashed,dash pattern={on 4pt off 2pt}] ([xshift=0.5\tabcolsep,yshift=7pt]c12.north) -- ([xshift=0.5\tabcolsep,yshift=-2pt]d12.south);\tikz[remember picture,overlay]   \draw[dashed,dash pattern={on 4pt off 2pt}] ([xshift=0.5\tabcolsep,yshift=7pt]g12.north) -- ([xshift=0.5\tabcolsep,yshift=-2pt]h12.south);
\end{split}    
\end{equation}

\item[($W_4$)] $W_{4}:={\rm Spec}\,\mathbb Z\Bigl[\epsilon^+_2,\eta^1_{12},\cdots,\eta^1_{1s_3},\mathfrak h^1_{21},\mathfrak h^1_{23},\mathfrak h^1_{24},\cdots,\mathfrak h^1_{2s_3},\eta^2_{12},\cdots,\eta^2_{1s_4},\mathfrak h^2_{22},\cdots,\mathfrak h^2_{2s_4},\epsilon^-_1,z_2,\cdots,$  $\left.z_{s_2},\mathfrak v_3,\cdots,\mathfrak v_{s_2},\epsilon_3^+,y_2,\cdots,y_{s_1},\mathfrak u_3,\cdots,\mathfrak u_{s_1},\frac{1}{1-\eta_{12}^1\epsilon^+_2},\frac{1}{1-\eta_{12}^1\mathfrak h_{21}^1}\right]$. The associated morphism $L:W_{4}\longrightarrow G(2,n)$ is defined by
\begin{equation}\label{w4}
\begin{split}
&\left(\begin{matrix}1&y_2&y_3&\cdots\\
0&\epsilon_2^+\epsilon^+_3&\epsilon_2^+\epsilon^+_3\mathfrak u_3&\cdots\\
\end{matrix}\hspace{-0.12in}\begin{matrix} &\hfill\tikzmark{a12}\\&\hfill\tikzmark{b12}
\end{matrix}\,\,\,\begin{matrix}0&\epsilon^-_1&\epsilon^-_1\mathfrak v_3&\cdots\\
1&z_2&z_3&\cdots\\
\end{matrix}\hspace{-0.12in}\begin{matrix} &\hfill\tikzmark{g12}\\&\hfill\tikzmark{h12}\end{matrix}\,\,\,\begin{matrix}1&\eta_{12}^1&\eta_{13}^1&\cdots\\
\mathfrak{h}^1_{21}&1&\mathfrak{h}^1_{23}&\cdots\\
\end{matrix}\hspace{-0.12in}\begin{matrix} &\hfill\tikzmark{e12}\\&\hfill\tikzmark{f12}\end{matrix}\right.\\
&\left.\,\,\,\,\,\,\,\,\,\,\,\,\,\,\,\,\,\,\,\begin{matrix} &\hfill\tikzmark{c12}\\&\hfill\tikzmark{d12}\end{matrix}\,\,\,\begin{matrix}1&\eta_{12}^2&\eta_{13}^2&\cdots\\
\epsilon_2^+&\epsilon_2^+\left(\mathfrak{h}^2_{22}+\eta^2_{12}\right)&\epsilon_2^+\left(\mathfrak{h}^2_{22}\mathfrak{h}^2_{23}+\eta^2_{13}\right)&\cdots\\
\end{matrix}\right)=:\Theta.
\tikz[remember picture,overlay]   \draw[dashed,dash pattern={on 4pt off 2pt}] ([xshift=0.5\tabcolsep,yshift=7pt]a12.north) -- ([xshift=0.5\tabcolsep,yshift=-2pt]b12.south);\tikz[remember picture,overlay]   \draw[dashed,dash pattern={on 4pt off 2pt}] ([xshift=0.5\tabcolsep,yshift=7pt]c12.north) -- ([xshift=0.5\tabcolsep,yshift=-2pt]d12.south);\tikz[remember picture,overlay]   \draw[dashed,dash pattern={on 4pt off 2pt}] ([xshift=0.5\tabcolsep,yshift=7pt]e12.north) -- ([xshift=0.5\tabcolsep,yshift=-2pt]f12.south);\tikz[remember picture,overlay]   \draw[dashed,dash pattern={on 4pt off 2pt}] ([xshift=0.5\tabcolsep,yshift=7pt]g12.north) -- ([xshift=0.5\tabcolsep,yshift=-2pt]h12.south);
\end{split}    
\end{equation}

\item[($W_5$)] 
$W_{5}:={\rm Spec}\,\mathbb Z\Bigl[\epsilon^+_2,\eta^1_{12},\cdots,\eta^1_{1s_3},\mathfrak h^1_{22},\cdots,\mathfrak h^1_{2s_3},\eta^2_{12},\cdots,\eta^2_{1s_4},\mathfrak h^2_{21},\mathfrak h^2_{23},\mathfrak h^2_{24},\cdots,\mathfrak h^2_{2s_4},\epsilon^-_1,z_2,\cdots,$  $\left.z_{s_2},\mathfrak v_3,\cdots,\mathfrak v_{s_2},\epsilon_3^+,y_2,\cdots,y_{s_1},\mathfrak u_3,\cdots,\mathfrak u_{s_1},\frac{1}{1-\epsilon^+_2\mathfrak h_{21}^2},\frac{1}{1-\eta_{12}^2\mathfrak h_{21}^2}\right]$. The associated morphism $L:W_{5}\longrightarrow G(2,n)$ is defined by
\begin{equation}\label{w5}
\begin{split}
&\left(\begin{matrix}1&y_2&y_3&\cdots\\
0&\epsilon_2^+\epsilon^+_3&\epsilon_2^+\epsilon^+_3\mathfrak u_3&\cdots\\
\end{matrix}\hspace{-0.12in}\begin{matrix} &\hfill\tikzmark{a12}\\&\hfill\tikzmark{b12}
\end{matrix}\,\,\,\begin{matrix}0&\epsilon^-_1&\epsilon^-_1\mathfrak v_3&\cdots\\
1&z_2&z_3&\cdots\\
\end{matrix}\hspace{-0.12in}\begin{matrix} &\hfill\tikzmark{g12}\\&\hfill\tikzmark{h12}\end{matrix}\,\,\,\begin{matrix}1&\eta_{12}^1&\eta_{13}^1&\cdots\\
1&\mathfrak{h}^1_{22}+\eta^1_{12}&\mathfrak{h}^1_{23}\mathfrak h^1_{22}+\eta^1_{13}&\cdots\\
\end{matrix}\hspace{-0.12in}\begin{matrix} &\hfill\tikzmark{e12}\\&\hfill\tikzmark{f12}\end{matrix}\right.\\
&\left.\,\,\,\,\,\,\,\,\,\,\,\,\,\,\,\,\,\,\,\begin{matrix} &\hfill\tikzmark{c12}\\&\hfill\tikzmark{d12}\end{matrix}\,\,\,\begin{matrix}1&\eta_{12}^2&\eta_{13}^2&\cdots\\
\epsilon_2^+\mathfrak{h}^2_{21}&\epsilon_2^+&\epsilon_2^+\mathfrak{h}^2_{23}&\cdots\\
\end{matrix}\right)=:\Theta.
\tikz[remember picture,overlay]   \draw[dashed,dash pattern={on 4pt off 2pt}] ([xshift=0.5\tabcolsep,yshift=7pt]a12.north) -- ([xshift=0.5\tabcolsep,yshift=-2pt]b12.south);\tikz[remember picture,overlay]   \draw[dashed,dash pattern={on 4pt off 2pt}] ([xshift=0.5\tabcolsep,yshift=7pt]c12.north) -- ([xshift=0.5\tabcolsep,yshift=-2pt]d12.south);\tikz[remember picture,overlay]   \draw[dashed,dash pattern={on 4pt off 2pt}] ([xshift=0.5\tabcolsep,yshift=7pt]e12.north) -- ([xshift=0.5\tabcolsep,yshift=-2pt]f12.south);\tikz[remember picture,overlay]   \draw[dashed,dash pattern={on 4pt off 2pt}] ([xshift=0.5\tabcolsep,yshift=7pt]g12.north) -- ([xshift=0.5\tabcolsep,yshift=-2pt]h12.south);
\end{split}    
\end{equation}

\item[($W_6$)] 
$W_{6}:={\rm Spec}\,\mathbb Z\Bigl[\epsilon^+_2,\eta^1_{12},\cdots,\eta^1_{1s_3},\mathfrak h^1_{21},\mathfrak h^1_{23},\mathfrak h^1_{24},\cdots,\mathfrak h^1_{2s_3},\eta^2_{12},\cdots,\eta^2_{1s_4},\mathfrak h^2_{21},\mathfrak h^2_{23},\mathfrak h^2_{24},\cdots,\mathfrak h^2_{2s_4},\epsilon^-_1,$  $\left.z_2,\cdots,z_{s_2},\mathfrak v_3,\cdots,\mathfrak v_{s_2},\epsilon_3^+,y_2,\cdots,y_{s_1},\mathfrak u_3,\cdots,\mathfrak u_{s_1},\frac{1}{1-\eta_{12}^1\mathfrak h_{21}^1},\frac{1}{1-\eta_{12}^2\mathfrak h_{21}^2},\frac{1}{1-\eta_{12}^1\epsilon^+_2\mathfrak h_{21}^2}\right]$. The associated morphism $L:W_{6}\longrightarrow G(2,n)$ is defined by

\begin{equation}\label{w6}
\begin{split}
&\left(\begin{matrix}1&y_2&y_3&\cdots\\
0&\epsilon_2^+\epsilon^+_3&\epsilon_2^+\epsilon^+_3\mathfrak u_3&\cdots\\
\end{matrix}\hspace{-0.12in}\begin{matrix} &\hfill\tikzmark{a12}\\&\hfill\tikzmark{b12}
\end{matrix}\,\,\,\begin{matrix}0&\epsilon^-_1&\epsilon^-_1\mathfrak v_3&\cdots\\
1&z_2&z_3&\cdots\\
\end{matrix}\hspace{-0.12in}\begin{matrix} &\hfill\tikzmark{g12}\\&\hfill\tikzmark{h12}\end{matrix}\,\,\,\begin{matrix}1&\eta_{12}^1&\eta_{13}^1&\cdots\\
\mathfrak{h}^1_{21}&1&\mathfrak{h}^1_{23}&\cdots\\
\end{matrix}\hspace{-0.12in}\begin{matrix} &\hfill\tikzmark{e12}\\&\hfill\tikzmark{f12}\end{matrix}\,\,\,\begin{matrix}1&\eta_{12}^2&\eta_{13}^2&\cdots\\
\epsilon_2^+\mathfrak{h}^2_{21}&\epsilon_2^+&\epsilon_2^+\mathfrak{h}^2_{23}&\cdots\\
\end{matrix}\right)=:\Theta.
\tikz[remember picture,overlay]   \draw[dashed,dash pattern={on 4pt off 2pt}] ([xshift=0.5\tabcolsep,yshift=7pt]a12.north) -- ([xshift=0.5\tabcolsep,yshift=-2pt]b12.south);\tikz[remember picture,overlay]   \draw[dashed,dash pattern={on 4pt off 2pt}] ([xshift=0.5\tabcolsep,yshift=7pt]e12.north) -- ([xshift=0.5\tabcolsep,yshift=-2pt]f12.south);\tikz[remember picture,overlay]   \draw[dashed,dash pattern={on 4pt off 2pt}] ([xshift=0.5\tabcolsep,yshift=7pt]g12.north) -- ([xshift=0.5\tabcolsep,yshift=-2pt]h12.south);
\end{split}    
\end{equation}

\item[($W_{7A}$)] $W_{7A}:={\rm Spec}\,\mathbb Z\Bigl[\epsilon^+_2,\eta^1_{12},\cdots,\eta^1_{1s_3},\mathfrak h^1_{22},\cdots,\mathfrak h^1_{2s_3},\eta^2_{12},\cdots,\eta^2_{1s_4},\mathfrak h^2_{21},\cdots,\mathfrak h^2_{2s_4},\epsilon^-_1,z_2,\cdots,z_{s_2},$  $\left.\mathfrak v_3,\cdots,\mathfrak v_{s_2},y_2,\cdots,y_{s_1},\mathfrak u_3,\cdots,\mathfrak u_{s_1},\frac{1}{1-\epsilon_{2}^+\mathfrak h_{21}^2},\frac{1}{\mathfrak h_{21}^2}\right]$. The associated morphism $L:W_{7A}\longrightarrow G(2,n)$ is defined by
\begin{equation}\label{w8a}
\begin{split}
&\left(\begin{matrix}1&y_2&y_3&\cdots\\
0&\epsilon_2^+&\epsilon_2^+\mathfrak u_3&\cdots\\
\end{matrix}\hspace{-0.12in}\begin{matrix} &\hfill\tikzmark{a12}\\&\hfill\tikzmark{b12}
\end{matrix}\,\,\,\begin{matrix}0&\epsilon^-_1&\epsilon^-_1\mathfrak v_3&\cdots\\
1&z_2&z_3&\cdots\\
\end{matrix}\hspace{-0.12in}\begin{matrix} &\hfill\tikzmark{g12}\\&\hfill\tikzmark{h12}\end{matrix}\,\,\,\begin{matrix}1&\eta_{12}^1&\eta_{13}^1&\cdots\\
1&\mathfrak{h}^1_{22}+\eta_{12}^1&\mathfrak{h}^1_{22}\mathfrak{h}^1_{23}+\eta_{13}^1&\cdots\\
\end{matrix}\hspace{-0.12in}\begin{matrix} &\hfill\tikzmark{e12}\\&\hfill\tikzmark{f12}\end{matrix}\right.\\
&\left.\,\,\,\,\,\,\,\,\,\,\,\,\,\,\,\,\,\,\,\begin{matrix} &\hfill\tikzmark{c12}\\&\hfill\tikzmark{d12}\end{matrix}\,\,\,\begin{matrix}1&\eta_{12}^2&\eta_{13}^2&\cdots\\
\epsilon_2^+\mathfrak{h}^2_{21}&\epsilon_2^+\left(\mathfrak{h}^2_{22}+\eta^2_{12}\mathfrak{h}^2_{21}\right)&\epsilon_2^+\left(\mathfrak{h}^2_{22}\mathfrak{h}^2_{23}+\eta^2_{13}\mathfrak{h}^2_{21}\right)&\cdots\\
\end{matrix}\right)=:\Theta.
\tikz[remember picture,overlay]   \draw[dashed,dash pattern={on 4pt off 2pt}] ([xshift=0.5\tabcolsep,yshift=7pt]a12.north) -- ([xshift=0.5\tabcolsep,yshift=-2pt]b12.south);\tikz[remember picture,overlay]   \draw[dashed,dash pattern={on 4pt off 2pt}] ([xshift=0.5\tabcolsep,yshift=7pt]c12.north) -- ([xshift=0.5\tabcolsep,yshift=-2pt]d12.south);\tikz[remember picture,overlay]   \draw[dashed,dash pattern={on 4pt off 2pt}] ([xshift=0.5\tabcolsep,yshift=7pt]e12.north) -- ([xshift=0.5\tabcolsep,yshift=-2pt]f12.south);\tikz[remember picture,overlay]   \draw[dashed,dash pattern={on 4pt off 2pt}] ([xshift=0.5\tabcolsep,yshift=7pt]g12.north) -- ([xshift=0.5\tabcolsep,yshift=-2pt]h12.south);
\end{split}    
\end{equation}

\item[($W_{7B}$)] $W_{7B}:={\rm Spec}\,\mathbb Z\Bigl[\epsilon^+_2,\eta^1_{12},\cdots,\eta^1_{1s_3},\mathfrak h^1_{22},\cdots,\mathfrak h^1_{2s_3},\eta^2_{12},\cdots,\eta^2_{1s_4},\mathfrak h^2_{21},\cdots,\mathfrak h^2_{2s_4},\epsilon^-_1,z_2,\cdots,z_{s_2},$  $\left.\mathfrak v_3,\cdots,\mathfrak v_{s_2},y_2,\cdots,y_{s_1},\mathfrak u_3,\cdots,\mathfrak u_{s_1},\frac{1}{1-\epsilon_{2}^+\mathfrak h_{21}^2},\frac{1}{1-y_2\mathfrak h_{21}^2}\right]$. The associated morphism $L:W_{7B}\longrightarrow G(2,n)$ is defined by  (\ref{w8a}).
\smallskip

\item[($W_{8A}$)] $W_{8A}:={\rm Spec}\,\mathbb Z\Bigl[\epsilon^+_2,\eta^1_{12},\cdots,\eta^1_{1s_3},\mathfrak h^1_{21},\mathfrak h^1_{23},\mathfrak h^1_{24},\cdots,\mathfrak h^1_{2s_3},\eta^2_{12},\cdots,\eta^2_{1s_4},\mathfrak h^2_{21},\cdots,\mathfrak h^2_{2s_4},\epsilon^-_1,z_2,$  $\left.\cdots,z_{s_2},\mathfrak v_3,\cdots,\mathfrak v_{s_2},y_2,\cdots,y_{s_1},\mathfrak u_3,\cdots,\mathfrak u_{s_1},\frac{1}{1-\eta_{12}^1\epsilon_{2}^+\mathfrak h_{21}^2},\frac{1}{1-\eta_{12}^1\mathfrak h_{21}^2},\frac{1}{\mathfrak h_{21}^2}\right]$. The associated morphism $L:W_{8A}\longrightarrow G(2,n)$ is defined by
\begin{equation}\label{w9a}
\begin{split}
&\left(\begin{matrix}1&y_2&y_3&\cdots\\
0&\epsilon_2^+&\epsilon_2^+\mathfrak u_3&\cdots\\
\end{matrix}\hspace{-0.12in}\begin{matrix} &\hfill\tikzmark{a12}\\&\hfill\tikzmark{b12}
\end{matrix}\,\,\,\begin{matrix}0&\epsilon^-_1&\epsilon^-_1\mathfrak v_3&\cdots\\
1&z_2&z_3&\cdots\\
\end{matrix}\hspace{-0.12in}\begin{matrix} &\hfill\tikzmark{g12}\\&\hfill\tikzmark{h12}\end{matrix}\,\,\,\begin{matrix}1&\eta_{12}^1&\eta_{13}^1&\cdots\\
\mathfrak{h}^1_{21}&1&\mathfrak{h}^1_{23}&\cdots\\
\end{matrix}\hspace{-0.12in}\begin{matrix} &\hfill\tikzmark{e12}\\&\hfill\tikzmark{f12}\end{matrix}\right.\\
&\left.\,\,\,\,\,\,\,\,\,\,\begin{matrix} &\hfill\tikzmark{c12}\\&\hfill\tikzmark{d12}\end{matrix}\,\,\,\begin{matrix}1&\eta_{12}^2&\eta_{13}^2&\cdots\\
\epsilon_2^+\mathfrak{h}^2_{21}&\epsilon_2^+\left(\mathfrak{h}^2_{22}+\eta^2_{12}\mathfrak{h}^2_{21}\right)&\epsilon_2^+\left(\mathfrak{h}^2_{22}\mathfrak{h}^2_{23}+\eta^2_{13}\mathfrak{h}^2_{21}\right)&\cdots\\
\end{matrix}\right)=:\Theta.
\tikz[remember picture,overlay]   \draw[dashed,dash pattern={on 4pt off 2pt}] ([xshift=0.5\tabcolsep,yshift=7pt]a12.north) -- ([xshift=0.5\tabcolsep,yshift=-2pt]b12.south);\tikz[remember picture,overlay]   \draw[dashed,dash pattern={on 4pt off 2pt}] ([xshift=0.5\tabcolsep,yshift=7pt]c12.north) -- ([xshift=0.5\tabcolsep,yshift=-2pt]d12.south);\tikz[remember picture,overlay]   \draw[dashed,dash pattern={on 4pt off 2pt}] ([xshift=0.5\tabcolsep,yshift=7pt]e12.north) -- ([xshift=0.5\tabcolsep,yshift=-2pt]f12.south);\tikz[remember picture,overlay]   \draw[dashed,dash pattern={on 4pt off 2pt}] ([xshift=0.5\tabcolsep,yshift=7pt]g12.north) -- ([xshift=0.5\tabcolsep,yshift=-2pt]h12.south);
\end{split}    
\end{equation}

\item[($W_{8B}$)] $W_{8B}:={\rm Spec}\,\mathbb Z\Bigl[\epsilon^+_2,\eta^1_{12},\cdots,\eta^1_{1s_3},\mathfrak h^1_{21},\mathfrak h^1_{23},\mathfrak h^1_{24},\cdots,\mathfrak h^1_{2s_3},\eta^2_{12},\cdots,\eta^2_{1s_4},\mathfrak h^2_{21},\cdots,\mathfrak h^2_{2s_4},\epsilon^-_1,z_2,$  $\left.\cdots,z_{s_2},\mathfrak v_3,\cdots,\mathfrak v_{s_2},y_2,\cdots,y_{s_1},\mathfrak u_3,\cdots,\mathfrak u_{s_1},\frac{1}{1-\eta_{12}^1\epsilon_{2}^+\mathfrak h_{21}^2},\frac{1}{1-\eta_{12}^1\mathfrak h_{21}^2},\frac{1}{1-y_2\mathfrak h_{21}^2}\right]$. The associated morphism $L:W_{8B}\longrightarrow G(2,n)$ is defined by  (\ref{w9a}).
\smallskip

\item[($W_{9A}$)]$W_{9A}:={\rm Spec}\,\mathbb Z\Bigl[\eta^1_{12},\cdots,\eta^1_{1s_3},\mathfrak h^1_{21},\cdots,\mathfrak h^1_{2s_3},\eta^2_{12},\cdots,\eta^2_{1s_4},\mathfrak h^2_{21},\cdots,\mathfrak h^2_{2s_4},\epsilon^-_1,z_2,\cdots,z_{s_2},\mathfrak v_3,$  $\left.\cdots,\mathfrak v_{s_2},y_2,\cdots,y_{s_1},\mathfrak u_3,\cdots,\mathfrak u_{s_1},\frac{1}{1-\eta_{12}^2\mathfrak h_{21}^2},\frac{1}{\mathfrak h_{21}^1},\frac{1}{\mathfrak{h}^2_{21}\left(\mathfrak{h}^2_{22}-\eta_{12}^2\mathfrak{h}^1_{21}\right)+\mathfrak{h}^1_{21}}\right]$. The associated morphism $L:W_{9A}\longrightarrow G(2,n)$ is defined by
\begin{equation}\label{w10a}
\begin{split}
&\left(\begin{matrix}1&y_2&y_3&\cdots\\
0&1&\mathfrak u_3&\cdots\\
\end{matrix}\hspace{-0.12in}\begin{matrix} &\hfill\tikzmark{a12}\\&\hfill\tikzmark{b12}
\end{matrix}\,\,\,\begin{matrix}0&\epsilon^-_1&\epsilon^-_1\mathfrak v_3&\cdots\\
1&z_2&z_3&\cdots\\
\end{matrix}\hspace{-0.12in}\begin{matrix} &\hfill\tikzmark{g12}\\&\hfill\tikzmark{h12}\end{matrix}\,\,\,\begin{matrix}1&\eta_{12}^1&\eta_{13}^1\\
\mathfrak{h}^1_{21}&\mathfrak{h}^1_{22}\left(\mathfrak{h}^2_{22}-\eta_{12}^2\mathfrak{h}^1_{21}\right)+\mathfrak{h}^1_{21}\eta_{12}^1&\mathfrak{h}^1_{22}\mathfrak{h}^1_{23}\left(\mathfrak{h}^2_{22}-\eta_{12}^2\mathfrak{h}^1_{21}\right)+\mathfrak{h}^1_{21}\eta_{13}^1\\
\end{matrix}\right.\\
&\left.\,\,\,\,\,\begin{matrix} \cdots\\
\cdots\\\end{matrix}\hspace{-0.12in}\begin{matrix} &\hfill\tikzmark{c12}\\&\hfill\tikzmark{d12}\end{matrix}\,\,\,\begin{matrix}1&\eta_{12}^2&\eta_{13}^2&\cdots\\
\mathfrak{h}^2_{21}\left(\mathfrak{h}^2_{22}-\eta_{12}^2\mathfrak{h}^1_{21}\right)+\mathfrak{h}^1_{21}&\mathfrak{h}^2_{22}&\mathfrak{h}^2_{23}\left(\mathfrak{h}^2_{22}-\eta_{12}^2\mathfrak{h}^1_{21}\right)+\mathfrak{h}^1_{21}\eta_{13}^2&\cdots\\
\end{matrix}\right)=:\Theta.
\tikz[remember picture,overlay]   \draw[dashed,dash pattern={on 4pt off 2pt}] ([xshift=0.5\tabcolsep,yshift=7pt]a12.north) -- ([xshift=0.5\tabcolsep,yshift=-2pt]b12.south);\tikz[remember picture,overlay]   \draw[dashed,dash pattern={on 4pt off 2pt}] ([xshift=0.5\tabcolsep,yshift=7pt]c12.north) -- ([xshift=0.5\tabcolsep,yshift=-2pt]d12.south);\tikz[remember picture,overlay]   \draw[dashed,dash pattern={on 4pt off 2pt}] ([xshift=0.5\tabcolsep,yshift=7pt]g12.north) -- ([xshift=0.5\tabcolsep,yshift=-2pt]h12.south);
\end{split}    
\end{equation}

\item[($W_{9B}$)] $W_{9B}:={\rm Spec}\,\mathbb Z\Bigl[\eta^1_{12},\cdots,\eta^1_{1s_3},\mathfrak h^1_{21},\cdots,\mathfrak h^1_{2s_3},\eta^2_{12},\cdots,\eta^2_{1s_4},\mathfrak h^2_{21},\cdots,\mathfrak h^2_{2s_4},\epsilon^-_1,z_2,\cdots,z_{s_2},\mathfrak v_3,$  $\left.\cdots,\mathfrak v_{s_2},y_2,\cdots,y_{s_1},\mathfrak u_3,\cdots,\mathfrak u_{s_1},\frac{1}{1-\eta_{12}^2\mathfrak h_{21}^2},\frac{1}{\mathfrak h_{21}^1},\frac{1}{1-y_2\left(\mathfrak{h}^2_{21}\left(\mathfrak{h}^2_{22}-\eta_{12}^2\mathfrak{h}^1_{21}\right)+\mathfrak{h}^1_{21}\right)}\right]$. The associated morphism $L:W_{9B}\longrightarrow G(2,n)$ is defined by  (\ref{w10a}).
\smallskip

\item[($W_{9C}$)] $W_{9C}:={\rm Spec}\,\mathbb Z\Bigl[\eta^1_{12},\cdots,\eta^1_{1s_3},\mathfrak h^1_{21},\cdots,\mathfrak h^1_{2s_3},\eta^2_{12},\cdots,\eta^2_{1s_4},\mathfrak h^2_{21},\cdots,\mathfrak h^2_{2s_4},\epsilon^-_1,z_2,\cdots,z_{s_2},\mathfrak v_3,$  $\left.\cdots,\mathfrak v_{s_2},y_2,\cdots,y_{s_1},\mathfrak u_3,\cdots,\mathfrak u_{s_1},\frac{1}{1-\eta_{12}^2\mathfrak h_{21}^2},\frac{1}{1-y_2\mathfrak h_{21}^1},\frac{1}{\mathfrak{h}^2_{21}\left(\mathfrak{h}^2_{22}-\eta_{12}^2\mathfrak{h}^1_{21}\right)+\mathfrak{h}^1_{21}}\right]$. The associated morphism $L:W_{9C}\longrightarrow G(2,n)$ is defined by (\ref{w10a}).
\smallskip

\item[($W_{9D}$)] $W_{9D}:={\rm Spec}\,\mathbb Z\Bigl[\eta^1_{12},\cdots,\eta^1_{1s_3},\mathfrak h^1_{21},\cdots,\mathfrak h^1_{2s_3},\eta^2_{12},\cdots,\eta^2_{1s_4},\mathfrak h^2_{21},\cdots,\mathfrak h^2_{2s_4},\epsilon^-_1,z_2,\cdots,z_{s_2},\mathfrak v_3,$  $\left.\cdots,\mathfrak v_{s_2},y_2,\cdots,y_{s_1},\mathfrak u_3,\cdots,\mathfrak u_{s_1},\frac{1}{1-\eta_{12}^2\mathfrak h_{21}^2},\frac{1}{1-y_2\mathfrak h_{21}^1},\frac{1}{1-y_2\left(\mathfrak{h}^2_{21}\left(\mathfrak{h}^2_{22}-\eta_{12}^2\mathfrak{h}^1_{21}\right)+\mathfrak{h}^1_{21}\right)}\right]$. The associated morphism $L:W_{9D}\longrightarrow G(2,n)$ is defined by  (\ref{w10a}).
\smallskip

\item[($W_{10A}$)] $W_{10A}:={\rm Spec}\,\mathbb Z\Bigl[\eta^1_{12},\cdots,\eta^1_{1s_3},\mathfrak h^1_{21},\cdots,\mathfrak h^1_{2s_3},\eta^2_{12},\cdots,\eta^2_{1s_4},\mathfrak h^2_{21},\cdots,\mathfrak h^2_{2s_4},\epsilon^-_1,z_2,\cdots,z_{s_2},\mathfrak v_3,$  $\left.\cdots,\mathfrak v_{s_2},y_2,\cdots,y_{s_1},\mathfrak u_3,\cdots,\mathfrak u_{s_1},\frac{1}{1-\eta_{12}^1\mathfrak h_{21}^2},\frac{1}{\mathfrak h_{21}^1},\frac{1}{\mathfrak{h}^2_{21}\left(\mathfrak{h}^1_{22}-\eta_{12}^1\mathfrak{h}^1_{21}\right)+\mathfrak{h}^1_{21}}\right]$. The associated morphism $L:W_{10A}\longrightarrow G(2,n)$ is defined by
\begin{equation}\label{w11a}
\begin{split}
&\left(\begin{matrix}1&y_2&y_3&\cdots\\
0&1&\mathfrak u_3&\cdots\\
\end{matrix}\hspace{-0.12in}\begin{matrix} &\hfill\tikzmark{a12}\\&\hfill\tikzmark{b12}
\end{matrix}\,\,\,\begin{matrix}0&\epsilon^-_1&\epsilon^-_1\mathfrak v_3&\cdots\\
1&z_2&z_3&\cdots\\
\end{matrix}\hspace{-0.12in}\begin{matrix} &\hfill\tikzmark{g12}\\&\hfill\tikzmark{h12}\end{matrix}\,\,\,\begin{matrix}1&\eta_{12}^1&\eta_{13}^1&\cdots\\
\mathfrak{h}^1_{21}&\mathfrak{h}^1_{22}&\mathfrak{h}^1_{23}\left(\mathfrak{h}^1_{22}-\eta_{12}^1\mathfrak{h}^1_{21}\right)+\eta_{13}^1\mathfrak{h}^1_{21}&\cdots\\
\end{matrix}\hspace{-0.12in}\begin{matrix} &\hfill\tikzmark{e12}\\&\hfill\tikzmark{f12}\end{matrix}\right.\\
&\left.\,\,\,\,\,\begin{matrix} &\hfill\tikzmark{c12}\\&\hfill\tikzmark{d12}\end{matrix}\,\,\,\begin{matrix}1&\eta_{12}^2\\
\mathfrak{h}^2_{21}\left(\mathfrak{h}^1_{22}-\eta_{12}^1\mathfrak{h}^1_{21}\right)+\mathfrak{h}^1_{21}&\left(\mathfrak{h}^2_{22}+\mathfrak h_{21}^2\eta_{12}^2\right)\left(\mathfrak{h}^1_{22}-\eta_{12}^1\mathfrak{h}^1_{21}\right)+\mathfrak{h}^1_{21}\eta_{12}^2\\
\end{matrix}\right.\\
&\left.\,\,\,\,\,\,\,\,\,\,\,\,\,\,\,\,\,\,\,\begin{matrix}\eta_{13}^2&\cdots\\
\left(\mathfrak{h}^2_{22}\mathfrak{h}^2_{23}+\mathfrak h_{21}^2\eta_{13}^2\right)\left(\mathfrak{h}^1_{22}-\eta_{12}^1\mathfrak{h}^1_{21}\right)+\mathfrak{h}^1_{21}\eta_{13}^2&\cdots\\
\end{matrix}\right)=:\Theta.
\tikz[remember picture,overlay]   \draw[dashed,dash pattern={on 4pt off 2pt}] ([xshift=0.5\tabcolsep,yshift=7pt]a12.north) -- ([xshift=0.5\tabcolsep,yshift=-2pt]b12.south);\tikz[remember picture,overlay]   \draw[dashed,dash pattern={on 4pt off 2pt}] ([xshift=0.5\tabcolsep,yshift=7pt]c12.north) -- ([xshift=0.5\tabcolsep,yshift=-2pt]d12.south);\tikz[remember picture,overlay]   \draw[dashed,dash pattern={on 4pt off 2pt}] ([xshift=0.5\tabcolsep,yshift=7pt]e12.north) -- ([xshift=0.5\tabcolsep,yshift=-2pt]f12.south);\tikz[remember picture,overlay]   \draw[dashed,dash pattern={on 4pt off 2pt}] ([xshift=0.5\tabcolsep,yshift=7pt]g12.north) -- ([xshift=0.5\tabcolsep,yshift=-2pt]h12.south);
\end{split}    
\end{equation}

\item[($W_{10B}$)] $W_{10B}:={\rm Spec}\,\mathbb Z\Bigl[\eta^1_{12},\cdots,\eta^1_{1s_3},\mathfrak h^1_{21},\cdots,\mathfrak h^1_{2s_3},\eta^2_{12},\cdots,\eta^2_{1s_4},\mathfrak h^2_{21},\cdots,\mathfrak h^2_{2s_4},\epsilon^-_1,z_2,\cdots,z_{s_2},\mathfrak v_3,$  $\left.\cdots,\mathfrak v_{s_2},y_2,\cdots,y_{s_1},\mathfrak u_3,\cdots,\mathfrak u_{s_1},\frac{1}{1-\eta_{12}^1\mathfrak h_{21}^2},\frac{1}{\mathfrak h_{21}^1},\frac{1}{1-y_2\left(\mathfrak{h}^2_{21}\left(\mathfrak{h}^1_{22}-\eta_{12}^1\mathfrak{h}^1_{21}\right)+\mathfrak{h}^1_{21}\right)}\right]$. The associated morphism $L:W_{10B}\longrightarrow G(2,n)$ is defined by (\ref{w11a}).
\smallskip

\item[($W_{10C}$)] $W_{10C}:={\rm Spec}\,\mathbb Z\Bigl[\eta^1_{12},\cdots,\eta^1_{1s_3},\mathfrak h^1_{21},\cdots,\mathfrak h^1_{2s_3},\eta^2_{12},\cdots,\eta^2_{1s_4},\mathfrak h^2_{21},\cdots,\mathfrak h^2_{2s_4},\epsilon^-_1,z_2,\cdots,z_{s_2},\mathfrak v_3,$  $\left.\cdots,\mathfrak v_{s_2},y_2,\cdots,y_{s_1},\mathfrak u_3,\cdots,\mathfrak u_{s_1},\frac{1}{1-\eta_{12}^1\mathfrak h_{21}^2},\frac{1}{1-y_2\mathfrak h_{21}^1},\frac{1}{\mathfrak{h}^2_{21}\left(\mathfrak{h}^1_{22}-\eta_{12}^1\mathfrak{h}^1_{21}\right)+\mathfrak{h}^1_{21}}\right]$. The associated morphism $L:W_{10C}\longrightarrow G(2,n)$ is defined by  (\ref{w11a}).
\smallskip

\item[($W_{10D}$)] $W_{10D}:={\rm Spec}\,\mathbb Z\Bigl[\eta^1_{12},\cdots,\eta^1_{1s_3},\mathfrak h^1_{21},\cdots,\mathfrak h^1_{2s_3},\eta^2_{12},\cdots,\eta^2_{1s_4},\mathfrak h^2_{21},\cdots,\mathfrak h^2_{2s_4},\epsilon^-_1,z_2,\cdots,z_{s_2},\mathfrak v_3,$  $\left.\cdots,\mathfrak v_{s_2},y_2,\cdots,y_{s_1},\mathfrak u_3,\cdots,\mathfrak u_{s_1},\frac{1}{1-\eta_{12}^1\mathfrak h_{21}^2},\frac{1}{1-y_2\mathfrak h_{21}^1},\frac{1}{1-y_2\left(\mathfrak{h}^2_{21}\left(\mathfrak{h}^1_{22}-\eta_{12}^1\mathfrak{h}^1_{21}\right)+\mathfrak{h}^1_{21}\right)}\right]$. The associated morphism $L:W_{10D}\longrightarrow G(2,n)$ is defined by (\ref{w11a}).
\smallskip

\item[($W_{10E}$)] $W_{10E}:={\rm Spec}\,\mathbb Z\Bigl[\eta^1_{12},\cdots,\eta^1_{1s_3},\mathfrak h^1_{21},\cdots,\mathfrak h^1_{2s_3},\eta^2_{12},\cdots,\eta^2_{1s_4},\mathfrak h^2_{21},\cdots,\mathfrak h^2_{2s_4},\epsilon^-_1,z_2,\cdots,z_{s_2},\mathfrak v_3,$  $\left.\cdots,\mathfrak v_{s_2},y_2,\cdots,y_{s_1},\mathfrak u_3,\cdots,\mathfrak u_{s_1},\frac{1}{\mathfrak h_{21}^2},\frac{1}{\mathfrak h_{21}^1},\frac{1}{\mathfrak{h}^2_{21}\left(\mathfrak{h}^1_{22}-\eta_{12}^1\mathfrak{h}^1_{21}\right)+\mathfrak{h}^1_{21}}\right]$. The associated morphism $L:W_{10E}\longrightarrow G(2,n)$ is defined by  (\ref{w11a}).
\smallskip

\item[($W_{10F}$)] $W_{10F}:={\rm Spec}\,\mathbb Z\Bigl[\eta^1_{12},\cdots,\eta^1_{1s_3},\mathfrak h^1_{21},\cdots,\mathfrak h^1_{2s_3},\eta^2_{12},\cdots,\eta^2_{1s_4},\mathfrak h^2_{21},\cdots,\mathfrak h^2_{2s_4},\epsilon^-_1,z_2,\cdots,z_{s_2},\mathfrak v_3,$  $\left.\cdots,\mathfrak v_{s_2},y_2,\cdots,y_{s_1},\mathfrak u_3,\cdots,\mathfrak u_{s_1},\frac{1}{\mathfrak h_{21}^2},\frac{1}{\mathfrak h_{21}^1},\frac{1}{1-y_2\left(\mathfrak{h}^2_{21}\left(\mathfrak{h}^1_{22}-\eta_{12}^1\mathfrak{h}^1_{21}\right)+\mathfrak{h}^1_{21}\right)}\right]$. The associated morphism $L:W_{10F}\longrightarrow G(2,n)$ is defined by  (\ref{w11a}).
\smallskip

\item[($W_{10G}$)] $W_{10G}:={\rm Spec}\,\mathbb Z\Bigl[\eta^1_{12},\cdots,\eta^1_{1s_3},\mathfrak h^1_{21},\cdots,\mathfrak h^1_{2s_3},\eta^2_{12},\cdots,\eta^2_{1s_4},\mathfrak h^2_{21},\cdots,\mathfrak h^2_{2s_4},\epsilon^-_1,z_2,\cdots,z_{s_2},\mathfrak v_3,$  $\left.\cdots,\mathfrak v_{s_2},y_2,\cdots,y_{s_1},\mathfrak u_3,\cdots,\mathfrak u_{s_1},\frac{1}{\mathfrak h_{21}^2},\frac{1}{1-y_2\mathfrak h_{21}^1},\frac{1}{\mathfrak{h}^2_{21}\left(\mathfrak{h}^1_{22}-\eta_{12}^1\mathfrak{h}^1_{21}\right)+\mathfrak{h}^1_{21}}\right]$. The associated morphism $L:W_{10G}\longrightarrow G(2,n)$ is defined by (\ref{w11a}).
\smallskip

\item[($W_{10H}$)] $W_{10H}:={\rm Spec}\,\mathbb Z\Bigl[\eta^1_{12},\cdots,\eta^1_{1s_3},\mathfrak h^1_{21},\cdots,\mathfrak h^1_{2s_3},\eta^2_{12},\cdots,\eta^2_{1s_4},\mathfrak h^2_{21},\cdots,\mathfrak h^2_{2s_4},\epsilon^-_1,z_2,\cdots,z_{s_2},\mathfrak v_3,$  $\left.\cdots,\mathfrak v_{s_2},y_2,\cdots,y_{s_1},\mathfrak u_3,\cdots,\mathfrak u_{s_1},\frac{1}{\mathfrak h_{21}^2},\frac{1}{1-y_2\mathfrak h_{21}^1},\frac{1}{1-y_2\left(\mathfrak{h}^2_{21}\left(\mathfrak{h}^1_{22}-\eta_{12}^1\mathfrak{h}^1_{21}\right)+\mathfrak{h}^1_{21}\right)}\right]$. The associated morphism $L:W_{10H}\longrightarrow G(2,n)$ is defined by  (\ref{w11a}).
\smallskip

\item[($W_{11A}$)] 
$W_{11A}:={\rm Spec}\,\mathbb Z\Bigl[\eta^1_{12},\cdots,\eta^1_{1s_3},\mathfrak h^1_{21},\cdots,\mathfrak h^1_{2s_3},\eta^2_{12},\cdots,\eta^2_{1s_4},\mathfrak h^2_{21},\cdots,\mathfrak h^2_{2s_4},\epsilon^-_1,z_2,\cdots,z_{s_2},\mathfrak v_3,$  $\left.\cdots,\mathfrak v_{s_2},y_2,\cdots,y_{s_1},\mathfrak u_3,\cdots,\mathfrak u_{s_1},\frac{1}{\mathfrak h_{21}^1},\frac{1}{\mathfrak h_{21}^2}\right]$. The associated morphism $L:W_{11A}\longrightarrow G(2,n)$ is defined by
\begin{equation}\label{w12a}
\begin{split}
&\left(\begin{matrix}1&y_2&y_3&\cdots\\
0&1&\mathfrak u_3&\cdots\\
\end{matrix}\hspace{-0.12in}\begin{matrix} &\hfill\tikzmark{a12}\\&\hfill\tikzmark{b12}
\end{matrix}\,\,\,\begin{matrix}0&\epsilon^-_1&\epsilon^-_1\mathfrak v_3&\cdots\\
1&z_2&z_3&\cdots\\
\end{matrix}\hspace{-0.12in}\begin{matrix} &\hfill\tikzmark{g12}\\&\hfill\tikzmark{h12}\end{matrix}\,\,\,\begin{matrix}1&\eta_{12}^1&\eta_{13}^1\\
\mathfrak{h}^1_{21}&\mathfrak{h}^1_{22}\left(\mathfrak{h}^2_{21}-\mathfrak{h}^1_{21}\right)+\mathfrak{h}^1_{21}\eta_{12}^1&\mathfrak{h}^1_{22}\mathfrak{h}^1_{23}\left(\mathfrak{h}^2_{21}-\mathfrak{h}^1_{21}\right)+\mathfrak{h}^1_{21}\eta_{13}^1\\
\end{matrix}\right.\\
&\left.\,\,\,\,\,\begin{matrix} \cdots\\
\cdots\\\end{matrix}\hspace{-0.12in}\begin{matrix} &\hfill\tikzmark{c12}\\&\hfill\tikzmark{d12}\end{matrix}\,\,\,\begin{matrix}1&\eta_{12}^2&\eta_{13}^2&\cdots\\
\mathfrak{h}^2_{21}&\left(\mathfrak{h}^2_{22}+\eta_{12}^2\right)\left(\mathfrak{h}^2_{21}-\mathfrak{h}^1_{21}\right)+\eta_{12}^2\mathfrak{h}^1_{21}&\left(\mathfrak{h}^2_{22}\mathfrak{h}^2_{23}+\eta_{13}^2\right)\left(\mathfrak{h}^2_{21}-\mathfrak{h}^1_{21}\right)+\eta_{13}^2\mathfrak{h}^1_{21}&\cdots\\
\end{matrix}\right)=:\Theta.
\tikz[remember picture,overlay]   \draw[dashed,dash pattern={on 4pt off 2pt}] ([xshift=0.5\tabcolsep,yshift=7pt]a12.north) -- ([xshift=0.5\tabcolsep,yshift=-2pt]b12.south);\tikz[remember picture,overlay]   \draw[dashed,dash pattern={on 4pt off 2pt}] ([xshift=0.5\tabcolsep,yshift=7pt]c12.north) -- ([xshift=0.5\tabcolsep,yshift=-2pt]d12.south);\tikz[remember picture,overlay]   \draw[dashed,dash pattern={on 4pt off 2pt}] ([xshift=0.5\tabcolsep,yshift=7pt]g12.north) -- ([xshift=0.5\tabcolsep,yshift=-2pt]h12.south);
\end{split}    
\end{equation}

\item[($W_{11B}$)] $W_{11B}:={\rm Spec}\,\mathbb Z\Bigl[\eta^1_{12},\cdots,\eta^1_{1s_3},\mathfrak h^1_{21},\cdots,\mathfrak h^1_{2s_3},\eta^2_{12},\cdots,\eta^2_{1s_4},\mathfrak h^2_{21},\cdots,\mathfrak h^2_{2s_4},\epsilon^-_1,z_2,\cdots,z_{s_2},\mathfrak v_3,$  $\left.\cdots,\mathfrak v_{s_2},y_2,\cdots,y_{s_1},\mathfrak u_3,\cdots,\mathfrak u_{s_1},\frac{1}{\mathfrak h_{21}^1},\frac{1}{1-y_2\mathfrak h_{21}^2}\right]$. The associated morphism $L:W_{11B}\longrightarrow G(2,n)$ is defined by  (\ref{w12a}).
\smallskip

\item[($W_{11C}$)] $W_{11C}:={\rm Spec}\,\mathbb Z\Bigl[\eta^1_{12},\cdots,\eta^1_{1s_3},\mathfrak h^1_{21},\cdots,\mathfrak h^1_{2s_3},\eta^2_{12},\cdots,\eta^2_{1s_4},\mathfrak h^2_{21},\cdots,\mathfrak h^2_{2s_4},\epsilon^-_1,z_2,\cdots,z_{s_2},\mathfrak v_3,$  $\left.\cdots,\mathfrak v_{s_2},y_2,\cdots,y_{s_1},\mathfrak u_3,\cdots,\mathfrak u_{s_1},\frac{1}{1-y_2\mathfrak h_{21}^1},\frac{1}{\mathfrak h_{21}^2}\right]$. The associated morphism $L:W_{11C}\longrightarrow G(2,n)$ is defined by  (\ref{w12a}).
\smallskip

\item[($W_{11D}$)] $W_{11D}:={\rm Spec}\,\mathbb Z\Bigl[\eta^1_{12},\cdots,\eta^1_{1s_3},\mathfrak h^1_{21},\cdots,\mathfrak h^1_{2s_3},\eta^2_{12},\cdots,\eta^2_{1s_4},\mathfrak h^2_{21},\cdots,\mathfrak h^2_{2s_4},\epsilon^-_1,z_2,\cdots,z_{s_2},\mathfrak v_3,$  $\left.\cdots,\mathfrak v_{s_2},y_2,\cdots,y_{s_1},\mathfrak u_3,\cdots,\mathfrak u_{s_1},\frac{1}{1-y_2\mathfrak h_{21}^1},\frac{1}{1-y_2\mathfrak h_{21}^2}\right]$. The associated morphism $L:W_{11D}$ $\longrightarrow G(2,n)$ is defined by  (\ref{w12a}).

\end{enumerate}

\section{Construction of local embeddings}\label{le4}

For convenience, in the following we will denote the regular functions $P_{(i_1,i_2)}\circ e_{1(s_1+1)}\circ L$ by $\mathcal P^{\Theta}_{(i_1,i_2)}$. Here the embedding $e_{1(s_1+1)}:U_{1(s_1+1)}\hookrightarrow G(2,n)$ is given by (\ref{ej1j2}), $L$ are the associated morphisms defined in Appendix \ref{ccm4}, $\Theta$ are the $2\times n$ matrices defining $L$.   Notice that,  $P_{(i_1,i_2)}\circ e_{1(s_1+1)}\circ L$ can be given by the determinants of $2\times 2$ submatrices of $\Theta$ consisting of the $i_1^{\rm th}$ and $i_2^{\rm th}$ columns.

Recall that $\tau:\{1,2,\cdots,n\}\longrightarrow\{1,2,3,4\}$ is the map defined by (\ref{intau}) such that $1+\sum_{t=1}^{\tau(i)-1}s_t\leq i\leq\sum_{t=1}^{\tau(i)}s_t$ for $1\leq i\leq n$.

\subsection{{\texorpdfstring{$s_1=s_2=s_3=1$ and $s_4\geq2$}{dd}}}\label{lx1}

The locally closed embeddings $f_{1B}:X_{1B}\times\mathbb G_m^{\mathbb V^{\underline s}}/\mathbb G_m\rightarrow\mathcal A_{13}\times \mathbb G_m{\Big\backslash\prod\nolimits_{\underline i\in \mathbb V^{\underline s}}\left(\wedge^{\underline i}E_{\bullet}-\{0\}\right)}$ and
$f_{4}:X_{4}\times\mathbb G_m^{\mathbb V^{\underline s}}/\mathbb G_m\rightarrow\mathcal A_{23}\times \mathbb G_m{\Big\backslash\prod\nolimits_{\underline i\in \mathbb V^{\underline s}}\left(\wedge^{\underline i}E_{\bullet}-\{0\}\right)}$ are induced by the ring homomorphisms
\begin{equation*}
\frac{\chi_{i_1j_1}\chi_{i_2j_2}}{\chi_{i_3j_3}\chi_{i_4j_4}}\mapsto \frac{\mathcal P^{\Theta}_{(i_1,j_1)}\mathcal P^{\Theta}_{(i_2,j_2)}}{\mathcal P^{\Theta}_{(i_3,j_3)}\mathcal P^{\Theta}_{(i_4,j_4)}}\frac{x_{i_3j_3}x_{i_4j_4}}{x_{i_1j_1}x_{i_2j_2}},\,\,\,\,\,\, z_{(i_1,i_2)}\mapsto\frac{\mathcal P^{\Theta}_{(i_1,i_2)}}{\mathcal P^{\Theta}_{(\tau(i_1),\tau(i_2))}}x_{\tau(i_1)\tau(i_2)},  
\end{equation*}
where the defining matrices $\Theta$ are given by (\ref{x1a}) and (\ref{x4}), respectively.

The locally closed embedding $f_{3}:X_{3}\times\mathbb G_m^{\mathbb V^{\underline s}}/\mathbb G_m\rightarrow\mathcal A_{13}\times \mathbb G_m{\Big\backslash\prod\nolimits_{\underline i\in \mathbb V^{\underline s}}\left(\wedge^{\underline i}E_{\bullet}-\{0\}\right)}$ is induced by the ring homomorhpism
\begin{equation}\label{yeben}
\frac{\chi_{i_1j_1}\chi_{i_2j_2}}{\chi_{i_3j_3}\chi_{i_4j_4}}\mapsto\frac{g_{i_1j_1}g_{i_2j_2}}{g_{i_3j_3}g_{i_4j_4}}\frac{x_{i_3j_3}x_{i_4j_4}}{x_{i_1j_1}x_{i_2j_2}},\,\,\,\,\,\, z_{(i_1,i_2)}\mapsto\frac{\mathcal P^{\Theta}_{(i_1,i_2)}}{g_{\tau(i_1)\tau(i_2)}}x_{\tau(i_1)\tau(i_2)},  
\end{equation} where 
\begin{equation*}
\begin{split}
&g_{12}:=\mathcal P^{\Theta}_{(1,2)},\,\,g_{13}:=\mathcal P^{\Theta}_{(1,3)},\,\,g_{14}:=\mathcal P^{\Theta}_{(1,5)},\,\,g_{23}:=\mathcal P^{\Theta}_{(2,3)},\\
&g_{24}:=\mathcal P^{\Theta}_{(2,4)},\,\,g_{34}:=\mathcal P^{\Theta}_{(3,5)},\,\,g_{44}:=\mathcal P^{\Theta}_{(4,5)},\\
\end{split}   
\end{equation*} 
with $\Theta$ given by (\ref{x3a}).


The locally closed embedding $f_{5}:X_{5}\times\mathbb G_m^{\mathbb V^{\underline s}}/\mathbb G_m\rightarrow\mathcal A_{23}\times \mathbb G_m{\Big\backslash\prod\nolimits_{\underline i\in \mathbb V^{\underline s}}\left(\wedge^{\underline i}E_{\bullet}-\{0\}\right)}$ is induced by the ring homomorphism (\ref{yeben}),
where 
\begin{equation*}
\begin{split}
&g_{12}:=\mathcal P^{\Theta}_{(1,2)},\,\,g_{13}:=\mathcal P^{\Theta}_{(1,3)},\,\,g_{14}:=\mathcal P^{\Theta}_{(1,5)},\,\,g_{23}:=\mathcal P^{\Theta}_{(2,3)},\\
&g_{24}:=\mathcal P^{\Theta}_{(2,4)},\,\,g_{34}:=\mathcal P^{\Theta}_{(3,4)},\,\,g_{44}:=\mathcal P^{\Theta}_{(4,5)},\\
\end{split}   
\end{equation*} 
with $\Theta$ given by (\ref{x5}).

\subsection{{\texorpdfstring{$s_1=s_2=1$ and $s_3,s_4\geq2$}{dd}}}\label{ly1}
Denote $t_{33}:=(0,0,2,0)$. Define functions  $l_{12}$, $l_{13}$, $l_{23}$, $l_{3}$, $l_{4}$ on $\mathbb V^{\underline s}$ as follows.
\begin{equation*}
\begin{split}
&l_{12}(t_{12})=l_{12}(t_{13})=l_{12}(t_{23})=l_{12}(t_{34})=0,\,\, \,l_{12}(t_{14})=l_{12}(t_{24})=l_{12}(t_{44})=-1,\,\, \,l_{12}(t_{33})=1;\\
&l_{13}(t_{12})=l_{13}(t_{13})=l_{13}(t_{14})=l_{13}(t_{23})=l_{13}(t_{34})=l_{13}(t_{33})=0,\,\, \,l_{13}(t_{24})=l_{13}(t_{44})=1;\\
&l_{23}(t_{12})=l_{23}(t_{13})=l_{23}(t_{23})=l_{23}(t_{24})=l_{23}(t_{34})=l_{23}(t_{33})=0,\,\, \,l_{23}(t_{14})=l_{23}(t_{44})=1;\\
&l_{3}(t_{12})=l_{3}(t_{13})=l_{3}(t_{23})=l_{3}(t_{34})=l_{3}(t_{14})=l_{3}(t_{24})=l_{3}(t_{44})=0,\,\, \,l_{3}(t_{33})=1;\\
&l_{4}(t_{12})=l_{4}(t_{13})=l_{4}(t_{23})=l_{4}(t_{14})=l_{4}(t_{24})=l_{4}(t_{34})=l_{4}(t_{33})=0,\,\, \,l_{4}(t_{44})=1.\\
\end{split}
\end{equation*}
Let $\mathcal C_{12}$ be the convex cone generated by $\mathcal C_{\emptyset}^{\mathbb V^{\widetilde{\underline S}}}$ and the functions $l_{12}$, $l_{3}$, $l_{4}$; let $\mathcal C_{13}$ be the convex cone generated by $\mathcal C_{\emptyset}^{\mathbb V^{\widetilde{\underline S}}}$ and the functions $l_{13}$, $l_{3}$, $l_{4}$; let $\mathcal C_{23}$ be the convex cone generated by $\mathcal C_{\emptyset}^{\mathbb V^{\widetilde{\underline S}}}$ and the functions $l_{23}$, $l_{3}$, $l_{4}$. One can show that $\mathcal A^{\mathbb V^{\underline s}}$ is covered by affine open subschemes \begin{equation*}
\mathcal A_{12}:={\rm Spec}\,\mathbb Z\left[\left(\overline{\mathcal C_{12}}\right)^*\cap M\right],\,\mathcal A_{13}:={\rm Spec}\,\mathbb Z\left[\left(\overline{\mathcal C_{13}}\right)^*\cap M\right],\,\mathcal A_{23}:={\rm Spec}\,\mathbb Z\left[\left(\overline{\mathcal C_{23}}\right)^*\cap M\right].   
\end{equation*} 

Write the coordinate ring of the torus $\mathbb G^{\mathbb V^{\underline s}}_m$ as
\begin{equation*}
\mathbb Z\left[x_{12},\frac{1}{x_{12}},x_{13},\frac{1}{x_{13}},x_{23},\frac{1}{x_{23}},x_{14},\frac{1}{x_{14}},x_{24},\frac{1}{x_{24}},x_{34},\frac{1}{x_{34}},x_{33},\frac{1}{x_{33}},x_{44},\frac{1}{x_{44}}\right].   
\end{equation*}
Define characters $\chi_{i_1i_2}$ of $\mathbb G^{\mathbb V^{\underline s}}_m$ by
$\chi_{i_1i_2}(x_{12},x_{13},x_{23},x_{14},x_{24},x_{34},x_{33},x_{44})=x_{i_1i_2}.$
Then
\begin{equation*}
\begin{split}
&\mathcal A_{12}:={\rm Spec}\,\mathbb Z\left[\frac{\chi_{14}\chi_{23}}{\chi_{13}\chi_{24}},\frac{\chi_{13}\chi_{24}}{\chi_{14}\chi_{23}},\frac{\chi_{12}\chi_{34}}{\chi_{23}\chi_{14}},\frac{\chi_{14}\chi_{33}}{\chi_{13}\chi_{34}},\frac{\chi_{13}\chi_{44}}{\chi_{14}\chi_{34}}\right],\\
&\mathcal A_{13}:={\rm Spec}\,\mathbb Z\left[\frac{\chi_{23}\chi_{14}}{\chi_{12}\chi_{34}},\frac{\chi_{12}\chi_{34}}{\chi_{23}\chi_{14}},\frac{\chi_{12}\chi_{33}}{\chi_{13}\chi_{23}},\frac{\chi_{13}\chi_{24}}{\chi_{12}\chi_{34}},\frac{\chi_{23}\chi_{44}}{\chi_{24}\chi_{34}}\right],\\
&\mathcal A_{23}:={\rm Spec}\,\mathbb Z\left[\frac{\chi_{12}\chi_{34}}{\chi_{13}\chi_{24}},\frac{\chi_{13}\chi_{24}}{\chi_{12}\chi_{34}},\frac{\chi_{23}\chi_{14}}{\chi_{12}\chi_{34}},\frac{\chi_{12}\chi_{33}}{\chi_{13}\chi_{23}},\frac{\chi_{13}\chi_{44}}{\chi_{14}\chi_{34}}\right].      
\end{split}  
\end{equation*} 
We view $\mathbb G_m^{\mathbb V^{\underline s}}/\mathbb G_m$ as a subscheme of ${\rm Proj}\,\mathbb Z\left[x_{12},x_{13},x_{23},x_{14},x_{24},x_{34},x_{33},x_{44}\right]$.

The locally closed embeddings $f_{1A}:Y_{1A}\times\mathbb G_m^{\mathbb V^{\underline s}}/\mathbb G_m\rightarrow\mathcal A_{23}\times \mathbb G_m{\Big\backslash\prod\nolimits_{\underline i\in \mathbb V^{\underline s}}\left(\wedge^{\underline i}E_{\bullet}-\{0\}\right)}$, $f_{1B}:Y_{1B}\times\mathbb G_m^{\mathbb V^{\underline s}}/\mathbb G_m\rightarrow\mathcal A_{12}\times \mathbb G_m{\Big\backslash\prod\nolimits_{\underline i\in \mathbb V^{\underline s}}\left(\wedge^{\underline i}E_{\bullet}-\{0\}\right)}$, $f_{2A}:Y_{2A}\times\mathbb G_m^{\mathbb V^{\underline s}}/\mathbb G_m\rightarrow\mathcal A_{12}\times \mathbb G_m{\Big\backslash\prod\nolimits_{\underline i\in \mathbb V^{\underline s}}\left(\wedge^{\underline i}E_{\bullet}-\{0\}\right)}$ are induced
by the ring homomorphisms 
\begin{equation}\label{yeben2}
\frac{\chi_{i_1j_1}\chi_{i_2j_2}}{\chi_{i_3j_3}\chi_{i_4j_4}}\mapsto\frac{g_{i_1j_1}g_{i_2j_2}}{g_{i_3j_3}g_{i_4j_4}}\frac{x_{i_3j_3}x_{i_4j_4}}{x_{i_1j_1}x_{i_2j_2}},\,\,\,\,\,\, z_{(i_1,i_2)}\mapsto\frac{\mathcal P^{\Theta}_{(i_1,i_2)}}{g_{\tau(i_1)\tau(i_2)}}x_{\tau(i_1)\tau(i_2)},  
\end{equation}
where 
\begin{equation*}
\begin{split}
&g_{12}:=\mathcal P^{\Theta}_{(1,2)},\,\,g_{13}:=\mathcal P^{\Theta}_{(1,3)},\,\,g_{14}:=\mathcal P^{\Theta}_{(1,3+s_3)},\,\,g_{23}:=\mathcal P^{\Theta}_{(2,3)},\\
&g_{24}:=\mathcal P^{\Theta}_{(2,3+s_3)},\,\,g_{33}:=\mathcal P^{\Theta}_{(3,4)},\,\,g_{34}:=\mathcal P^{\Theta}_{(3,3+s_3)},\,\,g_{44}:=\mathcal P^{\Theta}_{(3+s_3,4+s_3)},\\
\end{split}   
\end{equation*}
with the defining matrices $\Theta$ given by (\ref{didi}), (\ref{didi}), (\ref{y3a}), respectively.

The locally closed embedding $f_{2B}:Y_{2B}\times\mathbb G_m^{\mathbb V^{\underline s}}/\mathbb G_m\rightarrow\mathcal A_{12}\times \mathbb G_m{\Big\backslash\prod\nolimits_{\underline i\in \mathbb V^{\underline s}}\left(\wedge^{\underline i}E_{\bullet}-\{0\}\right)}$ is induced by the ring homomorphism (\ref{yeben2}), where 
\begin{equation*}
\begin{split}
&g_{12}:=\mathcal P^{\Theta}_{(1,2)},\,\,g_{13}:=\mathcal P^{\Theta}_{(1,3)},\,\,g_{14}:=\mathcal P^{\Theta}_{(1,3+s_3)},\,\,g_{23}:=\mathcal P^{\Theta}_{(2,3)},\\
&g_{24}:=\mathcal P^{\Theta}_{(2,3+s_3)},\,\,g_{33}:=\mathcal P^{\Theta}_{(3,4)},\,\,g_{34}:=\mathcal P^{\Theta}_{(4,3+s_3)},\,\,g_{44}:=\mathcal P^{\Theta}_{(3+s_3,4+s_3)},\\
\end{split}   
\end{equation*}  
with $\Theta$ given by (\ref{y3a}).

The locally closed embedding $f_{3}:Y_{3}\times\mathbb G_m^{\mathbb V^{\underline s}}/\mathbb G_m\rightarrow\mathcal A_{12}\times \mathbb G_m{\Big\backslash\prod\nolimits_{\underline i\in \mathbb V^{\underline s}}\left(\wedge^{\underline i}E_{\bullet}-\{0\}\right)}$ is induced by the ring homomorphism (\ref{yeben2}), where  
\begin{equation*}
\begin{split}
&g_{12}:=\mathcal P^{\Theta}_{(1,2)},\,\,g_{13}:=\mathcal P^{\Theta}_{(1,3)},\,\,g_{14}:=\mathcal P^{\Theta}_{(1,3+s_3)},\,\,g_{23}:=\mathcal P^{\Theta}_{(2,3)},\\
&g_{24}:=\mathcal P^{\Theta}_{(2,3+s_3)},\,\,g_{33}:=\mathcal P^{\Theta}_{(3,4)},\,\,g_{34}:=\mathcal P^{\Theta}_{(3,4+s_3)},\,\,g_{44}:=\mathcal P^{\Theta}_{(3+s_3,4+s_3)},\\
\end{split}   
\end{equation*}  
with $\Theta$ given by (\ref{y3}).

The locally closed embedding $f_{4}:Y_{4}\times\mathbb G_m^{\mathbb V^{\underline s}}/\mathbb G_m\rightarrow\mathcal A_{23}\times \mathbb G_m{\Big\backslash\prod\nolimits_{\underline i\in \mathbb V^{\underline s}}\left(\wedge^{\underline i}E_{\bullet}-\{0\}\right)}$ is induced by the ring homomorphism (\ref{yeben2}), where  
\begin{equation*}
\begin{split}
&g_{12}:=\mathcal P^{\Theta}_{(1,2)},\,\,g_{13}:=\mathcal P^{\Theta}_{(1,4)},\,\,g_{14}:=\mathcal P^{\Theta}_{(1,3+s_3)},\,\,g_{23}:=\mathcal P^{\Theta}_{(2,3)},\\
&g_{24}:=\mathcal P^{\Theta}_{(2,3+s_3)},\,\,g_{33}:=\mathcal P^{\Theta}_{(3,4)},\,\,g_{34}:=\mathcal P^{\Theta}_{(4,3+s_3)},\,\,g_{44}:=\mathcal P^{\Theta}_{(3+s_3,4+s_3)},\\
\end{split}   
\end{equation*}  
with $\Theta$ given by (\ref{y4}).

The locally closed embedding $f_{5}:Y_{5}\times\mathbb G_m^{\mathbb V^{\underline s}}/\mathbb G_m\rightarrow\mathcal A_{23}\times \mathbb G_m{\Big\backslash\prod\nolimits_{\underline i\in \mathbb V^{\underline s}}\left(\wedge^{\underline i}E_{\bullet}-\{0\}\right)}$ is induced by the ring homomorphism (\ref{yeben2}), where  
\begin{equation*}
\begin{split}
&g_{12}:=\mathcal P^{\Theta}_{(1,2)},\,\,g_{13}:=\mathcal P^{\Theta}_{(1,3)},\,\,g_{14}:=\mathcal P^{\Theta}_{(1,4+s_3)},\,\,g_{23}:=\mathcal P^{\Theta}_{(2,3)},\\
&g_{24}:=\mathcal P^{\Theta}_{(2,3+s_3)},\,\,g_{33}:=\mathcal P^{\Theta}_{(3,4)},\,\,g_{34}:=\mathcal P^{\Theta}_{(3,3+s_3)},\,\,g_{44}:=\mathcal P^{\Theta}_{(3+s_3,4+s_3)},\\
\end{split}   
\end{equation*} 
with $\Theta$ given by (\ref{y5}).


The locally closed embedding $f_{6}:Y_{6}\times\mathbb G_m^{\mathbb V^{\underline s}}/\mathbb G_m\rightarrow\mathcal A_{23}\times \mathbb G_m{\Big\backslash\prod\nolimits_{\underline i\in \mathbb V^{\underline s}}\left(\wedge^{\underline i}E_{\bullet}-\{0\}\right)}$ is induced by the ring homomorphism (\ref{yeben2}), where  
\begin{equation*}
\begin{split}
&g_{12}:=\mathcal P^{\Theta}_{(1,2)},\,\,g_{13}:=\mathcal P^{\Theta}_{(1,4)},\,\,g_{14}:=\mathcal P^{\Theta}_{(1,4+s_3)},\,\,g_{23}:=\mathcal P^{\Theta}_{(2,3)},\\
&g_{24}:=\mathcal P^{\Theta}_{(2,3+s_3)},\,\,g_{33}:=\mathcal P^{\Theta}_{(3,4)},\,\,g_{34}:=\mathcal P^{\Theta}_{(4,3+s_3)},\,\,g_{44}:=\mathcal P^{\Theta}_{(3+s_3,4+s_3)},\\
\end{split}   
\end{equation*} 
with $\Theta$ given by (\ref{y6}).


\subsection{{\texorpdfstring{$s_1=1$ and $s_2,s_3,s_4\geq2$}{dd}}}\label{lz1}
Write the coordinate ring of $\mathbb G^{\mathbb V^{\underline s}}_m$ as
\begin{equation*}
\mathbb Z\left[x_{12},\frac{1}{x_{12}},x_{13},\frac{1}{x_{13}},x_{23},\frac{1}{x_{23}},x_{14},\frac{1}{x_{14}},x_{24},\frac{1}{x_{24}},x_{34},\frac{1}{x_{34}},x_{22},\frac{1}{x_{22}},x_{33},\frac{1}{x_{33}},x_{44},\frac{1}{x_{44}}\right].   
\end{equation*}
Define characters $\chi_{i_1i_2}$ of $\mathbb G^{\mathbb V^{\underline s}}_m$ by
$\chi_{i_1i_2}(x_{12},x_{13},x_{23},x_{14},x_{24},x_{34},x_{22},x_{33},x_{44})=x_{i_1i_2}.$ One can show that $\mathcal A^{\mathbb V^{\underline s}}$ is covered by affine open subschemes 
\begin{equation*}
\begin{split}
&\mathcal A_{12}:={\rm Spec}\,\mathbb Z\left[\frac{\chi_{14}\chi_{23}}{\chi_{13}\chi_{24}},\frac{\chi_{13}\chi_{24}}{\chi_{14}\chi_{23}},\frac{\chi_{12}\chi_{34}}{\chi_{23}\chi_{14}},\frac{\chi_{14}\chi_{33}}{\chi_{13}\chi_{34}},\frac{\chi_{13}\chi_{44}}{\chi_{14}\chi_{34}},\frac{\chi_{13}\chi_{22}}{\chi_{12}\chi_{23}}\right],\\
&\mathcal A_{13}:={\rm Spec}\,\mathbb Z\left[\frac{\chi_{23}\chi_{14}}{\chi_{12}\chi_{34}},\frac{\chi_{12}\chi_{34}}{\chi_{23}\chi_{14}},\frac{\chi_{12}\chi_{33}}{\chi_{13}\chi_{23}},\frac{\chi_{13}\chi_{24}}{\chi_{12}\chi_{34}},\frac{\chi_{23}\chi_{44}}{\chi_{24}\chi_{34}},\frac{\chi_{22}\chi_{34}}{\chi_{23}\chi_{24}}\right],\\
&\mathcal A_{23}:={\rm Spec}\,\mathbb Z\left[\frac{\chi_{12}\chi_{34}}{\chi_{13}\chi_{24}},\frac{\chi_{13}\chi_{24}}{\chi_{12}\chi_{34}},\frac{\chi_{23}\chi_{14}}{\chi_{12}\chi_{34}},\frac{\chi_{12}\chi_{33}}{\chi_{13}\chi_{23}},\frac{\chi_{13}\chi_{44}}{\chi_{14}\chi_{34}},\frac{\chi_{13}\chi_{22}}{\chi_{12}\chi_{23}}\right].      
\end{split}  
\end{equation*} 
We view $\mathbb G_m^{\mathbb V^{\underline s}}/\mathbb G_m$ as a subscheme of ${\rm Proj}\,\mathbb Z\left[x_{12},x_{13},x_{23},x_{14},x_{24},x_{34},x_{22},x_{33},x_{44}\right]$.

The locally closed embeddings $f_{1A}:Z_{1A}\times\mathbb G_m^{\mathbb V^{\underline s}}/\mathbb G_m\rightarrow\mathcal A_{23}\times \mathbb G_m{\Big\backslash\prod\nolimits_{\underline i\in \mathbb V^{\underline s}}\left(\wedge^{\underline i}E_{\bullet}-\{0\}\right)}$, $f_{1B}:Z_{1B}\times\mathbb G_m^{\mathbb V^{\underline s}}/\mathbb G_m\rightarrow\mathcal A_{12}\times \mathbb G_m{\Big\backslash\prod\nolimits_{\underline i\in \mathbb V^{\underline s}}\left(\wedge^{\underline i}E_{\bullet}-\{0\}\right)}$, $f_{2A}:Z_{2A}\times\mathbb G_m^{\mathbb V^{\underline s}}/\mathbb G_m\rightarrow\mathcal A_{12}\times \mathbb G_m{\Big\backslash\prod\nolimits_{\underline i\in \mathbb V^{\underline s}}\left(\wedge^{\underline i}E_{\bullet}-\{0\}\right)}$  are induced
by the ring homomorphisms 
\begin{equation}\label{yeben3}
\frac{\chi_{i_1j_1}\chi_{i_2j_2}}{\chi_{i_3j_3}\chi_{i_4j_4}}\mapsto\frac{g_{i_1j_1}g_{i_2j_2}}{g_{i_3j_3}g_{i_4j_4}}\frac{x_{i_3j_3}x_{i_4j_4}}{x_{i_1j_1}x_{i_2j_2}},\,\,\,\,\,\, z_{(i_1,i_2)}\mapsto\frac{\mathcal P^{\Theta}_{(i_1,i_2)}}{g_{\tau(i_1)\tau(i_2)}}x_{\tau(i_1)\tau(i_2)},  
\end{equation}
where 
\begin{equation*}
\begin{split}
&g_{12}:=\mathcal P^{\Theta}_{(1,2)},\,\,g_{13}:=\mathcal P^{\Theta}_{(1,2+s_2)},\,\,g_{14}:=\mathcal P^{\Theta}_{(1,2+s_2+s_3)},\,\,g_{22}:=\mathcal P^{\Theta}_{(2,3)},\,\,g_{23}:=\mathcal P^{\Theta}_{(2,2+s_2)},\\
&g_{24}:=\mathcal P^{\Theta}_{(2,2+s_2+s_3)},\,\,g_{33}:=\mathcal P^{\Theta}_{(2+s_2,3+s_2)},\,\,g_{34}:=\mathcal P^{\Theta}_{(2+s_2,2+s_2+s_3)},\,\,g_{44}:=\mathcal P^{\Theta}_{(2+s_2+s_3,3+s_2+s_3)},\\
\end{split}   
\end{equation*}
with the defining matrices $\Theta$ given by (\ref{didi2}), (\ref{didi2}),  (\ref{z3a}), respectively.

The locally closed embedding $f_{2B}:Z_{2B}\times\mathbb G_m^{\mathbb V^{\underline s}}/\mathbb G_m\rightarrow\mathcal A_{12}\times \mathbb G_m{\Big\backslash\prod\nolimits_{\underline i\in \mathbb V^{\underline s}}\left(\wedge^{\underline i}E_{\bullet}-\{0\}\right)}$ is induced by the ring homomorphism (\ref{yeben3}), where 
\begin{equation*}
\begin{split}
&g_{12}:=\mathcal P^{\Theta}_{(1,2)},\,\,g_{13}:=\mathcal P^{\Theta}_{(1,2+s_2)},\,\,g_{14}:=\mathcal P^{\Theta}_{(1,2+s_2+s_3)},\,\,g_{22}:=\mathcal P^{\Theta}_{(2,3)},\,\,g_{23}:=\mathcal P^{\Theta}_{(2,2+s_2)},\\
&g_{24}:=\mathcal P^{\Theta}_{(2,2+s_2+s_3)},\,\,g_{33}:=\mathcal P^{\Theta}_{(2+s_2,3+s_2)},\,\,g_{34}:=\mathcal P^{\Theta}_{(3+s_2,2+s_2+s_3)},\,\,g_{44}:=\mathcal P^{\Theta}_{(2+s_2+s_3,3+s_2+s_3)},\\
\end{split}   
\end{equation*}
with $\Theta$ given by (\ref{z3a}).


The locally closed embedding $f_{3}:Z_{3}\times\mathbb G_m^{\mathbb V^{\underline s}}/\mathbb G_m\rightarrow\mathcal A_{12}\times \mathbb G_m{\Big\backslash\prod\nolimits_{\underline i\in \mathbb V^{\underline s}}\left(\wedge^{\underline i}E_{\bullet}-\{0\}\right)}$ is induced by the ring homomorphism (\ref{yeben3}), where 
\begin{equation*}
\begin{split}
&g_{12}:=\mathcal P^{\Theta}_{(1,2)},\,\,g_{13}:=\mathcal P^{\Theta}_{(1,2+s_2)},\,\,g_{14}:=\mathcal P^{\Theta}_{(1,2+s_2+s_3)},\,\,g_{22}:=\mathcal P^{\Theta}_{(2,3)},\,\,g_{23}:=\mathcal P^{\Theta}_{(2,2+s_2)},\\
&g_{24}:=\mathcal P^{\Theta}_{(2,2+s_2+s_3)},\,\,g_{33}:=\mathcal P^{\Theta}_{(2+s_2,3+s_2)},\,\,g_{34}:=\mathcal P^{\Theta}_{(2+s_2,3+s_2+s_3)},\,\,g_{44}:=\mathcal P^{\Theta}_{(2+s_2+s_3,3+s_2+s_3)},\\
\end{split}   
\end{equation*}
with $\Theta$ given by (\ref{z3}).


The locally closed embedding $f_{4}:Z_{4}\times\mathbb G_m^{\mathbb V^{\underline s}}/\mathbb G_m\rightarrow\mathcal A_{23}\times \mathbb G_m{\Big\backslash\prod\nolimits_{\underline i\in \mathbb V^{\underline s}}\left(\wedge^{\underline i}E_{\bullet}-\{0\}\right)}$ is induced by the ring homomorphism (\ref{yeben3}), where 
\begin{equation*}
\begin{split}
&g_{12}:=\mathcal P^{\Theta}_{(1,2)},\,\,g_{13}:=\mathcal P^{\Theta}_{(1,3+s_2)},\,\,g_{14}:=\mathcal P^{\Theta}_{(1,2+s_2+s_3)},\,\,g_{22}:=\mathcal P^{\Theta}_{(2,3)},\,\,g_{23}:=\mathcal P^{\Theta}_{(2,2+s_2)},\\
&g_{24}:=\mathcal P^{\Theta}_{(2,2+s_2+s_3)},\,\,g_{33}:=\mathcal P^{\Theta}_{(2+s_2,3+s_2)},\,\,g_{34}:=\mathcal P^{\Theta}_{(3+s_2,2+s_2+s_3)},\,\,g_{44}:=\mathcal P^{\Theta}_{(2+s_2+s_3,3+s_2+s_3)},\\
\end{split}   
\end{equation*}
with $\Theta$ given by (\ref{z4}).


The locally closed embeddings $f_{5}:Z_{5}\times\mathbb G_m^{\mathbb V^{\underline s}}/\mathbb G_m\rightarrow\mathcal A_{23}\times \mathbb G_m{\Big\backslash\prod\nolimits_{\underline i\in \mathbb V^{\underline s}}\left(\wedge^{\underline i}E_{\bullet}-\{0\}\right)}$, $f_{7}:Z_{7}\times\mathbb G_m^{\mathbb V^{\underline s}}/\mathbb G_m\rightarrow\mathcal A_{23}\times \mathbb G_m{\Big\backslash\prod\nolimits_{\underline i\in \mathbb V^{\underline s}}\left(\wedge^{\underline i}E_{\bullet}-\{0\}\right)}$ are induced by the ring homomorphisms (\ref{yeben3}), where 
\begin{equation*}
\begin{split}
&g_{12}:=\mathcal P^{\Theta}_{(1,2)},\,\,g_{13}:=\mathcal P^{\Theta}_{(1,2+s_2)},\,\,g_{14}:=\mathcal P^{\Theta}_{(1,3+s_2+s_3)},\,\,g_{22}:=\mathcal P^{\Theta}_{(2,3)},\,\,g_{23}:=\mathcal P^{\Theta}_{(2,2+s_2)},\\
&g_{24}:=\mathcal P^{\Theta}_{(2,2+s_2+s_3)},\,\,g_{33}:=\mathcal P^{\Theta}_{(2+s_2,3+s_2)},\,\,g_{34}:=\mathcal P^{\Theta}_{(2+s_2,2+s_2+s_3)},\,\,g_{44}:=\mathcal P^{\Theta}_{(2+s_2+s_3,3+s_2+s_3)},\\
\end{split}   
\end{equation*}
with $\Theta$ given by (\ref{z5}), (\ref{z7}), respectively.


The locally closed embedding $f_{6}:Z_{6}\times\mathbb G_m^{\mathbb V^{\underline s}}/\mathbb G_m\rightarrow\mathcal A_{23}\times \mathbb G_m{\Big\backslash\prod\nolimits_{\underline i\in \mathbb V^{\underline s}}\left(\wedge^{\underline i}E_{\bullet}-\{0\}\right)}$ is induced by the ring homomorphism (\ref{yeben3}), where 
\begin{equation*}
\begin{split}
&g_{12}:=\mathcal P^{\Theta}_{(1,2)},\,\,g_{13}:=\mathcal P^{\Theta}_{(1,3+s_2)},\,\,g_{14}:=\mathcal P^{\Theta}_{(1,3+s_2+s_3)},\,\,g_{22}:=\mathcal P^{\Theta}_{(2,3)},\,\,g_{23}:=\mathcal P^{\Theta}_{(2,2+s_2)},\\
&g_{24}:=\mathcal P^{\Theta}_{(2,2+s_2+s_3)},\,\,g_{33}:=\mathcal P^{\Theta}_{(2+s_2,3+s_2)},\,\,g_{34}:=\mathcal P^{\Theta}_{(3+s_2,2+s_2+s_3)},\,\,g_{44}:=\mathcal P^{\Theta}_{(2+s_2+s_3,3+s_2+s_3)},\\
\end{split}   
\end{equation*}
with $\Theta$ given by (\ref{z6}).


The locally closed embedding $f_{8}:Z_{8}\times\mathbb G_m^{\mathbb V^{\underline s}}/\mathbb G_m\rightarrow\mathcal A_{23}\times \mathbb G_m{\Big\backslash\prod\nolimits_{\underline i\in \mathbb V^{\underline s}}\left(\wedge^{\underline i}E_{\bullet}-\{0\}\right)}$ is induced by the ring homomorphism (\ref{yeben3}), where 
\begin{equation*}
\begin{split}
&g_{12}:=\mathcal P^{\Theta}_{(1,2)},\,\,g_{13}:=\mathcal P^{\Theta}_{(1,2+s_2)},\,\,g_{14}:=\mathcal P^{\Theta}_{(1,3+s_2+s_3)},\,\,g_{22}:=\mathcal P^{\Theta}_{(2,3)},\,\,g_{23}:=\mathcal P^{\Theta}_{(2,3+s_2)},\\
&g_{24}:=\mathcal P^{\Theta}_{(2,2+s_2+s_3)},\,\,g_{33}:=\mathcal P^{\Theta}_{(2+s_2,3+s_2)},\,\,g_{34}:=\mathcal P^{\Theta}_{(2+s_2,2+s_2+s_3)},\,\,g_{44}:=\mathcal P^{\Theta}_{(2+s_2+s_3,3+s_2+s_3)},\\
\end{split}   
\end{equation*}
with $\Theta$ given by (\ref{z8}).


The locally closed embedding $f_{9}:Z_{9}\times\mathbb G_m^{\mathbb V^{\underline s}}/\mathbb G_m\rightarrow\mathcal A_{23}\times \mathbb G_m{\Big\backslash\prod\nolimits_{\underline i\in \mathbb V^{\underline s}}\left(\wedge^{\underline i}E_{\bullet}-\{0\}\right)}$ is induced by the ring homomorphism (\ref{yeben3}), where 
\begin{equation*}
\begin{split}
&g_{12}:=\mathcal P^{\Theta}_{(1,2)},\,\,g_{13}:=\mathcal P^{\Theta}_{(1,2+s_2)},\,\,g_{14}:=\mathcal P^{\Theta}_{(1,3+s_2+s_3)},\,\,g_{22}:=\mathcal P^{\Theta}_{(2,3)},\,\,g_{23}:=\mathcal P^{\Theta}_{(3,2+s_2)},\\
&g_{24}:=\mathcal P^{\Theta}_{(2,2+s_2+s_3)},\,\,g_{33}:=\mathcal P^{\Theta}_{(2+s_2,3+s_2)},\,\,g_{34}:=\mathcal P^{\Theta}_{(2+s_2,2+s_2+s_3)},\,\,g_{44}:=\mathcal P^{\Theta}_{(2+s_2+s_3,3+s_2+s_3)},\\
\end{split}   
\end{equation*}
with $\Theta$ given by (\ref{z9}).


\subsection{{\texorpdfstring{$s_1,s_2,s_3,s_4\geq2$}{dd}}}\label{lw1} 
Write the coordinate ring of $\mathbb G^{\mathbb V^{\underline s}}_m$ as
\begin{equation*}
\mathbb Z\left[x_{12},\frac{1}{x_{12}},x_{13},\frac{1}{x_{13}},x_{23},\frac{1}{x_{23}},x_{14},\frac{1}{x_{14}},x_{24},\frac{1}{x_{24}},x_{34},\frac{1}{x_{34}},x_{11},\frac{1}{x_{11}},x_{22},\frac{1}{x_{22}},x_{33},\frac{1}{x_{33}},x_{44},\frac{1}{x_{44}}\right].   
\end{equation*}
Define characters $\chi_{i_1i_2}$ of $\mathbb G^{\mathbb V^{\underline s}}_m$ by
$\chi_{i_1i_2}(x_{12},x_{13},x_{23},x_{14},x_{24},x_{34},x_{11},x_{22},x_{33},x_{44})=x_{i_1i_2}.$ One can show that $\mathcal A^{\mathbb V^{\underline s}}$ is covered by affine open subschemes
\begin{equation}
\begin{split}
&\mathcal A_{12}:={\rm Spec}\,\mathbb Z\left[\frac{\chi_{14}\chi_{23}}{\chi_{13}\chi_{24}},\frac{\chi_{13}\chi_{24}}{\chi_{14}\chi_{23}},\frac{\chi_{12}\chi_{34}}{\chi_{23}\chi_{14}},\frac{\chi_{14}\chi_{33}}{\chi_{13}\chi_{34}},\frac{\chi_{13}\chi_{44}}{\chi_{14}\chi_{34}},\frac{\chi_{13}\chi_{22}}{\chi_{12}\chi_{23}},\frac{\chi_{11}\chi_{23}}{\chi_{12}\chi_{13}}\right],\\
&\mathcal A_{13}:={\rm Spec}\,\mathbb Z\left[\frac{\chi_{23}\chi_{14}}{\chi_{12}\chi_{34}},\frac{\chi_{12}\chi_{34}}{\chi_{23}\chi_{14}},\frac{\chi_{12}\chi_{33}}{\chi_{13}\chi_{23}},\frac{\chi_{13}\chi_{24}}{\chi_{12}\chi_{34}},\frac{\chi_{23}\chi_{44}}{\chi_{24}\chi_{34}},\frac{\chi_{22}\chi_{34}}{\chi_{23}\chi_{24}},\frac{\chi_{11}\chi_{23}}{\chi_{12}\chi_{13}}\right],\\
&\mathcal A_{23}:={\rm Spec}\,\mathbb Z\left[\frac{\chi_{12}\chi_{34}}{\chi_{13}\chi_{24}},\frac{\chi_{13}\chi_{24}}{\chi_{12}\chi_{34}},\frac{\chi_{23}\chi_{14}}{\chi_{12}\chi_{34}},\frac{\chi_{12}\chi_{33}}{\chi_{13}\chi_{23}},\frac{\chi_{13}\chi_{44}}{\chi_{14}\chi_{34}},\frac{\chi_{13}\chi_{22}}{\chi_{12}\chi_{23}},\frac{\chi_{11}\chi_{34}}{\chi_{13}\chi_{14}}\right].      
\end{split}  
\end{equation} 
We view $\mathbb G_m^{\mathbb V^{\underline s}}/\mathbb G_m$ as a subscheme of ${\rm Proj}\,\mathbb Z\left[x_{12},x_{13},x_{23},x_{14},x_{24},x_{34},x_{11},x_{22},x_{33},x_{44}\right]$.

The locally closed embeddings $f_{1A}:W_{1A}\times\mathbb G_m^{\mathbb V^{\underline s}}/\mathbb G_m\rightarrow\mathcal A_{23}\times \mathbb G_m{\Big\backslash\prod\nolimits_{\underline i\in \mathbb V^{\underline s}}\left(\wedge^{\underline i}E_{\bullet}-\{0\}\right)}$, $f_{1B}:W_{1B}\times\mathbb G_m^{\mathbb V^{\underline s}}/\mathbb G_m\rightarrow\mathcal A_{12}\times \mathbb G_m{\Big\backslash\prod\nolimits_{\underline i\in \mathbb V^{\underline s}}\left(\wedge^{\underline i}E_{\bullet}-\{0\}\right)}$, $f_{2A}:W_{2A}\times\mathbb G_m^{\mathbb V^{\underline s}}/\mathbb G_m\rightarrow\mathcal A_{12}\times \mathbb G_m{\Big\backslash\prod\nolimits_{\underline i\in \mathbb V^{\underline s}}\left(\wedge^{\underline i}E_{\bullet}-\{0\}\right)}$, $f_{7A}:W_{7A}\times\mathbb G_m^{\mathbb V^{\underline s}}/\mathbb G_m\rightarrow\mathcal A_{23}\times \mathbb G_m{\Big\backslash\prod\nolimits_{\underline i\in \mathbb V^{\underline s}}\left(\wedge^{\underline i}E_{\bullet}-\{0\}\right)}$, $f_{10E}:W_{10E}\times\mathbb G_m^{\mathbb V^{\underline s}}/\mathbb G_m\rightarrow\mathcal A_{12}\times \mathbb G_m{\Big\backslash\prod\nolimits_{\underline i\in \mathbb V^{\underline s}}\left(\wedge^{\underline i}E_{\bullet}-\{0\}\right)}$, $f_{11A}:W_{11A}\times\mathbb G_m^{\mathbb V^{\underline s}}/\mathbb G_m\rightarrow\mathcal A_{12}\times \mathbb G_m{\Big\backslash\prod\nolimits_{\underline i\in \mathbb V^{\underline s}}\left(\wedge^{\underline i}E_{\bullet}-\{0\}\right)}$  
 are induced
by the ring homomorphisms 
\begin{equation}\label{yeben4}
\frac{\chi_{i_1j_1}\chi_{i_2j_2}}{\chi_{i_3j_3}\chi_{i_4j_4}}\mapsto\frac{g_{i_1j_1}g_{i_2j_2}}{g_{i_3j_3}g_{i_4j_4}}\frac{x_{i_3j_3}x_{i_4j_4}}{x_{i_1j_1}x_{i_2j_2}},\,\,\,\,\,\, z_{(i_1,i_2)}\mapsto\frac{\mathcal P^{\Theta}_{(i_1,i_2)}}{g_{\tau(i_1)\tau(i_2)}}x_{\tau(i_1)\tau(i_2)},  
\end{equation}
where 
\begin{equation*}
\begin{split}
&g_{11}:=\mathcal P^{\Theta}_{(1,2)},\,\,g_{12}:=\mathcal P^{\Theta}_{(1,1+s_1)},\,\,g_{13}:=\mathcal P^{\Theta}_{(1,1+s_1+s_2)},\,\,g_{14}:=\mathcal P^{\Theta}_{(1,1+s_1+s_2+s_3)},\\
&g_{22}:=\mathcal P^{\Theta}_{(1+s_1,2+s_1)},\,\,g_{23}:=\mathcal P^{\Theta}_{(1+s_1,1+s_1+s_2)},\,\,g_{24}:=\mathcal P^{\Theta}_{(1+s_1,1+s_1+s_2+s_3)},\\
&g_{33}:=\mathcal P^{\Theta}_{(1+s_1+s_2,2+s_1+s_2)},\,\,g_{34}:=\mathcal P^{\Theta}_{(1+s_1+s_2,1+s_1+s_2+s_3)},\,\,g_{44}:=\mathcal P^{\Theta}_{(1+s_1+s_2+s_3,2+s_1+s_2+s_3)},\\
\end{split}   
\end{equation*}
with the defining matrices $\Theta$ given by (\ref{didi3}), (\ref{didi3}), (\ref{w3a}),  (\ref{w8a}), (\ref{w11a}), (\ref{w12a}), respectively.

The locally closed embeddings $f_{2B}:W_{2B}\times\mathbb G_m^{\mathbb V^{\underline s}}/\mathbb G_m\rightarrow\mathcal A_{12}\times \mathbb G_m{\Big\backslash\prod\nolimits_{\underline i\in \mathbb V^{\underline s}}\left(\wedge^{\underline i}E_{\bullet}-\{0\}\right)}$, $f_{10A}:W_{10A}\times\mathbb G_m^{\mathbb V^{\underline s}}/\mathbb G_m\rightarrow\mathcal A_{12}\times \mathbb G_m{\Big\backslash\prod\nolimits_{\underline i\in \mathbb V^{\underline s}}\left(\wedge^{\underline i}E_{\bullet}-\{0\}\right)}$ are induced by the ring homomorphisms (\ref{yeben4}), where \begin{equation*}
\begin{split}
&g_{11}:=\mathcal P^{\Theta}_{(1,2)},\,\,g_{12}:=\mathcal P^{\Theta}_{(1,1+s_1)},\,\,g_{13}:=\mathcal P^{\Theta}_{(1,1+s_1+s_2)},\,\,g_{14}:=\mathcal P^{\Theta}_{(1,1+s_1+s_2+s_3)},\\
&g_{22}:=\mathcal P^{\Theta}_{(1+s_1,2+s_1)},\,\,g_{23}:=\mathcal P^{\Theta}_{(1+s_1,1+s_1+s_2)},\,\,g_{24}:=\mathcal P^{\Theta}_{(1+s_1,1+s_1+s_2+s_3)},\\
&g_{33}:=\mathcal P^{\Theta}_{(1+s_1+s_2,2+s_1+s_2)},\,\,g_{34}:=\mathcal P^{\Theta}_{(2+s_1+s_2,1+s_1+s_2+s_3)},\,\,g_{44}:=\mathcal P^{\Theta}_{(1+s_1+s_2+s_3,2+s_1+s_2+s_3)},\\
\end{split}   
\end{equation*}
with $\Theta$ given by (\ref{w3a}), (\ref{w11a}), respectively.


The locally closed embeddings $f_{3}:W_{3}\times\mathbb G_m^{\mathbb V^{\underline s}}/\mathbb G_m\rightarrow\mathcal A_{12}\times \mathbb G_m{\Big\backslash\prod\nolimits_{\underline i\in \mathbb V^{\underline s}}\left(\wedge^{\underline i}E_{\bullet}-\{0\}\right)}$, $f_{9A}:W_{9A}\times\mathbb G_m^{\mathbb V^{\underline s}}/\mathbb G_m\rightarrow\mathcal A_{12}\times \mathbb G_m{\Big\backslash\prod\nolimits_{\underline i\in \mathbb V^{\underline s}}\left(\wedge^{\underline i}E_{\bullet}-\{0\}\right)}$ are induced by the ring homomorphisms (\ref{yeben4}), where \begin{equation*}
\begin{split}
&g_{11}:=\mathcal P^{\Theta}_{(1,2)},\,\,g_{12}:=\mathcal P^{\Theta}_{(1,1+s_1)},\,\,g_{13}:=\mathcal P^{\Theta}_{(1,1+s_1+s_2)},\,\,g_{14}:=\mathcal P^{\Theta}_{(1,1+s_1+s_2+s_3)},\\
&g_{22}:=\mathcal P^{\Theta}_{(1+s_1,2+s_1)},\,\,g_{23}:=\mathcal P^{\Theta}_{(1+s_1,1+s_1+s_2)},\,\,g_{24}:=\mathcal P^{\Theta}_{(1+s_1,1+s_1+s_2+s_3)},\\
&g_{33}:=\mathcal P^{\Theta}_{(1+s_1+s_2,2+s_1+s_2)},\,\,g_{34}:=\mathcal P^{\Theta}_{(1+s_1+s_2,2+s_1+s_2+s_3)},\,\,g_{44}:=\mathcal P^{\Theta}_{(1+s_1+s_2+s_3,2+s_1+s_2+s_3)},\\
\end{split}   
\end{equation*}
with $\Theta$ given by (\ref{w3}), (\ref{w10a}), respectively.


The locally closed embeddings $f_{4}:W_{4}\times\mathbb G_m^{\mathbb V^{\underline s}}/\mathbb G_m\rightarrow\mathcal A_{23}\times \mathbb G_m{\Big\backslash\prod\nolimits_{\underline i\in \mathbb V^{\underline s}}\left(\wedge^{\underline i}E_{\bullet}-\{0\}\right)}$, $f_{8A}:W_{8A}\times\mathbb G_m^{\mathbb V^{\underline s}}/\mathbb G_m\rightarrow\mathcal A_{23}\times \mathbb G_m{\Big\backslash\prod\nolimits_{\underline i\in \mathbb V^{\underline s}}\left(\wedge^{\underline i}E_{\bullet}-\{0\}\right)}$ are induced by the ring homomorphisms (\ref{yeben4}), where \begin{equation*}
\begin{split}
&g_{11}:=\mathcal P^{\Theta}_{(1,2)},\,\,g_{12}:=\mathcal P^{\Theta}_{(1,1+s_1)},\,\,g_{13}:=\mathcal P^{\Theta}_{(1,2+s_1+s_2)},\,\,g_{14}:=\mathcal P^{\Theta}_{(1,1+s_1+s_2+s_3)},\\
&g_{22}:=\mathcal P^{\Theta}_{(1+s_1,2+s_1)},\,\,g_{23}:=\mathcal P^{\Theta}_{(1+s_1,1+s_1+s_2)},\,\,g_{24}:=\mathcal P^{\Theta}_{(1+s_1,1+s_1+s_2+s_3)},\\
&g_{33}:=\mathcal P^{\Theta}_{(1+s_1+s_2,2+s_1+s_2)},\,\,g_{34}:=\mathcal P^{\Theta}_{(2+s_1+s_2,1+s_1+s_2+s_3)},\,\,g_{44}:=\mathcal P^{\Theta}_{(1+s_1+s_2+s_3,2+s_1+s_2+s_3)},\\
\end{split}   
\end{equation*}
with $\Theta$ given by (\ref{w4}), (\ref{w9a}), respectively.


The locally closed embedding $f_{5}:W_{5}\times\mathbb G_m^{\mathbb V^{\underline s}}/\mathbb G_m\rightarrow\mathcal A_{23}\times \mathbb G_m{\Big\backslash\prod\nolimits_{\underline i\in \mathbb V^{\underline s}}\left(\wedge^{\underline i}E_{\bullet}-\{0\}\right)}$  is induced by the ring homomorphism (\ref{yeben4}), where \begin{equation*}
\begin{split}
&g_{11}:=\mathcal P^{\Theta}_{(1,2)},\,\,g_{12}:=\mathcal P^{\Theta}_{(1,1+s_1)},\,\,g_{13}:=\mathcal P^{\Theta}_{(1,1+s_1+s_2)},\,\,g_{14}:=\mathcal P^{\Theta}_{(1,2+s_1+s_2+s_3)},\\
&g_{22}:=\mathcal P^{\Theta}_{(1+s_1,2+s_1)},\,\,g_{23}:=\mathcal P^{\Theta}_{(1+s_1,1+s_1+s_2)},\,\,g_{24}:=\mathcal P^{\Theta}_{(1+s_1,1+s_1+s_2+s_3)},\\
&g_{33}:=\mathcal P^{\Theta}_{(1+s_1+s_2,2+s_1+s_2)},\,\,g_{34}:=\mathcal P^{\Theta}_{(1+s_1+s_2,1+s_1+s_2+s_3)},\,\,g_{44}:=\mathcal P^{\Theta}_{(1+s_1+s_2+s_3,2+s_1+s_2+s_3)},\\
\end{split}   
\end{equation*}
with $\Theta$ given by (\ref{w5}).


The locally closed embedding $f_{6}:W_{6}\times\mathbb G_m^{\mathbb V^{\underline s}}/\mathbb G_m\rightarrow\mathcal A_{23}\times \mathbb G_m{\Big\backslash\prod\nolimits_{\underline i\in \mathbb V^{\underline s}}\left(\wedge^{\underline i}E_{\bullet}-\{0\}\right)}$  is induced by the ring homomorphism (\ref{yeben4}), where \begin{equation*}
\begin{split}
&g_{11}:=\mathcal P^{\Theta}_{(1,2)},\,\,g_{12}:=\mathcal P^{\Theta}_{(1,1+s_1)},\,\,g_{13}:=\mathcal P^{\Theta}_{(1,2+s_1+s_2)},\,\,g_{14}:=\mathcal P^{\Theta}_{(1,2+s_1+s_2+s_3)},\\
&g_{22}:=\mathcal P^{\Theta}_{(1+s_1,2+s_1)},\,\,g_{23}:=\mathcal P^{\Theta}_{(1+s_1,1+s_1+s_2)},\,\,g_{24}:=\mathcal P^{\Theta}_{(1+s_1,1+s_1+s_2+s_3)},\\
&g_{33}:=\mathcal P^{\Theta}_{(1+s_1+s_2,2+s_1+s_2)},\,\,g_{34}:=\mathcal P^{\Theta}_{(2+s_1+s_2,1+s_1+s_2+s_3)},\,\,g_{44}:=\mathcal P^{\Theta}_{(1+s_1+s_2+s_3,2+s_1+s_2+s_3)},\\
\end{split}   
\end{equation*}
with $\Theta$ given by (\ref{w6}).


The locally closed embeddings $f_{7B}:W_{7B}\times\mathbb G_m^{\mathbb V^{\underline s}}/\mathbb G_m\rightarrow\mathcal A_{23}\times \mathbb G_m{\Big\backslash\prod\nolimits_{\underline i\in \mathbb V^{\underline s}}\left(\wedge^{\underline i}E_{\bullet}-\{0\}\right)}$, $f_{10F}:W_{10F}\times\mathbb G_m^{\mathbb V^{\underline s}}/\mathbb G_m\rightarrow\mathcal A_{12}\times \mathbb G_m{\Big\backslash\prod\nolimits_{\underline i\in \mathbb V^{\underline s}}\left(\wedge^{\underline i}E_{\bullet}-\{0\}\right)}$, $f_{11B}:W_{11B}\times\mathbb G_m^{\mathbb V^{\underline s}}/\mathbb G_m\rightarrow\mathcal A_{12}\times \mathbb G_m{\Big\backslash\prod\nolimits_{\underline i\in \mathbb V^{\underline s}}\left(\wedge^{\underline i}E_{\bullet}-\{0\}\right)}$ are induced by the ring homomorphisms (\ref{yeben4}), where \begin{equation*}
\begin{split}
&g_{11}:=\mathcal P^{\Theta}_{(1,2)},\,\,g_{12}:=\mathcal P^{\Theta}_{(1,1+s_1)},\,\,g_{13}:=\mathcal P^{\Theta}_{(1,1+s_1+s_2)},\,\,g_{14}:=\mathcal P^{\Theta}_{(2,1+s_1+s_2+s_3)},\\
&g_{22}:=\mathcal P^{\Theta}_{(1+s_1,2+s_1)},\,\,g_{23}:=\mathcal P^{\Theta}_{(1+s_1,1+s_1+s_2)},\,\,g_{24}:=\mathcal P^{\Theta}_{(1+s_1,1+s_1+s_2+s_3)},\\
&g_{33}:=\mathcal P^{\Theta}_{(1+s_1+s_2,2+s_1+s_2)},\,\,g_{34}:=\mathcal P^{\Theta}_{(1+s_1+s_2,1+s_1+s_2+s_3)},\,\,g_{44}:=\mathcal P^{\Theta}_{(1+s_1+s_2+s_3,2+s_1+s_2+s_3)},\\
\end{split}   
\end{equation*}
with $\Theta$ given by (\ref{w8a}), (\ref{w11a}), (\ref{w12a}), respectively.


The locally closed embedding $f_{8B}:W_{8B}\times\mathbb G_m^{\mathbb V^{\underline s}}/\mathbb G_m\rightarrow\mathcal A_{23}\times \mathbb G_m{\Big\backslash\prod\nolimits_{\underline i\in \mathbb V^{\underline s}}\left(\wedge^{\underline i}E_{\bullet}-\{0\}\right)}$  is induced by the ring homomorphism (\ref{yeben4}), where \begin{equation*}
\begin{split}
&g_{11}:=\mathcal P^{\Theta}_{(1,2)},\,\,g_{12}:=\mathcal P^{\Theta}_{(1,1+s_1)},\,\,g_{13}:=\mathcal P^{\Theta}_{(1,2+s_1+s_2)},\,\,g_{14}:=\mathcal P^{\Theta}_{(2,1+s_1+s_2+s_3)},\\
&g_{22}:=\mathcal P^{\Theta}_{(1+s_1,2+s_1)},\,\,g_{23}:=\mathcal P^{\Theta}_{(1+s_1,1+s_1+s_2)},\,\,g_{24}:=\mathcal P^{\Theta}_{(1+s_1,1+s_1+s_2+s_3)},\\
&g_{33}:=\mathcal P^{\Theta}_{(1+s_1+s_2,2+s_1+s_2)},\,\,g_{34}:=\mathcal P^{\Theta}_{(2+s_1+s_2,1+s_1+s_2+s_3)},\,\,g_{44}:=\mathcal P^{\Theta}_{(1+s_1+s_2+s_3,2+s_1+s_2+s_3)},\\
\end{split}   
\end{equation*}
with $\Theta$ given by (\ref{w9a}).

The locally closed embedding $f_{9B}:W_{9B}\times\mathbb G_m^{\mathbb V^{\underline s}}/\mathbb G_m\rightarrow\mathcal A_{12}\times \mathbb G_m{\Big\backslash\prod\nolimits_{\underline i\in \mathbb V^{\underline s}}\left(\wedge^{\underline i}E_{\bullet}-\{0\}\right)}$  is induced by the ring homomorphism (\ref{yeben4}), where \begin{equation*}
\begin{split}
&g_{11}:=\mathcal P^{\Theta}_{(1,2)},\,\,g_{12}:=\mathcal P^{\Theta}_{(1,1+s_1)},\,\,g_{13}:=\mathcal P^{\Theta}_{(1,1+s_1+s_2)},\,\,g_{14}:=\mathcal P^{\Theta}_{(2,1+s_1+s_2+s_3)},\\
&g_{22}:=\mathcal P^{\Theta}_{(1+s_1,2+s_1)},\,\,g_{23}:=\mathcal P^{\Theta}_{(1+s_1,1+s_1+s_2)},\,\,g_{24}:=\mathcal P^{\Theta}_{(1+s_1,1+s_1+s_2+s_3)},\\
&g_{33}:=\mathcal P^{\Theta}_{(1+s_1+s_2,2+s_1+s_2)},\,\,g_{34}:=\mathcal P^{\Theta}_{(1+s_1+s_2,2+s_1+s_2+s_3)},\,\,g_{44}:=\mathcal P^{\Theta}_{(1+s_1+s_2+s_3,2+s_1+s_2+s_3)},\\
\end{split}   
\end{equation*}
with $\Theta$ given by (\ref{w10a}).


The locally closed embedding $f_{9C}:W_{9C}\times\mathbb G_m^{\mathbb V^{\underline s}}/\mathbb G_m\rightarrow\mathcal A_{12}\times \mathbb G_m{\Big\backslash\prod\nolimits_{\underline i\in \mathbb V^{\underline s}}\left(\wedge^{\underline i}E_{\bullet}-\{0\}\right)}$ is induced by the ring homomorphism (\ref{yeben4}), where \begin{equation*}
\begin{split}
&g_{11}:=\mathcal P^{\Theta}_{(1,2)},\,\,g_{12}:=\mathcal P^{\Theta}_{(1,1+s_1)},\,\,g_{13}:=\mathcal P^{\Theta}_{(2,1+s_1+s_2)},\,\,g_{14}:=\mathcal P^{\Theta}_{(1,1+s_1+s_2+s_3)},\\
&g_{22}:=\mathcal P^{\Theta}_{(1+s_1,2+s_1)},\,\,g_{23}:=\mathcal P^{\Theta}_{(1+s_1,1+s_1+s_2)},\,\,g_{24}:=\mathcal P^{\Theta}_{(1+s_1,1+s_1+s_2+s_3)},\\
&g_{33}:=\mathcal P^{\Theta}_{(1+s_1+s_2,2+s_1+s_2)},\,\,g_{34}:=\mathcal P^{\Theta}_{(1+s_1+s_2,2+s_1+s_2+s_3)},\,\,g_{44}:=\mathcal P^{\Theta}_{(1+s_1+s_2+s_3,2+s_1+s_2+s_3)},\\
\end{split}   
\end{equation*}
with $\Theta$ given by (\ref{w10a}).


The locally closed embedding $f_{9D}:W_{9D}\times\mathbb G_m^{\mathbb V^{\underline s}}/\mathbb G_m\rightarrow\mathcal A_{12}\times \mathbb G_m{\Big\backslash\prod\nolimits_{\underline i\in \mathbb V^{\underline s}}\left(\wedge^{\underline i}E_{\bullet}-\{0\}\right)}$ is induced by the ring homomorphism (\ref{yeben4}), where \begin{equation*}
\begin{split}
&g_{11}:=\mathcal P^{\Theta}_{(1,2)},\,\,g_{12}:=\mathcal P^{\Theta}_{(1,1+s_1)},\,\,g_{13}:=\mathcal P^{\Theta}_{(2,1+s_1+s_2)},\,\,g_{14}:=\mathcal P^{\Theta}_{(2,1+s_1+s_2+s_3)},\\
&g_{22}:=\mathcal P^{\Theta}_{(1+s_1,2+s_1)},\,\,g_{23}:=\mathcal P^{\Theta}_{(1+s_1,1+s_1+s_2)},\,\,g_{24}:=\mathcal P^{\Theta}_{(1+s_1,1+s_1+s_2+s_3)},\\
&g_{33}:=\mathcal P^{\Theta}_{(1+s_1+s_2,2+s_1+s_2)},\,\,g_{34}:=\mathcal P^{\Theta}_{(1+s_1+s_2,2+s_1+s_2+s_3)},\,\,g_{44}:=\mathcal P^{\Theta}_{(1+s_1+s_2+s_3,2+s_1+s_2+s_3)},\\
\end{split}   
\end{equation*}
with $\Theta$ given by (\ref{w10a}).

The locally closed embedding $f_{10B}:W_{10B}\times\mathbb G_m^{\mathbb V^{\underline s}}/\mathbb G_m\rightarrow\mathcal A_{12}\times \mathbb G_m{\Big\backslash\prod\nolimits_{\underline i\in \mathbb V^{\underline s}}\left(\wedge^{\underline i}E_{\bullet}-\{0\}\right)}$ is induced by the ring homomorphism (\ref{yeben4}), where \begin{equation*}
\begin{split}
&g_{11}:=\mathcal P^{\Theta}_{(1,2)},\,\,g_{12}:=\mathcal P^{\Theta}_{(1,1+s_1)},\,\,g_{13}:=\mathcal P^{\Theta}_{(1,1+s_1+s_2)},\,\,g_{14}:=\mathcal P^{\Theta}_{(2,1+s_1+s_2+s_3)},\\
&g_{22}:=\mathcal P^{\Theta}_{(1+s_1,2+s_1)},\,\,g_{23}:=\mathcal P^{\Theta}_{(1+s_1,1+s_1+s_2)},\,\,g_{24}:=\mathcal P^{\Theta}_{(1+s_1,1+s_1+s_2+s_3)},\\
&g_{33}:=\mathcal P^{\Theta}_{(1+s_1+s_2,2+s_1+s_2)},\,\,g_{34}:=\mathcal P^{\Theta}_{(2+s_1+s_2,1+s_1+s_2+s_3)},\,\,g_{44}:=\mathcal P^{\Theta}_{(1+s_1+s_2+s_3,2+s_1+s_2+s_3)},\\
\end{split}   
\end{equation*}
with $\Theta$ given by (\ref{w11a}).


The locally closed embedding $f_{10C}:W_{10C}\times\mathbb G_m^{\mathbb V^{\underline s}}/\mathbb G_m\rightarrow\mathcal A_{12}\times \mathbb G_m{\Big\backslash\prod\nolimits_{\underline i\in \mathbb V^{\underline s}}\left(\wedge^{\underline i}E_{\bullet}-\{0\}\right)}$ is induced by the ring homomorphism (\ref{yeben4}), where \begin{equation*}
\begin{split}
&g_{11}:=\mathcal P^{\Theta}_{(1,2)},\,\,g_{12}:=\mathcal P^{\Theta}_{(1,1+s_1)},\,\,g_{13}:=\mathcal P^{\Theta}_{(2,1+s_1+s_2)},\,\,g_{14}:=\mathcal P^{\Theta}_{(1,1+s_1+s_2+s_3)},\\
&g_{22}:=\mathcal P^{\Theta}_{(1+s_1,2+s_1)},\,\,g_{23}:=\mathcal P^{\Theta}_{(1+s_1,1+s_1+s_2)},\,\,g_{24}:=\mathcal P^{\Theta}_{(1+s_1,1+s_1+s_2+s_3)},\\
&g_{33}:=\mathcal P^{\Theta}_{(1+s_1+s_2,2+s_1+s_2)},\,\,g_{34}:=\mathcal P^{\Theta}_{(2+s_1+s_2,1+s_1+s_2+s_3)},\,\,g_{44}:=\mathcal P^{\Theta}_{(1+s_1+s_2+s_3,2+s_1+s_2+s_3)},\\
\end{split}   
\end{equation*}
with $\Theta$ given by (\ref{w11a}).


The locally closed embedding $f_{10D}:W_{10D}\times\mathbb G_m^{\mathbb V^{\underline s}}/\mathbb G_m\rightarrow\mathcal A_{12}\times \mathbb G_m{\Big\backslash\prod\nolimits_{\underline i\in \mathbb V^{\underline s}}\left(\wedge^{\underline i}E_{\bullet}-\{0\}\right)}$
is induced by the ring homomorphism (\ref{yeben4}), where \begin{equation*}
\begin{split}
&g_{11}:=\mathcal P^{\Theta}_{(1,2)},\,\,g_{12}:=\mathcal P^{\Theta}_{(1,1+s_1)},\,\,g_{13}:=\mathcal P^{\Theta}_{(2,1+s_1+s_2)},\,\,g_{14}:=\mathcal P^{\Theta}_{(2,1+s_1+s_2+s_3)},\\
&g_{22}:=\mathcal P^{\Theta}_{(1+s_1,2+s_1)},\,\,g_{23}:=\mathcal P^{\Theta}_{(1+s_1,1+s_1+s_2)},\,\,g_{24}:=\mathcal P^{\Theta}_{(1+s_1,1+s_1+s_2+s_3)},\\
&g_{33}:=\mathcal P^{\Theta}_{(1+s_1+s_2,2+s_1+s_2)},\,\,g_{34}:=\mathcal P^{\Theta}_{(2+s_1+s_2,1+s_1+s_2+s_3)},\,\,g_{44}:=\mathcal P^{\Theta}_{(1+s_1+s_2+s_3,2+s_1+s_2+s_3)},\\
\end{split}   
\end{equation*}
with $\Theta$ given by (\ref{w11a}).

The locally closed embeddings $f_{10G}:W_{10G}\times\mathbb G_m^{\mathbb V^{\underline s}}/\mathbb G_m\rightarrow\mathcal A_{12}\times \mathbb G_m{\Big\backslash\prod\nolimits_{\underline i\in \mathbb V^{\underline s}}\left(\wedge^{\underline i}E_{\bullet}-\{0\}\right)}$, $f_{11C}:W_{11C}\times\mathbb G_m^{\mathbb V^{\underline s}}/\mathbb G_m\rightarrow\mathcal A_{12}\times \mathbb G_m{\Big\backslash\prod\nolimits_{\underline i\in \mathbb V^{\underline s}}\left(\wedge^{\underline i}E_{\bullet}-\{0\}\right)}$ are induced by the ring homomorphisms (\ref{yeben4}), where \begin{equation*}
\begin{split}
&g_{11}:=\mathcal P^{\Theta}_{(1,2)},\,\,g_{12}:=\mathcal P^{\Theta}_{(1,1+s_1)},\,\,g_{13}:=\mathcal P^{\Theta}_{(2,1+s_1+s_2)},\,\,g_{14}:=\mathcal P^{\Theta}_{(1,1+s_1+s_2+s_3)},\\
&g_{22}:=\mathcal P^{\Theta}_{(1+s_1,2+s_1)},\,\,g_{23}:=\mathcal P^{\Theta}_{(1+s_1,1+s_1+s_2)},\,\,g_{24}:=\mathcal P^{\Theta}_{(1+s_1,1+s_1+s_2+s_3)},\\
&g_{33}:=\mathcal P^{\Theta}_{(1+s_1+s_2,2+s_1+s_2)},\,\,g_{34}:=\mathcal P^{\Theta}_{(1+s_1+s_2,1+s_1+s_2+s_3)},\,\,g_{44}:=\mathcal P^{\Theta}_{(1+s_1+s_2+s_3,2+s_1+s_2+s_3)},\\
\end{split}   
\end{equation*}
with $\Theta$ given by (\ref{w11a}), (\ref{w12a}), respectively.


The locally closed embeddings $f_{10H}:W_{10H}\times\mathbb G_m^{\mathbb V^{\underline s}}/\mathbb G_m\rightarrow\mathcal A_{12}\times \mathbb G_m{\Big\backslash\prod\nolimits_{\underline i\in \mathbb V^{\underline s}}\left(\wedge^{\underline i}E_{\bullet}-\{0\}\right)}$, $f_{11D}:W_{11D}\times\mathbb G_m^{\mathbb V^{\underline s}}/\mathbb G_m\rightarrow\mathcal A_{12}\times \mathbb G_m{\Big\backslash\prod\nolimits_{\underline i\in \mathbb V^{\underline s}}\left(\wedge^{\underline i}E_{\bullet}-\{0\}\right)}$  are induced by the ring homomorphisms (\ref{yeben4}), where \begin{equation*}
\begin{split}
&g_{11}:=\mathcal P^{\Theta}_{(1,2)},\,\,g_{12}:=\mathcal P^{\Theta}_{(1,1+s_1)},\,\,g_{13}:=\mathcal P^{\Theta}_{(2,1+s_1+s_2)},\,\,g_{14}:=\mathcal P^{\Theta}_{(2,1+s_1+s_2+s_3)},\\
&g_{22}:=\mathcal P^{\Theta}_{(1+s_1,2+s_1)},\,\,g_{23}:=\mathcal P^{\Theta}_{(1+s_1,1+s_1+s_2)},\,\,g_{24}:=\mathcal P^{\Theta}_{(1+s_1,1+s_1+s_2+s_3)},\\
&g_{33}:=\mathcal P^{\Theta}_{(1+s_1+s_2,2+s_1+s_2)},\,\,g_{34}:=\mathcal P^{\Theta}_{(1+s_1+s_2,1+s_1+s_2+s_3)},\,\,g_{44}:=\mathcal P^{\Theta}_{(1+s_1+s_2+s_3,2+s_1+s_2+s_3)},\\
\end{split}   
\end{equation*}
with $\Theta$ given by (\ref{w11a}), (\ref{w12a}), respectively.


\bigskip

\noindent H. Fang, School of Mathematical Sciences, Peking University, Beijing, Beijing, 100871, China. (hlfang$@$pku.edu.cn)

\noindent M. Zhang, School of Mathematical Sciences, Peking University, Beijing, Beijing, 100871, China. (zhangmy1990$@$gmail.com)

\end{document}